\pdfoutput=1
\RequirePackage{silence}
\PassOptionsToPackage{table,svgnames}{xcolor}
\documentclass[a4paper,svgnames]{amsart}
\usepackage[hmarginratio={1:1},vmarginratio={1:1},lmargin=60.0pt,tmargin=60.0pt]{geometry}

\usepackage{booktabs}
\allowdisplaybreaks

\usepackage[T1]{fontenc}
\usepackage{latexsym,exscale,amsfonts,amssymb,mathtools}
\usepackage{amsmath,amsthm,amsfonts,amssymb,amscd,textcomp,bbm}
\usepackage{mathrsfs,stackrel}
\usepackage{etoolbox}
\usepackage{tikz,tikz-cd}
\usetikzlibrary{matrix,arrows.meta,decorations.markings,shapes.geometric,shapes.multipart,decorations.pathreplacing,decorations.pathmorphing}
\usepackage{enumitem}
\setlist[enumerate]{itemsep=0.15cm,label=\emph{\upshape(\alph*)}}
\setlist[enumerate,2]{itemsep=0.15cm,label=\emph{\upshape(\roman*)}}
\usepackage{aliascnt}
\usepackage{ytableau}
\usepackage{xparse}
\usepackage{cite}

\usepackage{dynkin-diagrams}

\definecolor{strokecol}{rgb}{0.0,0.0,0.0}
\pgfsetstrokecolor{strokecol}

\usepackage{array}
\newcolumntype{C}{>{$}c<{$}}

\newcommand\floor[2]{\lfloor\tfrac{#1}{#2}\rfloor}

\usepackage{textcomp}
\usepackage{listings}
\lstdefinestyle{sagemath}{
    language=python,
    stringstyle=\color{Coral},
    tabsize=4,
    showtabs=false,
    showspaces=false,
    showstringspaces=false,
    inputencoding=utf8,
    extendedchars=true,
    upquote=true,
    backgroundcolor=\color{Ivory},
    numbers=none,
    breaklines=true,
    commentstyle=\color{Brown},
    resetmargins=true,
    keywordstyle={[1]\color{FireBrick}},
    keywordstyle={[2]\color{Blue}},
    keywordstyle={[3]\color{ForestGreen}},
    morekeywords={
        sage,
    },
    classoffset=1,
    morekeywords={
        RootSystem,
        crystals,
        basis,
        weight_space,
        subcrystal,
        LSPaths,
        digraph,
        crystals,
    },
    classoffset=2,
    morekeywords={
        CheckForRepeatedPositiveRoots,
    },
}
\lstnewenvironment{sagemath}{\lstset{style=sagemath}}{}

\usepackage{colortbl}
\usepackage{xcolor}

\definecolor{mygray}{gray}{0.6}
\definecolor{mygraydark}{gray}{0.4}
\definecolor{mygraylight}{gray}{0.85}
\definecolor{spinach}{RGB}{46,139,87}
\definecolor{tomato}{RGB}{255,99,71}
\definecolor{orchid}{RGB}{143,40,194}
\definecolor{neon}{RGB}{77,77,255}
\definecolor{pumpkin}{RGB}{224,180,80}
\definecolor{citron}{RGB}{190,180,90}

\definecolor{lava}{RGB}{207,16,32}
\definecolor{cream}{RGB}{255,253,208}
\definecolor{verdigris}{RGB}{67,179,174}
\definecolor{Black}{RGB}{0,0,0}
\definecolor{mydarkblue}{RGB}{10,10,170}
\definecolor{darkspinach}{RGB}{20,70,20}
\definecolor{darktomato}{RGB}{155,40,30}
\definecolor{darkorchid}{RGB}{50,10,100}
\definecolor{darklava}{RGB}{150,8,16}

\newcommand{\leftsquigarrow}{\,\reflectbox{$\rightsquigarrow$}\,}

\usepackage{todonotes}
\newcommand\AM[2][]{\todo[color=blue!20,inline,#1]{Andrew: #2}}
\newcommand\DT[2][]{\todo[color=orchid!20,inline,#1]{Dani: #2}}

\let\emph\relax
\DeclareTextFontCommand{\emph}{\bfseries\em}

\newcommand{\placeholder}{{}_{-}}

\newcommand{\mystrut}{\rule[-0.2\baselineskip]{0pt}{0.9\baselineskip}}

\let\<=\langle
\let\>=\rangle

\def\map#1#2{\,{:}\,#1\!\longrightarrow\!#2}

\newcommand\Luf{\mathsf{f}}

\newcommand\qbinom[3][]{\genfrac{[}{]}{0pt}{}{#2}{#3}_{#1}}
\DeclareMathOperator{\Proj}{Proj}
\DeclareMathOperator{\Rep}{Rep}
\DeclareMathOperator{\End}{End}
\DeclareMathOperator{\Hom}{Hom}

\newcommand{\ie}{\text{i.e.}}
\newcommand{\eg}{\text{e.g.}}
\newcommand{\cf}{\text{cf.}}
\newcommand{\etc}{\text{etc.}}

\newcommand{\ver}{\text{verbatim}}

\newcommand{\muta}{\text{mutatis mutandis}}
\newcommand{\loccit}{\text{loc. cit.}}

\newcommand{\R}{\mathbb{R}}
\newcommand{\N}{\mathbb{Z}_{\geq 0}}
\newcommand{\Q}{\mathbb{Q}}
\newcommand{\Z}{\mathbb{Z}}
\newcommand{\ring}{R}

\newcommand{\typea}[1][e]{A_{#1}}
\newcommand{\typeb}[1][e]{B_{#1}}
\newcommand{\typec}[1][e]{C_{#1}}
\newcommand{\typed}[1][e]{D_{#1}}
\newcommand{\typee}[1][6]{E_{#1}}
\newcommand{\typef}[1][4]{F_{#1}}
\newcommand{\typeg}[1][2]{G_{#1}}

\newcommand{\aonetypea}[1][e]{A_{#1}^{(1)}}
\newcommand{\atwotypea}[1][e>1]{A_{2\cdot #1}^{(2)}}
\newcommand{\athreetypea}[1][e>3]{A_{2\cdot #1-1}^{(2)}}
\newcommand{\atypeb}[1][e>1]{B_{#1}^{(1)}}
\newcommand{\atypec}[1][e>1]{C_{#1}^{(1)}}
\newcommand{\aonetyped}[1][e>3]{D_{#1}^{(1)}}
\newcommand{\atwotyped}[1][e>1]{D_{#1+1}^{(2)}}
\newcommand{\athreetyped}[1][4]{D_{#1}^{(3)}}
\newcommand{\atypee}[1][6]{E_{#1}^{(1)}}
\newcommand{\atwotypee}[1][6]{E_{#1}^{(2)}}
\newcommand{\atypef}[1][4]{F_{#1}^{(1)}}
\newcommand{\atypeg}[1][2]{G_{#1}^{(1)}}

\newcommand{\dynkinaut}{w_{0}}

\newcommand{\quiver}{\Gamma}
\newcommand\CrystalVertices[1][n]{V^{\bfweight}_{#1}}
\newcommand\aCrystalVertices[1][n]{V^{\abfweight}_{#1}}
\newcommand\dCrystalVertices[1][n]{\dot V^{\bfweight}_{#1}}
\newcommand\daCrystalVertices[1][n]{\dot V^{\abfweight}_{#1}}
\newcommand{\vertices}{I}
\newcommand{\edges}{E}
\newcommand{\sym}[1][n]{{\mathfrak{S}_{#1}}}

\newcommand{\sroot}[1][i]{\alpha_{#1}}
\newcommand{\scoroot}[1][i]{h_{#1}}
\newcommand{\sroots}{\Pi}

\newcommand{\fweight}[1][i]{\Lambda_{#1}}
\newcommand{\fweights}{\Lambda}

\newcommand{\crystal}{\mathcal{B}}
\newcommand\Crystalgraph[1][\abfweight]{\mathcal{G}_{\otimes}(#1)}
\NewDocumentCommand{\crystalgraph}{sO{\fweight}}{\mathcal{G}\IfBooleanF{#1}{(#2)}}

\newcommand{\dist}[1][\lambda]{\mathrm{d}(#1)}

\newcommand{\Parts}[2]{\mathrm{P}_{#1}^{#2}}
\newcommand{\fParts}[2]{\mathrm{P}_{#1}^{#2}}
\newcommand{\Res}[1][\blam]{\mathrm{res}(#1)}

\newcommand{\ppath}[1][\blam]{\mathbbm{p}_{#1}}
\newcommand{\qpath}[1][\blam]{\mathrm{p}_{#1}}

\newcommand\bpath{\boldsymbol{\mathrm{p}}}
\newcommand\cpath{\boldsymbol{\mathrm{q}}}

\newcommand{\affine}[1]{\underline{#1}}

\newcommand{\bu}{\boldsymbol{u}}
\newcommand{\bx}{\boldsymbol{x}}
\newcommand{\bi}{\boldsymbol{i}}
\newcommand{\bj}{\boldsymbol{j}}

\newcommand{\by}{\boldsymbol{y}}
\newcommand{\ba}{\boldsymbol{a}}
\newcommand{\bb}{\boldsymbol{b}}
\newcommand{\bfi}{\boldsymbol{f}}
\newcommand{\brho}{\boldsymbol{\rho}}
\newcommand{\charge}{\boldsymbol{\kappa}}
\newcommand{\bsig}{\boldsymbol{\sigma}}
\newcommand{\blam}{\boldsymbol{\lambda}}
\newcommand{\bbeta}{\boldsymbol{\beta}}
\newcommand{\bmu}{\boldsymbol{\mu}}
\newcommand{\bnu}{\boldsymbol{\nu}}

\newcommand{\bS}{\boldsymbol{S}}
\newcommand{\bT}{\boldsymbol{T}}

\newcommand{\bfweight}{\boldsymbol{\Lambda}}
\newcommand{\abfweight}{\affine{\boldsymbol{\Lambda}}}

\newcommand\idem[1][\blam]{\mathbf{1}_{#1}}
\newcommand\dotidem[1][\blam]{\mathbf{1}_{#1}^{y}}
\newcommand\sandwich[3]{h_{#1}^{#2,#3}}
\newcommand\fsandwich[2]{h_{#1}^{#2}}
\newcommand\daffine[1][\ba]{y^{#1}}
\newcommand\zeetwo[1][\bfi]{y^{#1}}

\newcommand{\aWA}[1][n]{\mathscr{W}_{#1}^{\affine{\brho}}(X)}
\newcommand{\WA}[1][n]{\mathscr{W}_{#1}^{\affine\brho}}
\newcommand{\WAc}[1][n]{\mathscr{R}_{#1}^{\brho}}

\newcommand\RLlam[1][\blam]{L^{\mathscr{R}}_{#1}}

\newcommand\RDlam[1][\blam]{\Delta^{\mathscr{R}}_{#1}}
\newcommand\WPlam[1][\blam]{P^{\mathscr{W}}_{#1}}
\newcommand\RPlam[1][\blam]{P^{\mathscr{R}}_{#1}}

\newcommand{\TA}[1][n]{\mathcal{W}_{#1}^{\brho}}
\newcommand{\TAc}[1][n]{\mathcal{R}_{#1}^{\brho}}

\newcommand{\gWA}[1][n]{\mathscr{W}_{#1}^{\bbeta}}
\newcommand{\gWAc}[1][n]{\mathscr{R}_{#1}^{\bbeta}}

\newcommand{\uideal}[1][n]{\mathscr{U}}

\newcommand{\hell}{\affine{\ell}}

\newcommand\Affch[1][\blam]{A^{y}(#1)}
\newcommand\Finch[1][\blam]{F^{y}(#1)}
\newcommand\Sandch[1][\blam]{H^{y}(#1)}
\newcommand{\fsand}[1][\lambda]{\mathscr{F}_{#1}}

\newcommand\BX[1][{\WA[n](X)}]{\mathrm{D}_{#1}}
\newcommand\BXc[1][{\WAc[n](X)}]{\mathrm{D}_{#1}}
\newcommand\EX[1][{\TA[n](X)}]{\mathrm{E}_{#1}}
\newcommand\EXc[1][{\TAc[n](X)}]{\mathrm{E}_{#1}}

\newcommand\WABasis[1][{\WA[n](X)}]{\mathcal{B}_{#1}}

\newcommand\Face[1][\abfweight]{\mathop{\rm Face}\nolimits_{#1}}
\newcommand\fFace[1][\bfweight]{\mathop{\rm Face}\nolimits_{#1}}
\newcommand\Detour[1][]{\mathop{\rm Detour}\nolimits_{#1}}

\renewcommand{\dim}{\mathrm{rk}}
\newcommand{\vpar}{v}
\newcommand{\grdim}[1][\vpar]{\mathrm{rk}^{#1}_{R}}

\newcommand{\Pcal}{\mathcal{P}}

\newcommand{\sandorder}[1][\R]{<_{#1}}
\newcommand{\rsandorder}[1][\R]{>_{#1}}
\newcommand{\gsandorder}[1][\R]{\geq_{#1}}
\newcommand{\lsandorder}[1][\R]{\leq_{#1}}
\newcommand{\sand}[1][\lambda]{\mathscr{H}_{#1}}
\newcommand{\asandbasis}[1][\lambda]{\affine{B}_{#1}^{\mathscr{H}}}
\newcommand{\sandbasis}[1][\lambda]{B_{#1}^{\mathscr{H}}}

\tikzset{
anchorbase/.style={baseline={([yshift=#1]current bounding box.center)}},
anchorbase/.default={-0.5ex},
dot colour/.initial=black,
dot colour/.default=black,
dynkin node/.style args={#1,#2,#3}{circle,inner sep=1.8pt,fill=#1,label={[#2=1mm]$#3$}},
dynkin/.style={dynkin node={DarkBlue,below,#1}},
dynkin/.default=,
tomato dynkin/.style={dynkin node={tomato,#1,0}},
tomato dynkin/.default=below,
tinynodes/.style={font=\tiny,text height=0.25ex,text depth=0.05ex},
smallnodes/.style={font=\scriptsize,text height=0.75ex,text depth=0.15ex},
mor/.style={line width=0.75,color=black,fill=cream},
dots/.style={line width=1pt,line cap=round, gray, dash pattern=on 0pt off 2\pgflinewidth},
redstring/.style = {draw=red!50,fill=none,line width=0.35mm,preaction={draw=red,line width=2.5pt,-},nodes={color=red}},
affine/.style= {draw=citron!50,fill=none,
line width=0.35mm,preaction={draw=citron,line width=2.5pt,-},nodes={color=citron}},
solid/.style = {draw=blue,fill=none,dot colour=blue,line width=0.4mm,nodes={color=blue}},
ghost/.style = {draw=darkgray,fill=none,dot colour=darkgray,
densely dashed,line width=0.4mm,nodes={color=darkgray}},
dghost/.style = {draw=darkgray,double,fill=none,dot colour=darkgray,
densely dashed,line width=0.4mm,nodes={color=darkgray}},
crossline/.style={preaction={draw=white,line width=4.75pt,-},preaction={draw=black,line width=0.9pt,-}},
dot/.style = {decoration={markings,
post length=0.25mm,
pre length=0.25mm,
mark=at position #1 with {\node[circle,radius=0.3cm,inner sep=-2.2pt,color=\pgfkeysvalueof{/tikz/dot colour},fill=\pgfkeysvalueof{/tikz/dot colour}]{};}
},
postaction={decorate}
},
dot/.default=0.5,
rootlyndon/.style={circle,text width=0.35cm,align=center,draw=spinach,very thick,fill=spinach!50},
stoplyndon/.style={circle,text width=0.35cm,align=center,draw=black,very thick,fill=white},
nostoplyndon/.style={circle,text width=0.35cm,align=center,draw=black,very thick,fill=tomato!50},
directed/.style={postaction={decorate,decoration={markings,
mark=at position #1 with {\arrow[line width=0.5mm, black]{>}}}}},
}

\DeclarePairedDelimiterX{\set}[1]{\{}{\}}{\setargs{#1}}
\NewDocumentCommand{\setargs}{>{\SplitArgument{1}{|}}m}{\setargsaux#1}
\NewDocumentCommand{\setargsaux}{mm}
{\IfNoValueTF{#2}{#1} {#1\,\delimsize|\,\mathopen{}#2}}%{#1\:;\:#2}

\usepackage[hypertexnames=false]{hyperref}
\usepackage{bookmark}
\hypersetup{
pdftoolbar=true,
pdfmenubar=true,
pdffitwindow=false,
pdfstartview={FitH},
pdftitle={Cellularity of KLR and wKLRW algebras via crystals},
pdfauthor={Andrew Mathas and Daniel Tubbenhauer},
pdfsubject={},
pdfcreator={Andrew Mathas and Daniel Tubbenhauer},
pdfproducer={Andrew Mathas and Daniel Tubbenhauer},
pdfkeywords={},
pdfnewwindow=true,
colorlinks=true,
linkcolor=mydarkblue,
citecolor=teal,
filecolor=magenta,
urlcolor=orchid,
linkbordercolor=lava,
citebordercolor=teal,
urlbordercolor=orchid,
linktocpage=true
}

\def\makeautorefname#1#2{\csdef{#1autorefname}{#2}}
\makeautorefname{section}{Section}%
\makeautorefname{subsection}{Section}%
\makeautorefname{subsubsection}{Section}%

\newcommand\NewTheorem[2][\relax]{%
\newaliascnt{#2}{equation}%
\newtheorem{#2}[#2]{#2}%
\aliascntresetthe{#2}%
\expandafter\def\csname#2autorefname\endcsname{#2}%
\ifx#1\relax\else
\AtEndEnvironment{#2}{\null\hfill$#1$}%
\newaliascnt{#2*}{equation}%
\newtheorem{#2*}[#2*]{#2}%
\aliascntresetthe{#2*}%
\expandafter\def\csname#2*autorefname\endcsname{#2}%
\fi
}

\def\equationautorefname~#1\null{(#1)\null}

\numberwithin{equation}{subsection}

\NewTheorem{Proposition}
\NewTheorem{Theorem}
\NewTheorem[\square]{Corollary}
\NewTheorem{Lemma}
\theoremstyle{definition}
\NewTheorem{Definition}
\NewTheorem{Notation}
\NewTheorem[\Diamond]{Example}
\NewTheorem[\Diamond]{Examples}
\theoremstyle{remark}
\NewTheorem{Remark}
\NewTheorem{Assumption}

\newcounter{Case}
\counterwithin*{Case}{equation}
\newcommand\Case[1]{%
\medskip\noindent\refstepcounter{Case}%
\textbf{Case \theCase:}\textit{ #1}%
}

\begin{document}
\title[Cellularity of KLR and wKLRW algebras via crystals]{Cellularity of KLR and wKLRW algebras via crystals}
\author[A. Mathas and D. Tubbenhauer]{Andrew Mathas and Daniel Tubbenhauer}

\address{A.M.: The University of Sydney, School of Mathematics and Statistics F07, Office Carslaw 718, NSW 2006, Australia,
\href{http://www.maths.usyd.edu.au/u/mathas/}{www.maths.usyd.edu.au/u/mathas/}, \href{https://orcid.org/0000-0001-7565-5798}{ORCID 0000-0001-7565-5798}}
\email{andrew.mathas@sydney.edu.au}

\address{D.T.: The University of Sydney, School of Mathematics and Statistics F07, Office Carslaw 827, NSW 2006, Australia, \href{http://www.dtubbenhauer.com}{www.dtubbenhauer.com}, \href{https://orcid.org/0000-0001-7265-5047}{ORCID 0000-0001-7265-5047}}
\email{daniel.tubbenhauer@sydney.edu.au}

\begin{abstract}
We prove that the wKLRW algebras of finite type, and their cyclotomic quotients, are sandwich cellular algebras. The sandwich cellular bases are explicitly described using crystal graphs. As a special case, this proves that the KLR algebras of finite type are sandwich cellular. As one of our applications, we give explicit formulas for some graded decomposition numbers of the cyclotomic algebras in level one.
\end{abstract}


\addtocontents{toc}{\protect\setcounter{tocdepth}{1}}

\maketitle

\tableofcontents

\ytableausetup{centertableaux,mathmode,boxsize=0.4cm}
{\global\arrayrulewidth=0.5mm}

%%%%%%%%%%%%%%%%%%%%%%%%%%%%%%%%%%%%%%%%%

\section{Introduction}\label{S:Introduction}

%%%%%%%%%%%%%%%%%%%%%%%%%%%%%%%%%%%%%%%%%

This introduction is organized as follows.
We begin with a historical overview in \autoref{S:History}, followed by a technicality-free summary in \autoref{S:Lay}; \autoref{S:Details} then provides the full details for experts.

\subsection{A historical overview}\label{S:History}

In 1996, Lascoux, Leclerc and Thibon~\cite{LLT} made their now-famous
LLT conjecture that predicted deep connections between the
representation theory of Iwahori-Hecke algebras of type~$A$ and the
canonical bases Kac-Moody algebras. Soon afterwards,
Ariki~\cite{Ariki:can} proved their conjecture, not just for the Hecke
algebras of type~$A$ but for the Ariki-Koike algebras at complex roots
of unity.

Rather than resolving the LLT conjecture, Ariki's work opened new
questions because the LLT conjecture says that the decomposition numbers
of the Iwahori--Hecke algebras of the symmetric groups can be computed
by evaluating certain parabolic Kazhdan--Lusztig polynomials at~$1$.
Interpreting the polynomial variable as a grading parameter, this is
precisely what one would expect if the decomposition numbers were
shadows of graded multiplicities. This strongly suggests that there
should be a nontrivial grading on these algebras, but no such grading
was known.

The first significant step towards finding such a grading was made by
Chuang and Rouquier~\cite{ChuangRouq:sl2}, who introduced
$\mathfrak{sl}_2$-categorifications, which are  a family of derived equivalences
that categorify the action of the simple reflections in the affine Weyl
group on the Fock spaces underlying the LLT conjecture. As an
application, Chuang and Rouquier showed that Brou{\'e}’s abelian defect
group conjecture~\cite{Broue:perfect} holds for the symmetric groups.
Motivated in part by the LLT conjecture and Ariki’s work, they further
suggested the existence of stronger Kac--Moody categorifications.

Independently,
Khovanov-Lauda~\cite{KhLa-cat-quantum-sln-first,KhLa-cat-quantum-sln-second}
and Rouquier~\cite{Ro-2-kac-moody} 
introduced the KLR algebras, which are a family of $\mathbb{Z}$-graded
algebras attached to a Kac-Moody quiver. These algebras can be defined 
both by generators and relations and diagrammatically, with the grading oh these
algebras being determined by the underlying quiver and its Cartan
matrix. They proved that the (affine) KLR algebras category the positive
part of the corresponding quantised Kac--Moody algebra.

The cyclotomic KLR algebras are finite dimensional quotients of the KLR
algebras, indexed by dominant weights of the corresponding
Kac-Moody algebra, which are defined in exactly the same way
as the cyclotomic quotients of the affine Hecke algebras. For affine
quivers of type $A$, Lauda and Vazirani~\cite{LaVa-crystals-cat} proved
that the simple modules of these algebras categorify the crystal graph
of the associated highest weight module.  Varagnolo--Vasserot
\cite{VaVa-canonical-bases-klr} and Rouquier~\cite{Ro-quiver-hecke}
proved that the cyclotomic KLR algebras of finite type categorify the
highest weight modules and their crystal graphs. For the simply laced
quivers, Varagnolo-Vasserot~\cite{VaVa-canonical-bases-klr} and
Kang-Kashiwara~\cite{KaKa-categorification-via-klr} proved that the
simple modules of the cyclotomic KLR algebras categorify the dual
canonical bases, and the projective indecomposable modules category the canonical bases, of the corresponding highest weight modules.

In the special case of affine quivers of type~$A$, Brundan and Kleshchev~\cite{BrKl-hecke-klr} and Rouquier~\cite{Ro-quiver-hecke} proved that the cyclotomic KLR algebras are isomorphic to the cyclotomic Hecke algebras of type~$A$. Indeed, before the work of Varagnolo--Vasserot and Kang--Kashiwara mentioned above, \cite{BrKl-graded-decomposition-hecke,Ro-quiver-hecke} proved that the canonical bases categorification theorem in type A. Consequently, these results for the KLR algebras can be viewed a vast generalisation of Ariki's categorification theorem. In particular, the grading on the KLR algebras fully explain the appearance of $q$ in the LLT conjecture.

In 2012, Kleshchev, Loubert and Miemietz
\cite{KlLoMi-KLR-affine-cellular-typea} conjectured that the KLR
algebras are graded (affine) cellular algebras. Some progress was made
towards this conjecture for the infinite dimensional KLR algebras
\cite{KlLo-klr-affine-cellular-finite-type}, but the conjecture for the
cyclotomic quotients proved to be much harder, except in type
$\aonetypea$~\cite{HuMa-klr-basis}.

In spite of the wealth of strong categorification theorems
describing the representation theory of the cyclotomic KLR algebras,
surprisingly little is known about these algebras, except in affine type~$A$
where the Brundan-Kleshchev isomorphism to the ungraded Hecke algebras
made it possible to classify the blocks of these
algebras~\cite{BrKl-hecke-klr,LyMa} and construct cellular
bases for them~\cite{HuMa-klr-basis}.

In order to understand categorifications of tensor products, Webster
\cite{We-knot-invariants,We-weighted-klr,We-rouquier-dia-algebra}
introduced the \emph{weighted KLRW algebras}, which include the KLR
algebras as a special case. These algebras admit a diagrammatic
presentation and Bowman~\cite{Bo-many-cellular-structures} was able to
use these algebras to construct new graded cellular bases for the
cyclotomic KLR and weighted KLRW algebras of affine type~$A$. Extending these
ideas, the authors proved that the weighted KLRW algebras and their
cyclotomic quotients are graded (affine) cellular algebras in types
$\atypec$, $\atypeb[\N]$, $\atwotypea$ and $\atwotyped$ in
\cite{MaTu-klrw-algebras} and \cite{MaTu-klrw-algebras-bad}.

In this paper we apply the KLRW machinery to construct, in principle,
explicit bases for the cyclotomic KLR and weighted KLRW algebras for
quivers of finite type. To construct bases for these algebras we use the
general approach that we developed in
\cite{MaTu-klrw-algebras,MaTu-klrw-algebras-bad}. In these papers we
described these bases using explicit tableau combinatorics, which
partially existed in the literature for describing Fock spaces for the
corresponding highest weight modules. Similar combinatorics exists for
the KLR algebras of simply laced types, but in this paper we take a
different approach and use the crystal graphs to construct bases for
these algebras. By using the crystal graphs, we have a uniform
combinatorics for describing bases of the cyclotomic KLR and weighted
KLRW algebras for all quivers of finite type.  This is very natural
because Lauda--Vazirani \cite{LaVa-crystals-cat} and Varagnolo--Vasserot
\cite{VaVa-canonical-bases-klr} have shown that the cyclotomic KLR
algebras of finite type categorify the highest weight modules and their
crystal graphs. In this way we are able to directly link the crystal
theory and the weighted KLRW algebras. To do this we consider a special
family of weighted KLRW algebras, which are closely related to the KLR
algebras.

\subsection{A layperson’s summary}\label{S:Lay}

One way to organize a complicated algebra is to give it a good set of building blocks and a recipe for how the blocks multiply.
For matrix algebras, Gaussian elimination plays this role; for symmetric groups, Young diagrams and tableaux do.
This paper shows that a much larger and more modern family of diagrammatic algebras, those arising from the Khovanov--Lauda--Rouquier (KLR) \cite{KhLa-cat-quantum-sln-first,KhLa-cat-quantum-sln-second,Ro-2-kac-moody} framework and their (weighted) variants \cite{We-knot-invariants,We-weighted-klr,We-rouquier-dia-algebra}, also admit such an organizing principle in finite type.
In plain terms: we prove that in the finite-type setting there are bases and a compatible structure that let you ``triangularize'' the algebra and read off its representations in a controlled and combinatorial way. (Finite type here refers to the case controlled by classical Dynkin diagrams, the regime where the underlying combinatorics is sharpest and well-understood.)

Why do people care about KLR algebras? Well, KLR algebras sit at a busy crossroads. We summarized part of their history in \autoref{S:History}, with the focus on the LLT conjecture part. Additionally, they categorify quantum groups, turning the combinatorics of canonical/crystal bases into concrete, computable objects; they provide a diagrammatic home for link-homology theories in low-dimensional topology (via braid and tangle actions); and they connect to geometric representation theory through quiver varieties and sheaf-theoretic models.
On the algebra/combinatorics side, they unify and extend familiar diagram algebras (Hecke, Temperley–Lieb, Brauer, Schur-like quotients) and inform questions about decomposition numbers and modular phenomena.
In short: understanding KLR algebras pays dividends across topology, geometry, and algebra, and having a cellular/triangular organization turns these conceptual links into practical tools.
Since their inception in the mid-2000s there has been enormous progress, yet basic structural features, most notably ``nice'' bases that behave well under multiplication, have remained elusive.
Providing such bases in finite type is the starting point and main contribution of this paper.

Let us now explain the results of this paper in more detail while keeping the details under the rug. First, we use a slightly enhanced algebra that we call wKLRW algebra. The difference between these algebras, for our purpose, is mostly cosmetic: the combinatorics is easier to see in the wKLRW picture.
Readers unfamiliar with KLR-style algebras can think of them as algebras presented by pictures: one draws colored labeled strings in a strip, allowed to cross or carry dots, and composes pictures by stacking. An example is:
\begin{gather*}
\begin{tikzpicture}[anchorbase,smallnodes,rounded corners]
\draw[ghost](1.5,1)--++(1,1)node[above,yshift=-1pt]{$2$};
\draw[ghost](2.5,1)--++(-1,1)node[above,yshift=-1pt]{$0$};
\draw[ghost](6,1)node[below]{$\phantom{i}$}--++(-1,1)node[above,yshift=-1pt]{$5$};
\draw[solid](1,1)node[below]{$1$}--++(1,1)node[above,yshift=-1pt]{$\phantom{i}$};
\draw[solid](2,1)node[below]{$1$}--++(-1,1);
\draw[solid](3.5,1)node[below]{$3$}--++(2,1)node[above,yshift=-1pt]{$\phantom{i}$};
\draw[solid](4.5,1)node[below]{$4$}--++(-1,1)node[above,yshift=-1pt]{$\phantom{i}$};
\draw[solid](5.5,1)node[below]{$0$}--++(-1,1)node[above,yshift=-1pt]{$\phantom{i}$};
\draw[redstring](2.25,1)node[below]{$1$}--++(0,1);
\draw[redstring](4,1)node[below]{$2$}--++(0,1);
\end{tikzpicture}
\circ
\begin{tikzpicture}[anchorbase,smallnodes,rounded corners]
\draw[ghost](1.5,1)node[above,yshift=-1pt]{$2$}--++(1,-1);
\draw[ghost](2.5,1)node[above,yshift=-1pt]{$0$}--++(-1,-1);
\draw[ghost](6,1)node[above,yshift=-1pt]{$5$}--++(-1,-1)node[below]{$\phantom{i}$};
\draw[solid](1,1)node[above,yshift=-1pt]{$\phantom{i}$}--++(1,-1)node[below]{$2$};
\draw[solid](2,1)--++(-1,-1)node[below]{$1$};
\draw[solid](3.5,1)node[above,yshift=-1pt]{$\phantom{i}$}--++(2,-1)node[below]{$3$};
\draw[solid](4.5,1)node[above,yshift=-1pt]{$\phantom{i}$}--++(-1,-1)node[below]{$4$};
\draw[solid](5.5,1)node[above,yshift=-1pt]{$\phantom{i}$}--++(-1,-1)node[below]{$0$};
\draw[redstring](2.25,1)--++(0,-1)node[below]{$1$};
\draw[redstring](4,1)--++(0,-1)node[below]{$2$};
\end{tikzpicture}
\\
=
\begin{tikzpicture}[anchorbase,smallnodes,rounded corners]
\draw[ghost](1.5,1)--++(1,-1);
\draw[ghost](2.5,1)--++(-1,-1);
\draw[ghost](6,1)--++(-1,-1)node[below]{$\phantom{i}$};
\draw[solid](1,1)--++(1,-1)node[below]{$2$};
\draw[solid](2,1)--++(-1,-1)node[below]{$1$};
\draw[solid](3.5,1)--++(2,-1)node[below]{$3$};
\draw[solid](4.5,1)--++(-1,-1)node[below]{$4$};
\draw[solid](5.5,1)--++(-1,-1)node[below]{$0$};
\draw[redstring](2.25,1)--++(0,-1)node[below]{$1$};
\draw[redstring](4,1)--++(0,-1)node[below]{$2$};
\draw[ghost](1.5,1)--++(1,1)node[above,yshift=-1pt]{$2$};
\draw[ghost](2.5,1)--++(-1,1)node[above,yshift=-1pt]{$0$};
\draw[ghost](6,1)--++(-1,1)node[above,yshift=-1pt]{$5$};
\draw[solid](1,1)--++(1,1);
\draw[solid](2,1)--++(-1,1);
\draw[solid](3.5,1)--++(2,1);
\draw[solid](4.5,1)--++(-1,1);
\draw[solid](5.5,1)--++(-1,1);
\draw[redstring](2.25,1)--++(0,1);
\draw[redstring](4,1)--++(0,1);
\end{tikzpicture}
\end{gather*}
Familiar older cousins are Temperley--Lieb and Brauer algebras, where noncrossing matchings or pairings govern multiplication.
The power of a diagram algebra comes from having a basis of diagrams for which multiplying and simplifying is predictable.
Our contribution is to show that, for finite type, the relevant KLR algebras admit exactly this kind of predictable organization, a sandwich cellular structure as e.g. in \cite{Tu-sandwich-cellular}, and that the indexing data are supplied by the crystal graphs from Lie theory.

Without technicalities this means the following.
A cellular algebra (in the sense of Graham--Lehrer \cite{GrLe-cellular} and generalizations) is one that comes with:
(i) a poset indexing of basis pieces (``cells''),
(ii) standard modules attached to these indices, and
(iii) a multiplication rule that is upper-triangular with respect to the poset.
Practically, this yields a clean bookkeeping device for simple modules and decomposition numbers, and it behaves well under standard operations (duality, quotients, specializations).
We work in a slightly more flexible setting, so-called sandwich cellularity, which isolates a straightforward piece in the middle of a diagram but keeps the same spirit: multiplication respects a transparent order, so representation-theoretic information can be read off in stages.

Before giving more details for the experts, here is what we prove informally.
\begin{itemize}
  \item \textit{Cellularity from crystals.}
  wKLRW algebras of finite type admit sandwich cellular structures with bases built from crystal paths. (More on crystals momentarily.)
  \item \textit{Same representation theory, different models.}
  wKLRW and KLR algebras in finite type are Morita equivalent; their module categories are the same up to equivalence.
  \item \textit{Structural corollaries.}
  We record criteria for when these algebras are nonzero, semisimple, quasi-hereditary, or indecomposable, and we point out situations where cellularity fails, clarifying the boundary of our results.
\end{itemize}

Although we were surprised by the deep connection between crystals and sandwich cellularity, there are conceptual reasons why crystals appear:
Crystals are the combinatorial ``skeletons'' of representations of semisimple Lie algebras: directed graphs whose vertices label basis elements and whose arrows record how simple raising/lowering operations move you around.
In finite type, it is a guiding principle that the representation theory of KLR algebras mirrors these crystal graphs \cite{LaVa-crystals-cat,VaVa-canonical-bases-klr}.
Our results make this principle concrete for bases and cellularity:
basis elements are built by following paths in the crystal and drawing the corresponding diagrams;
triangularity reflects a certain partial order coming from the crystal;
and the standard modules one expects in a cellular setting line up with crystal data in a hands-on way.

\textit{Takeaway.}
Diagrams, crystals, and sandwich cellularity fit together in finite type to provide a concrete, diagrammatic blueprint for the representation theory of KLR/wKLRW algebras.
The results put the combinatorics of crystals to work as an organizing spine for bases, modules, and decomposition numbers in a broad class of modern diagrammatic algebras.

\subsection{More details}\label{S:Details}

Using this framework, our main results are:

\begin{enumerate}[label=(\Alph*)]

\item \autoref{T:ConstructionGeneralMain} shows that the weighted KLRW algebras, and their cyclotomic quotients, are graded sandwich cellular algebras. This implies the corresponding results for the KLR algebras.

\item In \autoref{P:ConsequencesMorita} we prove that the cyclotomic wKLRW algebras and the cyclotomic KLR algebras of finite type are graded Morita equivalent.

\item In \autoref{P:ConsequencesSimples} we prove that simple modules of these algebras are indexed by the vertices of the corresponding crystal graph.

\item In level one, under certain technical assumptions, we give a basis for the projective indecomposable modules of the cyclotomic wKLRW algebras and compute the graded decomposition numbers of the cell modules in \autoref{SS:ConsequencesProj}. Except in type $\typef$,
\autoref{T:ConsequencesCanonicalBasis} shows that the projective indecomposable modules categorify the canonical basis of the highest weight module.
In this case, we construct the simple modules explicitly and determine their dimensions in \autoref{P:ConsequencesDimensions}.

\item In \autoref{S:Consequences} we list several more consequences of sandwich cellularity, including results on when these algebras are nonzero, semisimple, quasi-hereditary and indecomposable.

\item We also talk about the limitations of our results, see \autoref{P:ConsequencesNotCellular}.
\end{enumerate}
All of these results are independent of the characteristic.

We briefly outline the contents of this paper. \autoref{S:Sandwich} recalls the basic machinery of sandwich cellular algebras. \autoref{S:Recollection} contains the definition of wKLRW algebras and gives the basic rewriting results from \cite{MaTu-klrw-algebras} that we use to pull strings and jump dots to the right that underpins our basis theorems and hence all of the results in this paper. To explain (almost all) of the combinatorics used in this paper, \autoref{S:MainExample} is an extended running example of our construction of cellular bases for the wKLRW algebras. \autoref{S:FiniteTypes} starts with some reminders from the theory of crystal graphs, which is then extended and applied to the wKLRW algebras in \autoref{S:Construction}. It is striking that our crystal basis is built from idempotent diagrams, which correspond to paths in the crystal, together with permutations of the labels along these paths. \autoref{S:Proofs} proves our sandwich cellular basis theorems and \autoref{S:Consequences} contains the applications of these results to the representation theory of the KLR and wKLRW algebras.

\begin{Remark}\label{R:IntroductionSageMath}
Many of the results in this paper were inspired by extensive calculations with crystals using SageMath~\cite{sage} and some results for exceptional type rely on these calculations. The code that we used can be found at \cite{MaTu-sagemath-finite-type-klrw}, which also includes some additional features, such as code for drawing some of the wKLRW diagrams used in this paper.
\end{Remark}

%%%%%%%%%%%%%%%%
% Acknowledgments
%%%%%%%%%%%%%%%%

\noindent\textbf{Acknowledgments.}
We thank Anton Evseev and Sasha Kleshchev for sharing a sketch of a proof
that the affine type D cyclotomic KLR algebras are not cellular;
Robert Muth for helpful discussions that gave birth to \autoref{R:ConsequencesGraphsNotCellular}; and Travis Scrimshaw for invaluable discussions about SageMath and crystals.
We also thank Chris Bowman, Joe Chuang, Dinushi Munasinghe, Tao Qin, the referee, Liron Speyer
and Ben Webster for helpful discussions about
(weighted) KLR(W) algebras and cellularity.
All of this is very much appreciated.

Both authors were supported, in part, by the Australian Research Council.
D.T. thanks additionally their depression, burnout and usage of outdated ideas for long-term support.

%%%%%%%%%%%%%%%%%%%%%%%%%%%%%%%%%%%%%%%%%

\section{Sandwich cellularity}\label{S:Sandwich}

%%%%%%%%%%%%%%%%%%%%%%%%%%%%%%%%%%%%%%%%%

The main results of this paper show that the wKLRW algebras of
finite type are sandwich cellular algebra. This section defines sandwich cellular algebras and sets our notation.

%%%%%%%%%%%%%%%%%%%%%%%%%%%%%%%%%%%%%%%%%

\subsection{Sandwich cellular algebras}\label{SS:SandwichOne}

%%%%%%%%%%%%%%%%%%%%%%%%%%%%%%%%%%%%%%%%%

\begin{Notation}\label{N:SandwichReadingDiagrams}
\leavevmode
\begin{enumerate}

\item A grading will always mean a $\Z$ grading.

\item Throughout we will work with diagram algebras and left actions and left modules, which are given by acting from the top:
\begin{gather*}
E\circ D
=
\begin{tikzpicture}[anchorbase,smallnodes,rounded corners]
\node[rectangle,draw,minimum width=0.5cm,minimum height=0.5cm,ultra thick] at(0,0){\raisebox{-0.05cm}{$D$}};
\node[rectangle,draw,minimum width=0.5cm,minimum height=0.5cm,ultra thick] at(0,0.5){\raisebox{-0.05cm}{$E$}};
\end{tikzpicture}
.
\end{gather*}
Modules will always be left modules.

\item Let $\ring$ be an integral domain, such as $\ring=\Z$ or a field, which is our ground ring throughout.
(\autoref{SS:SandwichOne} works in more generality but this is not important for this paper.)

\item In this paper all $\ring$-modules are $\ring$-free, ranks are always considered over $\ring$
and we will use $\vpar$
to indicate the grading. More precisely, if $V$ is a graded $\ring$-module, then
the \emph{graded rank}
$\grdim(V)=\sum_{i\in\Z}\dim V_{i}\vpar^{i}\in\N[\vpar,\vpar^{-1}]$ of~$V$
is defined by the property that there is a homogeneous degree $0$ isomorphism $V\cong\bigoplus_{i\in\Z}R\langle i\rangle^{\oplus a_{i}}$, where $R\langle i\rangle$ is concentrated in degree $i$. Specializing
$\vpar\to 1$ recovers the usual (ungraded) rank $\grdim[1](V)$ of $V$.

\end{enumerate}
\end{Notation}

The following is copied, almost {\ver}, from \cite[Definition 2A.2]{MaTu-klrw-algebras-bad} and this, in turn, is a mild reformulation of \cite[Definition 2.2]{TuVa-handlebody} and
\cite[Definition 2A.3]{Tu-sandwich-cellular}:

\begin{Definition}\label{D:SandwichCellularAlgebra}
Let $A$ be a locally unital graded $\ring$-algebra.
A \emph{graded sandwich cell datum} for $A$ is a tuple $\big(\Pcal,(S,T),(\sand[],\sandbasis[]),C,\deg\big)$, where:
\begin{itemize}

\item $\Pcal=(\Pcal,<)$ is a poset (the \emph{middle poset} with \emph{sandwich order $<=<_{\Pcal}$}),

\item $S=\coprod_{\lambda\in\Pcal}S(\lambda)$ and $T=\coprod_{\lambda\in\Pcal}T(\lambda)$ are collections of finite sets
(the \emph{bottom/top sets}),

\item $\sand[]=\bigoplus_{\lambda\in\Pcal}\sand[\lambda]$ is a direct sum of
graded (not necessarily unital) algebras $\sand[\lambda]$
(the \emph{sandwiched algebras}) and $\sandbasis[]=\coprod{\lambda\in\Pcal}\sandbasis$ (the \emph{sandwiched basis}) such that $\sandbasis[\lambda]$ is a homogeneous
basis of $\sand[\lambda]$ (we write $\deg$ for the degree function on $\sand[\lambda]$),

\item $C\map{\coprod_{\lambda\in\Pcal}S(\lambda)\times\sandbasis[\lambda]\times T(\lambda)}{A};(S,b,T)\mapsto C_{ST}^{b}$ is an injective map
(the \emph{basis}),

\item $\deg\map{\coprod_{\lambda\in\Pcal}S(\lambda)\cup T(\lambda)}{\Z}$ is a function
(the \emph{degree}),

\end{itemize}
such that:
\begin{enumerate}[label=\upshape(AC${}_{\arabic*}$\upshape)]

\item For $S\in S(\lambda)$, $T\in T(\lambda)$ and $b\in\sandbasis[\lambda]$, $\lambda\in\Pcal$ the element $C_{ST}^{b}$
is homogeneous of degree $\deg(S)+\deg(b)+\deg(T)$.

\item The set $\set{C^{b}_{ST}|\lambda\in\Pcal,S\in S(\lambda),T\in T(\lambda),b\in\sandbasis[\lambda]}$
is an $\ring$-basis of $A$.

\item For all $x\in A$ there exist scalars $r_{SU}=r_{SU}(x)\in R$ that do not depend on $T$ or on $b$, such that
\begin{gather*}
xC_{ST}^{b}\equiv
\sum_{U\in T(\lambda),c\in\sandbasis[\lambda]}r_{SU}C_{UT}^{c}\pmod{A^{>\lambda}},
\end{gather*}
where $A^{>\lambda}=\bigcup_{\mu>\lambda}A^{\geq\mu}$ and $A^{\geq\mu}$ is the $\ring$-submodule of $A$ spanned by
$\set{C^{c}_{UV}|\nu\in\Pcal,\nu\geq\mu,U\in S(\nu),V\in T(\nu),c\in\sandbasis[\nu]}$.

\item Let $A(\lambda)=A^{\geq\lambda}/A^{>\lambda}$. Then there exist free graded right and left
$\sand[\lambda]$-modules $\Delta(\lambda)$ and $\nabla(\lambda)$, respectively, such that $A(\lambda)\cong\Delta(\lambda)\otimes_{\sand[\lambda]}\nabla(\lambda)$ for all $\lambda\in\Pcal$.

\end{enumerate}
The algebra $A$ is a graded \emph{sandwich cellular algebra} if it
has a graded sandwich cell datum.

Assume that $S(\lambda)=T(\lambda)$ for all $\lambda\in\Pcal$, and that there is a antiinvolution $(\placeholder)^{\star}\colon A\to A$
and a homogeneous bijection $(\placeholder)^{\star}\colon\sandbasis[\lambda]\to\sandbasis[\lambda]$ of order two such that:
\begin{enumerate}[label=\upshape(AC${}_{5}$\upshape)]

\item We have $(C^{b}_{ST})^{\star}\equiv C_{TS}^{b^{\star}}\pmod{A^{>\lambda}}$, for all $S,T\in T(\lambda)$ and $b\in\sandbasis[\lambda]$.

\end{enumerate}
In this case, we write $\big(\Pcal,S=T,(\sand[],\sandbasis[]), C,\deg,(\placeholder)^{\star}\big)$ and call this datum \emph{involutive}.
\end{Definition}

We omit the brackets around $\sand[],\sandbasis[]$ for readability. We do not distinguish
$(\placeholder)^{\star}\colon A\to A$
and the map $(\placeholder)^{\star}\colon\sandbasis[\lambda]
\to\sandbasis[\lambda]$.

\begin{Remark}\label{R:SandwichCellularAlgebra}
Sandwich cellularity is a common generalization of cellular
algebras as in \cite{GrLe-cellular} and affine cellular
algebras as in \cite{KoXi-affine-cellular}.
(Strictly speaking we use a weakened condition on the antiinvolution
as in \cite[Definition 2.3]{GoGr-cellularity-jones-basic} and we always
use this definition when we refer to cellular or affine algebras.)
The notion of sandwiched cellularity originates
in work of Brown \cite{Br-gen-matrix-algebras} on the Brauer
algebra, see also \cite{Tu-sandwich-cellular} for more historical comments.

For reader familiar with diagram algebras (as in the main body of this paper),
the diagrammatic incarnation of sandwich cellularity is:
\begin{gather*}
\text{cellular: }
C^{1}_{ST}
\leftrightsquigarrow
\scalebox{0.8}{$\begin{tikzpicture}[anchorbase,scale=1]
\draw[mor] (0,-0.5) to (0.25,0) to (0.75,0) to (1,-0.5) to (0,-0.5);
\node at (0.5,-0.25){$T$};
\draw[mor] (0,0.5) to (0.25,0) to (0.75,0) to (1,0.5) to (0,0.5);
\node at (0.5,0.25){$S$};
\draw[thick,->] (1.5,0.4)node[right]{top part} to (1,0.3);
\draw[thick,->] (1.5,-0.4)node[right]{bottom part} to (1,-0.3);
\draw[thick,->] (1.5,0)node[right]{$\sand\cong R$} to (1,0);
\end{tikzpicture}$}
,\quad
\text{affine cellular: }
C^{b}_{ST}
\leftrightsquigarrow
\scalebox{0.8}{$\begin{tikzpicture}[anchorbase,scale=1]
\draw[mor] (0,-0.5) to (0.25,0) to (0.75,0) to (1,-0.5) to (0,-0.5);
\node at (0.5,-0.25){$T$};
\draw[mor] (0,1) to (0.25,0.5) to (0.75,0.5) to (1,1) to (0,1);
\node at (0.5,0.75){$S$};
\draw[mor] (0.25,0) to (0.25,0.5) to (0.75,0.5) to (0.75,0) to (0.25,0);
\node at (0.5,0.25){$b$};
\draw[thick,->] (1.5,0.75)node[right]{top part} to (1,0.75);
\draw[thick,->] (1.5,0.25)node[right]{commutative $\sand$} to (1,0.25);
\draw[thick,->] (1.5,-0.25)node[right]{bottom part} to (1,-0.25);
\end{tikzpicture}$}
,
\end{gather*}
\begin{gather*}
\text{sandwich cellular: }
C^{b}_{ST}
\leftrightsquigarrow
\scalebox{0.8}{$\begin{tikzpicture}[anchorbase,scale=1]
\draw[mor] (0,-0.5) to (0.25,0) to (0.75,0) to (1,-0.5) to (0,-0.5);
\node at (0.5,-0.25){$T$};
\draw[mor] (0,1) to (0.25,0.5) to (0.75,0.5) to (1,1) to (0,1);
\node at (0.5,0.75){$S$};
\draw[mor] (0.25,0) to (0.25,0.5) to (0.75,0.5) to (0.75,0) to (0.25,0);
\node at (0.5,0.25){$b$};
\draw[thick,->] (1.5,0.75)node[right]{top part} to (1,0.75);
\draw[thick,->] (1.5,0.25)node[right]{general $\sand$} to (1,0.25);
\draw[thick,->] (1.5,-0.25)node[right]{bottom part} to (1,-0.25);
\end{tikzpicture}$}
.
\end{gather*}
In the first two pictures, the (affine) cellular algebras have involutive cell data with $\sand[\lambda]\cong R$, for $\lambda\in\Pcal$.
\end{Remark}

\begin{Notation}\label{N:SandwichCellularAlgebra}
We usually abuse language and simply say that an algebra \emph{is} sandwich cellular. We also use \emph{affine sandwich cellular} to highlight that the
sandwiched algebras are polynomial rings, even though
affine sandwich cellular algebras are just
special cases of sandwich cellular algebras.
\end{Notation}

%%%%%%%%%%%%%%%%%%%%%%%%%%%%%%%%%%%%%%%%%

\subsection{Some statements, including \texorpdfstring{$H$}{H}-reduction}\label{SS:SandwichTwo}

%%%%%%%%%%%%%%%%%%%%%%%%%%%%%%%%%%%%%%%%%

Sandwich cellularity essentially generalizes
the machinery of cellular algebras to a larger class of algebras, see {\eg} \cite[Section 2]{TuVa-handlebody}
and \cite{Tu-sandwich-cellular} for details. Moreover, incorporating a grading generalizes \cite[Section 2]{HuMa-klr-basis} to this larger class of algebras.

Instead of filling the next few pages with definitions, all of which are not very different
from the cellular case and can be found in \cite{Tu-sandwich-cellular}, we stress the statements that we will actually use.

\begin{Theorem}\label{T:SandwichMain}
Let $A$ be a (graded) sandwich cellular algebra.
\begin{enumerate}

\item Let $e\in A$ be a (homogeneous) idempotent.
Then $eAe$ is a (graded) sandwich cellular algebra.

\end{enumerate}

Assume now that $\ring$ is a field.

\begin{enumerate}[resume]
\item All (graded) simple $A$-modules are uniquely associated to a $\lambda\in\Pcal$, called their \emph{apex}.
Denote the set of apexes by
$\Pcal^{\neq 0}\subset\Pcal$ and
assume that $\sand[\lambda]$ is unital Artinian or unital commutative.
Then there exist 1:1-correspondences
\begin{gather*}
\set[\big]{\text{(graded) simple $A$-modules with apex $\lambda$}}/{\cong}
\;\xleftrightarrow{1:1}
\set[\big]{\text{(graded) simple $\sand[\lambda]$-modules}}/{\cong}
.
\end{gather*}

\item Assume that $\dim_{R}(A)<\infty$.
If $A$ is involutive, $\sand[\lambda]\cong R$ for $\lambda\in\Pcal$
and all cell modules are simple, then $A$ is semisimple.
Conversely, if at least one $\sand[\lambda]$ is not semisimple,
then $A$ is not semisimple.

\end{enumerate}
\end{Theorem}

\begin{proof}
\textit{(a).} This has not appeared in the literature,
but is an easy adaption of the standard arguments in
this setting, see
{\eg} \cite[Proposition 4.3]{KoXi-cellular}
or \cite[Proposition 2.8.(a)]{EhTu-relcell}. Details are omitted.

\textit{(b)+(c).}
See \cite[Theorem 2.16 and Proposition 2.9]{TuVa-handlebody}
and \cite[Proposition 2B.23]{Tu-sandwich-cellular}.
\end{proof}

\begin{Remark}\label{R:SandwichMain}
\leavevmode

\begin{enumerate}

\item Part (b) of \autoref{T:SandwichMain} is the
celebrated \emph{Clifford--Munn--Ponizovski\u{\i}} theorem
or \emph{$H$-reduction}. The assumption
that $\sand[\lambda]$ is unital Artinian or unital commutative
can be relaxed, see \cite[Section 2]{TuVa-handlebody},
but these conditions are automatically satisfied in this paper.

\item \autoref{T:SandwichMain}(b) can be strengthened at
the cost of requiring more notation. The above formulation
is sufficient for this paper and we refer the reader to
\cite[Proposition 2B.23]{Tu-sandwich-cellular} for a
stronger statement.

\end{enumerate}

\end{Remark}

%%%%%%%%%%%%%%%%%%%%%%%%%%%%%%%%%%%%%%%%%

\section{wKLRW algebras}\label{S:Recollection}

%%%%%%%%%%%%%%%%%%%%%%%%%%%%%%%%%%%%%%%%%

We now recall the main algebras considered in this paper, following \cite{MaTu-klrw-algebras}. These algebras, the \emph{wKLRW algebras} or \emph{wKLRW algebras} for short, were originally introduced in \cite{We-weighted-klr},
\cite{We-rouquier-dia-algebra}. The wKLRW algebras are diagram algebras, but not in the traditional sense. At first sight, their definition looks more complicated than it actually is. We have tried to give a self-contained account but we encourage the reader to consult these three papers for more details.

%%%%%%%%%%%%%%%%%%%%%%%%%%%%%%%%%%%%%%%%%

\subsection{Kac--Moody datum}\label{SS:RecollectionKacMoody}

%%%%%%%%%%%%%%%%%%%%%%%%%%%%%%%%%%%%%%%%%

In this paper, a \emph{symmetrizable Kac--Moody datum}
is a set $\vertices=\set{1,\dots,e}$
together with a bilinear form $(\placeholder,\placeholder)$ on $\Z[I]$
such that $(i,i)=2$ and $(i,j)\in\set{0,-1,-2,-3}$ for all $i\neq j$. Let $a_{ij}=(i,j)$. Then $\boldsymbol{A}=(a_{ij})_{i,j=1}^{e}$ is the \emph{Cartan matrix} of a symmetrizable Kac--Moody algebra.
The Cartan matrix has an associated symmetrizer $\boldsymbol{d}=(d_{1},\dots,d_{e})\in\N^{e}$, which is the minimal
tuple such that $(d_{i}a_{ij})_{i,j=1}^{e}$ is symmetric and positive
definite. We also associate a quiver $\quiver=(\vertices,\edges)$, with vertices
$\vertices=\set{1,\dots,e}$
and edges $\edges$. We fix an orientation on the simply laced edges, which will not affect the wKLRW algebras up to isomorphism (see
\cite[Proposition 3A.1]{MaTu-klrw-algebras} for a precise statement).
Following \cite[Sections 1.1, 1.2 and 1.3]{Ka-infdim-lie},
we have simple roots
$\sroots=\set{\sroot[1],\dots,\sroot[e]}$, coroots $\set{\scoroot[1],\dots,\scoroot[e]}$, fundamental weights
$\fweights=\set{\fweight[1],\dots,\fweight[d]}$, the \emph{positive root lattice} $Q^{+}=\bigoplus_{i\in\vertices}\N\alpha_{i}$, the \emph{dominant weight lattice } $P^{+}=\bigoplus_{i\in\vertices}\N\fweight[i]$, {\etc} that we
will use throughout.

\begin{Notation}\label{N:RecollectionArrows}
We write $i\rightsquigarrow j$ if there is an edge from $i$ to $j$.
If we need to specify the multiplicity, then we will write
$i\rightarrow j$, $i\Rightarrow j$ or $i\Rrightarrow j$ for
$d_{i}a_{ij}=-1$, $d_{i}a_{ij}=-2$ and $d_{i}a_{ij}=-3$, respectively.
In particular, if $d_{i}a_{ij}=-1$, then we a fix an orientation
$i\rightarrow j$ and, by convention, there is no edge from $j$ to $i$.
\end{Notation}

We will eventually focus on quivers of finite type,
with the conventions specified in \autoref{SS:FiniteTypesConventions} below. For now we allow arbitrary symmetrizable Kac--Moody data because we will need this for \autoref{S:Consequences}.

%%%%%%%%%%%%%%%%%%%%%%%%%%%%%%%%%%%%%%%%%

\subsection{wKLRW algebras}\label{SS:RecollectionDiagrams}

%%%%%%%%%%%%%%%%%%%%%%%%%%%%%%%%%%%%%%%%%

This section defines the \emph{wKLRW algebras} $\WA[n](X)$,
and their finite dimensional
\emph{cyclotomic} quotients $\WAc[n](X)$.

The definition of the wKLRW algebras involves two crucial ingredients:

\begin{enumerate}

\item The wKLRW algebras are graded diagram algebras.
That is, they are defined using isotopy equivalence classes of \emph{string diagrams}
in $\R^{2}$. These diagrams have a degree, which is determined by the Cartan datum, and they are subject to certain diagrammatic relations, which we recall in
\autoref{SS:RecollectionRelations} and \autoref{SS:RecollectionDegree} below. The strings in wKLRW diagrams are of three types:
\emph{solid}, \emph{ghost} and \emph{red}.
In \autoref{SS:ConstructionAffine}, we introduce \emph{affine red} strings, but these strings are used only as a notational convenience and are not part of the wKLRW diagrams.
These strings are illustrated as:
\begin{gather*}
\text{solid string}:
\begin{tikzpicture}[anchorbase,smallnodes,rounded corners]
\draw[solid] (0,0)node[below]{$i$} to (0,0.5)node[above,yshift=-1pt]{$\phantom{i}$};
\end{tikzpicture}
,\quad
\text{ghost string}:
\begin{tikzpicture}[anchorbase,smallnodes,rounded corners]
\draw[ghost] (0,0)node[below]{$\phantom{i}$} to (0,0.5)node[above,yshift=-1pt]{$i$};
\end{tikzpicture}
,\quad
\text{red string}:
\begin{tikzpicture}[anchorbase,smallnodes,rounded corners]
\draw[redstring] (0,0)node[below]{$i$} to (0,0.5)node[above,yshift=-1pt]{$\phantom{i}$};
\end{tikzpicture}
,\quad
\text{affine red string}:
\begin{tikzpicture}[anchorbase,smallnodes,rounded corners]
\draw[affine] (0,0)node[below]{$i$} to (0,0.5)node[above,yshift=-1pt]{$\phantom{i}$};
\end{tikzpicture}
.
\end{gather*}
These strings are labeled by \emph{residues} $i\in\vertices$,
which are written under
the solid and (affine) red strings and over the ghost strings.
Solid and ghost strings are
also allowed to carry finitely many
\emph{dots}. We list all possible local configurations in \autoref{Eq:RecollectionDegree}, where we also define the degrees of these configurations.

\item The various strings, their labels and positions are
determined as follows.

\begin{enumerate}

\item We fix $n\in\N$, the number of solid strings, the \emph{level} $\ell\in\Z_{>0}$, which is the number of red strings, and $\brho\in\vertices^{\ell}$. The $n$ solid strings
are labeled by resides in $\vertices$, and the $\ell$ red strings are labeled by the fixed choices $\brho_{1},\dots,\brho_{\ell}$.

For example, taking $n=3$, $\ell=2$ and $\brho=(1,0)$:
\begin{gather}\label{Eq:RecollectionExampleStrings}
\begin{tikzpicture}[anchorbase]
\node[dynkin=1] (0) at (0,0){};
\node[dynkin=2] (1) at (1,0){};
\node[dynkin=3] (2) at (2,0){};
\draw[directed=0.5](1)--(0)node[above,xshift=0.5cm,yshift=0.1cm]{$0.75$};
\draw[directed=0.5](1)--(2)node[above,xshift=-0.5cm,yshift=0.1cm]{$1.25$};
\end{tikzpicture}
\quad\rightsquigarrow\quad
\begin{tikzpicture}[anchorbase,smallnodes,rounded corners]
\draw[ghost](2.75,0)node[below]{$\phantom{i}$}--++(0,1)node[above,yshift=-0.05cm]{$2$};
\draw[ghost](3.25,0)node[below]{$\phantom{i}$}--++(0,1)node[above,yshift=-0.05cm]{$2$};
\draw[solid](0,0)node[below]{$1$}--++(0,1)node[above,yshift=-1pt]{$\phantom{i}$};
\draw[solid](1.2,0)node[below]{$3$}--++(0,1)node[above,yshift=-1pt]{$\phantom{i}$};
\draw[solid](2,0)node[below]{$2$}--++(0,1)node[above,yshift=-1pt]{$\phantom{i}$};
\draw[redstring] (-0.5,0)node[below]{$1$} to (-0.5,1)node[above,yshift=-1pt]{$\phantom{i}$};
\draw[redstring] (1,0)node[below]{$0$} to (1,1)node[above,yshift=-1pt]{$\phantom{i}$};
\end{tikzpicture}
.
\end{gather}
On the left-hand side, we have drawn an example quiver $\Gamma$. For the wKLRW diagram on the right-hand side we could use any quiver with no edges leaving $1$ and $3$, and two edges leaving $2$.

\item The number of ghost strings and their labels is determined
as follows. Every edge $\epsilon\colon i\rightsquigarrow j$ is given a positive real number $\sigma_{\epsilon}>0$ called its \emph{ghost shift}.
(More general constructions are allowed in see \cite[Section 2]{MaTu-klrw-algebras}, but we do not need the most general setting here.)
Then each solid $i$-string has a ghost $i$-string that is shifted $\sigma_{\epsilon}$ units to the right. The ghost string mimics its solid string. In particular, for each dot on the solid string there is a corresponding dot on the ghost string that is shifted $\sigma_{\epsilon}$ units to the right.

An example is:
\begin{gather*}
\begin{tikzpicture}[anchorbase]
\node[dynkin=1] (0) at (0,0){};
\node[dynkin=2] (1) at (1,0){};
\node[dynkin=3] (2) at (2,0){};
\draw[directed=0.5](1)--(0)node[above,xshift=0.5cm,yshift=0.1cm]{$0.75$};
\draw[directed=0.5](1)--(2)node[above,xshift=-0.5cm,yshift=0.1cm]{$1.25$};
\end{tikzpicture}
\quad\rightsquigarrow\quad
\begin{tikzpicture}[anchorbase,smallnodes,rounded corners]
\draw[ghost](0.75,0)node[below]{$\phantom{i}$}--++(0,1)node[above,yshift=-0.05cm]{$2$};
\draw[ghost](1.25,0)node[below]{$\phantom{i}$}--++(0,1)node[above,yshift=-0.05cm]{$2$};
\draw[solid](0,0)node[below]{$2$}--++(0,1)node[above,yshift=-1pt]{$\phantom{i}$};
\end{tikzpicture}
\text{ and }
\begin{tikzpicture}[anchorbase,smallnodes,rounded corners]
\draw[ghost,dot](0.75,0)node[below]{$\phantom{i}$}--++(0,1)node[above,yshift=-0.05cm]{$2$};
\draw[ghost,dot](1.25,0)node[below]{$\phantom{i}$}--++(0,1)node[above,yshift=-0.05cm]{$2$};
\draw[solid,dot](0,0)node[below]{$2$}--++(0,1)node[above,yshift=-1pt]{$\phantom{i}$};
\end{tikzpicture}
.
\end{gather*}

\item To position the solid and ghost strings we
fix a \emph{solid positioning} $\bx=(x_{1},\dots,x_{n})\in\R^{n}$,
and a \emph{ghost shift} $\bsig=(\sigma_{\epsilon}\in\R_{>0})_{\epsilon\in\edges}$, which is graph theoretic \emph{weighting} of the quiver $\Gamma$. The solid positioning gives the
\emph{coordinates} of the solid strings on the real line,
with the corresponding ghost strings shifted $\sigma_{\epsilon}$ units to the right.
The positions of the red strings
are determined by a \emph{charge}
$\charge=(\kappa_{1},\dots,\kappa_{\ell})\in\R^{\ell}$ with
$\kappa_{1}<\dots<\kappa_{\ell}$. The solid positioning, ghost shifts and charge are chosen so that no two strings have the same endpoints. All diagrams use the same ghost shifts and charge but we allow different diagrams to have different solid positionings.

For example, in \autoref{Eq:RecollectionExampleStrings} we have
$\bx=(0,1.2,2)$ and $\charge=(-0.5,1)$. In this diagram, the solid $2$-string has two ghost strings, with ghost shifts $0.75$ and $1.25$, because the corresponding quiver has two edges starting at $2$.

\end{enumerate}

\end{enumerate}

\begin{Notation}\label{N:RecollectionOmit}
\leavevmode

\begin{enumerate}

\item For clarity, we often omit strings from diagrams if they are not relevant to the features of the diagram that we are describing.

\item Care needs to be taken when drawing diagrams because they are very sensitive
with respect to the positions of the strings. However, we can rescale diagrams without harm by \cite[Section 5]{MaTu-klrw-algebras}, so we often do this to improve our exposition.

\end{enumerate}
\end{Notation}

%%%%%%%%%%%%%%%%%%%%%%%%%%%%%%%%%%%%%%%%%

\subsection{Q-polynomials}\label{SS:RecollectionQPoly}

%%%%%%%%%%%%%%%%%%%%%%%%%%%%%%%%%%%%%%%%%

Except in \autoref{SS:Affine}, we only consider quivers $\quiver$ with $a_{ij}a_{ji}\leq 3$.
This includes all quivers of finite and affine type except
$A^{(2)}_2$.
We will always use the following choice of
$Q$-polynomials (although the relations in
\autoref{SS:RecollectionRelations} work in greater
generality):
\begin{gather*}
Q_{ij}(u,v)=
\left\{
\begin{array}{CC@{\quad}|@{\quad}CC}
$u-v^{\phantom{1}}$ & if $i\rightarrow j$, & $v-u^{\phantom{1}}$ & if $i\leftarrow j$,
\\
$u-v^{2}$ & if $i\Rightarrow j$, & $v-u^{2}$ & if $i\Leftarrow j$,
\\
$u-v^{3}$ & if $i\Rrightarrow j$, & $v-u^{3}$ & if $i\Lleftarrow j$,
\\
$0^{\phantom{1}}$ & if $i=j$, & $1^{\phantom{1}}$ & otherwise.
\\
\end{array}
\right.
\end{gather*}
With these choices the values for $Q_{i,j,i}(u,v,w)=\frac{Q_{ij}(u,v)-Q_{ij}(u,w)}{w-v}$ are
\begin{gather*}
Q_{i,j,i}(u,v,w)=
\left\{
\begin{array}{CC@{\quad}|@{\quad}CC}
$1$ & if $i\rightarrow j$, & $-1$ & if $i\leftarrow j$,
\\
$v+w$ & if $i\Rightarrow j$, & $-v-w$ & if $i\Leftarrow j$,
\\
$v^{2}+vw+w^{2}$ & if $i\Rrightarrow j$, & $-v^{2}-vw-w^{2}$ & if $i\Lleftarrow j$,
\\
$0^{\phantom{1}}$ & if $i=j$, & $0$ & otherwise.
\\
\end{array}
\right.
\end{gather*}
To use these polynomials in relations we make the following
definition.
Let $\bu=(u_{1},\dots,u_{n})$ and let $D$ be a diagram with
$n$ solid strings.
For a monomial $f(\bu)=u_{1}^{a_{1}}\dots u_{n}^{a_{n}}\in\N[u_{1},\dots,u_{n}]$ the diagram $f(\bu)D$ is the diagram obtained from $D$ by putting $a_{k}$ dots at the top of the $k$th solid string, and $a_{k}$ dots on all of its ghosts. Here we
number the solid strings from left to right along the top of the diagram so that the $k$th solid string is the $k$th string from the left.
(From \autoref{N:ConstructionRepeated} onwards, we use
a different numbering.)
More generally, given an arbitrary polynomial $f(\bu)\in\N[u_{1},\dots,u_{n}]$ let $f(\bu)D$ be the corresponding linear combination of diagrams.

%%%%%%%%%%%%%%%%%%%%%%%%%%%%%%%%%%%%%%%%%

\subsection{The relations of wKLRW algebras}\label{SS:RecollectionRelations}

%%%%%%%%%%%%%%%%%%%%%%%%%%%%%%%%%%%%%%%%%

We repeat the multilocal relations of wKLRW algebras from \cite[Definition 2C.7]{MaTu-klrw-algebras}. The reader may consult that paper for examples and general results about the wKLRW algebras.

\begin{Definition}\label{D:RecollectionClose}
Two strings are \emph{close} if,  using isotopies, they can be pulled arbitrarily
close together in a neighborhood that does not contain any other strings.
\end{Definition}

It is often necessary to pull strings close together in order to apply the relations below.

\begin{Example}\label{E:RecollectionBilocal}
Let $i,j\in I$ and consider the following two diagrams:
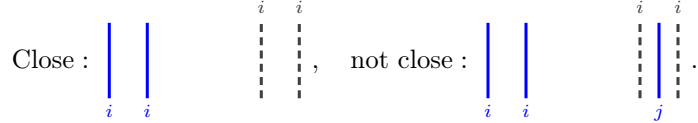
\begin{gather*}
\text{Close}:
\begin{tikzpicture}[anchorbase,smallnodes]
\draw[ghost](2,0)--++(0,1)node[above,yshift=-1pt]{$i$};
\draw[ghost](2.5,0)--++(0,1)node[above,yshift=-1pt]{$i$};
\draw[solid](0,0)node[below]{$i$}--++(0,1)node[above,yshift=-1pt]{$\phantom{i}$};
\draw[solid](0.5,0)node[below]{$i$}--++(0,1)node[above,yshift=-1pt]{$\phantom{i}$};
\end{tikzpicture}
,\quad
\text{not close}:
\begin{tikzpicture}[anchorbase,smallnodes]
\draw[ghost](2,0)--++(0,1)node[above,yshift=-1pt]{$i$};
\draw[ghost](2.5,0)--++(0,1)node[above,yshift=-1pt]{$i$};
\draw[solid](0,0)node[below]{$i$}--++(0,1)node[above,yshift=-1pt]{$\phantom{i}$};
\draw[solid](0.5,0)node[below]{$i$}--++(0,1)node[above,yshift=-1pt]{$\phantom{i}$};
\draw[solid](2.25,0)node[below]{$j$}--++(0,1)node[above,yshift=-1pt]{$\phantom{i}$};
\end{tikzpicture}
.
\end{gather*}
The two solid (and ghost) $i$-strings on the left are close, but the
two solid $i$-strings on the right are not because the solid $j$-string prevents the ghost $i$-strings from being close.
This implies that the relation \autoref{Eq:RecollectionDotCrossing}
below can be applied to the left-hand diagram. In contrast, \autoref{Eq:RecollectionDotCrossing} cannot be applied to the right-hand diagram even though it appears that it can be applied if we look only in a local neighborhood of the two solid $i$-strings.
\end{Example}

Following \cite{MaTu-klrw-algebras}, a \emph{multilocal} relation
is a relation that needs to be applied simultaneously
in two or more local neighborhoods. We use multilocal relations because we want to simultaneously apply relations to the solids strings and their ghosts.
For simplicity, we often only illustrate one of the local neighborhoods where the relations are applied.
Multilocal relations can only be applied if the strings in all relevant neighborhoods are close.

Recall from \autoref{SS:RecollectionDiagrams}.(b).(iii) that we have fixed a charge $\charge=(\kappa_1,\dots,\kappa_\ell)\in\mathbb{R}^\ell$. Fix $\varepsilon\ll1$ such that
\begin{gather}\label{Eq:ConstructionEpsilon}
0<\varepsilon<\tfrac{1}{n^{2}}\min\set[\big]{\tfrac{1}{4n(\ell+ne)},\kappa_{l+1}-\kappa_{l}|1\le i<\ell},
\end{gather}
and define the set $X$ of \emph{positions} to be
\begin{gather}\label{Eq:X}
X=X_{n}
=\set{\kappa_{1}-m\varepsilon-m|1\leq m\leq n}.
\end{gather}
The set $X$ gives the possible solid positionings for the solid strings. That is, the top and bottom coordinates of strings are given by the $n$ coordinates in $X$. The definition of $\varepsilon$ in \autoref{Eq:ConstructionEpsilon} ensures that there will always be enough parking positions for our strings, in the sense of \autoref{D:ConstructionParking} below. The positions in~$X$ ensure that the endpoints solid strings are always to the left of the red strings, which implies that the wKLRW algebras defined below are isomorphic to KLR algebras.

Recall that we fixed  $\brho\in I^\ell$ in \autoref{SS:RecollectionDiagrams}. The red strings in the wKLRW diagrams have residues $\rho_1,\dots,\rho_\ell$, when read from left to right.

\begin{Definition}\label{D:hell}
The \emph{affine level} is
$\hell=\ell+ne$.
Define the \emph{affine charge}
$\affine{\charge}=
(\affine{\kappa}_{1},\dots,\affine{\kappa}_{\hell})\in\Z^{\hell}$
and the \emph{affine red labels}
$\affine{\brho}=(\affine{\rho}_{1},\dots,\affine{\rho}_{\hell})\in \vertices^{\hell}$ by
\begin{gather*}
\affine{\kappa}_{m}=
\begin{cases*}
\kappa_{m} & if $1\leq m\leq\ell$,\\
\kappa_{\ell}+(2n+1)(m-\ell) & otherwise,
\end{cases*}
\quad\text{and}\quad
\affine{\rho}_{m}=
\begin{cases*}
\rho_{m} & if $1\leq m\leq\ell$,\\
\floor{m-\ell-1}{n}+e\Z & otherwise,
\end{cases*}
\end{gather*}
and set
\begin{gather*}
\abfweight=(\fweight[\affine{\rho}_{1}],\dots,
\fweight[\affine{\rho}_{\hell}])
=(\fweight[{\rho_{1}}],\dots,\fweight[{\rho_{\ell}}],
\underbrace{\fweight[1],\dots,\fweight[1]}_{n\text{ times}},\dots,\underbrace{\fweight[e],\dots,\fweight[e]}_{n\text{ times}}).
\end{gather*}
\end{Definition}

%Let $\aWA$ be the wKLRW algebra defined using $\affine{\brho}$.
The \emph{affine wKLRW algebra} $\WA=\WA[n](X)$
is the graded unital associative $\ring$-algebra spanned by
these string diagrams, with \emph{isotopy relations}
in $\R^{2}$ and subject to the multilocal relations:

\begin{enumerate}

\item The \emph{(honest) dot sliding relations} hold (that is, solid and ghost dots can pass through any crossing)
except in the following cases:
\begin{gather}\label{Eq:RecollectionDotCrossing}
\begin{tikzpicture}[anchorbase,smallnodes,rounded corners]
\draw[solid](0.5,0.5)node[above,yshift=-1pt]{$\phantom{i}$}--(0,0) node[below]{$i$};
\draw[solid,dot=0.25](0,0.5)--(0.5,0) node[below]{$i$};
\end{tikzpicture}
-
\begin{tikzpicture}[anchorbase,smallnodes,rounded corners]
\draw[solid](0.5,0.5)node[above,yshift=-1pt]{$\phantom{i}$}--(0,0) node[below]{$i$};
\draw[solid,dot=0.75](0,0.5)--(0.5,0) node[below]{$i$};
\end{tikzpicture}
=
\begin{tikzpicture}[anchorbase,smallnodes,rounded corners]
\draw[solid](0,0.5)node[above,yshift=-1pt]{$\phantom{i}$}--(0,0) node[below]{$i$};
\draw[solid](0.5,0.5)--(0.5,0) node[below]{$i$};
\end{tikzpicture}
=
\begin{tikzpicture}[anchorbase,smallnodes,rounded corners]
\draw[solid,dot=0.75](0.5,0.5)node[above,yshift=-1pt]{$\phantom{i}$}--(0,0) node[below]{$i$};
\draw[solid](0,0.5)--(0.5,0) node[below]{$i$};
\end{tikzpicture}
-
\begin{tikzpicture}[anchorbase,smallnodes,rounded corners]
\draw[solid,dot=0.25](0.5,0.5)node[above,yshift=-1pt]{$\phantom{i}$}--(0,0) node[below]{$i$};
\draw[solid](0,0.5)--(0.5,0) node[below]{$i$};
\end{tikzpicture}
.
\end{gather}

\item The \emph{(honest) Reidemeister II relations} hold except in the following cases:
\begin{gather}\label{Eq:RecollectionReidemeisterII}
\begin{gathered}
\begin{tikzpicture}[anchorbase,smallnodes,rounded corners]
\draw[solid](0,1)--++(0.5,-0.5)--++(-0.5,-0.5) node[below]{$i$};
\draw[solid](0.5,1)node[above,yshift=-1pt]{$\phantom{i}$}--++(-0.5,-0.5)--++(0.5,-0.5) node[below]{$i$};
\end{tikzpicture}
=0
,\quad
\begin{tikzpicture}[anchorbase,smallnodes,rounded corners]
\draw[ghost](0,1)node[above,yshift=-1pt]{$i$}--++(0.5,-0.5)--++(-0.5,-0.5) node[below]{$\phantom{i}$};
\draw[solid](0.5,1)--++(-0.5,-0.5)--++(0.5,-0.5) node[below]{$j$};
\end{tikzpicture}
=Q_{ij}(\bu)
\begin{tikzpicture}[anchorbase,smallnodes,rounded corners]
\draw[ghost](0,1)node[above,yshift=-1pt]{$i$}--++(0,-1)node[below]{$\phantom{i}$};
\draw[solid](0.5,1)--++(0,-1)node[below]{$j$};
\end{tikzpicture}
\quad\text{or}\quad
\begin{tikzpicture}[anchorbase,smallnodes,rounded corners]
\draw[ghost](0.5,1)node[above,yshift=-1pt]{$i$}--++(-0.5,-0.5)--++(0.5,-0.5) node[below]{$\phantom{i}$};
\draw[solid](0,1)--++(0.5,-0.5)--++(-0.5,-0.5) node[below]{$j$};
\end{tikzpicture}
=Q_{ji}(\bu)
\begin{tikzpicture}[anchorbase,smallnodes,rounded corners]
\draw[ghost](0.5,1)node[above,yshift=-1pt]{$i$}--++(0,-1)node[below]{$\phantom{i}$};
\draw[solid](0,1)--++(0,-1)node[below]{$j$};
\end{tikzpicture}
\quad
\text{if $i\rightsquigarrow j$}
,\\
\begin{tikzpicture}[anchorbase,smallnodes,rounded corners]
\draw[solid](0.5,1)node[above,yshift=-1pt]{$\phantom{i}$}--++(-0.5,-0.5)--++(0.5,-0.5) node[below]{$i$};
\draw[redstring](0,1)--++(0.5,-0.5)--++(-0.5,-0.5) node[below]{$i$};
\end{tikzpicture}
=
\begin{tikzpicture}[anchorbase,smallnodes,rounded corners]
\draw[solid,dot](0.5,0)node[below]{$i$}--++(0,1)node[above,yshift=-1pt]{$\phantom{i}$};
\draw[redstring](0,0)node[below]{$i$}--++(0,1);
\end{tikzpicture}
,\quad
\begin{tikzpicture}[anchorbase,smallnodes,rounded corners]
\draw[solid](0,1)node[above,yshift=-1pt]{$\phantom{i}$}--++(0.5,-0.5)--++(-0.5,-0.5) node[below]{$i$};
\draw[redstring](0.5,1)--++(-0.5,-0.5)--++(0.5,-0.5) node[below]{$i$};
\end{tikzpicture}
=
\begin{tikzpicture}[anchorbase,smallnodes,rounded corners]
\draw[solid,dot](0,0)node[below]{$i$}--++(0,1)node[above,yshift=-1pt]{$\phantom{i}$};
\draw[redstring](0.5,0)node[below]{$i$}--++(0,1);
\end{tikzpicture}.
\end{gathered}
\end{gather}

\item The \emph{(honest) Reidemeister III relations} hold except in the following cases:
\begin{gather}\label{Eq:RecollectionReidemeisterIII}
\begin{gathered}
\begin{tikzpicture}[anchorbase,smallnodes,rounded corners]
\draw[ghost](1,1)node[above,yshift=-1pt]{$i$}--++(1,-1)node[below]{$\phantom{i}$};
\draw[ghost](2,1)node[above,yshift=-1pt]{$i$}--++(-1,-1)node[below]{$\phantom{i}$};
\draw[solid,smallnodes,rounded corners](1.5,1)--++(-0.5,-0.5)--++(0.5,-0.5)node[below]{$j$};
\end{tikzpicture}
{=}
\begin{tikzpicture}[anchorbase,smallnodes,rounded corners]
\draw[ghost](3,1)node[above,yshift=-1pt]{$\phantom{i}$}--++(1,-1)node[below]{$\phantom{i}$};
\draw[ghost](4,1)node[above,yshift=-1pt]{$i$}--++(-1,-1)node[below]{$\phantom{i}$};
\draw[solid,smallnodes,rounded corners](3.5,1)--++(0.5,-0.5)--++(-0.5,-0.5)node[below]{$j$};
\end{tikzpicture}
-Q_{iji}(\bu)
\begin{tikzpicture}[anchorbase,smallnodes,rounded corners]
\draw[ghost](6.2,1)node[above,yshift=-1pt]{$i$}--++(0,-1)node[below]{$\phantom{i}$};
\draw[ghost](7.2,1)node[above,yshift=-1pt]{$i$}--++(0,-1)node[below]{$\phantom{i}$};
\draw[solid](6.7,1)--++(0,-1)node[below]{$j$};
\end{tikzpicture}
\;\text{or}\;
\begin{tikzpicture}[anchorbase,smallnodes,rounded corners]
\draw[solid](1,1)node[above,yshift=-1pt]{$\phantom{i}$}--++(1,-1)node[below]{$j$};
\draw[solid](2,1)--++(-1,-1)node[below]{$j$};
\draw[ghost,smallnodes,rounded corners](1.5,1)node[above,yshift=-1pt]{$i$}--++(-0.5,-0.5)--++(0.5,-0.5)node[below]{$\phantom{i}$};
\end{tikzpicture}
{=}
\begin{tikzpicture}[anchorbase,smallnodes,rounded corners]
\draw[solid](3,1)node[above,yshift=-1pt]{$\phantom{i}$}--++(1,-1)node[below]{$j$};
\draw[solid](4,1)--++(-1,-1)node[below]{$j$};
\draw[ghost,smallnodes,rounded corners](3.5,1)node[above,yshift=-1pt]{$i$}--++(0.5,-0.5)--++(-0.5,-0.5)node[below]{$\phantom{i}$};
\end{tikzpicture}
+Q_{jij}(\bu)
\begin{tikzpicture}[anchorbase,smallnodes,rounded corners]
\draw[solid](7.2,1)node[above,yshift=-1pt]{$\phantom{i}$}--++(0,-1)node[below]{$j$};
\draw[solid](8.2,1)--++(0,-1)node[below]{$j$};
\draw[ghost](7.7,1)node[above,yshift=-1pt]{$i$}--++(0,-1)node[below]{$\phantom{i}$};
\end{tikzpicture}
\quad
\text{if $i\rightsquigarrow j$}
,\\
\begin{tikzpicture}[anchorbase,smallnodes,rounded corners]
\draw[solid](1,1)node[above,yshift=-1pt]{$\phantom{i}$}--++(1,-1)node[below]{$i$};
\draw[solid](2,1)--++(-1,-1)node[below]{$i$};
\draw[redstring](1.5,1)--++(-0.5,-0.5)--++(0.5,-0.5)node[below]{$i$};
\end{tikzpicture}
=
\begin{tikzpicture}[anchorbase,smallnodes,rounded corners]
\draw[solid](3,1)node[above,yshift=-1pt]{$\phantom{i}$}--++(1,-1)node[below]{$i$};
\draw[solid](4,1)--++(-1,-1)node[below]{$i$};
\draw[redstring](3.5,1)--++(0.5,-0.5)--++(-0.5,-0.5)node[below]{$i$};
\end{tikzpicture}
-
\begin{tikzpicture}[anchorbase,smallnodes,rounded corners]
\draw[solid](5,1)node[above,yshift=-1pt]{$\phantom{i}$}--++(0,-1)node[below]{$i$};
\draw[solid](6,1)--++(0,-1)node[below]{$i$};
\draw[redstring](5.5,1)--++(0,-1)node[below]{$i$};
\end{tikzpicture}
.
\end{gathered}
\end{gather}
\end{enumerate}

The algebra $\WA$ is graded because all of the relations are homogeneous with respect to the degree function given in \autoref{Eq:RecollectionDegree} below.

As in \autoref{SS:RecollectionQPoly},
the right-hand side of relations \autoref{Eq:RecollectionReidemeisterII} and
\autoref{Eq:RecollectionReidemeisterIII} carry dots determined by the fixed choices of $Q$-polynomials. We emphasize that in these relations we have drawn only a solid string or its ghost, but the other string is still implicitly part of the relation. The relations can only be applied if they can be simultaneously applied in local neighborhoods around both the solid and ghost strings in the relation.

\begin{Example}\label{E:RecollectionReidemeisterII}
Note that the relations \autoref{Eq:RecollectionReidemeisterII} and \autoref{Eq:RecollectionReidemeisterIII} are not symmetric with respect to the solid and ghost strings.
For example, take the quiver $\dynkin A3$
of type $\typea[3]$ numbered and
oriented left to right.
On the one hand, \autoref{Eq:RecollectionReidemeisterII}
gives a nontrivial relation between the ghost $1$-string and
the solid $2$-string. On the other hand, the
ghost $2$-string and the solid $1$-string satisfy an honest
Reidemeister II relation. That is, we have:
\begin{gather*}
\begin{tikzpicture}[anchorbase,smallnodes,rounded corners]
\draw[ghost](0,1)node[above,yshift=-1pt]{$1$}--++(0.5,-0.5)--++(-0.5,-0.5) node[below]{$\phantom{i}$};
\draw[solid](0.5,1)--++(-0.5,-0.5)--++(0.5,-0.5) node[below]{$2$};
\end{tikzpicture}
=Q_{ij}(\bu)
\begin{tikzpicture}[anchorbase,smallnodes,rounded corners]
\draw[ghost](0,1)node[above,yshift=-1pt]{$1$}--++(0,-1)node[below]{$\phantom{i}$};
\draw[solid](0.5,1)--++(0,-1)node[below]{$2$};
\end{tikzpicture}
\text{ and }
\begin{tikzpicture}[anchorbase,smallnodes,rounded corners]
\draw[ghost](0,1)node[above,yshift=-1pt]{$2$}--++(0.5,-0.5)--++(-0.5,-0.5) node[below]{$\phantom{i}$};
\draw[solid](0.5,1)--++(-0.5,-0.5)--++(0.5,-0.5) node[below]{$1$};
\end{tikzpicture}
=
\begin{tikzpicture}[anchorbase,smallnodes,rounded corners]
\draw[ghost](0,1)node[above,yshift=-1pt]{$2$}--++(0,-1)node[below]{$\phantom{i}$};
\draw[solid](0.5,1)--++(0,-1)node[below]{$1$};
\end{tikzpicture}
.
\end{gather*}
This asymmetry is an important feature of the wKLRW algebras.
\end{Example}

%%%%%%%%%%%%%%%%%%%%%%%%%%%%%%%%%%%%%%%%%

\subsection{Degrees of diagrams}\label{SS:RecollectionDegree}

%%%%%%%%%%%%%%%%%%%%%%%%%%%%%%%%%%%%%%%%%

Diagrams are given a grading by using the following local rules (not multilocal).
The \emph{degree} of a fixed diagram is the sum of its local degrees. The local degrees are:
\begin{gather}\label{Eq:RecollectionDegree}
\begin{gathered}
\deg\begin{tikzpicture}[anchorbase,smallnodes,rounded corners]
\draw[solid,dot] (0,0)node[below]{$i$} to (0,0.5)node[above,yshift=-1pt]{$\phantom{i}$};
\end{tikzpicture}
=2d_{i}
,\quad
\deg\begin{tikzpicture}[anchorbase,smallnodes,rounded corners]
\draw[ghost,dot] (0,0)node[below]{$\phantom{i}$} to (0,0.5)node[above,yshift=-1pt]{$i$};
\end{tikzpicture}=0,
\quad
\deg\begin{tikzpicture}[anchorbase,smallnodes,rounded corners]
\draw[solid] (0,0)node[below]{$i$} to (0.5,0.5);
\draw[solid] (0.5,0)node[below]{$j$} to (0,0.5)node[above,yshift=-1pt]{$\phantom{i}$};
\end{tikzpicture}
=-\delta_{i,j}2d_{i}
,\quad
\deg\begin{tikzpicture}[anchorbase,smallnodes,rounded corners]
\draw[ghost] (0,0)node[below]{$\phantom{i}$} to (0.5,0.5)node[above,yshift=-1pt]{$i$};
\draw[solid] (0.5,0)node[below]{$j$} to (0,0.5)node[above,yshift=-1pt]{$\phantom{i}$};
\end{tikzpicture}
=
\deg\begin{tikzpicture}[anchorbase,smallnodes,rounded corners]
\draw[ghost] (0.5,0)node[below]{$\phantom{i}$} to (0,0.5)node[above,yshift=-1pt]{$i$};
\draw[solid] (0,0)node[below]{$j$} to (0.5,0.5)node[above,yshift=-1pt]{$\phantom{i}$};
\end{tikzpicture}
=
\begin{cases}
-(i,j)
&\text{if $i\rightsquigarrow j$},
\\
0&\text{else},
\end{cases}
\\
\deg\begin{tikzpicture}[anchorbase,smallnodes,rounded corners]
\draw[ghost] (0,0)node[below]{$\phantom{i}$} to (0.5,0.5)node[above,yshift=-1pt]{$i$};
\draw[ghost] (0.5,0)node[below]{$\phantom{i}$} to (0,0.5)node[above,yshift=-1pt]{$j$};
\end{tikzpicture}
=0
,\quad
\deg\begin{tikzpicture}[anchorbase,smallnodes,rounded corners]
\draw[solid] (0,0)node[below]{$i$} to (0.5,0.5)node[above,yshift=-1pt]{$\phantom{i}$};
\draw[redstring] (0.5,0)node[below]{$j$} to (0,0.5);
\end{tikzpicture}
=
\deg\begin{tikzpicture}[anchorbase,smallnodes,rounded corners]
\draw[solid] (0.5,0)node[below]{$i$} to (0,0.5)node[above,yshift=-1pt]{$\phantom{i}$};
\draw[redstring] (0,0)node[below]{$j$} to (0.5,0.5);
\end{tikzpicture}
=\tfrac{1}{2}\delta_{i,j}(i,i)
,\quad
\deg\begin{tikzpicture}[anchorbase,smallnodes,rounded corners]
\draw[ghost] (0,0)node[below]{$\phantom{i}$} to (0.5,0.5)node[above,yshift=-1pt]{$i$};
\draw[redstring] (0.5,0)node[below]{$j$} to (0,0.5)node[above,yshift=-1pt]{$\phantom{i}$};
\end{tikzpicture}
=
\deg\begin{tikzpicture}[anchorbase,smallnodes,rounded corners]
\draw[ghost] (0.5,0)node[below]{$\phantom{i}$} to (0,0.5)node[above,yshift=-1pt]{$i$};
\draw[redstring] (0,0)node[below]{$j$} to (0.5,0.5)node[above,yshift=-1pt]{$\phantom{i}$};
\end{tikzpicture}
=0.
\end{gathered}
\end{gather}
Here, $(\placeholder,\placeholder)$ is the Cartan pairing
associated to $\Gamma$, which gives the entries of the
Cartan matrix, and $(d_{1},\dots,d_{e})$ is the symmetrizer. It is straightforward to check that the relations of~$\WA(X)$ are homogeneous with respect to these local rules, so $\WA(X)$ is a graded algebra with degree function $\deg$.

%%%%%%%%%%%%%%%%%%%%%%%%%%%%%%%%%%%%%%%%%

\subsection{Cyclotomic quotients}\label{SS:RecollectionSteady}

%%%%%%%%%%%%%%%%%%%%%%%%%%%%%%%%%%%%%%%%%

We now define of the finite dimensional
quotients of $\WA(X)$ that play a key role in this paper.

\begin{Example}\label{E:RecollectionIsotopy}
All diagrams in this paper are defined up to isotopy, so isotopy relations play a crucial role when working with wKLRW algebras. Let us give an example and a nonexample of isotopy:
\begin{gather*}
\begin{tikzpicture}[anchorbase,smallnodes,rounded corners]
\draw[ghost](1,0)--++(0,1)node[above,yshift=-1pt]{$i$};
\draw[solid](0.5,0)node[below]{$i$}--++(0,1)node[above,yshift=-1pt]{$\phantom{i}$};
\draw[redstring](0.25,0)node[below]{$\rho$}--++(0,1)node[above,yshift=-1pt]{$\phantom{i}$};
\end{tikzpicture}
\xrightarrow{\text{isotopy}}
\begin{tikzpicture}[anchorbase,smallnodes,rounded corners]
\draw[ghost](1,0)--++(-0.5,0.25)--++(0.5,0.25)--++(0.5,0.25)--++(-0.5,0.25)node[above,yshift=-1pt]{$i$};
\draw[solid](0.5,0)node[below]{$i$}--++(-0.5,0.25)--++(0.5,0.25)--++(0.5,0.25)--++(-0.5,0.25)node[above,yshift=-1pt]{$\phantom{i}$};
\draw[redstring](0.25,0)node[below]{$\rho$}--++(-1,0.5)--++(1,0.5)node[above,yshift=-1pt]{$\phantom{i}$};
\end{tikzpicture}
,\quad
\text{not an isotopy of the left-hand diagram}:
\begin{tikzpicture}[anchorbase,smallnodes,rounded corners]
\draw[ghost](1,0)--++(-0.5,0.5)--++(0.5,0.5)node[above,yshift=-1pt]{$i$};
\draw[solid](0.5,0)node[below]{$i$}--++(-0.5,0.5)--++(0.5,0.5)node[above,yshift=-1pt]{$\phantom{i}$};
\draw[redstring](0.25,0)node[below]{$\rho$}--++(0,1)node[above,yshift=-1pt]{$\phantom{i}$};
\end{tikzpicture}
.
\end{gather*}
Any move that changes the topological nature of
a diagram is not an isotopy. All isotopies need to deform solid strings and their ghosts in the same way and dots must move with the strings that they are on. Isotopies never introduce new crossings between strings.
\end{Example}

A solid string is \emph{unsteady} if it can be pulled, using isotopies, arbitrarily far to the right while the
red strings stay inside a bounded region. A diagram is \emph{unsteady} if it contains an unsteady string. A string or diagram is \emph{steady} if it is not unsteady.

Note that being unsteady is biased towards the right-hand side. A unsteady diagram can contain solid strings that are steady.

\begin{Example}\label{E:RecollectionUnsteady}
The following diagrams are examples of an unsteady and a steady diagram.
\begin{gather*}
\text{Unsteady}\colon
\begin{tikzpicture}[anchorbase,smallnodes]
\draw[solid](0.5,0)node[below]{$i$}--++(0,1)node[above,yshift=-1pt]{$\phantom{i}$};
\draw[redstring](0.25,0)node[below]{$\rho$}--++(0,1)node[above,yshift=-1pt]{$\phantom{i}$};
\draw[->,decorate,decoration={snake,amplitude=0.4mm,segment length=2mm,post length=1mm},spinach] (0.5,0.5) to (1.5,0.5)node[right]{pulls freely};
\end{tikzpicture}
,\quad
\text{steady}\colon
\begin{tikzpicture}[anchorbase,smallnodes]
\draw[solid](0,0)node[below]{$i$}--++(0,1)node[above,yshift=-1pt]{$\phantom{i}$};
\draw[redstring](0.25,0)node[below]{$i$}--++(0,1)node[above,yshift=-1pt]{$\phantom{i}$};
\draw[->,decorate,decoration={snake,amplitude=0.4mm,segment length=2mm,post length=1mm},spinach] (0,0.5) to (1,0.5)node[right]{blocked by the red string};
\end{tikzpicture}
.
\end{gather*}
The relations imply that the solid $i$-string
cannot be pulled through the red $i$-string, so the right-hand diagram is steady.
Note that the assumption that the red strings
stays inside a bounded region is necessary because otherwise
the solid $i$-string and the red $i$-string
in the right-hand diagram can be simultaneously pulled arbitrarily far to the right under isotopy.
\end{Example}

The \emph{cyclotomic wKLRW algebra} $\WAc[n](X)$ is the quotient
of $\WA[n](X)$ by the two-sided ideal $\uideal$ generated by the unsteady diagrams. As the unsteady ideal is homogeneous, the cyclotomic wKLRW algebra $\WAc(X)$ is also graded.

By \cite[Proposition 3B.12]{MaTu-klrw-algebras}, the affine wKLRW algebras are infinite dimensional and by \cite[Proposition 3D.4]{MaTu-klrw-algebras} the cyclotomic wKLRW algebras are finite dimensional algebras.

\begin{Remark}\label{R:RecollectionUnsteady}
One of the key features of the homogeneous affine sandwich cellular
bases constructed in \cite{MaTu-klrw-algebras} and \cite{MaTu-klrw-algebras-bad} is that they restrict to give homogeneous sandwich cellular bases of $\WAc[n](X)$.
We will see that the same is true in this paper.
\end{Remark}

%%%%%%%%%%%%%%%%%%%%%%%%%%%%%%%%%%%%%%%%%

\subsection{Pulling strings and jumping dots}\label{SS:RecollectionMoveRight}

%%%%%%%%%%%%%%%%%%%%%%%%%%%%%%%%%%%%%%%%%

Note that the relations
\autoref{Eq:RecollectionDotCrossing}--\autoref{Eq:RecollectionReidemeisterIII} come in mirrored pairs, where the roles of the solid and ghost strings. We use this as follows:

\begin{Notation}\label{N:RecollectionPartner}
Up to scalars, the set of relations is invariant under reflecting diagrams in their vertical axis and, in some circumstances, swapping the roles of solid and ghost strings (in the sense of \autoref{Eq:RecollectionPartner} below).
Many of the relations we will use have this type of duality, which we call \emph{partner relations}.
\end{Notation}

Related diagrams appearing in a pair of partner relations can have different scalars.
Nevertheless, as we will see, we use the symmetry in the partner relations to reduce the number of relations that we need to consider.

\begin{Example}\label{E:RecollectionPartner}
Let $i\rightsquigarrow j$.
Then
\begin{gather}\label{Eq:RecollectionPartner}
\begin{tikzpicture}[anchorbase,smallnodes,rounded corners]
\draw[ghost](1,1)node[above,yshift=-1pt]{$i$}--++(1,-1)node[below]{$\phantom{i}$};
\draw[ghost](2,1)node[above,yshift=-1pt]{$i$}--++(-1,-1)node[below]{$\phantom{i}$};
\draw[solid,smallnodes,rounded corners](1.5,1)--++(-0.5,-0.5)--++(0.5,-0.5)node[below]{$j$};
\end{tikzpicture}
{=}
\begin{tikzpicture}[anchorbase,smallnodes,rounded corners]
\draw[ghost](3,1)node[above,yshift=-1pt]{$i$}--++(1,-1)node[below]{$\phantom{i}$};
\draw[ghost](4,1)node[above,yshift=-1pt]{$i$}--++(-1,-1)node[below]{$\phantom{i}$};
\draw[solid,smallnodes,rounded corners](3.5,1)--++(0.5,-0.5)--++(-0.5,-0.5)node[below]{$j$};
\end{tikzpicture}
-Q_{iji}(\bu)
\begin{tikzpicture}[anchorbase,smallnodes,rounded corners]
\draw[ghost](6.2,1)node[above,yshift=-1pt]{$i$}--++(0,-1)node[below]{$\phantom{i}$};
\draw[ghost](7.2,1)node[above,yshift=-1pt]{$i$}--++(0,-1)node[below]{$\phantom{i}$};
\draw[solid](6.7,1)--++(0,-1)node[below]{$j$};
\end{tikzpicture}
\xleftrightarrow[\text{relations}]{\text{partner}}
\begin{tikzpicture}[anchorbase,smallnodes,rounded corners]
\draw[solid](1,1)node[above,yshift=-1pt]{$\phantom{i}$}--++(1,-1)node[below]{$j$};
\draw[solid](2,1)--++(-1,-1)node[below]{$j$};
\draw[ghost,smallnodes,rounded corners](1.5,1)node[above,yshift=-1pt]{$i$}--++(-0.5,-0.5)--++(0.5,-0.5)node[below]{$\phantom{i}$};
\end{tikzpicture}
{=}
\begin{tikzpicture}[anchorbase,smallnodes,rounded corners]
\draw[solid](3,1)node[above,yshift=-1pt]{$\phantom{i}$}--++(1,-1)node[below]{$j$};
\draw[solid](4,1)--++(-1,-1)node[below]{$j$};
\draw[ghost,smallnodes,rounded corners](3.5,1)node[above,yshift=-1pt]{$i$}--++(0.5,-0.5)--++(-0.5,-0.5)node[below]{$\phantom{i}$};
\end{tikzpicture}
+Q_{jij}(\bu)
\begin{tikzpicture}[anchorbase,smallnodes,rounded corners]
\draw[solid](7.2,1)node[above,yshift=-1pt]{$\phantom{i}$}--++(0,-1)node[below]{$j$};
\draw[solid](8.2,1)--++(0,-1)node[below]{$j$};
\draw[ghost](7.7,1)node[above,yshift=-1pt]{$i$}--++(0,-1)node[below]{$\phantom{i}$};
\end{tikzpicture}
,
\end{gather}
is an example of partner relations. Note that the scalar in front of the identity diagram changes from $-Q_{iji}(\bu)$ in the left-hand relation to $+Q_{jij}(\bu)$ in the right-hand relation.
\end{Example}

\begin{Remark}\label{R:RecollectionPullingStrings}
One of the key diagrammatic techniques in this paper is
\emph{pulling strings} and
\emph{jumping dots to the right}. As the next lemma shows, this allows us to rewrite some diagrams as a linear combination of diagrams with certain strings or dots moved further to the right. To make it clearer which strings or dots are moving, we mark the string that is being pulled, or the string that carries a jumping dot, in green.

In all the relations listed below, the diagrams on the right-hand side of each relation are obtained from the diagram on the left-hand side by pulling a (green) string, or the dot on a string, further to the right. In \autoref{D:ConstructionChoice} below, we introduce a total order $\sandorder$
that keeps track of how far strings are placed to the right.
\autoref{L:RecollectionMovingStringsDots} can be viewed as saying that diagrams become bigger when strings are pulled to the right. To give our (homogeneous affine) sandwich cellular bases we also want to jump dots as far to the right as possible.
\end{Remark}

As in \autoref{D:RecollectionClose}, all of the multilocal relations in
\autoref{L:RecollectionMovingStringsDots}, and \autoref{L:RecollectionReidemeisterIII}, are close configurations. The relations in \autoref{L:RecollectionMovingStringsDots}.(a) hold for any quiver and they give a stronger form of the failure of the
honest Reidemeister II relations. The relations in \autoref{L:RecollectionMovingStringsDots}.(b) depend on the quiver (and our choice of $Q$-polynomials).

Recall that partner relations are defined in \autoref{N:RecollectionPartner}.

\begin{Lemma}\label{L:RecollectionMovingStringsDots}
Suppose that $i,j\in I$.
\begin{enumerate}
\item The following relations, and their partner relations, hold:
\begin{align*}
\begin{tikzpicture}[anchorbase,smallnodes,rounded corners]
\draw[solid,dot,spinach](0,0)node[below]{$i$}--++(0,1)node[above,yshift=-1pt]{$\phantom{i}$};
\draw[redstring](0.5,0)node[below]{$i$}--++(0,1);
\end{tikzpicture}
=&
\begin{tikzpicture}[anchorbase,smallnodes,rounded corners]
\draw[solid](0,1)node[above,yshift=-1pt]{$\phantom{i}$}--++(0.8,-0.5)--++(-0.8,-0.5)node[below]{$i$};
\draw[redstring](0.5,1)--++(0,-1) node[below]{$i$};
\end{tikzpicture}
,\\
\begin{tikzpicture}[anchorbase,smallnodes,rounded corners]
\draw[solid,spinach](0,1)node[above,yshift=-1pt]{$\phantom{i}$}--++(0,-1)node[below]{$i$};
\draw[solid](0.5,1)--++(0,-1)node[below]{$i$};
\end{tikzpicture}
=&
\begin{tikzpicture}[anchorbase,smallnodes,rounded corners]
\draw[solid,dot](0,1)node[above,yshift=-1pt]{$\phantom{i}$}--++(0.8,-0.5)--++(-0.8,-0.5)node[below]{$i$};
\draw[solid,dot=0.1](0.5,1)--++(0,-1) node[below]{$i$};
\end{tikzpicture}
-
\begin{tikzpicture}[anchorbase,smallnodes,rounded corners]
\draw[solid,dot,dot=0.9](0,1)node[above,yshift=-1pt]{$\phantom{i}$}--++(0.8,-0.5)--++(-0.8,-0.5)node[below]{$i$};
\draw[solid](0.5,1)--++(0,-1) node[below]{$i$};
\end{tikzpicture}
,\\
\begin{tikzpicture}[anchorbase,smallnodes,rounded corners]
\draw[solid,spinach](0,1)node[above,yshift=-1pt]{$\phantom{i}$}--++(0,-1)node[below]{$i$};
\draw[solid](0.5,1)--++(0,-1)node[below]{$i$};
\draw[solid](1,1)--++(0,-1)node[below]{$i$};
\end{tikzpicture}
=&
\begin{tikzpicture}[anchorbase,smallnodes,rounded corners]
\draw[solid,dot=0.425,dot=0.575](0,1)node[above,yshift=-1pt]{$\phantom{i}$}--++(0,-0.2)--++(1.3,-0.15)--++(0,-0.3)--++(-1.3,-0.15)--++(0,-0.2)node[below]{$i$};
\draw[solid,dot=0.1](0.5,1)--++(0,-1) node[below]{$i$};
\draw[solid,dot=0.1](1,1)--++(0,-1)node[below]{$i$};
\end{tikzpicture}
-
\begin{tikzpicture}[anchorbase,smallnodes,rounded corners]
\draw[solid,dot=0.425,dot=0.575,dot=0.725](0,1)node[above,yshift=-1pt]{$\phantom{i}$}--++(0,-0.2)--++(1.3,-0.15)--++(0,-0.3)--++(-1.3,-0.15)--++(0,-0.2)node[below]{$i$};
\draw[solid,dot=0.1](0.5,1)--++(0,-1) node[below]{$i$};
\draw[solid](1,1)--++(0,-1)node[below]{$i$};
\end{tikzpicture}
-
\begin{tikzpicture}[anchorbase,smallnodes,rounded corners]
\draw[solid,dot=0.425,dot=0.575,dot=0.9](0,1)node[above,yshift=-1pt]{$\phantom{i}$}--++(0,-0.2)--++(1.3,-0.15)--++(0,-0.3)--++(-1.3,-0.15)--++(0,-0.2)node[below]{$i$};
\draw[solid](0.5,1)--++(0,-1) node[below]{$i$};
\draw[solid,dot=0.1](1,1)--++(0,-1)node[below]{$i$};
\end{tikzpicture}
+
\begin{tikzpicture}[anchorbase,smallnodes,rounded corners]
\draw[solid,dot=0.425,dot=0.575,dot=0.725,dot=0.9](0,1)node[above,yshift=-1pt]{$\phantom{i}$}--++(0,-0.2)--++(1.3,-0.15)--++(0,-0.3)--++(-1.3,-0.15)--++(0,-0.2)node[below]{$i$};
\draw[solid](0.5,1)--++(0,-1) node[below]{$i$};
\draw[solid](1,1)--++(0,-1)node[below]{$i$};
\end{tikzpicture}
,
\\
\begin{tikzpicture}[anchorbase,smallnodes,rounded corners]
\draw[solid,spinach](0,1)node[above,yshift=-1pt]{$\phantom{i}$}--++(0,-1)node[below]{$i$};
\draw[solid](0.5,1)--++(0,-1)node[below]{$i$};
\draw[solid](1,1)--++(0,-1)node[below]{$i$};
\draw[solid](1.5,1)--++(0,-1)node[below]{$i$};
\end{tikzpicture}
=&
\begin{tikzpicture}[anchorbase,smallnodes,rounded corners]
\draw[solid,dot=0.45,dot=0.5,dot=0.55](0,1)node[above,yshift=-1pt]{$\phantom{i}$}--++(0,-0.2)--++(1.8,-0.15)--++(0,-0.3)--++(-1.8,-0.15)--++(0,-0.2)node[below]{$i$};
\draw[solid,dot=0.1](0.5,1)--++(0,-1) node[below]{$i$};
\draw[solid,dot=0.1](1,1)--++(0,-1)node[below]{$i$};
\draw[solid,dot=0.1](1.5,1)--++(0,-1)node[below]{$i$};
\end{tikzpicture}
-
\begin{tikzpicture}[anchorbase,smallnodes,rounded corners]
\draw[solid,dot=0.45,dot=0.5,dot=0.55,dot=0.66](0,1)node[above,yshift=-1pt]{$\phantom{i}$}--++(0,-0.2)--++(1.8,-0.15)--++(0,-0.3)--++(-1.8,-0.15)--++(0,-0.2)node[below]{$i$};
\draw[solid,dot=0.1](0.5,1)--++(0,-1) node[below]{$i$};
\draw[solid,dot=0.1](1,1)--++(0,-1)node[below]{$i$};
\draw[solid](1.5,1)--++(0,-1)node[below]{$i$};
\end{tikzpicture}
-
\begin{tikzpicture}[anchorbase,smallnodes,rounded corners]
\draw[solid,dot=0.45,dot=0.5,dot=0.55,dot=0.79](0,1)node[above,yshift=-1pt]{$\phantom{i}$}--++(0,-0.2)--++(1.8,-0.15)--++(0,-0.3)--++(-1.8,-0.15)--++(0,-0.2)node[below]{$i$};
\draw[solid,dot=0.1](0.5,1)--++(0,-1) node[below]{$i$};
\draw[solid](1,1)--++(0,-1)node[below]{$i$};
\draw[solid,dot=0.1](1.5,1)--++(0,-1)node[below]{$i$};
\end{tikzpicture}
-
\begin{tikzpicture}[anchorbase,smallnodes,rounded corners]
\draw[solid,dot=0.45,dot=0.5,dot=0.55,dot=0.9](0,1)node[above,yshift=-1pt]{$\phantom{i}$}--++(0,-0.2)--++(1.8,-0.15)--++(0,-0.3)--++(-1.8,-0.15)--++(0,-0.2)node[below]{$i$};
\draw[solid](0.5,1)--++(0,-1) node[below]{$i$};
\draw[solid,dot=0.1](1,1)--++(0,-1)node[below]{$i$};
\draw[solid,dot=0.1](1.5,1)--++(0,-1)node[below]{$i$};
\end{tikzpicture}
\\
&+
\begin{tikzpicture}[anchorbase,smallnodes,rounded corners]
\draw[solid,dot=0.45,dot=0.5,dot=0.55,dot=0.66,dot=0.79](0,1)node[above,yshift=-1pt]{$\phantom{i}$}--++(0,-0.2)--++(1.8,-0.15)--++(0,-0.3)--++(-1.8,-0.15)--++(0,-0.2)node[below]{$i$};
\draw[solid,dot=0.1](0.5,1)--++(0,-1) node[below]{$i$};
\draw[solid](1,1)--++(0,-1)node[below]{$i$};
\draw[solid](1.5,1)--++(0,-1)node[below]{$i$};
\end{tikzpicture}
+
\begin{tikzpicture}[anchorbase,smallnodes,rounded corners]
\draw[solid,dot=0.45,dot=0.5,dot=0.55,dot=0.66,dot=0.9](0,1)node[above,yshift=-1pt]{$\phantom{i}$}--++(0,-0.2)--++(1.8,-0.15)--++(0,-0.3)--++(-1.8,-0.15)--++(0,-0.2)node[below]{$i$};
\draw[solid](0.5,1)--++(0,-1) node[below]{$i$};
\draw[solid,dot=0.1](1,1)--++(0,-1)node[below]{$i$};
\draw[solid](1.5,1)--++(0,-1)node[below]{$i$};
\end{tikzpicture}
+
\begin{tikzpicture}[anchorbase,smallnodes,rounded corners]
\draw[solid,dot=0.45,dot=0.5,dot=0.55,dot=0.79,dot=0.9](0,1)node[above,yshift=-1pt]{$\phantom{i}$}--++(0,-0.2)--++(1.8,-0.15)--++(0,-0.3)--++(-1.8,-0.15)--++(0,-0.2)node[below]{$i$};
\draw[solid](0.5,1)--++(0,-1) node[below]{$i$};
\draw[solid](1,1)--++(0,-1)node[below]{$i$};
\draw[solid,dot=0.1](1.5,1)--++(0,-1)node[below]{$i$};
\end{tikzpicture}
-
\begin{tikzpicture}[anchorbase,smallnodes,rounded corners]
\draw[solid,dot=0.45,dot=0.5,dot=0.55,dot=0.66,dot=0.79,dot=0.9](0,1)node[above,yshift=-1pt]{$\phantom{i}$}--++(0,-0.2)--++(1.8,-0.15)--++(0,-0.3)--++(-1.8,-0.15)--++(0,-0.2)node[below]{$i$};
\draw[solid](0.5,1)--++(0,-1) node[below]{$i$};
\draw[solid](1,1)--++(0,-1)node[below]{$i$};
\draw[solid](1.5,1)--++(0,-1)node[below]{$i$};
\end{tikzpicture}
.
\end{align*}

\item The following relations, and their partner relations, hold:
\begin{gather*}
i\rightarrow j\colon
\begin{tikzpicture}[anchorbase,smallnodes,rounded corners]
\draw[ghost,dot,spinach](0,1)node[above,yshift=-1pt]{$i$}--++(0,-1);
\draw[solid](0.5,1)--++(0,-1)node[below]{$j$};
\end{tikzpicture}
=
\begin{tikzpicture}[anchorbase,smallnodes,rounded corners]
\draw[ghost](0,1)node[above,yshift=-1pt]{$i$}--++(0.8,-0.5)--++(-0.8,-0.5);
\draw[solid](0.5,1)--++(0,-1) node[below]{$j$};
\end{tikzpicture}
+
\begin{tikzpicture}[anchorbase,smallnodes,rounded corners]
\draw[ghost](0,1)node[above,yshift=-1pt]{$i$}--++(0,-1);
\draw[solid,dot](0.5,1)--++(0,-1)node[below]{$j$};
\end{tikzpicture}
,\quad
i\Rightarrow j\colon
\begin{tikzpicture}[anchorbase,smallnodes,rounded corners]
\draw[ghost,dot,spinach](0,1)node[above,yshift=-1pt]{$i$}--++(0,-1);
\draw[solid](0.5,1)--++(0,-1)node[below]{$j$};
\end{tikzpicture}
=
\begin{tikzpicture}[anchorbase,smallnodes,rounded corners]
\draw[ghost](0,1)node[above,yshift=-1pt]{$i$}--++(0.8,-0.5)--++(-0.8,-0.5);
\draw[solid](0.5,1)--++(0,-1) node[below]{$j$};
\end{tikzpicture}
+
\begin{tikzpicture}[anchorbase,smallnodes,rounded corners]
\draw[ghost](0,1)node[above,yshift=-1pt]{$i$}--++(0,-1);
\draw[solid,dot=0.4,dot=0.6](0.5,1)--++(0,-1)node[below]{$j$};
\end{tikzpicture}
,\quad
i\Rrightarrow j\colon
\begin{tikzpicture}[anchorbase,smallnodes,rounded corners]
\draw[ghost,dot,spinach](0,1)node[above,yshift=-1pt]{$i$}--++(0,-1);
\draw[solid](0.5,1)--++(0,-1)node[below]{$j$};
\end{tikzpicture}
=
\begin{tikzpicture}[anchorbase,smallnodes,rounded corners]
\draw[ghost](0,1)node[above,yshift=-1pt]{$i$}--++(0.8,-0.5)--++(-0.8,-0.5);
\draw[solid](0.5,1)--++(0,-1) node[below]{$j$};
\end{tikzpicture}
+
\begin{tikzpicture}[anchorbase,smallnodes,rounded corners]
\draw[ghost](0,1)node[above,yshift=-1pt]{$i$}--++(0,-1);
\draw[solid,dot=0.3,dot=0.5,dot=0.7](0.5,1)--++(0,-1)node[below]{$j$};
\end{tikzpicture}
.
\end{gather*}

\end{enumerate}
\end{Lemma}

\begin{proof}
The first relation in (a), and all of part (b), follows immediately from the defining relations. The second relation in (a) is proved in
\cite[Lemma 6D.1]{MaTu-klrw-algebras}, and the other two relations follow by using this relation, together with \autoref{Eq:RecollectionDotCrossing} and
\autoref{Eq:RecollectionReidemeisterII}, and then cancelling some terms.
\end{proof}

\autoref{L:RecollectionMovingStringsDots} is really just a reformulation of some of the defining relations in $\WA(X)$. As in \cite{MaTu-klrw-algebras} and
\cite{MaTu-klrw-algebras-bad}, this result should be viewed as a mechanism for pulling strings and jumping dots to the right, which is crucial for this paper.

The following relations are inspired by \cite[Equation (5.2)]{Bo-many-cellular-structures}. We call them \emph{plactic relations} for the reasons explained in
\autoref{R:RecollectionPlactic} below.

\begin{Lemma}\label{L:RecollectionReidemeisterIII}
The following relations hold, together with their partner relations:
\begin{align*}
i\not\!\rightsquigarrow j\colon
\begin{tikzpicture}[anchorbase,smallnodes,rounded corners]
\draw[ghost,spinach](1,0)--++(0,1)node[above,yshift=-1pt]{$i$};
\draw[solid](1.25,0)node[below]{$j$}--++(0,1);
\end{tikzpicture}
=&
\begin{tikzpicture}[anchorbase,smallnodes,rounded corners]
\draw[ghost](1,0)--++(0.5,0.5)--++(-0.5,0.5)node[above,yshift=-1pt]{$i$};
\draw[solid](1.25,0)node[below]{$j$}--++(0,1);
\end{tikzpicture}
,\\
i\rightarrow j\colon
\begin{tikzpicture}[anchorbase,smallnodes,rounded corners]
\draw[ghost,spinach](1,0)--++(0,1)node[above,yshift=-1pt]{$i$};
\draw[ghost](1.5,0)--++(0,1)node[above,yshift=-1pt]{$i$};
\draw[solid](1.25,0)node[below]{$j$}--++(0,1);
\end{tikzpicture}
=&
-
\begin{tikzpicture}[anchorbase,smallnodes,rounded corners]
\draw[ghost](1,0)--++(0.75,0.5)--++(-0.75,0.5)node[above,yshift=-1pt]{$i$};
\draw[ghost,dot](1.5,0)--++(0,1)node[above,yshift=-1pt]{$i$};
\draw[solid](1.25,0)node[below]{$j$}--++(0,1);
\end{tikzpicture}
-
\begin{tikzpicture}[anchorbase,smallnodes,rounded corners]
\draw[ghost,dot](1,0)--++(0.75,0.5)--++(-0.75,0.5)node[above,yshift=-1pt]{$i$};
\draw[ghost](1.5,0)--++(0,1)node[above,yshift=-1pt]{$i$};
\draw[solid](1.25,0)node[below]{$j$}--++(0.75,0.5)--++(-0.75,0.5);
\end{tikzpicture}
,
\\
i\Rightarrow j\colon
\begin{tikzpicture}[anchorbase,smallnodes,rounded corners]
\draw[ghost,spinach](0.75,0)--++(0,1)node[above,yshift=-1pt]{$i$};
\draw[ghost](1,0)--++(0,1)node[above,yshift=-1pt]{$i$};
\draw[ghost](1.5,0)--++(0,1)node[above,yshift=-1pt]{$i$};
\draw[solid](1.25,0)node[below]{$j$}--++(0,1);
\end{tikzpicture}
=&
-
\begin{tikzpicture}[anchorbase,smallnodes,rounded corners]
\draw[ghost,dot](0.75,0)--++(0,0.2)--++(1,0.2)--++(0,0.2)--++(-1,0.2)--++(0,0.2)node[above,yshift=-1pt]{$i$};
\draw[ghost,dot=0.95](1,0)--++(0,1)node[above,yshift=-1pt]{$i$};
\draw[ghost](1.5,0)--++(0,1)node[above,yshift=-1pt]{$i$};
\draw[solid](1.25,0)node[below]{$j$}--++(0.75,0.5)--++(-0.75,0.5);
\end{tikzpicture}
+
\begin{tikzpicture}[anchorbase,smallnodes,rounded corners]
\draw[ghost,dot=0.05,dot](0.75,0)--++(0,0.2)--++(1,0.2)--++(0,0.2)--++(-1,0.2)--++(0,0.2)node[above,yshift=-1pt]{$i$};
\draw[ghost](1,0)--++(0,1)node[above,yshift=-1pt]{$i$};
\draw[ghost](1.5,0)--++(0,1)node[above,yshift=-1pt]{$i$};
\draw[solid](1.25,0)node[below]{$j$}--++(0.75,0.5)--++(-0.75,0.5);
\end{tikzpicture}
-
\begin{tikzpicture}[anchorbase,smallnodes,rounded corners]
\draw[ghost](0.75,0)--++(0,0.2)--++(1,0.2)--++(0,0.2)--++(-1,0.2)--++(0,0.2)node[above,yshift=-1pt]{$i$};
\draw[ghost,dot=0.95](1,0)--++(0,1)node[above,yshift=-1pt]{$i$};
\draw[ghost,dot](1.5,0)--++(0,1)node[above,yshift=-1pt]{$i$};
\draw[solid](1.25,0)node[below]{$j$}--++(0,1);
\end{tikzpicture}
+
\begin{tikzpicture}[anchorbase,smallnodes,rounded corners]
\draw[ghost,dot=0.05](0.75,0)--++(0,0.2)--++(1,0.2)--++(0,0.2)--++(-1,0.2)--++(0,0.2)node[above,yshift=-1pt]{$i$};
\draw[ghost](1,0)--++(0,1)node[above,yshift=-1pt]{$i$};
\draw[ghost,dot](1.5,0)--++(0,1)node[above,yshift=-1pt]{$i$};
\draw[solid](1.25,0)node[below]{$j$}--++(0,1);
\end{tikzpicture}
,
\\
i\Rrightarrow j\colon
\begin{tikzpicture}[anchorbase,smallnodes,rounded corners]
\draw[ghost,spinach](0.5,0)--++(0,1)node[above,yshift=-1pt]{$i$};
\draw[ghost](0.75,0)--++(0,1)node[above,yshift=-1pt]{$i$};
\draw[ghost](1,0)--++(0,1)node[above,yshift=-1pt]{$i$};
\draw[ghost](1.5,0)--++(0,1)node[above,yshift=-1pt]{$i$};
\draw[solid](1.25,0)node[below]{$j$}--++(0,1);
\end{tikzpicture}
=&
-
\begin{tikzpicture}[anchorbase,smallnodes,rounded corners]
\draw[ghost,dot=0.08,dot=0.165,dot](0.5,0)--++(0,0.2)--++(1.25,0.2)--++(0,0.2)--++(-1.25,0.2)--++(0,0.2)node[above,yshift=-1pt]{$i$};
\draw[ghost](0.75,0)--++(0,1)node[above,yshift=-1pt]{$i$};
\draw[ghost](1,0)--++(0,1)node[above,yshift=-1pt]{$i$};
\draw[ghost](1.5,0)--++(0,1)node[above,yshift=-1pt]{$i$};
\draw[solid](1.25,0)node[below]{$j$}--++(0.75,0.5)--++(-0.75,0.5);
\end{tikzpicture}
+
\begin{tikzpicture}[anchorbase,smallnodes,rounded corners]
\draw[ghost,dot=0.08,dot](0.5,0)--++(0,0.2)--++(1.25,0.2)--++(0,0.2)--++(-1.25,0.2)--++(0,0.2)node[above,yshift=-1pt]{$i$};
\draw[ghost](0.75,0)--++(0,1)node[above,yshift=-1pt]{$i$};
\draw[ghost,dot=0.9](1,0)--++(0,1)node[above,yshift=-1pt]{$i$};
\draw[ghost](1.5,0)--++(0,1)node[above,yshift=-1pt]{$i$};
\draw[solid](1.25,0)node[below]{$j$}--++(0.75,0.5)--++(-0.75,0.5);
\end{tikzpicture}
+
\begin{tikzpicture}[anchorbase,smallnodes,rounded corners]
\draw[ghost,dot=0.165,dot](0.5,0)--++(0,0.2)--++(1.25,0.2)--++(0,0.2)--++(-1.25,0.2)--++(0,0.2)node[above,yshift=-1pt]{$i$};
\draw[ghost](0.75,0)--++(0,1)node[above,yshift=-1pt]{$i$};
\draw[ghost,dot=0.9](1,0)--++(0,1)node[above,yshift=-1pt]{$i$};
\draw[ghost](1.5,0)--++(0,1)node[above,yshift=-1pt]{$i$};
\draw[solid](1.25,0)node[below]{$j$}--++(0.75,0.5)--++(-0.75,0.5);
\end{tikzpicture}
-
\begin{tikzpicture}[anchorbase,smallnodes,rounded corners]
\draw[ghost,dot=0.08,dot=0.165,dot](0.5,0)--++(0,0.2)--++(1.25,0.2)--++(0,0.2)--++(-1.25,0.2)--++(0,0.2)node[above,yshift=-1pt]{$i$};
\draw[ghost](0.75,0)--++(0,1)node[above,yshift=-1pt]{$i$};
\draw[ghost](1,0)--++(0,1)node[above,yshift=-1pt]{$i$};
\draw[ghost](1.5,0)--++(0,1)node[above,yshift=-1pt]{$i$};
\draw[solid](1.25,0)node[below]{$j$}--++(0.75,0.5)--++(-0.75,0.5);
\end{tikzpicture}
\\
&-
\begin{tikzpicture}[anchorbase,smallnodes,rounded corners]
\draw[ghost](0.5,0)--++(0,0.2)--++(1.25,0.2)--++(0,0.2)--++(-1.25,0.2)--++(0,0.2)node[above,yshift=-1pt]{$i$};
\draw[ghost,dot=0.9](0.75,0)--++(0,1)node[above,yshift=-1pt]{$i$};
\draw[ghost,dot=0.9](1,0)--++(0,1)node[above,yshift=-1pt]{$i$};
\draw[ghost,dot](1.5,0)--++(0,1)node[above,yshift=-1pt]{$i$};
\draw[solid](1.25,0)node[below]{$j$}--++(0,1);
\end{tikzpicture}
+
\begin{tikzpicture}[anchorbase,smallnodes,rounded corners]
\draw[ghost,dot=0.08](0.5,0)--++(0,0.2)--++(1.25,0.2)--++(0,0.2)--++(-1.25,0.2)--++(0,0.2)node[above,yshift=-1pt]{$i$};
\draw[ghost](0.75,0)--++(0,1)node[above,yshift=-1pt]{$i$};
\draw[ghost,dot=0.9](1,0)--++(0,1)node[above,yshift=-1pt]{$i$};
\draw[ghost,dot](1.5,0)--++(0,1)node[above,yshift=-1pt]{$i$};
\draw[solid](1.25,0)node[below]{$j$}--++(0,1);
\end{tikzpicture}
+
\begin{tikzpicture}[anchorbase,smallnodes,rounded corners]
\draw[ghost,dot=0.165](0.5,0)--++(0,0.2)--++(1.25,0.2)--++(0,0.2)--++(-1.25,0.2)--++(0,0.2)node[above,yshift=-1pt]{$i$};
\draw[ghost,dot=0.9](0.75,0)--++(0,1)node[above,yshift=-1pt]{$i$};
\draw[ghost](1,0)--++(0,1)node[above,yshift=-1pt]{$i$};
\draw[ghost,dot](1.5,0)--++(0,1)node[above,yshift=-1pt]{$i$};
\draw[solid](1.25,0)node[below]{$j$}--++(0,1);
\end{tikzpicture}
-
\begin{tikzpicture}[anchorbase,smallnodes,rounded corners]
\draw[ghost,dot=0.08,dot=0.165](0.5,0)--++(0,0.2)--++(1.25,0.2)--++(0,0.2)--++(-1.25,0.2)--++(0,0.2)node[above,yshift=-1pt]{$i$};
\draw[ghost](0.75,0)--++(0,1)node[above,yshift=-1pt]{$i$};
\draw[ghost](1,0)--++(0,1)node[above,yshift=-1pt]{$i$};
\draw[ghost,dot](1.5,0)--++(0,1)node[above,yshift=-1pt]{$i$};
\draw[solid](1.25,0)node[below]{$j$}--++(0,1);
\end{tikzpicture}
.
\end{align*}
The dotted plactic relations, and their partner relations,
\begin{gather*}%\label{L:RecollectionReidemeisterIII}
i\Rightarrow j\colon
\begin{tikzpicture}[anchorbase,smallnodes,rounded corners]
\draw[ghost,spinach,dot](1,0)--++(0,1)node[above,yshift=-1pt]{$i$};
\draw[ghost](1.5,0)--++(0,1)node[above,yshift=-1pt]{$i$};
\draw[solid](1.25,0)node[below]{$j$}--++(0,1);
\end{tikzpicture}
=
\begin{tikzpicture}[anchorbase,smallnodes,rounded corners]
\draw[ghost](1,0)--++(0,1)node[above,yshift=-1pt]{$i$};
\draw[ghost,dot](1.5,0)--++(0,1)node[above,yshift=-1pt]{$i$};
\draw[solid](1.25,0)node[below]{$j$}--++(0,1);
\end{tikzpicture}
-
\begin{tikzpicture}[anchorbase,smallnodes,rounded corners]
\draw[ghost](1,0)--++(0.75,0.5)--++(-0.75,0.5)node[above,yshift=-1pt]{$i$};
\draw[ghost,dot](1.5,0)--++(0,1)node[above,yshift=-1pt]{$i$};
\draw[solid](1.25,0)node[below]{$j$}--++(0,1);
\end{tikzpicture}
-
\begin{tikzpicture}[anchorbase,smallnodes,rounded corners]
\draw[ghost,dot](1,0)--++(0.75,0.5)--++(-0.75,0.5)node[above,yshift=-1pt]{$i$};
\draw[ghost](1.5,0)--++(0,1)node[above,yshift=-1pt]{$i$};
\draw[solid](1.25,0)node[below]{$j$}--++(0.75,0.5)--++(-0.75,0.5);
\end{tikzpicture}
,\quad
i\Rrightarrow j\colon
\begin{tikzpicture}[anchorbase,smallnodes,rounded corners]
\draw[ghost,spinach,dot=0.33,dot=0.66](1,0)--++(0,1)node[above,yshift=-1pt]{$i$};
\draw[ghost](1.5,0)--++(0,1)node[above,yshift=-1pt]{$i$};
\draw[solid](1.25,0)node[below]{$j$}--++(0,1);
\end{tikzpicture}
=
\begin{tikzpicture}[anchorbase,smallnodes,rounded corners]
\draw[ghost](1,0)--++(0,1)node[above,yshift=-1pt]{$i$};
\draw[ghost,dot=0.33,dot=0.66](1.5,0)--++(0,1)node[above,yshift=-1pt]{$i$};
\draw[solid](1.25,0)node[below]{$j$}--++(0,1);
\end{tikzpicture}
+
\begin{tikzpicture}[anchorbase,smallnodes,rounded corners]
\draw[ghost,dot](1,0)--++(0,1)node[above,yshift=-1pt]{$i$};
\draw[ghost,dot](1.5,0)--++(0,1)node[above,yshift=-1pt]{$i$};
\draw[solid](1.25,0)node[below]{$j$}--++(0,1);
\end{tikzpicture}
-
\begin{tikzpicture}[anchorbase,smallnodes,rounded corners]
\draw[ghost](1,0)--++(0.75,0.5)--++(-0.75,0.5)node[above,yshift=-1pt]{$i$};
\draw[ghost,dot](1.5,0)--++(0,1)node[above,yshift=-1pt]{$i$};
\draw[solid](1.25,0)node[below]{$j$}--++(0,1);
\end{tikzpicture}
-
\begin{tikzpicture}[anchorbase,smallnodes,rounded corners]
\draw[ghost,dot](1,0)--++(0.75,0.5)--++(-0.75,0.5)node[above,yshift=-1pt]{$i$};
\draw[ghost](1.5,0)--++(0,1)node[above,yshift=-1pt]{$i$};
\draw[solid](1.25,0)node[below]{$j$}--++(0.75,0.5)--++(-0.75,0.5);
\end{tikzpicture}
,
\end{gather*}
hold as well.
\end{Lemma}

\begin{proof}
We only prove one of these identities, namely the plactic
relation for a doubly laced edge. All of the other
relations can be proven similarly.
Using \autoref{L:RecollectionMovingStringsDots} we compute
\begin{gather*}
\begin{tikzpicture}[anchorbase,smallnodes,rounded corners]
\draw[ghost,spinach](0.75,0)--++(0,1)node[above,yshift=-1pt]{$i$};
\draw[ghost](1,0)--++(0,1)node[above,yshift=-1pt]{$i$};
\draw[ghost](1.5,0)--++(0,1)node[above,yshift=-1pt]{$i$};
\draw[solid](1.25,0)node[below]{$j$}--++(0,1);
\end{tikzpicture}
=
\begin{tikzpicture}[anchorbase,smallnodes,rounded corners]
\draw[ghost,dot](0.75,0)--++(0.4,0.5)--++(-0.4,0.5)node[above,yshift=-1pt]{$i$};
\draw[ghost,dot=0.9](1,0)--++(0,1)node[above,yshift=-1pt]{$i$};
\draw[ghost](1.5,0)--++(0,1)node[above,yshift=-1pt]{$i$};
\draw[solid](1.25,0)node[below]{$j$}--++(0,1);
\end{tikzpicture}
-
\begin{tikzpicture}[anchorbase,smallnodes,rounded corners]
\draw[ghost,dot=0.1,dot](0.75,0)--++(0.4,0.5)--++(-0.4,0.5)node[above,yshift=-1pt]{$i$};
\draw[ghost](1,0)--++(0,1)node[above,yshift=-1pt]{$i$};
\draw[ghost](1.5,0)--++(0,1)node[above,yshift=-1pt]{$i$};
\draw[solid](1.25,0)node[below]{$j$}--++(0,1);
\end{tikzpicture}
.
\end{gather*}
Now, the Reidemeister III relation, {\ie}
\autoref{Eq:RecollectionReidemeisterIII}, for $i\Rightarrow j$
and our choice of $Q$-polynomials gives
the following left-hand equation, while
\autoref{Eq:RecollectionDotCrossing} gives the right-hand relation, again including their partner relations:
\begin{gather*}
\begin{tikzpicture}[anchorbase,smallnodes,rounded corners]
\draw[ghost,dot](6.2,1)node[above,yshift=-1pt]{$i$}--++(0,-1)node[below]{$\phantom{i}$};
\draw[ghost](7.2,1)node[above,yshift=-1pt]{$i$}--++(0,-1)node[below]{$\phantom{i}$};
\draw[solid](6.7,1)--++(0,-1)node[below]{$j$};
\end{tikzpicture}
=
\begin{tikzpicture}[anchorbase,smallnodes,rounded corners]
\draw[ghost](3,1)node[above,yshift=-1pt]{$\phantom{i}$}--++(1,-1)node[below]{$\phantom{i}$};
\draw[ghost](4,1)node[above,yshift=-1pt]{$i$}--++(-1,-1)node[below]{$\phantom{i}$};
\draw[solid,smallnodes,rounded corners](3.5,1)--++(0.5,-0.5)--++(-0.5,-0.5)node[below]{$j$};
\end{tikzpicture}
-
\begin{tikzpicture}[anchorbase,smallnodes,rounded corners]
\draw[ghost](1,1)node[above,yshift=-1pt]{$i$}--++(1,-1)node[below]{$\phantom{i}$};
\draw[ghost](2,1)node[above,yshift=-1pt]{$i$}--++(-1,-1)node[below]{$\phantom{i}$};
\draw[solid,smallnodes,rounded corners](1.5,1)--++(-0.5,-0.5)--++(0.5,-0.5)node[below]{$j$};
\end{tikzpicture}
-
\begin{tikzpicture}[anchorbase,smallnodes,rounded corners]
\draw[ghost](6.2,1)node[above,yshift=-1pt]{$i$}--++(0,-1)node[below]{$\phantom{i}$};
\draw[ghost,dot](7.2,1)node[above,yshift=-1pt]{$i$}--++(0,-1)node[below]{$\phantom{i}$};
\draw[solid](6.7,1)--++(0,-1)node[below]{$j$};
\end{tikzpicture}
,\quad
\begin{tikzpicture}[anchorbase,smallnodes,rounded corners]
\draw[ghost](0,1)node[above,yshift=-1pt]{$i$}--++(0.5,-1) node[below]{\phantom{$i$}};
\draw[ghost](0.5,1)node[above,yshift=-1pt]{$i$}--++(-0.5,-1);
\end{tikzpicture}
=
-
\begin{tikzpicture}[anchorbase,smallnodes,rounded corners]
\draw[ghost,dot](0,1)node[above,yshift=-1pt]{$i$}--++(0.5,-0.5)--++(-0.5,-0.5) node[below]{\phantom{$i$}};
\draw[ghost](0.5,1)node[above,yshift=-1pt]{$i$}--++(-0.5,-0.5)--++(0.5,-0.5);
\end{tikzpicture}
=
\begin{tikzpicture}[anchorbase,smallnodes,rounded corners]
\draw[ghost](0,1)node[above,yshift=-1pt]{$i$}--++(0.5,-0.5)--++(-0.5,-0.5) node[below]{\phantom{$i$}};
\draw[ghost,dot](0.5,1)node[above,yshift=-1pt]{$i$}--++(-0.5,-0.5)--++(0.5,-0.5) node[below]{$\phantom{i}$};
\end{tikzpicture}
.
\end{gather*}
Using these two relations to rewrite the relation in the preceding displayed equation completes the proof.
\end{proof}

\begin{Remark}\label{R:RecollectionPlactic}
In Lie theory, the plactic relations are the image of the Serre relations in the crystals, and the relations in \autoref{L:RecollectionReidemeisterIII} mimic these as follows.
Cutting the pictures in \autoref{L:RecollectionReidemeisterIII} along their equator gives linear relations between different words in the residues $\vertices$.
More precisely, we obtain relations between the following words:
\begin{gather*}
i\not\!\rightsquigarrow j\colon\underline{ij},ji
,\quad
i\rightarrow j\colon\overline{iij},\underline{iji},jii
,\quad
i\Rightarrow j\colon\overline{iiij},\underline{iiji},ijii,jiii
,\quad
i\Rrightarrow j\colon\overline{iiiij},\underline{iiiji},iijii,ijiii,jiiii.
\end{gather*}
These relations express the underlined words in
terms of the other words, all of which have $i$-strings further to the right. It is clear that the $i$-strings are further to the right in all of these words, except for the overlined
words for which we can use \autoref{L:RecollectionMovingStringsDots} to create dots. The extra dots again allow us to pull an $i$-string further to the right.
\end{Remark}

%%%%%%%%%%%%%%%%%%%%%%%%%%%%%%%%%%%%%%%%%

\section{Main example}\label{S:MainExample}

%%%%%%%%%%%%%%%%%%%%%%%%%%%%%%%%%%%%%%%%%

We now give an example that emphasizes the ideas underpinning our  construction of the sandwich cellular bases.
Some of the terms we use have not been defined yet but, hopefully, the meaning will be clear. We encourage the reader to come back to this example while reading the definitions in \autoref{S:FiniteTypes}.

%%%%%%%%%%%%%%%%%%%%%%%%%%%%%%%%%%%%%%%%%

\subsection{Strategic interlude}\label{SS:MainExampleOne}

%%%%%%%%%%%%%%%%%%%%%%%%%%%%%%%%%%%%%%%%%

%\begin{Remark}\label{R:MainExamplePicture}
Similarly to \cite{MaTu-klrw-algebras} and
\cite{MaTu-klrw-algebras-bad}, and also
\autoref{R:SandwichCellularAlgebra}, the picture
to keep in mind for the sandwich cellular
basis elements $D_{\bS\bT}^{\ba,\bfi}$ that we construct is:
\begin{gather}\label{Eq:MainExampleSandwich}
D_{\bS\bT}^{\ba,\bfi}
\leftrightsquigarrow
\begin{tikzpicture}[anchorbase,scale=1]
\draw[line width=0.75,color=black,fill=cream] (0,-0.5) to (0.25,0) to (0.75,0) to (1,-0.5) to (0,-0.5);
\node at (0.5,-0.25){$\bT$};
\draw[line width=0.75,color=black,fill=cream] (0,1) to (0.25,0.5) to (0.75,0.5) to (1,1) to (0,1);
\node at (0.5,0.75){$\bS$};
\draw[line width=0.75,color=black,fill=cream] (0.25,0) to (0.25,0.5) to (0.75,0.5) to (0.75,0) to (0.25,0);
\node at (0.51,0.25){\scalebox{0.65}{$\sandwich{\blam}{\ba}{\bfi}$}};
\end{tikzpicture}
=
\begin{tikzpicture}[anchorbase,scale=1]
\draw[line width=0.75,color=black,fill=cream] (0,-0.5) to (0.25,0) to (0.75,0) to (1,-0.5) to (0,-0.5);
\node at (0.5,-0.25){$\bT$};
\draw[line width=0.75,color=black,fill=cream] (0,2) to (0.25,1.5) to (0.75,1.5) to (1,2) to (0,2);
\node at (0.5,1.75){$\bS$};
\draw[line width=0.75,color=black,fill=cream] (0.25,0) to (0.25,0.5) to (0.75,0.5) to (0.75,0) to (0.25,0);
\node at (0.5,0.25){$\dotidem$};
\draw[line width=0.75,color=black,fill=cream] (0.25,0.5) to (0.25,1) to (0.75,1) to (0.75,0.5) to (0.25,0.5);
\node at (0.5,0.75){$\zeetwo$};
\draw[line width=0.75,color=black,fill=cream] (0.25,1) to (0.25,1.5) to (0.75,1.5) to (0.75,1) to (0.25,1);
\node at (0.5,1.25){$\daffine$};
\end{tikzpicture},
\quad\text{where}\quad
\begin{aligned}
\begin{tikzpicture}[anchorbase,scale=1]
\draw[line width=0.75,color=black,fill=cream] (0,1) to (0.25,0.5) to (0.75,0.5) to (1,1) to (0,1);
\node at (0.5,0.75){$\bS$};
\end{tikzpicture}
&\text{ a permutation diagram,}
\\
\begin{tikzpicture}[anchorbase,scale=1]
\draw[white,ultra thin] (0,0) to (1,0);
\draw[line width=0.75,color=black,fill=cream] (0.25,0) to (0.25,0.5) to (0.75,0.5) to (0.75,0) to (0.25,0);
\node at (0.5,0.25){$\daffine$};
\end{tikzpicture}
&\text{ an ``affine/unsteady'' dot placement,}
\\
\begin{tikzpicture}[anchorbase,scale=1]
\draw[white,ultra thin] (0,0) to (1,0);
\draw[line width=0.75,color=black,fill=cream] (0.25,0) to (0.25,0.5) to (0.75,0.5) to (0.75,0) to (0.25,0);
\node at (0.5,0.25){$\zeetwo$};
\end{tikzpicture}
&\text{ a ``finite/steady'' dot placement,}
\\
\begin{tikzpicture}[anchorbase,scale=1]
\draw[white,ultra thin] (0,0) to (1,0);
\draw[line width=0.75,color=black,fill=cream] (0.25,0) to (0.25,0.5) to (0.75,0.5) to (0.75,0) to (0.25,0);
\node at (0.5,0.25){$\dotidem$};
\end{tikzpicture}
&\text{ a dotted idempotent,}
\\
\begin{tikzpicture}[anchorbase,scale=1]
\draw[line width=0.75,color=black,fill=cream] (0,-0.5) to (0.25,0) to (0.75,0) to (1,-0.5) to (0,-0.5);
\node at (0.5,-0.25){$\bT$};
\end{tikzpicture}
&\text{ a permutation diagram.}
\end{aligned}
\end{gather}
As our terminology suggests, the steady and unsteady dot placements will correspond to steady and unsteady diagrams. The middle of the diagram is given by a \emph{(dotted) idempotent diagram}, possibly with sandwiched dots $\daffine\zeetwo$,
and the bottom and top are given by
\emph{permutation diagrams}. Examples of these types of diagrams are:
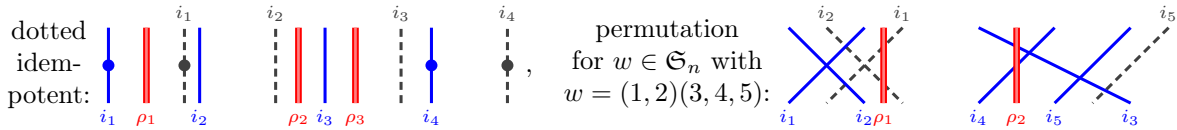
\begin{gather*}
\begin{gathered}
\text{dotted}
\\[-0.1cm]
\text{idem-}
\\[-0.1cm]
\text{potent:}
\end{gathered}
\begin{tikzpicture}[anchorbase,smallnodes,rounded corners]
\draw[ghost,dot](0,0)node[below]{$\phantom{i}$}--++(0,1)node[above,yshift=-1pt]{$i_{1}$};
\draw[ghost](1.2,0)node[below]{$\phantom{i}$}--++(0,1)node[above,yshift=-1pt]{$i_{2}$};
\draw[ghost](2.85,0)node[below]{$\phantom{i}$}--++(0,1)node[above,yshift=-1pt]{$i_{3}$};
\draw[ghost,dot](4.25,0)node[below]{$\phantom{i}$}--++(0,1)node[above,yshift=-1pt]{$i_{4}$};
\draw[solid,dot](-1,0)node[below]{$i_{1}$}--++(0,1) node[above,yshift=-1pt]{$\phantom{i_{1}}$};
\draw[solid](0.2,0)node[below]{$i_{2}$}--++(0,1);
\draw[solid](1.85,0)node[below]{$i_{3}$}--++(0,1);
\draw[solid,dot](3.25,0)node[below]{$i_{4}$}--++(0,1);
\draw[redstring](-0.5,0)node[below]{$\rho_{1}$}--++(0,1);
\draw[redstring](1.5,0)node[below]{$\rho_{2}$}--++(0,1);
\draw[redstring](2.25,0)node[below]{$\rho_{3}$}--++(0,1);
\end{tikzpicture}
,\quad
\begin{gathered}
\text{permutation}
\\[-0.1cm]
\text{for $w\in\sym$ with}
\\[-0.1cm]
\text{$w=(1,2)(3,4,5)$:}
\end{gathered}
\begin{tikzpicture}[anchorbase,smallnodes,rounded corners]
\draw[ghost](1.5,1)node[above,yshift=-1pt]{$i_{2}$}--++(1,-1);
\draw[ghost](2.5,1)node[above,yshift=-1pt]{$i_{1}$}--++(-1,-1);
\draw[ghost](6,1)node[above,yshift=-1pt]{$i_{5}$}--++(-1,-1)node[below]{$\phantom{i}$};
\draw[solid](1,1)node[above,yshift=-1pt]{$\phantom{i}$}--++(1,-1)node[below]{$i_{2}$};
\draw[solid](2,1)--++(-1,-1)node[below]{$i_{1}$};
\draw[solid](3.5,1)node[above,yshift=-1pt]{$\phantom{i}$}--++(2,-1)node[below]{$i_{3}$};
\draw[solid](4.5,1)node[above,yshift=-1pt]{$\phantom{i}$}--++(-1,-1)node[below]{$i_{4}$};
\draw[solid](5.5,1)node[above,yshift=-1pt]{$\phantom{i}$}--++(-1,-1)node[below]{$i_{5}$};
\draw[redstring](2.25,1)--++(0,-1)node[below]{$\rho_{1}$};
\draw[redstring](4,1)--++(0,-1)node[below]{$\rho_{2}$};
\end{tikzpicture}
\end{gather*}
Here, and throughout, we let $\sym=\mathrm{Aut}(\set{1,\dots,n})$ be the symmetric group on $\set{1,\dots,n}$, where we use the
cycle notation for its elements.
%\end{Remark}

%\begin{Remark}\label{SS:MainExampleOne}
As we have said, our basic strategy is to
\emph{pull strings and jump dots to the right} using the results of \autoref{SS:RecollectionMoveRight}.

As in \cite{MaTu-klrw-algebras} and \cite{MaTu-klrw-algebras-bad},
an important observation is that certain
Reidemeister II relations do not hold, which stops us from pulling strings to the right. For example,
assume that $i\Rightarrow j$. Then \autoref{L:RecollectionMovingStringsDots}.(b)
tells us that we can pull
a ghost $i$-string with one dot to the right,
or jump the dot on a ghost $i$-string to the right. However, none of the relations allow us to pull a ghost $i$-string without dot,
or a solid $j$-string with one dot, to the right, so these strings are \emph{blocked}. It is possible that some of the other relations allow us to pull these strings further to the right but we eventually see that it is not possible to pull these strings further to right.

Generalizing the simple observation from the last paragraph to all residues leads us to our main strategy
for constructing idempotents: we place the solid strings and their ghosts inductively by putting the new strings on the left of the diagram and then pull them to the right until they are blocked
by another string for which there is no
Reidemeister II relation. Using categorification, we show that the residue sequences for the idempotent diagrams are given by paths in the corresponding crystal graph, so the idempotent diagrams are naturally indexed by vertices of the crystal graph.
We will see that the idempotent diagrams in our sandwich cellular basis are maximal with respect to
a total order $\sandorder$, which measures how far strings are placed to the right.

The inductive process used to construct idempotent diagrams can be thought of as follows.
For the idempotent diagram that we are constructing, we assume that the first $(k-1)$ solid strings, and their ghosts, have already been placed. We then put the $k$th
solid string, and any ghosts, on the far left of the diagram and pull these strings to the right using isotopies and honest Reidemeister II moves. We emphasize that the solid strings and their ghosts have to be pulled at the same time, because the ghost strings are shifts of the solid string. These strings
will eventually become blocked, by either an affine red string or by a previously placed solid or ghost string, which becomes the final position of these strings. This process
ensures that idempotents in the cellular basis are maximal with respect to $\sandorder$. Pulling strings to the right in this way also allows us to argue by induction by working modulo higher order terms with respect to $\sandorder$.

After the idempotent diagram is constructed we need to decorate the idempotent
with dots whenever there are close strings
with repeated residues. This is necessary because we want to flank
diagrams with crossings, {\cf} \autoref{Eq:MainExampleSandwich}, but \autoref{Eq:RecollectionReidemeisterII}
annihilates two close solid $i$-strings whenever they are flanked with crossings.

In addition, it is sometimes possible to add more dots to a
dotted idempotent diagram in such a way that the strings are still blocked. For example, the partner relations of the relations in
\autoref{L:RecollectionMovingStringsDots}.(b) only allow strings to be pulled to the right when they have enough dots. This means that these strings are blocked if we add one or two dots, with the number of allowed dots depending on the quiver. These additional dot placements define the middle, or the sandwiched part, of our basis. In contrast, for the strings in the unsteady part of these diagrams we allow arbitrarily many dots, which corresponds to the sandwich algebra being a polynomial algebra.

Permutations diagrams can then be added to the bottom and top of diagrams, as
it is very common in the KLR and wKLRW world.
%\end{Remark}

%%%%%%%%%%%%%%%%%%%%%%%%%%%%%%%%%%%%%%%%%

\subsection{The main example}\label{SS:MainExampleTwo}

%%%%%%%%%%%%%%%%%%%%%%%%%%%%%%%%%%%%%%%%%

\begin{Example}\label{E:MainExampleTheBeastItself}
We will now combine all of the ingredients above in one example.
In this example we take a quiver $\quiver$ of type $\typeb[3]$, that is
$\dynkin B3$ with vertices $1$, $2$ and $3$, read left to right,
with the additional orientation $1\rightarrow 2$. The symmetrizer for this quiver is $\mathbf{d}=(2,2,1)$. We take the crystal graph
$\crystalgraph[{\fweight[2]}]$ for the second fundamental weight.
Fix $n=5$, $\ell=1$, $\charge=(0)$ and $\brho=(2)$, and set all
ghost shifts to be $1$. We will describe how we will construct
a homogeneous sandwich cellular basis for $\WAc(X)$.

\begin{enumerate}

\item \textbf{Residue sequences and idempotents.}
The root of the crystal graph $\crystalgraph[{\fweight[2]}]$ is the highest weight vertex of weight
$\fweight[2]$. The graph on the left-hand side of \autoref{Eq:MainExample} below
is the crystal graph of $\fweight[2]$.

For each vertex $\sigma$ of distance $n$ in $\crystalgraph[{\fweight[2]}]$ we fix a rooted path $\qpath[\sigma]$ from $\fweight[2]$ to $\sigma$. The residue sequence $\Res[{\qpath[\sigma]}]=(i_{1},\dots,i_{n})\in\vertices^{n}$ of
such a path $\qpath[\sigma]$ is the sequence of colors on the edges of the path. By standard properties of crystal graphs, a path is uniquely determined by its residue sequence. To each
vertex $\sigma$ we associate an idempotent diagram
$\idem[\sigma]$ using $\Res[{\qpath[\sigma]}]$ by pulling the strings to the right, by the inductive procedure outlined in \autoref{SS:MainExampleOne}.

For this example, we have chosen the highlighted paths $\ppath[\lambda]$, $\ppath[\mu]$ and $\ppath[\nu]$ as our \emph{preferred paths}. Therefore,
$\Res[{\ppath[\lambda]}]=23312$, which is shorthand notation for
$\Res[{\ppath[\lambda]}]=(2,3,3,1,2)$, $\Res[{\ppath[\mu]}]=23321$ and $\Res[{\ppath[\nu]}]=21323$,
as illustrated in \autoref{Eq:MainExample} below. Note that there are several paths in the crystal to each of these vertices and, for our purposes it does not matter which paths we choose.

Let us consider $\mu$ first. We order the five solid strings, from left to right, using $\Res[{\ppath[\mu]}]$. That is, in order, we have a solid
$2$-string, two solid $3$-strings, a solid $2$-string and a solid $1$-string.
Our strategy from \autoref{SS:MainExampleOne} for placing these strings in order, and then pulling them to the right, gives the following sequence of diagrams:
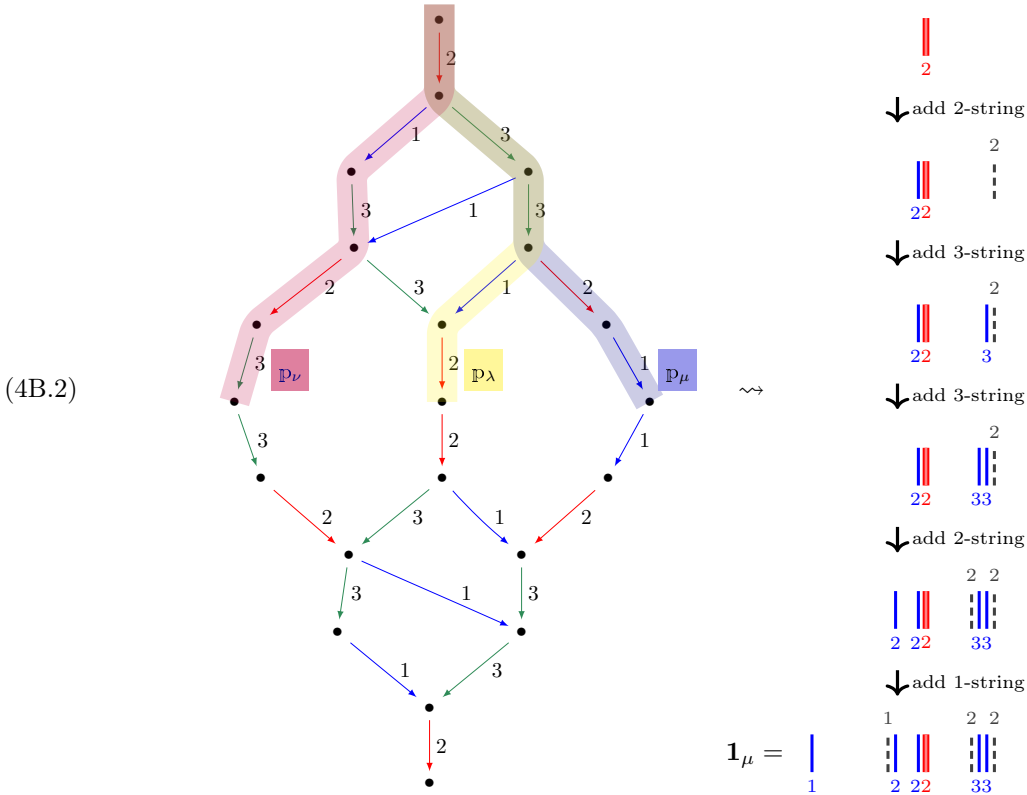
\begin{gather}\label{Eq:MainExample}
\scalebox{0.75}{$\begin{tikzpicture}[anchorbase,>=latex,line join=bevel,xscale=0.65,yscale=0.5]
\node (node_0) at (69.5bp,439.0bp) [draw,draw=none] {$\bullet$};
\node (node_6) at (53.5bp,367.0bp) [draw,draw=none] {$\bullet$};
\node (node_1) at (137.5bp,583.0bp) [draw,draw=none] {$\bullet$};
\node (node_9) at (139.5bp,511.0bp) [draw,draw=none] {$\bullet$};
\node (node_2) at (127.5bp,151.0bp) [draw,draw=none] {$\bullet$};
\node (node_8) at (193.5bp,79.5bp) [draw,draw=none] {$\bullet$};
\node (node_3) at (264.5bp,511.0bp) [draw,draw=none] {$\bullet$};
\node (node_17) at (320.5bp,439.0bp) [draw,draw=none] {$\bullet$};
\node (node_18) at (202.5bp,439.0bp) [draw,draw=none] {$\bullet$};
\node (node_4) at (259.5bp,223.0bp) [draw,draw=none] {$\bullet$};
\node (node_16) at (259.5bp,151.0bp) [draw,draw=none] {$\bullet$};
\node (node_5) at (72.5bp,295.0bp) [draw,draw=none] {$\bullet$};
\node (node_14) at (135.5bp,223.0bp) [draw,draw=none] {$\bullet$};
\node (node_7) at (321.5bp,295.0bp) [draw,draw=none] {$\bullet$};
\node (node_13) at (193.5bp,8.5bp) [draw,draw=none] {$\bullet$};
\node (node_10) at (200.5bp,725.5bp) [draw,draw=none] {$\bullet$};
\node (node_12) at (200.5bp,654.5bp) [draw,draw=none] {$\bullet$};
\node (node_11) at (202.5bp,295.0bp) [draw,draw=none] {$\bullet$};
\node (node_19) at (264.5bp,583.0bp) [draw,draw=none] {$\bullet$};
\node (node_15) at (351.5bp,367.0bp) [draw,draw=none] {$\bullet$};
\node (node_20) at (202.5bp,367.0bp) [draw,draw=none] {$\bullet$};
\draw [spinach,->] (node_0) ..controls (65.255bp,419.43bp) and (60.953bp,400.61bp)  .. (node_6);
\draw (72.0bp,403.0bp) node {$3$};
\draw [spinach,->] (node_1) ..controls (138.02bp,563.68bp) and (138.57bp,544.46bp)  .. (node_9);
\draw (148.0bp,547.0bp) node {$3$};
\draw [blue,->] (node_2) ..controls (145.48bp,131.07bp) and (165.5bp,109.99bp)  .. (node_8);
\draw (176.0bp,115.0bp) node {$1$};
\draw [red,->] (node_3) ..controls (279.67bp,491.04bp) and (296.42bp,470.1bp)  .. (node_17);
\draw (307.0bp,475.0bp) node {$2$};
\draw [blue,->] (node_3) ..controls (247.61bp,490.93bp) and (228.8bp,469.7bp)  .. (node_18);
\draw (249.0bp,475.0bp) node {$1$};
\draw [spinach,->] (node_4) ..controls (259.5bp,203.68bp) and (259.5bp,184.46bp)  .. (node_16);
\draw (268.0bp,187.0bp) node {$3$};
\draw [red,->] (node_5) ..controls (89.778bp,274.8bp) and (108.19bp,254.34bp)  .. (node_14);
\draw (120.0bp,259.0bp) node {$2$};
\draw [spinach,->] (node_6) ..controls (58.569bp,347.32bp) and (63.751bp,328.23bp)  .. (node_5);
\draw (74.0bp,331.0bp) node {$3$};
\draw [red,->] (node_7) ..controls (304.61bp,274.93bp) and (285.8bp,253.7bp)  .. (node_4);
\draw (306.0bp,259.0bp) node {$2$};
\draw [red,->] (node_8) ..controls (193.5bp,60.442bp) and (193.5bp,41.496bp)  .. (node_13);
\draw (202.0bp,44.0bp) node {$2$};
\draw [red,->] (node_9) ..controls (120.2bp,490.7bp) and (99.457bp,469.96bp)  .. (node_0);
\draw (122.0bp,475.0bp) node {$2$};
\draw [spinach,->] (node_9) ..controls (156.99bp,490.56bp) and (175.99bp,469.46bp)  .. (node_18);
\draw (186.0bp,475.0bp) node {$3$};
\draw [red,->] (node_10) ..controls (200.5bp,706.44bp) and (200.5bp,687.5bp)  .. (node_12);
\draw (209.0bp,690.0bp) node {$2$};
\draw [blue,->] (node_11) ..controls (214.23bp,277.42bp) and (225.23bp,262.33bp)  .. (235.5bp,250.0bp) .. controls (238.51bp,246.39bp) and (241.86bp,242.61bp)  .. (node_4);
\draw (244.0bp,259.0bp) node {$1$};
\draw [spinach,->] (node_11) ..controls (184.26bp,274.95bp) and (163.98bp,253.75bp)  .. (node_14);
\draw (185.0bp,259.0bp) node {$3$};
\draw [blue,->] (node_12) ..controls (183.34bp,634.57bp) and (164.23bp,613.49bp)  .. (node_1);
\draw (184.0bp,619.0bp) node {$1$};
\draw [spinach,->] (node_12) ..controls (217.71bp,634.81bp) and (236.53bp,614.37bp)  .. (node_19);
\draw (248.0bp,619.0bp) node {$3$};
\draw [spinach,->] (node_14) ..controls (133.35bp,203.2bp) and (131.13bp,183.79bp)  .. (node_2);
\draw (141.0bp,187.0bp) node {$3$};
\draw [blue,->] (node_14) ..controls (170.82bp,202.06bp) and (210.64bp,179.58bp)  .. (node_16);
\draw (220.0bp,187.0bp) node {$1$};
\draw [blue,->] (node_15) ..controls (343.35bp,346.99bp) and (334.79bp,327.02bp)  .. (node_7);
\draw (348.0bp,331.0bp) node {$1$};
\draw [spinach,->] (node_16) ..controls (240.94bp,130.45bp) and (221.46bp,109.94bp)  .. (node_8);
\draw (242.0bp,115.0bp) node {$3$};
\draw [blue,->] (node_17) ..controls (328.71bp,419.47bp) and (337.46bp,399.69bp)  .. (node_15);
\draw (348.0bp,403.0bp) node {$1$};
\draw [red,->] (node_18) ..controls (202.5bp,419.68bp) and (202.5bp,400.46bp)  .. (node_20);
\draw (211.0bp,403.0bp) node {$2$};
\draw [spinach,->] (node_19) ..controls (264.5bp,563.2bp) and (264.5bp,543.79bp)  .. (node_3);
\draw (273.0bp,547.0bp) node {$3$};
\draw [blue,->] (node_19) ..controls (228.89bp,562.06bp) and (188.75bp,539.58bp)  .. (node_9);
\draw (225.0bp,547.0bp) node {$1$};
\draw [red,->] (node_20) ..controls (202.5bp,347.2bp) and (202.5bp,327.79bp)  .. (node_11);
\draw (211.0bp,331.0bp) node {$2$};
\draw[black!50!blue,line width=0.5cm,opacity=0.2,rounded corners] (200.5bp,740.5bp) to (200.5bp,654.5bp) to (264.5bp,583.0bp) to (264.5bp,511.0bp) to (320.5bp,439.0bp) to (351.5bp,367.0bp);
\node[black] at (371.5bp,397.0bp) {\colorbox{black!20!blue!40!white}{\mystrut$\ppath[\mu]$}};
\draw[yellow,line width=0.5cm,opacity=0.2,rounded corners] (200.5bp,740.5bp) to (200.5bp,654.5bp) to (264.5bp,583.0bp) to (264.5bp,511.0bp) to (202.5bp,439.0bp) to (202.5bp,367.0bp);
\node[black] at (232.5bp,397.0bp) {\colorbox{yellow!50}{\mystrut$\ppath[\lambda]$}};
\draw[purple,line width=0.5cm,opacity=0.2,rounded corners] (200.5bp,740.5bp) to (200.5bp,654.5bp) to (137.5bp,583.0bp) to (139.5bp,511.0bp) to (69.5bp,439.0bp) to (53.5bp,367.0bp);
\node[black!50!blue] at (93.5bp,397.0bp) {\colorbox{purple!50}{\mystrut$\ppath[\nu]$}};
\end{tikzpicture}$}
\quad\rightsquigarrow
\scalebox{0.9}{$\begin{tikzcd}[ampersand replacement=\&,column sep=0.15em,row sep=1em]
\kern-2mm
\begin{tikzpicture}[anchorbase,smallnodes,yscale=0.5]
\draw[white](-1.65,0)--++(0,1)node[above,yshift=-1pt]{\phantom{$1$}};
\draw[white](0.7,0)--++(0,1)node[above,yshift=-1pt]{\phantom{$1$}};
\draw[redstring](0,0)node[below]{$2$}--++(0,1);
\end{tikzpicture}
\ar[d,very thick,->,"\text{ add $2$-string}"]
\\
\begin{tikzpicture}[anchorbase,smallnodes,yscale=0.5]
\draw[white](-1.65,0)--++(0,1)node[above,yshift=-1pt]{\phantom{$1$}};
\draw[white](0.7,0)--++(0,1)node[above,yshift=-1pt]{\phantom{$1$}};
\draw[ghost](0.9,0)--++(0,1)node[above,yshift=-1pt]{$2$};
\draw[solid](-0.1,0)node[below]{$2\;$}--++(0,1);
\draw[redstring](0,0)node[below]{$2$}--++(0,1);
\end{tikzpicture}
\ar[d,very thick,->,"\text{ add $3$-string}"]
\\
\begin{tikzpicture}[anchorbase,smallnodes,yscale=0.5]
\draw[white](-1.65,0)--++(0,1)node[above,yshift=-1pt]{\phantom{$1$}};
\draw[white](0.7,0)--++(0,1)node[above,yshift=-1pt]{\phantom{$1$}};
\draw[ghost](0.9,0)--++(0,1)node[above,yshift=-1pt]{$2$};
\draw[solid](-0.1,0)node[below]{$2\;$}--++(0,1);
\draw[solid](0.8,0)node[below]{$3$}--++(0,1);
\draw[redstring](0,0)node[below]{$2$}--++(0,1);
\end{tikzpicture}
\ar[d,very thick,->,"\text{ add $3$-string}"]
\\
\begin{tikzpicture}[anchorbase,smallnodes,yscale=0.5]
\draw[white](-1.65,0)--++(0,1)node[above,yshift=-1pt]{\phantom{$1$}};
\draw[white](0.7,0)--++(0,1)node[above,yshift=-1pt]{\phantom{$1$}};
\draw[ghost](0.9,0)--++(0,1)node[above,yshift=-1pt]{$2$};
\draw[solid](-0.1,0)node[below]{$2\;$}--++(0,1);
\draw[solid](0.7,0)node[below]{$3\;$}--++(0,1);
\draw[solid](0.8,0)node[below]{$3$}--++(0,1);
\draw[redstring](0,0)node[below]{$2$}--++(0,1);
\end{tikzpicture}
\ar[d,very thick,->,"\text{ add $2$-string}"]
\\
\begin{tikzpicture}[anchorbase,smallnodes,yscale=0.5]
\draw[white](-1.65,0)--++(0,1)node[above,yshift=-1pt]{\phantom{$1$}};
\draw[white](0.7,0)--++(0,1)node[above,yshift=-1pt]{\phantom{$1$}};
\draw[ghost](0.6,0)--++(0,1)node[above,yshift=-1pt]{$2$};
\draw[ghost](0.9,0)--++(0,1)node[above,yshift=-1pt]{$2$};
\draw[solid](-0.4,0)node[below]{$2$}--++(0,1);
\draw[solid](-0.1,0)node[below]{$2\;$}--++(0,1);
\draw[solid](0.7,0)node[below]{$3\;$}--++(0,1);
\draw[solid](0.8,0)node[below]{$3$}--++(0,1);
\draw[redstring](0,0)node[below]{$2$}--++(0,1);
\end{tikzpicture}
\ar[d,very thick,->,"\text{ add $1$-string}"]
\\
\llap{$\idem[\mu]=$}\begin{tikzpicture}[anchorbase,smallnodes,yscale=0.5]
\draw[white](-1.65,0)--++(0,1)node[above,yshift=-1pt]{\phantom{$1$}};
\draw[white](0.7,0)--++(0,1)node[above,yshift=-1pt]{\phantom{$1$}};
\draw[ghost](-0.5,0)--++(0,1)node[above,yshift=-1pt]{$1$};
\draw[ghost](0.6,0)--++(0,1)node[above,yshift=-1pt]{$2$};
\draw[ghost](0.9,0)--++(0,1)node[above,yshift=-1pt]{$2$};
\draw[solid](-1.5,0)node[below]{$1$}--++(0,1);
\draw[solid](-0.4,0)node[below]{$2$}--++(0,1);
\draw[solid](-0.1,0)node[below]{$2\;$}--++(0,1);
\draw[solid](0.7,0)node[below]{$3\;$}--++(0,1);
\draw[solid](0.8,0)node[below]{$3$}--++(0,1);
\draw[redstring](0,0)node[below]{$2$}--++(0,1);
\end{tikzpicture}
\end{tikzcd}$}
\quad.
\end{gather}
The final picture on the right-hand side is the idempotent diagram $\idem[\mu]$. Note that the first placed solid $2$-string
is blocked by the red $2$-string, since in this case
the Reidemeister II move \autoref{Eq:RecollectionReidemeisterII} can only be applied if there a dot. Moreover, the second string
that we place is a solid $3$-string, which is blocked by the ghost $2$-string,
again by \autoref{Eq:RecollectionReidemeisterII}.
The second solid $3$-string is blocked similarly whereas the second solid $2$-string, and the solid $1$-string, have their ghost strings blocked by the solid strings on their right. That is, each string in $\idem[\mu]$ is blocked in the sense that none of the defining relations allow any string to be pulled the string further to the right. A priori, it is possible that there is a some cunning way to apply the relations to pull strings further to the right but we will eventually see that this is not the case.

In the same way, we construct the idempotent diagrams for $\lambda$ and $\nu$, using the residue sequences $\Res[{\ppath[\lambda]}]=23312$ and $\Res[{\ppath[\nu]}]=21323$. This gives the three idempotent diagrams:
\begin{gather}\label{Eq:MainExampleOrder}
\begin{aligned}
\idem[\lambda]
=&
\begin{tikzpicture}[anchorbase,smallnodes,yscale=0.5]
\draw[white](-1.65,0)--++(0,1)node[above,yshift=-1pt]{\phantom{$1$}};
\draw[white](0.7,0)--++(0,1)node[above,yshift=-1pt]{\phantom{$1$}};
\draw[ghost](-0.2,0)--++(0,1)node[above,yshift=-1pt]{$1$};
\draw[ghost](0.9,0)--++(0,1)node[above,yshift=-1pt]{$2$};
\draw[ghost](0.6,0)--++(0,1)node[above,yshift=-1pt]{$2$};
\draw[solid](-1.2,0)node[below]{$1$}--++(0,1);
\draw[solid](-0.4,0)node[below]{$2$}--++(0,1);
\draw[solid](-0.1,0)node[below]{$2\;$}--++(0,1);
\draw[solid](0.7,0)node[below]{$3\;$}--++(0,1);
\draw[solid](0.8,0)node[below]{$3$}--++(0,1);
\draw[redstring](0,0)node[below]{$2$}--++(0,1);
\end{tikzpicture}
,
\\[-0.1cm]
\idem[\mu]
=&
\begin{tikzpicture}[anchorbase,smallnodes,yscale=0.5]
\draw[white](-1.65,0)--++(0,1)node[above,yshift=-1pt]{\phantom{$1$}};
\draw[white](0.7,0)--++(0,1)node[above,yshift=-1pt]{\phantom{$1$}};
\draw[ghost](-0.5,0)--++(0,1)node[above,yshift=-1pt]{$1$};
\draw[ghost](0.6,0)--++(0,1)node[above,yshift=-1pt]{$2$};
\draw[ghost](0.9,0)--++(0,1)node[above,yshift=-1pt]{$2$};
\draw[solid](-1.5,0)node[below]{$1$}--++(0,1);
\draw[solid](-0.4,0)node[below]{$2$}--++(0,1);
\draw[solid](-0.1,0)node[below]{$2\;$}--++(0,1);
\draw[solid](0.7,0)node[below]{$3\;$}--++(0,1);
\draw[solid](0.8,0)node[below]{$3$}--++(0,1);
\draw[redstring](0,0)node[below]{$2$}--++(0,1);
\end{tikzpicture}
,
\\[-0.1cm]
\idem[\nu]
=&
\begin{tikzpicture}[anchorbase,smallnodes,yscale=0.5]
\draw[white](-1.65,0)--++(0,1)node[above,yshift=-1pt]{\phantom{$1$}};
\draw[white](0.7,0)--++(0,1)node[above,yshift=-1pt]{\phantom{$1$}};
\draw[ghost](-0.2,0)--++(0,1)node[above,yshift=-1pt]{$1$};
\draw[ghost](0.7,0)--++(0,1)node[above,yshift=-1pt]{$2$};
\draw[ghost](0.9,0)--++(0,1)node[above,yshift=-1pt]{$2$};
\draw[solid](-1.2,0)node[below]{$1$}--++(0,1);
\draw[solid](-0.3,0)node[below]{$2$}--++(0,1);
\draw[solid](-0.1,0)node[below]{$2\;$}--++(0,1);
\draw[solid](0.6,0)node[below]{$3$}--++(0,1);
\draw[solid](0.8,0)node[below]{$3$}--++(0,1);
\draw[redstring](0,0)node[below]{$2$}--++(0,1);
\end{tikzpicture}
.
\end{aligned}
\end{gather}
Reading left to right, the sequences of coordinates of the solid strings in these idempotent diagrams are
$(-1.3,-0.4,-0.1,0.7,0.8)$ for $\lambda$,
$(-1.5,-0.4,-0.1,0.7,0.8)$ for $\mu$ and
$(-1.2,-0.3,-0.1,0.6,0.8)$ for $\nu$.
Reading these coordinate sequences lexicographically defines an order $\sandorder$ on the vertices. In this case, $\nu\sandorder\mu\sandorder\lambda$.
Note that the order of the strings, from left to right, in the idempotent diagrams is, in general, different from the order coming from reading along the residue sequences.

By definition, the idempotent diagrams $\idem[\lambda]$, $\idem[\mu]$ and $\idem[\nu]$ do not belong to the algebra $\WAc[5](X)$ because not all of the $x$-coordinates of the strings in these diagrams belong to $X$. However, we will see that we can use these diagrams to give a basis for $\WAc[5](X)$ because the steady diagrams in this algebra factor through these idempotents.

\item \textbf{Dots on idempotents.}
Note the two close solid $3$-strings in $\idem[\lambda]$ and
$\idem[\mu]$. We want
to flank these two strings with crossings,  however, doing this annihilates these
diagrams by \autoref{Eq:RecollectionReidemeisterII}. To avoid
this we put a dot on the right-hand solid $3$-string in each diagram, which gives the dotted idempotents $\dotidem[\lambda]$ and $\dotidem[\mu]$, associated to $\lambda$ and $\mu$, respectively. The diagram
for $\nu$ does not need any extra dots because it has no repeated residues.
Let $y_{k}$ be the operation of putting a dot on the $k$th
solid string when reading along the residue sequence. Then:
\begin{gather*}
\dotidem[\lambda]=y_{2}\idem[\lambda]
=
\begin{tikzpicture}[anchorbase,smallnodes,yscale=0.5]
\draw[ghost](-0.3,0)--++(0,1)node[above,yshift=-1pt]{$1$};
\draw[ghost](0.6,0)--++(0,1)node[above,yshift=-1pt]{$2$};
\draw[ghost](0.9,0)--++(0,1)node[above,yshift=-1pt]{$2$};
\draw[solid](-1.3,0)node[below]{$1$}--++(0,1);
\draw[solid](-0.4,0)node[below]{$2$}--++(0,1);
\draw[solid](-0.1,0)node[below]{$2\;$}--++(0,1);
\draw[solid](0.7,0)node[below]{$3\;$}--++(0,1);
\draw[solid,dot](0.8,0)node[below]{$3$}--++(0,1);
\draw[redstring](0,0)node[below]{$2$}--++(0,1);
\end{tikzpicture}
,\quad
\dotidem[\mu]=y_{2}\idem[\mu]
=
\begin{tikzpicture}[anchorbase,smallnodes,yscale=0.5]
\draw[ghost](-0.5,0)--++(0,1)node[above,yshift=-1pt]{$1$};
\draw[ghost](0.6,0)--++(0,1)node[above,yshift=-1pt]{$2$};
\draw[ghost](0.9,0)--++(0,1)node[above,yshift=-1pt]{$2$};
\draw[solid](-1.5,0)node[below]{$1$}--++(0,1);
\draw[solid](-0.4,0)node[below]{$2$}--++(0,1);
\draw[solid](-0.1,0)node[below]{$2\;$}--++(0,1);
\draw[solid](0.7,0)node[below]{$3\;$}--++(0,1);
\draw[solid,dot](0.8,0)node[below]{$3$}--++(0,1);
\draw[redstring](0,0)node[below]{$2$}--++(0,1);
\end{tikzpicture}
\quad\text{and}\quad
\dotidem[\nu]
=\idem[\nu]
.
\end{gather*}
Note that our convention for multiplying idempotent diagrams by polynomials $f(\bu)\in\N[\bu]$ is based on the residue sequence, which is different to the initial convention that we used in \autoref{SS:RecollectionQPoly} but more convenient for what follows. These dotted idempotents have positive degree if they carry dots. For example, above $\deg(\dotidem[\lambda])=\deg(\dotidem[\mu])=2$ and
$\deg(\dotidem[\nu])=0$.

\item \textbf{Sandwiched dots.}
Our basis will require sandwich dots on some dotted idempotents $\dotidem[\lambda]$, which will give the basis $\sandbasis[\lambda]$ for the corresponding sandwich algebra $\sand[\lambda]$. If there are no sandwich dots then, up to shift, $\sand[\lambda]=\ring$. The idempotent $\dotidem[\nu]$ has
two solid $3$-strings that are each blocked
by a ghost $2$-string. By \autoref{L:RecollectionMovingStringsDots}.(b)
they are still blocked when we put a dot on them,
so we allow putting a dot on either these strings.
We will see that
$\sandbasis[\nu]=\set{\dotidem[\nu],y_{3}\dotidem[\nu],y_{5}\dotidem[\nu],y_{3}y_{5}\dotidem[\nu]}$ is a homogeneous basis of the sandwich algebra
$\sand[\nu]$. The other two sandwiched bases are trivial, so
$\sandbasis[\lambda]=\set{\dotidem[\lambda]}$,
$\sandbasis[\mu]=\set{\dotidem[\mu]}$ and $\sand[\lambda]\cong\sand[\mu]\cong R$.
Hence, the graded ranks of the sandwiched
algebras are $\grdim(\sand[\lambda])=\grdim(\sand[\mu])=\vpar^{2}$
and $\grdim(\sand[\nu])=(1+\vpar^{2})^{2}$.

\item \textbf{Face permutations.}
A permutation $\sigma\in\sym[5]$ acts on a residue sequence of length $5$
by permuting entries.
Let $D_{w}$ be the associated permutation diagram that
permutes the solid strings and their ghosts. We want to find all permutations of the residue sequences that are compatible with the dotted idempotent diagrams for the rooted paths in $\crystalgraph$. We consider those  permutation for which composition of diagrams $\dotidem[\lambda_{1}]D_{w}\dotidem[\lambda_{2}]$ is nonzero.
As we will see later, these permutations are
precisely the \emph{face permutations}, which permute the paths in the crystal around faces.

Note that there are three paths in the crystal graph to $\lambda$, which have residue sequences $23321$, $231323$ and $21332$. It turns out that these three paths correspond to the same wKLRW idempotent diagram, so we do not need to consider them. The point is that these paths differ only by the direction in which we travel around a square with edges labeled by $1$ and $3$, and these strings commute in the wKLRW algebra (as do the corresponding generators of the quantum group). Similarly, the different paths to $\mu$ and $\nu$ in the crystal graph correspond to the same wKLRW diagrams.

The complete list of face permutations is given by the following table:
\begin{gather*}
\raisebox{-0.25cm}{$\begin{array}{ccc}
\toprule
\lambda_{1} & \lambda_{2}  &  \set{\sigma} \\ \midrule
\lambda & \lambda & \set{1,(2,3)} \\
\lambda & \mu & \set{(2,3)(4,5)} \\
\lambda & \nu & \set{(2,5,4)}\\
\mu & \mu & \set{1,(2,3)} \\
\nu & \nu & \set{1} \\
\bottomrule
\end{array}$}
\;.
\end{gather*}

\item \textbf{The basis.} The basis we construct starts with the
smallest idempotent diagram $\dotidem[\nu]$ and its sandwiched
basis $\sandbasis[\nu]$, which is flanked with
\begin{enumerate}[label=$\bullet$]

\item $D_{(2,4,5)}=D_{(2,5,4)^{-1}}$ going from $\nu$ to $\lambda$;

\item $D_{(2,3)}$ going from $\nu$ to $\nu_{2}$.

\end{enumerate}
The next biggest is $\mu$ and its sandwiched
basis $\sandbasis[\mu]$ is flanked with
\begin{enumerate}[label=$\bullet$]

\item $D_{(2,3)}$ going from $\mu$ to $\mu$ itself;

\item $D_{(2,3)(4,5)}$  going from $\mu$ to $\nu$.

\end{enumerate}
Finally, the sandwiched basis
$\sandbasis[\lambda]$ is flanked with
\begin{enumerate}[label=$\bullet$]

\item $D_{(2,3)}$ going from $\lambda$ to $\lambda$ itself;

\item $D_{(2,4,3)}$ and
$D_{(2,4)}$ going from $\lambda$ to $\lambda_{2}$ and $\lambda_{3}$, respectively.

\end{enumerate}
Hence, the graded rank is $\grdim(\WAc)=
6\vpar^{-2}+19+26\vpar^{2}+19\vpar^{4}+6\vpar^{6}$, so
that $\grdim[1](\WAc)=76$. These numbers are obtained by computing the degrees of the basis elements listed above. More explicitly, the nonzero ranks
$\grdim(\idem[\lambda_{1}]\WAc\idem[\lambda_{2}])=\grdim(\idem[\lambda_{2}]\WAc\idem[\lambda_{1}])$ are given by:
\begin{gather*}
\begin{array}{ccc}\toprule
\lambda_{1} & \lambda_{2}  &  \grdim(\idem[\lambda_{1}]\WAc\idem[\lambda_{2}]) \\ \midrule
\lambda & \lambda & (\vpar^{4}+1)(\vpar^{2}+1)^{2}\vpar^{-2} \\
\lambda & \mu & (\vpar^{2}+1)^{2} \\
\lambda & \nu & (\vpar^{2}+1)^{2} \\
\mu & \mu & (\vpar^{4}+1)(\vpar^{2}+1)^{2}\vpar^{-2} \\
\nu & \nu & (\vpar^{2}+1)^{2} \\
\bottomrule
\end{array}
\;.
\end{gather*}
These numbers agree with the general formula for the graded ranks of the cyclotomic KLR algebras given by Hu--Shi \cite[Theorem 1.1]{HuSh-monomial-klr-basis}.
Note that the graded rank is palindromic up to a shift, which is
expected since the cyclotomic KLR algebras are graded symmetric algebras by \cite[Remark 3.19]{We-knot-invariants}.

\item \textbf{Detour permutations.} Finally, consider the two paths
$\ppath[\lambda^{\prime}]$ and $\qpath[\lambda^{\prime}]$
below, which end at the same vertex $\lambda^{\prime}$ of the crystal:
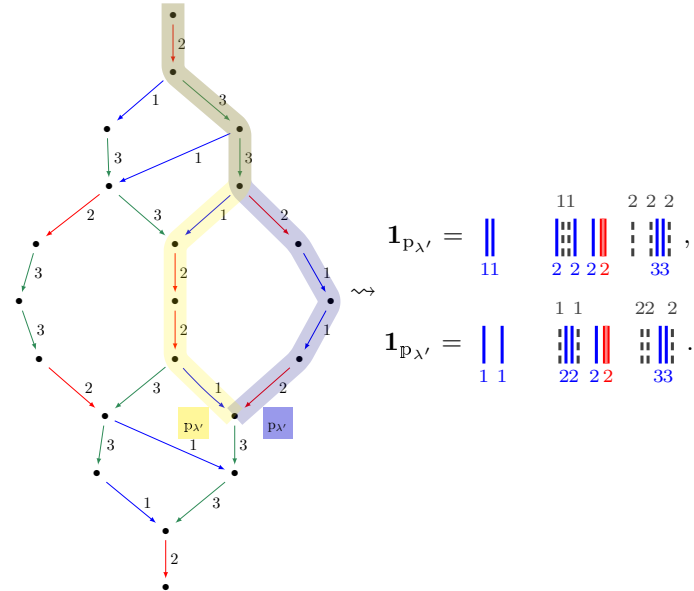
\begin{gather*}
\scalebox{0.6}{$\begin{tikzpicture}[anchorbase,>=latex,line join=bevel,xscale=0.65,yscale=0.5]
\node (node_0) at (69.5bp,439.0bp) [draw,draw=none] {$\bullet$};
\node (node_6) at (53.5bp,367.0bp) [draw,draw=none] {$\bullet$};
\node (node_1) at (137.5bp,583.0bp) [draw,draw=none] {$\bullet$};
\node (node_9) at (139.5bp,511.0bp) [draw,draw=none] {$\bullet$};
\node (node_2) at (127.5bp,151.0bp) [draw,draw=none] {$\bullet$};
\node (node_8) at (193.5bp,79.5bp) [draw,draw=none] {$\bullet$};
\node (node_3) at (264.5bp,511.0bp) [draw,draw=none] {$\bullet$};
\node (node_17) at (320.5bp,439.0bp) [draw,draw=none] {$\bullet$};
\node (node_18) at (202.5bp,439.0bp) [draw,draw=none] {$\bullet$};
\node (node_4) at (259.5bp,223.0bp) [draw,draw=none] {$\bullet$};
\node (node_16) at (259.5bp,151.0bp) [draw,draw=none] {$\bullet$};
\node (node_5) at (72.5bp,295.0bp) [draw,draw=none] {$\bullet$};
\node (node_14) at (135.5bp,223.0bp) [draw,draw=none] {$\bullet$};
\node (node_7) at (321.5bp,295.0bp) [draw,draw=none] {$\bullet$};
\node (node_13) at (193.5bp,8.5bp) [draw,draw=none] {$\bullet$};
\node (node_10) at (200.5bp,725.5bp) [draw,draw=none] {$\bullet$};
\node (node_12) at (200.5bp,654.5bp) [draw,draw=none] {$\bullet$};
\node (node_11) at (202.5bp,295.0bp) [draw,draw=none] {$\bullet$};
\node (node_19) at (264.5bp,583.0bp) [draw,draw=none] {$\bullet$};
\node (node_15) at (351.5bp,367.0bp) [draw,draw=none] {$\bullet$};
\node (node_20) at (202.5bp,367.0bp) [draw,draw=none] {$\bullet$};
\draw [spinach,->] (node_0) ..controls (65.255bp,419.43bp) and (60.953bp,400.61bp)  .. (node_6);
\draw (72.0bp,403.0bp) node {$3$};
\draw [spinach,->] (node_1) ..controls (138.02bp,563.68bp) and (138.57bp,544.46bp)  .. (node_9);
\draw (148.0bp,547.0bp) node {$3$};
\draw [blue,->] (node_2) ..controls (145.48bp,131.07bp) and (165.5bp,109.99bp)  .. (node_8);
\draw (176.0bp,115.0bp) node {$1$};
\draw [red,->] (node_3) ..controls (279.67bp,491.04bp) and (296.42bp,470.1bp)  .. (node_17);
\draw (307.0bp,475.0bp) node {$2$};
\draw [blue,->] (node_3) ..controls (247.61bp,490.93bp) and (228.8bp,469.7bp)  .. (node_18);
\draw (249.0bp,475.0bp) node {$1$};
\draw [spinach,->] (node_4) ..controls (259.5bp,203.68bp) and (259.5bp,184.46bp)  .. (node_16);
\draw (268.0bp,187.0bp) node {$3$};
\draw [red,->] (node_5) ..controls (89.778bp,274.8bp) and (108.19bp,254.34bp)  .. (node_14);
\draw (120.0bp,259.0bp) node {$2$};
\draw [spinach,->] (node_6) ..controls (58.569bp,347.32bp) and (63.751bp,328.23bp)  .. (node_5);
\draw (74.0bp,331.0bp) node {$3$};
\draw [red,->] (node_7) ..controls (304.61bp,274.93bp) and (285.8bp,253.7bp)  .. (node_4);
\draw (306.0bp,259.0bp) node {$2$};
\draw [red,->] (node_8) ..controls (193.5bp,60.442bp) and (193.5bp,41.496bp)  .. (node_13);
\draw (202.0bp,44.0bp) node {$2$};
\draw [red,->] (node_9) ..controls (120.2bp,490.7bp) and (99.457bp,469.96bp)  .. (node_0);
\draw (122.0bp,475.0bp) node {$2$};
\draw [spinach,->] (node_9) ..controls (156.99bp,490.56bp) and (175.99bp,469.46bp)  .. (node_18);
\draw (186.0bp,475.0bp) node {$3$};
\draw [red,->] (node_10) ..controls (200.5bp,706.44bp) and (200.5bp,687.5bp)  .. (node_12);
\draw (209.0bp,690.0bp) node {$2$};
\draw [blue,->] (node_11) ..controls (214.23bp,277.42bp) and (225.23bp,262.33bp)  .. (235.5bp,250.0bp) .. controls (238.51bp,246.39bp) and (241.86bp,242.61bp)  .. (node_4);
\draw (244.0bp,259.0bp) node {$1$};
\draw [spinach,->] (node_11) ..controls (184.26bp,274.95bp) and (163.98bp,253.75bp)  .. (node_14);
\draw (185.0bp,259.0bp) node {$3$};
\draw [blue,->] (node_12) ..controls (183.34bp,634.57bp) and (164.23bp,613.49bp)  .. (node_1);
\draw (184.0bp,619.0bp) node {$1$};
\draw [spinach,->] (node_12) ..controls (217.71bp,634.81bp) and (236.53bp,614.37bp)  .. (node_19);
\draw (248.0bp,619.0bp) node {$3$};
\draw [spinach,->] (node_14) ..controls (133.35bp,203.2bp) and (131.13bp,183.79bp)  .. (node_2);
\draw (141.0bp,187.0bp) node {$3$};
\draw [blue,->] (node_14) ..controls (170.82bp,202.06bp) and (210.64bp,179.58bp)  .. (node_16);
\draw (220.0bp,187.0bp) node {$1$};
\draw [blue,->] (node_15) ..controls (343.35bp,346.99bp) and (334.79bp,327.02bp)  .. (node_7);
\draw (348.0bp,331.0bp) node {$1$};
\draw [spinach,->] (node_16) ..controls (240.94bp,130.45bp) and (221.46bp,109.94bp)  .. (node_8);
\draw (242.0bp,115.0bp) node {$3$};
\draw [blue,->] (node_17) ..controls (328.71bp,419.47bp) and (337.46bp,399.69bp)  .. (node_15);
\draw (348.0bp,403.0bp) node {$1$};
\draw [red,->] (node_18) ..controls (202.5bp,419.68bp) and (202.5bp,400.46bp)  .. (node_20);
\draw (211.0bp,403.0bp) node {$2$};
\draw [spinach,->] (node_19) ..controls (264.5bp,563.2bp) and (264.5bp,543.79bp)  .. (node_3);
\draw (273.0bp,547.0bp) node {$3$};
\draw [blue,->] (node_19) ..controls (228.89bp,562.06bp) and (188.75bp,539.58bp)  .. (node_9);
\draw (225.0bp,547.0bp) node {$1$};
\draw [red,->] (node_20) ..controls (202.5bp,347.2bp) and (202.5bp,327.79bp)  .. (node_11);
\draw (211.0bp,331.0bp) node {$2$};
\draw[black!50!blue,line width=0.5cm,opacity=0.2,rounded corners] (200.5bp,740.5bp) to (200.5bp,654.5bp) to (264.5bp,583.0bp) to (264.5bp,511.0bp) to (320.5bp,439.0bp) to (351.5bp,367.0bp) to (321.5bp,295.0bp) to (259.5bp,223.0bp);
\node[black] at (300.5bp,215.0bp) {\colorbox{black!20!blue!40!white}{\mystrut\scalebox{0.9}{$\qpath[\lambda^{\prime}]$}}};
\draw[yellow,line width=0.5cm,opacity=0.2,rounded corners] (200.5bp,740.5bp) to (200.5bp,654.5bp) to (264.5bp,583.0bp) to (264.5bp,511.0bp) to (202.5bp,439.0bp) to (202.5bp,367.0bp) to (202.5bp,295.0bp) to (259.5bp,223.0bp);
\node[black] at (220.5bp,215.0bp) {\colorbox{yellow!50}{\mystrut\scalebox{0.9}{$\ppath[\lambda^{\prime}]$}}};
\end{tikzpicture}$}
\rightsquigarrow
\begin{gathered}
\idem[{\ppath[\lambda^{\prime}]}]=
\begin{tikzpicture}[anchorbase,smallnodes,yscale=0.5]
\draw[redstring](0.05,0)node[below]{\,$2$}--+(0,1);
\draw[solid](-0.08,0)node[below]{$2$\,}--+(0,1);
\draw[ghost](0.92,0)--+(0,1)node[above]{$2$};
\draw[solid](0.84,0)node[below]{\,$3$}--+(0,1);
\draw[solid](0.76,0)node[below]{$3$\,}--+(0,1);
\draw[solid](-0.32,0)node[below]{$2$}--+(0,1);
\draw[ghost](0.68,0)--+(0,1)node[above]{$2$};
\draw[solid](-1.40,0)node[below]{\,$1$}--+(0,1);
\draw[ghost](-0.40,0)--+(0,1)node[above]{\,$1$};
\draw[solid](-1.48,0)node[below]{$1$\,}--+(0,1);
\draw[ghost](-0.48,0)--+(0,1)node[above]{$1$\,};
\draw[solid](-0.56,0)node[below]{$2$}--+(0,1);
\draw[ghost](0.44,0)--+(0,1)node[above]{$2$};
\end{tikzpicture}
,
\\
\idem[{\qpath[\lambda^{\prime}]}]=
\begin{tikzpicture}[anchorbase,smallnodes,yscale=0.5]
\draw[redstring](0.05,0)node[below]{\,$2$}--+(0,1);
\draw[solid](-0.08,0)node[below]{$2$\,}--+(0,1);
\draw[ghost](0.92,0)--+(0,1)node[above]{$2$};
\draw[solid](0.84,0)node[below]{\,$3$}--+(0,1);
\draw[solid](0.76,0)node[below]{$3$\,}--+(0,1);
\draw[solid](-1.32,0)node[below]{$1$}--+(0,1);
\draw[ghost](-0.32,0)--+(0,1)node[above]{$1$};
\draw[solid](-0.40,0)node[below]{\,$2$}--+(0,1);
\draw[ghost](0.60,0)--+(0,1)node[above]{\,$2$};
\draw[solid](-0.48,0)node[below]{$2$\,}--+(0,1);
\draw[ghost](0.52,0)--+(0,1)node[above]{$2$\,};
\draw[solid](-1.56,0)node[below]{$1$}--+(0,1);
\draw[ghost](-0.56,0)--+(0,1)node[above]{$1$};
\end{tikzpicture}
.
\end{gathered}
\end{gather*}
The paths $\ppath[\lambda^{\prime}]$ and $\qpath[\lambda^{\prime}]$ differ by the direction in which they travel around the octagon, which has edges labeled by noncommuting residues $1$ and $2$. As a result, these paths correspond to different idempotent diagrams in the wKLRW algebra, so both of them need to be considered. For our basis, we fix the path that is minimal with respect to the total order $\lsandorder$, which in this case is $\ppath[\lambda^{\prime}]$. We also need to include the face permutation $(4,6)(5,7)$ from $\ppath[\lambda^{\prime}]$ to $\qpath[\lambda^{\prime}]$. We call such permutations \emph{detour permutations} because they give different paths to the same vertex in the crystal graph.
\item \textbf{Sandwich cellularity.}
By \autoref{T:ConstructionMain}.(b) below, this basis is a sandwich cellular basis of $\WAc(X)$. As the sandwich algebra $\sand[\nu]$ is nontrivial, it turns out that this basis cannot be refined to a cellular basis. In fact, by results of Robert Muth, $\WAc(X)$ is not cellular, cf. \autoref{R:ConsequencesGraphsNotCellular}.
\end{enumerate}
We use this example as a running example throughout the paper.
\end{Example}

%%%%%%%%%%%%%%%%%%%%%%%%%%%%%%%%%%%%%%%%%

\section{Finite types and crystal graphs}\label{S:FiniteTypes}

%%%%%%%%%%%%%%%%%%%%%%%%%%%%%%%%%%%%%%%%%

We now specify the combinatorial notions we use.

%%%%%%%%%%%%%%%%%%%%%%%%%%%%%%%%%%%%%%%%%

\subsection{Finite type combinatorics}\label{SS:FiniteTypesConventions}

%%%%%%%%%%%%%%%%%%%%%%%%%%%%%%%%%%%%%%%%%

We use the following conventions for the \emph{finite types}.
Let $e\in\Z_{\geq 1}$.
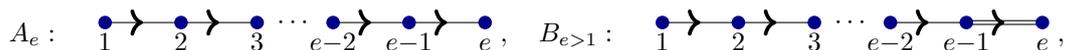
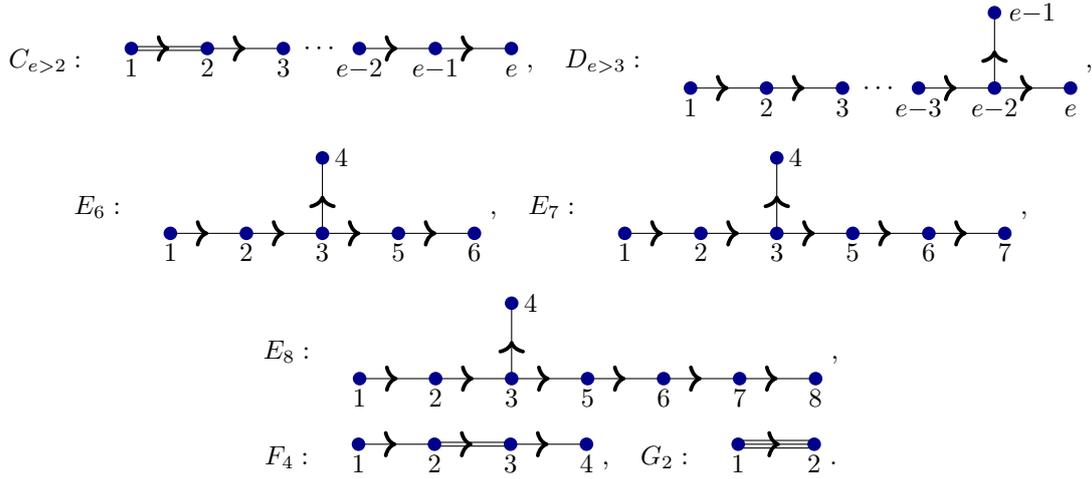
\begin{gather*}
\typea:\quad
\begin{tikzpicture}[anchorbase]
\draw[directed=0.5](1,0)--(2,0);
\draw[directed=0.5](2,0)--(3,0);
\draw[directed=0.5](4,0)--(5,0);
\draw[directed=0.5](5,0)--(6,0);
\node at (3.5,0){$\cdots$};
\foreach \x in {1,...,3} {
\node[dynkin=\x] (\x) at (\x,0){};
}
\foreach \x [evaluate=\x as \c using {int(6-\x)}] in {4,5} {
\node[dynkin=e{-}\c] (\x) at (\x,0){};
}
\node[dynkin=e] (6) at (6,0){};
\end{tikzpicture}
,\quad
\typeb[e>1]:\quad
\begin{tikzpicture}[anchorbase]
\draw[directed=0.5](1,0)--(2,0);
\draw[directed=0.5](2,0)--(3,0);
\draw[directed=0.5](4,0)--(5,0);
\draw[directed=0.5,white](5,0)--(6,0);
\draw(5,0.025)--(6,0.025);
\draw(5,-0.025)--(6,-0.025);
\node at (3.5,0){$\cdots$};
\foreach \x in {1,...,3} {
\node[dynkin=\x] (\x) at (\x,0){};
}
\foreach \x [evaluate=\x as \c using {int(6-\x)}] in {4,5} {
\node[dynkin=e{-}\c] (\x) at (\x,0){};
}
\node[dynkin=e] (6) at (6,0){};
\end{tikzpicture}
,
\\
\typec[e>2]:\quad
\begin{tikzpicture}[anchorbase]
\draw[directed=0.5,white](1,0)--(2,0);
\draw(1,0.025)--+(1,0);
\draw(1,-0.025)--+(1,0);
\draw[directed=0.5](2,0)--+(1,0);
\draw[directed=0.5](4,0)--+(1,0);
\draw[directed=0.5](5,0)--+(1,0);
\node at (3.5,0){$\cdots$};
\foreach \x in {1,...,3} {
\node[dynkin=\x] (\x) at (\x,0){};
}
\foreach \x [evaluate=\x as \c using {int(6-\x)}] in {4,5} {
\node[dynkin=e{-}\c] (\x) at (\x,0){};
}
\node[dynkin=e] (6) at (6,0){};
\end{tikzpicture}
,\quad
\typed[e>3]:\quad
\begin{tikzpicture}[anchorbase]
\draw[directed=0.5](5,0)--(5,1);
\draw[directed=0.5](1,0)--(2,0);
\draw[directed=0.5](2,0)--(3,0);
\draw[directed=0.5](4,0)--(5,0);
\draw[directed=0.5](5,0)--(6,0);
\node at (3.5,0){$\cdots$};
\node[dynkin] (1) at (5,1){};
\node at (5.5,1) {$e{-}1$};
\node[dynkin=e] (1) at (6,0){};
\foreach \x in {1,...,3} {
\node[dynkin=\x] (\x) at (\x,0){};
}
\foreach \x [evaluate=\x as \c using {int(6-\x)}] in {3,4} {
\node[dynkin=e{-}\c] (\x) at (\x+1,0){};
}
\end{tikzpicture}
,
\\
\typee:\quad
\begin{tikzpicture}[anchorbase]
\draw[directed=0.5](3,0)--(3,1);
\draw[directed=0.5](1,0)--(2,0);
\draw[directed=0.5](2,0)--(3,0);
\draw[directed=0.5](3,0)--(4,0);
\draw[directed=0.5](4,0)--(5,0);
\foreach \x [evaluate=\x as \xx using {int(\x>3?\x+1:\x)}] in {1,...,5} {
\node[dynkin=\xx] (\x) at (\x,0){};
}
\node[dynkin] (6) at (3,1){};
\node at (3.25,1) {$4$};
\end{tikzpicture}
,\quad
\typee[7]:\quad
\begin{tikzpicture}[anchorbase]
\draw[directed=0.5](3,0)--(3,1);
\draw[directed=0.5](1,0)--(2,0);
\draw[directed=0.5](2,0)--(3,0);
\draw[directed=0.5](3,0)--(4,0);
\draw[directed=0.5](4,0)--(5,0);
\draw[directed=0.5](5,0)--(6,0);
\foreach \x [evaluate=\x as \xx using {int(\x>3?\x+1:\x)}] in {1,...,6} {
\node[dynkin=\xx] (\x) at (\x,0){};
}
\node[dynkin] (7) at (3,1){};
\node at (3.25,1) {$4$};
\end{tikzpicture}
,
\\
\typee[8]:\quad
\begin{tikzpicture}[anchorbase]
\draw[directed=0.5](3,0)--(3,1);
\draw[directed=0.5](1,0)--(2,0);
\draw[directed=0.5](2,0)--(3,0);
\draw[directed=0.5](3,0)--(4,0);
\draw[directed=0.5](4,0)--(5,0);
\draw[directed=0.5](5,0)--(6,0);
\draw[directed=0.5](6,0)--(7,0);
\foreach \x [evaluate=\x as \xx using {int(\x>3?\x+1:\x)}] in {1,...,7} {
\node[dynkin=\xx] (\x) at (\x,0){};
}
\node[dynkin] (8) at (3,1){};
\node at (3.25,1) {$4$};
\end{tikzpicture}
,
\\
\typef:\quad
\begin{tikzpicture}[anchorbase]
\draw[directed=0.5](1,0)--(2,0);
\draw[directed=0.5,white](2,0)--(3,0);
\draw (2,0.025)--(3,0.025);
\draw (2,-0.025)--(3,-0.025);
\draw[directed=0.5](3,0)--(4,0);
\foreach \x in {1,...,4} {
\node[dynkin=\x] (\x) at (\x,0){};
}
\end{tikzpicture}
,\quad
\typeg:\quad
\begin{tikzpicture}[anchorbase]
\draw[directed=0.55](1,0)--(2,0);
\draw (1,0.05)--(2,0.05);
\draw (1,-0.05)--(2,-0.05);
\foreach \x in {1,...,2} {
\node[dynkin=\x] (\x) at (\x,0){};
}
\end{tikzpicture}
.
\end{gather*}

The types $\typea$, $\typed$ and $\typee[i]$ are \emph{simply laced},
and their symmetrizers are $\boldsymbol{d}=(1,\dots,1)$. The \emph{doubly laced} types
are $\typeb$ (symmetrizer $\boldsymbol{d}=(2,\dots,2,1)$), $\typec$ (symmetrizer $\boldsymbol{d}=(2,1\dots,1)$) and $\typef[4]$ (symmetrizer $\boldsymbol{d}=(2,2,1,1)$).
The types $\typea$, $\typeb$, $\typec$ and $\typed$ are the \emph{classical types}, while $\typee[i]$, $\typef[4]$ and $\typeg[2]$
(symmetrizer $\boldsymbol{d}=(3,1)$) are the \emph{exceptional types}. The number of vertices in $\Gamma$ is the \emph{rank}.

The vertices $e-1$ and $e$ in type $\typed$, and the vertices $4$ and $5$ in types $\typee[i]$ are \emph{fishtail vertices}.
We call $i,j\in\vertices$ \emph{adjacent} if there is an edge between them in $\quiver$, otherwise they are \emph{nonadjacent}.

\begin{Remark}\label{R:FiniteTypesSageMath}
Our labeling conventions are almost the same as \cite[Pages 265--290, plates I--IX]{Bo-chapters-4-6} and SageMath \cite{sage} but differs from both. Specifically, compared with SageMath
we have reflected the diagrams in type $\typec[e>2]$ and $\typeg[2]$
and renumbered the vertices in type $\typee[i]$. We have chosen this labeling because it leads to slightly better combinatorics for the ghost strings.
\end{Remark}

%%%%%%%%%%%%%%%%%%%%%%%%%%%%%%%%%%%%%%%%%

\subsection{Crystal graphs}\label{SS:FiniteTypesCrystals}

%%%%%%%%%%%%%%%%%%%%%%%%%%%%%%%%%%%%%%%%%

We explain the combinatorial background on crystal graphs
that is needed for this paper.
Everything below is well-known and can be found, for example, in \cite{BuSc-crystal-bases}. We also fix our notation for paths in the crystal, which is not standard but is necessary in what follows because we need to keep track of the different paths in the crystal.

Let $\crystal$ be a (Lusztig--Kashiwara) crystal,
as defined in \cite[Definition 2.13]{BuSc-crystal-bases}.
As usual, we associate to $\crystal$ a labeled
directed graph $\crystalgraph[\crystal]$,
called the \textbf{crystal graph}. The vertices in $\crystal$ are labeled
by weights and the
edges labeled by residues $i\in\vertices$. More explicitly, $x,y\in\crystal$ the graph $\crystalgraph[\crystal]$
has an edge $x\to y$ labeled $i\in\vertices$
if and only if $f_{i}(x)=y$ in $\crystal$ or, equivalently, $e_{i}(y)=x$, where $e_{i}$ and $f_{i}$ are the Kashiwara operators. In addition, for each $i\in\vertices$ the crystal $\crystal$ comes equipped with maps
$\varepsilon_{i},\varphi_{i}\colon\crystal\to\N\cup\set{-\infty}$ such that if $\varepsilon_{i}(x)\neq-\infty$ then $\varepsilon_{i}(x)=\max\set{k\geq0|e_{i}^k(x)\neq0}$ and if $\varphi_{i}(x)\neq-\infty$ then $\varphi_{i}(x)=\max\set{k\geq0|f_{i}^k(x)\neq0}$.

In this paper we only consider crystal graphs $\crystalgraph$ associated to
fundamental weights $\fweight$ in finite types, for
$i\in\vertices$.

The following well-known result will be used silently throughout.

\begin{Lemma}\label{L:FiniteTypesCrystalsExist}
Suppose that $\quiver$ is a quiver of finite type and let
$\Lambda\in P^{+}$. There exists a unique crystal graph
$\crystalgraph[\Lambda]$ for the highest weight module $L(\Lambda)$ of
the quantum Lie algebra for $\quiver$.  The crystal graph
$\crystalgraph[\Lambda]$ has a unique source, the highest weight, and
sink, the lowest weight. In the crystal limit, the action of the
Kashiwara operations $e_{i}$ and $f_{i}$ on $\crystalgraph[\Lambda]$ correspond
to the action of the Chevalley generators $E_{i}$ and $F_{i}$ on
$L(\Lambda)$, for $i\in\vertices$.
\end{Lemma}

\begin{proof}
See, for example, \cite[Theorem 9.25]{Ja-lectures-qgroups}
for an explicit construction.
\end{proof}

We identify the vertices of $\crystalgraph[\Lambda]$ with the weights in
the associated highest weight module, written in terms of the
fundamental weight, with the source vertex of $\crystalgraph[\Lambda]$
labeled by $\Lambda$.  This said, as on the right-hand side of
\autoref{Eq:FiniteTypesSageMath} below, the vertex labels do not play an
important role in what follows, so we often replace the labels with a
$\bullet$ to increase readability.

\begin{Example}\label{E:FiniteTypesSageMath}
One can use SageMath \cite{sage} to
draw crystal graphs. The code below generated the \LaTeX\ TikZ code for the crystal graph below:
\begin{sagemath}
sage: R = RootSystem(['D',5])
sage: La = R.weight_space().basis()
sage: LS = crystals.LSPaths(La[5])
sage: C = LS.subcrystal()
sage: G = LS.digraph(subset = C)
sage: latex(G)
\end{sagemath}
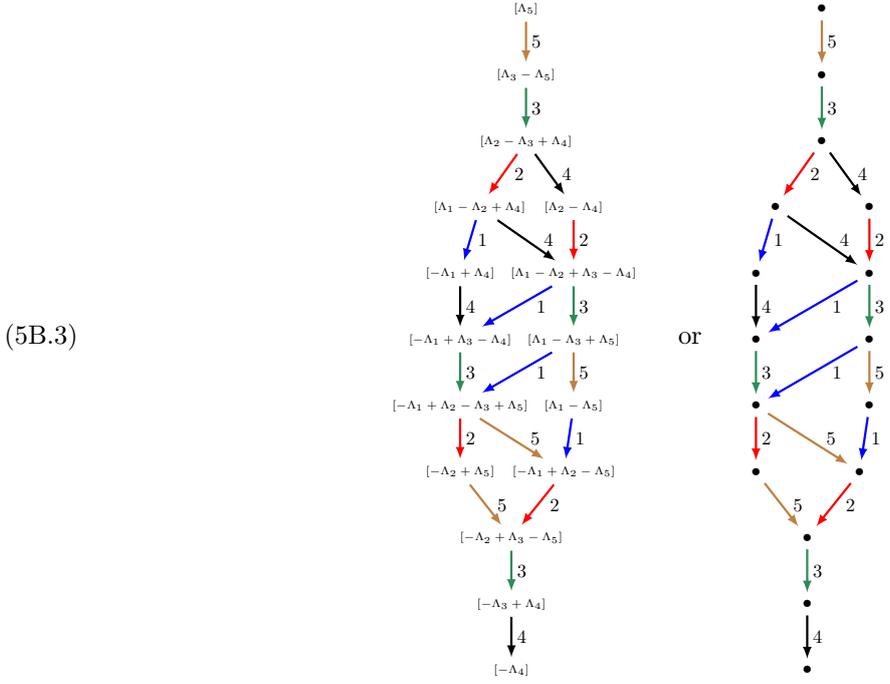
\begin{gather}\label{Eq:FiniteTypesSageMath}
\scalebox{0.9}{\begin{tikzpicture}[>=latex,line join=bevel,anchorbase,xscale=0.65,yscale=0.5,every path/.style={very thick},
every path/.style={very thick}]
\node (node_0) at (67.0bp,505.5bp) [draw,draw=none] {\tiny$[\Lambda_{1} - \Lambda_{2} + \Lambda_{4}]$};
\node (node_11) at (144.0bp,434.5bp) [draw,draw=none] {\tiny$[\Lambda_{1} - \Lambda_{2} + \Lambda_{3} - \Lambda_{4}]$};
\node (node_15) at (51.0bp,434.5bp) [draw,draw=none] {\tiny$[-\Lambda_{1} + \Lambda_{4}]$};
\node (node_1) at (93.0bp,150.5bp) [draw,draw=none] {\tiny$[-\Lambda_{2} + \Lambda_{3} - \Lambda_{5}]$};
\node (node_7) at (93.0bp,79.5bp) [draw,draw=none] {\tiny$[-\Lambda_{3} + \Lambda_{4}]$};
\node (node_2) at (51.0bp,221.5bp) [draw,draw=none] {\tiny$[-\Lambda_{2} + \Lambda_{5}]$};
\node (node_3) at (105.0bp,718.5bp) [draw,draw=none] {\tiny$[\Lambda_{5}]$};
\node (node_10) at (105.0bp,647.5bp) [draw,draw=none] {\tiny$[\Lambda_{3} - \Lambda_{5}]$};
\node (node_4) at (136.0bp,221.5bp) [draw,draw=none] {\tiny$[-\Lambda_{1} + \Lambda_{2} - \Lambda_{5}]$};
\node (node_5) at (93.0bp,8.5bp) [draw,draw=none] {\tiny$[-\Lambda_{4}]$};
\node (node_6) at (144.0bp,505.5bp) [draw,draw=none] {\tiny$[\Lambda_{2} - \Lambda_{4}]$};
\node (node_8) at (105.0bp,576.5bp) [draw,draw=none] {\tiny$[\Lambda_{2} - \Lambda_{3} + \Lambda_{4}]$};
\node (node_9) at (51.0bp,292.5bp) [draw,draw=none] {\tiny$[-\Lambda_{1} + \Lambda_{2} - \Lambda_{3} + \Lambda_{5}]$};
\node (node_12) at (144.0bp,363.5bp) [draw,draw=none] {\tiny$[\Lambda_{1} - \Lambda_{3} + \Lambda_{5}]$};
\node (node_14) at (51.0bp,363.5bp) [draw,draw=none] {\tiny$[-\Lambda_{1} + \Lambda_{3} - \Lambda_{4}]$};
\node (node_13) at (144.0bp,292.5bp) [draw,draw=none] {\tiny$[\Lambda_{1} - \Lambda_{5}]$};
\draw [black,->] (node_0) ..controls (88.072bp,485.62bp) and (111.7bp,464.44bp)  .. (node_11);
\draw (123.5bp,470.0bp) node {$4$};
\draw [blue,->] (node_0) ..controls (62.813bp,486.44bp) and (58.419bp,467.5bp)  .. (node_15);
\draw (69.5bp,470.0bp) node {$1$};
\draw [spinach,->] (node_1) ..controls (93.0bp,131.44bp) and (93.0bp,112.5bp)  .. (node_7);
\draw (101.5bp,115.0bp) node {$3$};
\draw [brown,->] (node_2) ..controls (62.18bp,202.13bp) and (74.214bp,182.36bp)  .. (node_1);
\draw (85.5bp,186.0bp) node {$5$};
\draw [brown,->] (node_3) ..controls (105.0bp,699.44bp) and (105.0bp,680.5bp)  .. (node_10);
\draw (113.5bp,683.0bp) node {$5$};
\draw [red,->] (node_4) ..controls (124.55bp,202.13bp) and (112.23bp,182.36bp)  .. (node_1);
\draw (128.5bp,186.0bp) node {$2$};
\draw [red,->] (node_6) ..controls (144.0bp,486.44bp) and (144.0bp,467.5bp)  .. (node_11);
\draw (152.5bp,470.0bp) node {$2$};
\draw [black,->] (node_7) ..controls (93.0bp,60.442bp) and (93.0bp,41.496bp)  .. (node_5);
\draw (101.5bp,44.0bp) node {$4$};
\draw [red,->] (node_8) ..controls (94.885bp,557.13bp) and (83.997bp,537.36bp)  .. (node_0);
\draw (99.5bp,541.0bp) node {$2$};
\draw [black,->] (node_8) ..controls (115.38bp,557.13bp) and (126.56bp,537.36bp)  .. (node_6);
\draw (138.5bp,541.0bp) node {$4$};
\draw [red,->] (node_9) ..controls (51.0bp,273.44bp) and (51.0bp,254.5bp)  .. (node_2);
\draw (59.5bp,257.0bp) node {$2$};
\draw [brown,->] (node_9) ..controls (74.462bp,272.45bp) and (101.11bp,250.82bp)  .. (node_4);
\draw (112.5bp,257.0bp) node {$5$};
\draw [spinach,->] (node_10) ..controls (105.0bp,628.44bp) and (105.0bp,609.5bp)  .. (node_8);
\draw (113.5bp,612.0bp) node {$3$};
\draw [spinach,->] (node_11) ..controls (144.0bp,415.44bp) and (144.0bp,396.5bp)  .. (node_12);
\draw (152.5bp,399.0bp) node {$3$};
\draw [blue,->] (node_11) ..controls (118.19bp,414.35bp) and (88.643bp,392.43bp)  .. (node_14);
\draw (117.5bp,399.0bp) node {$1$};
\draw [blue,->] (node_12) ..controls (118.19bp,343.35bp) and (88.643bp,321.43bp)  .. (node_9);
\draw (117.5bp,328.0bp) node {$1$};
\draw [brown,->] (node_12) ..controls (144.0bp,344.44bp) and (144.0bp,325.5bp)  .. (node_13);
\draw (152.5bp,328.0bp) node {$5$};
\draw [blue,->] (node_13) ..controls (141.91bp,273.44bp) and (139.71bp,254.5bp)  .. (node_4);
\draw (149.5bp,257.0bp) node {$1$};
\draw [spinach,->] (node_14) ..controls (51.0bp,344.44bp) and (51.0bp,325.5bp)  .. (node_9);
\draw (59.5bp,328.0bp) node {$3$};
\draw [black,->] (node_15) ..controls (51.0bp,415.44bp) and (51.0bp,396.5bp)  .. (node_14);
\draw (59.5bp,399.0bp) node {$4$};
\end{tikzpicture}}
\quad\text{ or }\quad\;
\scalebox{0.9}{\begin{tikzpicture}[>=latex,line join=bevel,anchorbase,xscale=0.65,yscale=0.5,every path/.style={very thick}]
\node (node_0) at (67.0bp,505.5bp) [draw,draw=none] {$\bullet$};
\node (node_11) at (144.0bp,434.5bp) [draw,draw=none] {$\bullet$};
\node (node_15) at (51.0bp,434.5bp) [draw,draw=none] {$\bullet$};
\node (node_1) at (93.0bp,150.5bp) [draw,draw=none] {$\bullet$};
\node (node_7) at (93.0bp,79.5bp) [draw,draw=none] {$\bullet$};
\node (node_2) at (51.0bp,221.5bp) [draw,draw=none] {$\bullet$};
\node (node_3) at (105.0bp,718.5bp) [draw,draw=none] {$\bullet$};
\node (node_10) at (105.0bp,647.5bp) [draw,draw=none] {$\bullet$};
\node (node_4) at (136.0bp,221.5bp) [draw,draw=none] {$\bullet$};
\node (node_5) at (93.0bp,8.5bp) [draw,draw=none] {$\bullet$};
\node (node_6) at (144.0bp,505.5bp) [draw,draw=none] {$\bullet$};
\node (node_8) at (105.0bp,576.5bp) [draw,draw=none] {$\bullet$};
\node (node_9) at (51.0bp,292.5bp) [draw,draw=none] {$\bullet$};
\node (node_12) at (144.0bp,363.5bp) [draw,draw=none] {$\bullet$};
\node (node_14) at (51.0bp,363.5bp) [draw,draw=none] {$\bullet$};
\node (node_13) at (144.0bp,292.5bp) [draw,draw=none] {$\bullet$};
\draw [black,->] (node_0) ..controls (88.072bp,485.62bp) and (111.7bp,464.44bp)  .. (node_11);
\draw (123.5bp,470.0bp) node {$4$};
\draw [blue,->] (node_0) ..controls (62.813bp,486.44bp) and (58.419bp,467.5bp)  .. (node_15);
\draw (69.5bp,470.0bp) node {$1$};
\draw [spinach,->] (node_1) ..controls (93.0bp,131.44bp) and (93.0bp,112.5bp)  .. (node_7);
\draw (101.5bp,115.0bp) node {$3$};
\draw [brown,->] (node_2) ..controls (62.18bp,202.13bp) and (74.214bp,182.36bp)  .. (node_1);
\draw (85.5bp,186.0bp) node {$5$};
\draw [brown,->] (node_3) ..controls (105.0bp,699.44bp) and (105.0bp,680.5bp)  .. (node_10);
\draw (113.5bp,683.0bp) node {$5$};
\draw [red,->] (node_4) ..controls (124.55bp,202.13bp) and (112.23bp,182.36bp)  .. (node_1);
\draw (128.5bp,186.0bp) node {$2$};
\draw [red,->] (node_6) ..controls (144.0bp,486.44bp) and (144.0bp,467.5bp)  .. (node_11);
\draw (152.5bp,470.0bp) node {$2$};
\draw [black,->] (node_7) ..controls (93.0bp,60.442bp) and (93.0bp,41.496bp)  .. (node_5);
\draw (101.5bp,44.0bp) node {$4$};
\draw [red,->] (node_8) ..controls (94.885bp,557.13bp) and (83.997bp,537.36bp)  .. (node_0);
\draw (99.5bp,541.0bp) node {$2$};
\draw [black,->] (node_8) ..controls (115.38bp,557.13bp) and (126.56bp,537.36bp)  .. (node_6);
\draw (138.5bp,541.0bp) node {$4$};
\draw [red,->] (node_9) ..controls (51.0bp,273.44bp) and (51.0bp,254.5bp)  .. (node_2);
\draw (59.5bp,257.0bp) node {$2$};
\draw [brown,->] (node_9) ..controls (74.462bp,272.45bp) and (101.11bp,250.82bp)  .. (node_4);
\draw (112.5bp,257.0bp) node {$5$};
\draw [spinach,->] (node_10) ..controls (105.0bp,628.44bp) and (105.0bp,609.5bp)  .. (node_8);
\draw (113.5bp,612.0bp) node {$3$};
\draw [spinach,->] (node_11) ..controls (144.0bp,415.44bp) and (144.0bp,396.5bp)  .. (node_12);
\draw (152.5bp,399.0bp) node {$3$};
\draw [blue,->] (node_11) ..controls (118.19bp,414.35bp) and (88.643bp,392.43bp)  .. (node_14);
\draw (117.5bp,399.0bp) node {$1$};
\draw [blue,->] (node_12) ..controls (118.19bp,343.35bp) and (88.643bp,321.43bp)  .. (node_9);
\draw (117.5bp,328.0bp) node {$1$};
\draw [brown,->] (node_12) ..controls (144.0bp,344.44bp) and (144.0bp,325.5bp)  .. (node_13);
\draw (152.5bp,328.0bp) node {$5$};
\draw [blue,->] (node_13) ..controls (141.91bp,273.44bp) and (139.71bp,254.5bp)  .. (node_4);
\draw (149.5bp,257.0bp) node {$1$};
\draw [spinach,->] (node_14) ..controls (51.0bp,344.44bp) and (51.0bp,325.5bp)  .. (node_9);
\draw (59.5bp,328.0bp) node {$3$};
\draw [black,->] (node_15) ..controls (51.0bp,415.44bp) and (51.0bp,396.5bp)  .. (node_14);
\draw (59.5bp,399.0bp) node {$4$};
\end{tikzpicture}}
.
\end{gather}
This is the crystal graph $\crystalgraph[\Lambda_5]$ in type $\typed[5]$.
\end{Example}

\begin{Example}\label{E:FiniteTypesCrystalList}
The crystal graphs $\crystalgraph[\Lambda_{1}]$ in types $\typea[4]$,
$\typeb[4]$ and $\typed[4]$, and the crystal graph $\crystalgraph[\Lambda_{4}]$ in type
$\typec[4]$, are:
\begin{gather}\label{Eq:ConstructionMainCrystals}
\typea[4]\colon
\scalebox{0.8}{$\begin{tikzpicture}[anchorbase,>=latex,line join=bevel,scale=0.5,
every path/.style={very thick}]
\node (node_0) at (27.5bp,150.5bp) [draw,draw=none] {$\bullet$};
\node (node_4) at (27.5bp,79.5bp) [draw,draw=none] {$\bullet$};
\node (node_1) at (27.5bp,221.5bp) [draw,draw=none] {$\bullet$};
\node (node_2) at (27.5bp,292.5bp) [draw,draw=none] {$\bullet$};
\node (node_3) at (27.5bp,8.5bp) [draw,draw=none] {$\bullet$};
\draw [spinach,->] (node_0) ..controls (27.5bp,131.44bp) and (27.5bp,112.5bp)  .. (node_4);
\draw (36.0bp,115.0bp) node {$3$};
\draw [red,->] (node_1) ..controls (27.5bp,202.44bp) and (27.5bp,183.5bp)  .. (node_0);
\draw (36.0bp,186.0bp) node {$2$};
\draw [blue,->] (node_2) ..controls (27.5bp,273.44bp) and (27.5bp,254.5bp)  .. (node_1);
\draw (36.0bp,257.0bp) node {$1$};
\draw [black,->] (node_4) ..controls (27.5bp,60.442bp) and (27.5bp,41.496bp)  .. (node_3);
\draw (36.0bp,44.0bp) node {$4$};
\end{tikzpicture}$}
,\quad
\typeb[4]\colon
\scalebox{0.8}{$\begin{tikzpicture}[anchorbase,>=latex,line join=bevel,scale=0.45,every path/.style={very thick}]
\node (node_0) at (53.5bp,435.5bp) [draw,draw=none] {$\bullet$};
\node (node_8) at (53.5bp,364.5bp) [draw,draw=none] {$\bullet$};
\node (node_1) at (53.5bp,506.5bp) [draw,draw=none] {$\bullet$};
\node (node_2) at (53.5bp,577.5bp) [draw,draw=none] {$\bullet$};
\node (node_3) at (53.5bp,150.5bp) [draw,draw=none] {$\bullet$};
\node (node_5) at (53.5bp,79.5bp) [draw,draw=none] {$\bullet$};
\node (node_4) at (53.5bp,293.0bp) [draw,draw=none] {$\bullet$};
\node (node_6) at (53.5bp,221.5bp) [draw,draw=none] {$\bullet$};
\node (node_7) at (53.5bp,8.5bp) [draw,draw=none] {$\bullet$};
\draw [spinach,->] (node_0) ..controls (53.5bp,416.44bp) and (53.5bp,397.5bp)  .. (node_8);
\draw (62.0bp,400.0bp) node {$3$};
\draw [red,->] (node_1) ..controls (53.5bp,487.44bp) and (53.5bp,468.5bp)  .. (node_0);
\draw (62.0bp,471.0bp) node {$2$};
\draw [blue,->] (node_2) ..controls (53.5bp,558.44bp) and (53.5bp,539.5bp)  .. (node_1);
\draw (62.0bp,542.0bp) node {$1$};
\draw [red,->] (node_3) ..controls (53.5bp,131.44bp) and (53.5bp,112.5bp)  .. (node_5);
\draw (62.0bp,115.0bp) node {$2$};
\draw [black,->] (node_4) ..controls (53.5bp,273.19bp) and (53.5bp,254.59bp)  .. (node_6);
\draw (62.0bp,257.0bp) node {$4$};
\draw [blue,->] (node_5) ..controls (53.5bp,60.442bp) and (53.5bp,41.496bp)  .. (node_7);
\draw (62.0bp,44.0bp) node {$1$};
\draw [spinach,->] (node_6) ..controls (53.5bp,202.44bp) and (53.5bp,183.5bp)  .. (node_3);
\draw (62.0bp,186.0bp) node {$3$};
\draw [black,->] (node_8) ..controls (53.5bp,345.42bp) and (53.5bp,326.63bp)  .. (node_4);
\draw (62.0bp,329.0bp) node {$4$};
\end{tikzpicture}$}
,\quad
\typec[4]\colon
\scalebox{0.8}{$\begin{tikzpicture}[anchorbase,>=latex,line join=bevel,scale=0.5,
every path/.style={very thick}]
\node (node_0) at (27.5bp,150.5bp) [draw,draw=none] {$\bullet$};
\node (node_5) at (27.5bp,79.5bp) [draw,draw=none] {$\bullet$};
\node (node_1) at (27.5bp,221.5bp) [draw,draw=none] {$\bullet$};
\node (node_2) at (27.5bp,505.5bp) [draw,draw=none] {$\bullet$};
\node (node_7) at (27.5bp,434.5bp) [draw,draw=none] {$\bullet$};
\node (node_3) at (27.5bp,8.5bp) [draw,draw=none] {$\bullet$};
\node (node_4) at (27.5bp,363.5bp) [draw,draw=none] {$\bullet$};
\node (node_6) at (27.5bp,292.5bp) [draw,draw=none] {$\bullet$};
\draw [spinach,->] (node_0) ..controls (27.5bp,131.44bp) and (27.5bp,112.5bp)  .. (node_5);
\draw (36.0bp,115.0bp) node {$3$};
\draw [red,->] (node_1) ..controls (27.5bp,202.44bp) and (27.5bp,183.5bp)  .. (node_0);
\draw (36.0bp,186.0bp) node {$2$};
\draw [black,->] (node_2) ..controls (27.5bp,486.44bp) and (27.5bp,467.5bp)  .. (node_7);
\draw (36.0bp,470.0bp) node {$4$};
\draw [red,->] (node_4) ..controls (27.5bp,344.44bp) and (27.5bp,325.5bp)  .. (node_6);
\draw (36.0bp,328.0bp) node {$2$};
\draw [black,->] (node_5) ..controls (27.5bp,60.442bp) and (27.5bp,41.496bp)  .. (node_3);
\draw (36.0bp,44.0bp) node {$4$};
\draw [blue,->] (node_6) ..controls (27.5bp,273.44bp) and (27.5bp,254.5bp)  .. (node_1);
\draw (36.0bp,257.0bp) node {$1$};
\draw [spinach,->] (node_7) ..controls (27.5bp,415.44bp) and (27.5bp,396.5bp)  .. (node_4);
\draw (36.0bp,399.0bp) node {$3$};
\end{tikzpicture}$}
,\quad
\typed[4]\colon
\scalebox{0.8}{$\begin{tikzpicture}[anchorbase,>=latex,line join=bevel,xscale=0.65,yscale=0.5,
every path/.style={very thick}]
\node (node_0) at (61.5bp,150.5bp) [draw,draw=none] {$\bullet$};
\node (node_5) at (61.5bp,79.5bp) [draw,draw=none] {$\bullet$};
\node (node_1) at (61.5bp,363.5bp) [draw,draw=none] {$\bullet$};
\node (node_3) at (61.5bp,292.5bp) [draw,draw=none] {$\bullet$};
\node (node_2) at (61.5bp,434.5bp) [draw,draw=none] {$\bullet$};
\node (node_4) at (27.5bp,221.5bp) [draw,draw=none] {$\bullet$};
\node (node_6) at (96.5bp,221.5bp) [draw,draw=none] {$\bullet$};
\node (node_7) at (61.5bp,8.5bp) [draw,draw=none] {$\bullet$};
\draw [red,->] (node_0) ..controls (61.5bp,131.44bp) and (61.5bp,112.5bp)  .. (node_5);
\draw (70.0bp,115.0bp) node {$2$};
\draw [red,->] (node_1) ..controls (61.5bp,344.44bp) and (61.5bp,325.5bp)  .. (node_3);
\draw (70.0bp,328.0bp) node {$2$};
\draw [blue,->] (node_2) ..controls (61.5bp,415.44bp) and (61.5bp,396.5bp)  .. (node_1);
\draw (70.0bp,399.0bp) node {$1$};
\draw [spinach,->] (node_3) ..controls (52.5bp,273.24bp) and (42.895bp,253.74bp)  .. (node_4);
\draw (58.0bp,257.0bp) node {$3$};
\draw [black,->] (node_3) ..controls (70.764bp,273.24bp) and (80.653bp,253.74bp)  .. (node_6);
\draw (91.0bp,257.0bp) node {$4$};
\draw [black,->] (node_4) ..controls (36.5bp,202.24bp) and (46.105bp,182.74bp)  .. (node_0);
\draw (58.0bp,186.0bp) node {$4$};
\draw [blue,->] (node_5) ..controls (61.5bp,60.442bp) and (61.5bp,41.496bp)  .. (node_7);
\draw (70.0bp,44.0bp) node {$1$};
\draw [spinach,->] (node_6) ..controls (87.236bp,202.24bp) and (77.347bp,182.74bp)  .. (node_0);
\draw (91.0bp,186.0bp) node {$3$};
\end{tikzpicture}$}
.
\end{gather}
The associated representations are the vector representations.
For the classical types, these examples for the crystals of the vector representations generalize in the obvious way to all ranks.
\end{Example}

The crystal graph $\crystalgraph[\Lambda]$ is a directed graph. A \emph{path} in $\crystalgraph[\Lambda]$ is a directed path. A path in $\crystalgraph[\Lambda]$ is \emph{rooted} if it starts at the highest weight vector. The \emph{distance} $\dist[\lambda]\in\N$ of a vertex $\lambda\in\crystalgraph[\Lambda]$ from the highest weight is the length of any rooted
path to $\lambda$, which is the same for all paths. The vertices with $\dist[\lambda]=n$
form the \emph{$n$th layer} of $\crystalgraph[\Lambda]$. For example, the source vertex forms the $0$th layer.

\begin{Definition}\label{D:FiniteTypesResidueSequence}
For $\lambda\in\crystalgraph[\Lambda]$,
a \emph{path to $\lambda$} is a rooted path $\qpath[\lambda]$ in $\crystalgraph[\Lambda]$
that ends in $\lambda$.
A \emph{path of length~$n$} is a rooted path to some $\lambda\in\crystalgraph[\Lambda]$
with $\dist[\lambda]=n$.
Let $\Parts{n}{\Lambda}$ be the set of all paths in $\crystalgraph[\Lambda]$ of length $n$.
For every $\qpath[\lambda]\in\Parts{n}{\Lambda}$ its \emph{residue
sequence} $\Res[{\qpath[\lambda]}]$ is the ordered tuple of the labels of the edges
in $\lambda$.
\end{Definition}

\begin{Remark}\label{R:FiniteTypesResidueSequenceTwo}
We will see that rooted paths play the same role in the representation theory of wKLRW algebras as partitions do in the representation theory of the symmetric groups. We use $\qpath[\lambda]$ for a generic path to the vertex $\lambda$. Later we will fix a \emph{preferred path} $\ppath[\lambda]$ for each vertex.
\end{Remark}

We let $\sym[n]$ act on residue sequences of length $n$ by place permutations. Crystal theory implies that a rooted path $\bpath$ in $\crystalgraph[\Lambda]$ is uniquely determined by its residue sequence $\Res[\bpath]$. Hence, we identify a path with its residue sequence. Note that not every sequence in $\vertices^{n}$ corresponds to a path in $\crystalgraph[\Lambda]$.

\begin{Lemma}\label{L:FiniteTypesResidueSequence}
Let $\lambda\in\crystalgraph[\Lambda]$ be a vertex with $\dist[\lambda]=n$.
Then any two paths to $\lambda$ are in the same $\sym[n]$-orbit. That is,
their residue sequences are the same up to permutation.
\end{Lemma}

\begin{proof}
This follows directly from the definitions, see \cite[Definition 2.13, Axiom A2]{BuSc-crystal-bases}.
\end{proof}

As in \autoref{E:MainExampleTheBeastItself}, the permutations appearing in \autoref{L:FiniteTypesResidueSequence} play an important role in this paper, so we give them a special name.

\begin{Definition}\label{D:FiniteTypesDetour}
Let $\bpath$ and $\cpath$ be two paths of length $n$ in $\crystalgraph[\Lambda]$ to the same vertex. Then $w\in\sym[n]$ is a \emph{detour permutation} if $\bpath=w\cpath$.
\end{Definition}

\begin{Remark}\label{R:FiniteTypesResidueSequence}
Strictly speaking $w$ is a detour permutation for the pair of paths $(\bpath,\cpath)$ to the same vertex. As we have already seen in \autoref{E:MainExampleTheBeastItself}, it also happens that $\qpath[\lambda]=w\qpath[\mu]$, for some $w\in\sym[n]$, where $\qpath[\lambda]$ and $\qpath[\mu]$ to are paths in $\crystalgraph[\Lambda]$ to different vertices $\lambda\neq\mu$. Such permutations are not detour permutations. See also \autoref{E:FiniteTypesResidueSequenceMain} for another example.
\end{Remark}

\begin{Example}\label{E:FiniteTypesResidueSequenceTwo}
We again consider the crystal graph as in \autoref{E:FiniteTypesSageMath}.
There are three paths to $\lambda=[-\Lambda_{1}+\Lambda_{3}-\Lambda_{4}]$, giving three different residue sequences $53214$,
$53241$ and $53421$. These residue sequences are all permutations of each other by detour permutations. The crystal graph $\crystalgraph[{\fweight[5]}]$
has ten layers.
\end{Example}

\begin{Example}\label{E:FiniteTypesResidueSequenceMain}
In the main example in \autoref{S:MainExample} there are
three vertices $\lambda$, $\mu$ and $\nu$
in $\Parts{5}{\fweight[2]}$ and the crystal graph $\crystalgraph[{\fweight[2]}]$ has ten layers. \autoref{E:MainExampleTheBeastItself} considers
three paths of length five:
$\ppath[\lambda]$ ,$\ppath[\mu]$ and $\ppath[\nu]$ for the vertices $\lambda$, $\mu$ and $\nu$, respectively. These paths are related by the permutations that are explicitly listed in \autoref{E:MainExampleTheBeastItself}. As these paths correspond to distinct vertices of $\crystalgraph[{\fweight[2]}]$, none of the permutations between distinct paths are detour permutations. \autoref{E:MainExampleTheBeastItself}.(f) also describes a detour permutation from $\ppath[\lambda^{\prime}]$ to $\qpath[\lambda^{\prime}]$.
\end{Example}

Recall that every Dynkin diagram of finite type admits
a Dynkin diagram automorphism $(\placeholder)^{\dynkinaut}$ of order at most two
induced by $\sroot[i]\mapsto -w_{0}\sroot$, for $i\in\vertices$, where $w_{0}$
is the longest element of the associated Weyl group.

\begin{Example}\label{E:FiniteTypesFlip}
The automorphism $(\placeholder)^{\dynkinaut}$ in classical types
is as follows:
\begin{gather*}
\begin{array}{cc}
\toprule
\text{Quiver} & (\placeholder)^{\dynkinaut} \\ \midrule
\typea[e] & \text{flip $i\leftrightarrow e-i$} \\
\typeb[e>1],\typec[e>2] & \text{identity} \\
\typed[e>3] & \text{flip $(e-1)\leftrightarrow e$} \\
\bottomrule
\end{array}
.
\end{gather*}
Note $(\placeholder)^{\dynkinaut}$ induces a relabeling of $\vertices$.
\end{Example}

\begin{Lemma}\label{L:FiniteTypesFlip}
The graph $\crystalgraph[\Lambda]$ has a symmetry given by changing the residues using
$(\placeholder)^{\dynkinaut}$, reversing the orientation
of the arrows, and multiplying the labels (=weights) of the vertices by $-1$.
\end{Lemma}

\begin{proof}
This follows because the associated highest weight module $L(\Lambda)$ has this symmetry.
\end{proof}

\begin{Example}
For example, in \autoref{Eq:FiniteTypesSageMath},
$(\placeholder)^{\dynkinaut}$ swaps the resides $4$ and $5$, rotates around the middle of the crystal graph, which is $\set[\big]{[-\fweight[1]+\fweight[3]-\fweight[4]],[\fweight[1]-\fweight[3]+\fweight[5]]}$, reverses the orientation of the arrows and finally multiplies the vertex labels by $-1$.
\end{Example}

\iffalse
For higher levels we use:

\begin{Definition}\label{D:FiniteTypesMulti}
Let $\ell\in\N$.
For $\bfweight=(\fweight[{i_{1}}],\dots,\fweight[{i_{\ell}}])$
a \emph{$\bfweight$-vertex of distance $n$} is a choice of $\ell$ vertices
$\blam=(\lambda_{1},\dots,\lambda_{\ell})$ with $\lambda_{k}\in\crystalgraph[{\fweight[i_{k}]}]$ and $\sum_{k=1}^{\ell}\dist[\lambda_{k}]=n$. A \emph{path} $\qpath[\blam]$ to $\blam$ is a choice of $\ell$ paths
$(\qpath[\lambda_{1}],\dots,\qpath[\lambda_{\ell}])$ to these vertices.
Let $\Parts{n}{\bfweight}$ be the set of all paths of length $n$.
For every path $\qpath[\blam]\in\Parts{n}{\bfweight}$ its \emph{residue
sequence} is $\Res[{\qpath[\blam]}]=\big(\Res[{\qpath[\lambda_{1}]}],\dots,\Res[{\qpath[\lambda_{\ell}]}]\big)$.
\end{Definition}
\fi

%%%%%%%%%%%%%%%%%%%%%%%%%%%%%%%%%%%%%%%%%

\subsection{Stembridge--Sternberg relations}\label{SS:FiniteTypesCrystalsSSR}

%%%%%%%%%%%%%%%%%%%%%%%%%%%%%%%%%%%%%%%%%

We need to identify certain faces in crystal graphs, where
\emph{faces} are defined respect to any embedding
of the graph into a surface of high enough genus. Combinatorially:

\begin{Definition}\label{D:FiniteTypesSSR}
A \emph{two color face} in a crystal graph
$\crystalgraph*$ is an equation $a=b$ between two (not necessarily rooted) paths $a,b$ that are of minimal length such that they have the same end points and their residue sequences only contain the symbols $i,j\in\vertices$ with $i\neq j$. The \emph{source} and \emph{sink} of a two color face are the vertices of minimal and maximal distance, respectively, from the highest weight vector in $\crystalgraph*$.
\end{Definition}

As we will see in \autoref{P:FiniteTypesSSR} below, the only two color faces that appear in the fundamental crystals
are \emph{squares}, \emph{octagons}, \emph{decagons} and \emph{tetradecagon}. Hence, keeping the face picture in mind,
some examples of two color faces are:
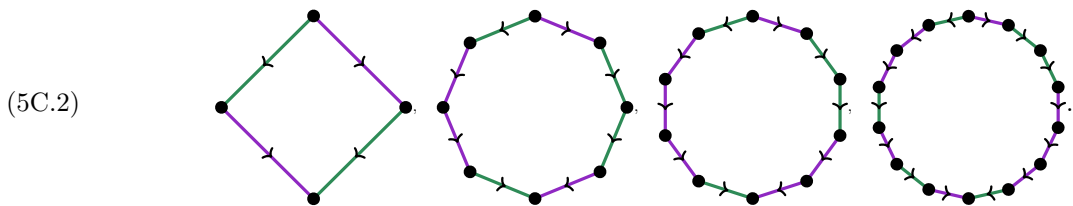
\begin{gather}\label{Eq:FiniteTypesSSR}
\scalebox{0.6}{$\begin{tikzpicture}[anchorbase,rotate=360/8,transform shape,
every path/.style={very thick}]
\node[regular polygon, regular polygon sides=4, minimum size=4cm,thick] (m) at (0,0) {};
\draw [spinach,line width=2.0,directed=0.55] (m.corner 1) -- (m.corner 2);
\draw [orchid,line width=2.0,directed=0.55] (m.corner 2) -- (m.corner 3);
\draw [orchid,line width=2.0,directed=0.55] (m.corner 1) -- (m.corner 4);
\draw [spinach,line width=2.0,directed=0.55] (m.corner 4) -- (m.corner 3);
\foreach \x in {1,...,4}{
\fill [black,thick] (m.corner \x) circle (4pt);
}
\end{tikzpicture}
,\quad
\begin{tikzpicture}[anchorbase,rotate=360/16,transform shape]
\node[regular polygon, regular polygon sides=8, minimum size=4cm,thick] (m) at (0,0) {};
\draw [spinach,line width=2.0,directed=0.55] (m.corner 1) -- (m.corner 2);
\draw [orchid,line width=2.0,directed=0.55] (m.corner 2) -- (m.corner 3);
\draw [orchid,line width=2.0,directed=0.55] (m.corner 3) -- (m.corner 4);
\draw [spinach,line width=2.0,directed=0.55] (m.corner 4) -- (m.corner 5);
\draw [orchid,line width=2.0,directed=0.55] (m.corner 1) -- (m.corner 8);
\draw [spinach,line width=2.0,directed=0.55] (m.corner 8) -- (m.corner 7);
\draw [spinach,line width=2.0,directed=0.55] (m.corner 7) -- (m.corner 6);
\draw [orchid,line width=2.0,directed=0.55] (m.corner 6) -- (m.corner 5);
\foreach \x in {1,...,8}{
\fill [black,thick] (m.corner \x) circle (4pt);
}
\end{tikzpicture}
,\quad
\begin{tikzpicture}[anchorbase,rotate=360/20,transform shape]
\node[regular polygon, regular polygon sides=10, minimum size=4cm,thick] (m) at (0,0) {};
\draw [spinach,line width=2.0,directed=0.55] (m.corner 1) -- (m.corner 2);
\draw [orchid,line width=2.0,directed=0.55] (m.corner 2) -- (m.corner 3);
\draw [orchid,line width=2.0,directed=0.55] (m.corner 3) -- (m.corner 4);
\draw [orchid,line width=2.0,directed=0.55] (m.corner 4) -- (m.corner 5);
\draw [spinach,line width=2.0,directed=0.55] (m.corner 5) -- (m.corner 6);
\draw [orchid,line width=2.0,directed=0.55] (m.corner 1) -- (m.corner 10);
\draw [spinach,line width=2.0,directed=0.55] (m.corner 10) -- (m.corner 9);
\draw [spinach,line width=2.0,directed=0.55] (m.corner 9) -- (m.corner 8);
\draw [orchid,line width=2.0,directed=0.55] (m.corner 8) -- (m.corner 7);
\draw [orchid,line width=2.0,directed=0.55] (m.corner 7) -- (m.corner 6);
\foreach \x in {1,...,10}{
\fill [black,thick] (m.corner \x) circle (4pt);
}
\end{tikzpicture}
,\quad
\begin{tikzpicture}[anchorbase,rotate=360/28,transform shape]
\node[regular polygon, regular polygon sides=14, minimum size=4cm,thick] (m) at (0,0) {};
\draw [spinach,line width=2.0,directed=0.55] (m.corner 1) -- (m.corner 2);
\draw [orchid,line width=2.0,directed=0.55] (m.corner 2) -- (m.corner 3);
\draw [orchid,line width=2.0,directed=0.55] (m.corner 3) -- (m.corner 4);
\draw [spinach,line width=2.0,directed=0.55] (m.corner 4) -- (m.corner 5);
\draw [orchid,line width=2.0,directed=0.55] (m.corner 5) -- (m.corner 6);
\draw [spinach,line width=2.0,directed=0.55] (m.corner 6) -- (m.corner 7);
\draw [orchid,line width=2.0,directed=0.55] (m.corner 7) -- (m.corner 8);
\draw [orchid,line width=2.0,directed=0.55] (m.corner 1) -- (m.corner 14);
\draw [spinach,line width=2.0,directed=0.55] (m.corner 14) -- (m.corner 13);
\draw [spinach,line width=2.0,directed=0.55] (m.corner 13) -- (m.corner 12);
\draw [orchid,line width=2.0,directed=0.55] (m.corner 12) -- (m.corner 11);
\draw [orchid,line width=2.0,directed=0.55] (m.corner 11) -- (m.corner 10);
\draw [orchid,line width=2.0,directed=0.55] (m.corner 10) -- (m.corner 9);
\draw [spinach,line width=2.0,directed=0.55] (m.corner 9) -- (m.corner 8);
\foreach \x in {1,...,14}{
\fill [black,thick] (m.corner \x) circle (4pt);
}
\end{tikzpicture}
$}.
\end{gather}
Other colorings of the edges around a two color face are allowed, with the two colors appearing in different orders. The crystal graph machinery ensures that the same number of edges of color $i$ and color $j$ appear in the two paths from the source to the sink.

\begin{Example}\label{E:FiniteTypesSSR}
In \autoref{Eq:FiniteTypesSageMath} the two color faces are $24=42$, $14=41$, $13=31$, $15=51$ and $25=52$.
\autoref{Eq:MainExample}(f) contains the two octogons $2332=3223$ and $1221=2112$, in addition to four squares.
\end{Example}

The next result is motivated by \cite{St-local-crystal-simply-laced} and \cite{St-local-crystal-doubly-laced}. In \autoref{SS:Plactic} below, we will impose additional constraints on the adjacent squares in types $\typee[i]$
and $\typef$.

\begin{Proposition}\label{P:FiniteTypesSSR}
Fix a crystal graph $\crystalgraph$, for some fundamental weight $\fweight[i]$.
\begin{enumerate}[nosep]
\item For type $\typea[e]$, the only two color faces are $ij=ji$ for nonadjacent $i,j\in\vertices$ (a nonadjacent square).

\item For types $\typeb$, $\typec$, $\typed$, the only two color faces are:

\begin{enumerate}

\item $ij=ji$ for nonadjacent $i,j\in\vertices$ (a nonadjacent square), or $ij=ji$ for adjacent $i,j\in\vertices$ (an adjacent square) with $i\not\!\Rightarrow j$ and $i\not\!\leftarrow j$ and there is a preceding $i$ or $j$, up to nonadjacent permutations.

\item $ijji=jiij$ for adjacent $i,j\in\vertices$ (an octagon).

\end{enumerate}

\item For type $\typee[i]$, the only two color faces are:

\begin{enumerate}

\item $ij=ji$ for nonadjacent $i,j\in\vertices$ (a nonadjacent square), or $ij=ji$ for adjacent $i,j\in\vertices$ (an adjacent square). We do not necessarily have a preceding $i$ or $j$, up to nonadjacent permutations.

\item $ijji=jiij$ for adjacent $i,j\in\vertices$ (an octagon).

\end{enumerate}

\item For type $\typef[4]$, the only two color faces are:

\begin{enumerate}

\item $ij=ji$ for nonadjacent $i,j\in\vertices$ (a nonadjacent square), or $ij=ji$ for adjacent $i,j\in\vertices$ (an adjacent square). We do not necessarily have a preceding $i$ or $j$, up to nonadjacent permutations.

\item $ijji=jiij$ for adjacent $i,j\in\vertices$ (an octagon).

\item $23332=32233=32323$ (a decagon).

\item $2332323=3223332=3232332=2333223$ (a tetradecagon).

\end{enumerate}

\item For type $\typeg[2]$, the only two color face is $1221=2112$ (an octagon).

\end{enumerate}
\end{Proposition}

\begin{proof}
\textit{Classical types.}
For minuscule weights the statement is clear
by their construction in, for example,
\cite[Section 5.4]{BuSc-crystal-bases}. So we only discuss the other fundamental weights.

We will use the Young diagrams combinatorics as
in {\eg} \cite[Section 6.2]{BuSc-crystal-bases}.
The combinatorics works as follows. (We do not need to consider
type $\typea$, since all weights are minuscule, but we list it for completeness.) We fix the
partially ordered \emph{filling sets}
\begin{gather*}
\typea[e]\colon
\set{1<\dots<e}
,\quad
\typeb[e>1]\colon
\set{1<\dots<e<0<\bar{e}<\dots<\bar{1}}
,
\\
\typec[e>2]\colon
\set{e<\dots<1<\bar{1}<\dots<\bar{e}}
,\quad
\typed[e>3]\colon
\set{1<\dots<e,\bar{e}<\dots<\bar{1}}
,
\end{gather*}
which we use to fill tableaux.
Note that $e$ and $\bar{e}$ are not comparable
in type $\typed[e>3]$.
A Young diagram of shape $(1^{k})$ is a column of height $k$, and we
fill these with numbers from the filling sets such that:
\begin{enumerate}[label=(\roman*)]

\item The entries are strictly increasing from top to bottom with the exception that the letter $0$ in $\typeb[e>1]$ can be repeated once, and
the letters $e$ and $\bar{e}$ in type $\typed[e>3]$ can alternate once.

\item If both letters $i$ and
$\bar{i}$ appear, and $i$ is in the $a$th node
and $\bar{i}$ is in the $b$th node from the top $a+b\leq i$.

\end{enumerate}
See \cite[Section 6.2]{BuSc-crystal-bases} for several examples (but note that our labeling of the Dynkin diagrams of type $\typec[e>2]$ is different).

The Young diagrams of shape $(1^{k})$ with the above fillings
correspond to the vertices of $\crystalgraph[{\fweight[k]}]$,
except in type $\typec[e>2]$ where this corresponds to $\crystalgraph[{\fweight[e-k+1]}]$.
The arrows in these graphs are determined by the following replacement rules, which generalize to larger $e$ accordingly):
\begin{gather}\label{Eq:ProofsCrystals}
\scalebox{0.6}{$\typea[5]\colon
\begin{tikzpicture}[anchorbase,>=latex,line join=bevel,scale=0.5,
every path/.style={very thick}]
\node (node_0) at (8.5bp,9.5bp) [draw,draw=none] {${\def\lr#1{\multicolumn{1}{|@{\hspace{.6ex}}c@{\hspace{.6ex}}|}{\raisebox{-.3ex}{$#1$}}}\raisebox{-.6ex}{$\begin{array}[b]{*{1}c}\cline{1-1}\lr{6}\\\cline{1-1}\end{array}$}}$};
\node (node_1) at (8.5bp,301.5bp) [draw,draw=none] {${\def\lr#1{\multicolumn{1}{|@{\hspace{.6ex}}c@{\hspace{.6ex}}|}{\raisebox{-.3ex}{$#1$}}}\raisebox{-.6ex}{$\begin{array}[b]{*{1}c}\cline{1-1}\lr{2}\\\cline{1-1}\end{array}$}}$};
\node (node_5) at (8.5bp,228.5bp) [draw,draw=none] {${\def\lr#1{\multicolumn{1}{|@{\hspace{.6ex}}c@{\hspace{.6ex}}|}{\raisebox{-.3ex}{$#1$}}}\raisebox{-.6ex}{$\begin{array}[b]{*{1}c}\cline{1-1}\lr{3}\\\cline{1-1}\end{array}$}}$};
\node (node_2) at (8.5bp,82.5bp) [draw,draw=none] {${\def\lr#1{\multicolumn{1}{|@{\hspace{.6ex}}c@{\hspace{.6ex}}|}{\raisebox{-.3ex}{$#1$}}}\raisebox{-.6ex}{$\begin{array}[b]{*{1}c}\cline{1-1}\lr{5}\\\cline{1-1}\end{array}$}}$};
\node (node_3) at (8.5bp,155.5bp) [draw,draw=none] {${\def\lr#1{\multicolumn{1}{|@{\hspace{.6ex}}c@{\hspace{.6ex}}|}{\raisebox{-.3ex}{$#1$}}}\raisebox{-.6ex}{$\begin{array}[b]{*{1}c}\cline{1-1}\lr{4}\\\cline{1-1}\end{array}$}}$};
\node (node_4) at (8.5bp,374.5bp) [draw,draw=none] {${\def\lr#1{\multicolumn{1}{|@{\hspace{.6ex}}c@{\hspace{.6ex}}|}{\raisebox{-.3ex}{$#1$}}}\raisebox{-.6ex}{$\begin{array}[b]{*{1}c}\cline{1-1}\lr{1}\\\cline{1-1}\end{array}$}}$};
\draw [red,->] (node_1) ..controls (8.5bp,281.04bp) and (8.5bp,262.45bp)  .. (node_5);
\draw (17.0bp,265.0bp) node {$2$};
\draw [brown,->] (node_2) ..controls (8.5bp,62.042bp) and (8.5bp,43.449bp)  .. (node_0);
\draw (17.0bp,46.0bp) node {$5$};
\draw [black,->] (node_3) ..controls (8.5bp,135.04bp) and (8.5bp,116.45bp)  .. (node_2);
\draw (17.0bp,119.0bp) node {$4$};
\draw [blue,->] (node_4) ..controls (8.5bp,354.04bp) and (8.5bp,335.45bp)  .. (node_1);
\draw (17.0bp,338.0bp) node {$1$};
\draw [spinach,->] (node_5) ..controls (8.5bp,208.04bp) and (8.5bp,189.45bp)  .. (node_3);
\draw (17.0bp,192.0bp) node {$3$};
\end{tikzpicture}
,\quad
\typeb[5]\colon
\begin{tikzpicture}[anchorbase,>=latex,line join=bevel,scale=0.5,every path/.style={very thick}]
\node (node_0) at (8.5bp,228.5bp) [draw,draw=none] {${\def\lr#1{\multicolumn{1}{|@{\hspace{.6ex}}c@{\hspace{.6ex}}|}{\raisebox{-.3ex}{$#1$}}}\raisebox{-.6ex}{$\begin{array}[b]{*{1}c}\cline{1-1}\lr{\overline{4}}\\\cline{1-1}\end{array}$}}$};
\node (node_8) at (8.5bp,155.5bp) [draw,draw=none] {${\def\lr#1{\multicolumn{1}{|@{\hspace{.6ex}}c@{\hspace{.6ex}}|}{\raisebox{-.3ex}{$#1$}}}\raisebox{-.6ex}{$\begin{array}[b]{*{1}c}\cline{1-1}\lr{\overline{3}}\\\cline{1-1}\end{array}$}}$};
\node (node_1) at (8.5bp,666.5bp) [draw,draw=none] {${\def\lr#1{\multicolumn{1}{|@{\hspace{.6ex}}c@{\hspace{.6ex}}|}{\raisebox{-.3ex}{$#1$}}}\raisebox{-.6ex}{$\begin{array}[b]{*{1}c}\cline{1-1}\lr{2}\\\cline{1-1}\end{array}$}}$};
\node (node_10) at (8.5bp,593.5bp) [draw,draw=none] {${\def\lr#1{\multicolumn{1}{|@{\hspace{.6ex}}c@{\hspace{.6ex}}|}{\raisebox{-.3ex}{$#1$}}}\raisebox{-.6ex}{$\begin{array}[b]{*{1}c}\cline{1-1}\lr{3}\\\cline{1-1}\end{array}$}}$};
\node (node_2) at (8.5bp,447.5bp) [draw,draw=none] {${\def\lr#1{\multicolumn{1}{|@{\hspace{.6ex}}c@{\hspace{.6ex}}|}{\raisebox{-.3ex}{$#1$}}}\raisebox{-.6ex}{$\begin{array}[b]{*{1}c}\cline{1-1}\lr{5}\\\cline{1-1}\end{array}$}}$};
\node (node_9) at (8.5bp,374.5bp) [draw,draw=none] {${\def\lr#1{\multicolumn{1}{|@{\hspace{.6ex}}c@{\hspace{.6ex}}|}{\raisebox{-.3ex}{$#1$}}}\raisebox{-.6ex}{$\begin{array}[b]{*{1}c}\cline{1-1}\lr{0}\\\cline{1-1}\end{array}$}}$};
\node (node_3) at (8.5bp,301.5bp) [draw,draw=none] {${\def\lr#1{\multicolumn{1}{|@{\hspace{.6ex}}c@{\hspace{.6ex}}|}{\raisebox{-.3ex}{$#1$}}}\raisebox{-.6ex}{$\begin{array}[b]{*{1}c}\cline{1-1}\lr{\overline{5}}\\\cline{1-1}\end{array}$}}$};
\node (node_4) at (8.5bp,82.5bp) [draw,draw=none] {${\def\lr#1{\multicolumn{1}{|@{\hspace{.6ex}}c@{\hspace{.6ex}}|}{\raisebox{-.3ex}{$#1$}}}\raisebox{-.6ex}{$\begin{array}[b]{*{1}c}\cline{1-1}\lr{\overline{2}}\\\cline{1-1}\end{array}$}}$};
\node (node_6) at (8.5bp,9.5bp) [draw,draw=none] {${\def\lr#1{\multicolumn{1}{|@{\hspace{.6ex}}c@{\hspace{.6ex}}|}{\raisebox{-.3ex}{$#1$}}}\raisebox{-.6ex}{$\begin{array}[b]{*{1}c}\cline{1-1}\lr{\overline{1}}\\\cline{1-1}\end{array}$}}$};
\node (node_5) at (8.5bp,520.5bp) [draw,draw=none] {${\def\lr#1{\multicolumn{1}{|@{\hspace{.6ex}}c@{\hspace{.6ex}}|}{\raisebox{-.3ex}{$#1$}}}\raisebox{-.6ex}{$\begin{array}[b]{*{1}c}\cline{1-1}\lr{4}\\\cline{1-1}\end{array}$}}$};
\node (node_7) at (8.5bp,739.5bp) [draw,draw=none] {${\def\lr#1{\multicolumn{1}{|@{\hspace{.6ex}}c@{\hspace{.6ex}}|}{\raisebox{-.3ex}{$#1$}}}\raisebox{-.6ex}{$\begin{array}[b]{*{1}c}\cline{1-1}\lr{1}\\\cline{1-1}\end{array}$}}$};
\draw [spinach,->] (node_0) ..controls (8.5bp,208.04bp) and (8.5bp,189.45bp)  .. (node_8);
\draw (17.0bp,192.0bp) node {$3$};
\draw [red,->] (node_1) ..controls (8.5bp,646.04bp) and (8.5bp,627.45bp)  .. (node_10);
\draw (17.0bp,630.0bp) node {$2$};
\draw [brown,->] (node_2) ..controls (8.5bp,427.04bp) and (8.5bp,408.45bp)  .. (node_9);
\draw (17.0bp,411.0bp) node {$5$};
\draw [black,->] (node_3) ..controls (8.5bp,281.04bp) and (8.5bp,262.45bp)  .. (node_0);
\draw (17.0bp,265.0bp) node {$4$};
\draw [blue,->] (node_4) ..controls (8.5bp,62.042bp) and (8.5bp,43.449bp)  .. (node_6);
\draw (17.0bp,46.0bp) node {$1$};
\draw [black,->] (node_5) ..controls (8.5bp,500.04bp) and (8.5bp,481.45bp)  .. (node_2);
\draw (17.0bp,484.0bp) node {$4$};
\draw [blue,->] (node_7) ..controls (8.5bp,719.04bp) and (8.5bp,700.45bp)  .. (node_1);
\draw (17.0bp,703.0bp) node {$1$};
\draw [red,->] (node_8) ..controls (8.5bp,135.04bp) and (8.5bp,116.45bp)  .. (node_4);
\draw (17.0bp,119.0bp) node {$2$};
\draw [brown,->] (node_9) ..controls (8.5bp,354.04bp) and (8.5bp,335.45bp)  .. (node_3);
\draw (17.0bp,338.0bp) node {$5$};
\draw [spinach,->] (node_10) ..controls (8.5bp,573.04bp) and (8.5bp,554.45bp)  .. (node_5);
\draw (17.0bp,557.0bp) node {$3$};
\end{tikzpicture}
,\quad
\typec[5]\colon
\begin{tikzpicture}[anchorbase,>=latex,line join=bevel,scale=0.5,every path/.style={very thick}]
\node (node_0) at (8.5bp,228.5bp) [draw,draw=none] {${\def\lr#1{\multicolumn{1}{|@{\hspace{.6ex}}c@{\hspace{.6ex}}|}{\raisebox{-.3ex}{$#1$}}}\raisebox{-.6ex}{$\begin{array}[b]{*{1}c}\cline{1-1}\lr{\overline{2}}\\\cline{1-1}\end{array}$}}$};
\node (node_8) at (8.5bp,155.5bp) [draw,draw=none] {${\def\lr#1{\multicolumn{1}{|@{\hspace{.6ex}}c@{\hspace{.6ex}}|}{\raisebox{-.3ex}{$#1$}}}\raisebox{-.6ex}{$\begin{array}[b]{*{1}c}\cline{1-1}\lr{\overline{3}}\\\cline{1-1}\end{array}$}}$};
\node (node_1) at (8.5bp,593.5bp) [draw,draw=none] {${\def\lr#1{\multicolumn{1}{|@{\hspace{.6ex}}c@{\hspace{.6ex}}|}{\raisebox{-.3ex}{$#1$}}}\raisebox{-.6ex}{$\begin{array}[b]{*{1}c}\cline{1-1}\lr{4}\\\cline{1-1}\end{array}$}}$};
\node (node_9) at (8.5bp,520.5bp) [draw,draw=none] {${\def\lr#1{\multicolumn{1}{|@{\hspace{.6ex}}c@{\hspace{.6ex}}|}{\raisebox{-.3ex}{$#1$}}}\raisebox{-.6ex}{$\begin{array}[b]{*{1}c}\cline{1-1}\lr{3}\\\cline{1-1}\end{array}$}}$};
\node (node_2) at (8.5bp,374.5bp) [draw,draw=none] {${\def\lr#1{\multicolumn{1}{|@{\hspace{.6ex}}c@{\hspace{.6ex}}|}{\raisebox{-.3ex}{$#1$}}}\raisebox{-.6ex}{$\begin{array}[b]{*{1}c}\cline{1-1}\lr{1}\\\cline{1-1}\end{array}$}}$};
\node (node_3) at (8.5bp,301.5bp) [draw,draw=none] {${\def\lr#1{\multicolumn{1}{|@{\hspace{.6ex}}c@{\hspace{.6ex}}|}{\raisebox{-.3ex}{$#1$}}}\raisebox{-.6ex}{$\begin{array}[b]{*{1}c}\cline{1-1}\lr{\overline{1}}\\\cline{1-1}\end{array}$}}$};
\node (node_4) at (8.5bp,82.5bp) [draw,draw=none] {${\def\lr#1{\multicolumn{1}{|@{\hspace{.6ex}}c@{\hspace{.6ex}}|}{\raisebox{-.3ex}{$#1$}}}\raisebox{-.6ex}{$\begin{array}[b]{*{1}c}\cline{1-1}\lr{\overline{4}}\\\cline{1-1}\end{array}$}}$};
\node (node_6) at (8.5bp,9.5bp) [draw,draw=none] {${\def\lr#1{\multicolumn{1}{|@{\hspace{.6ex}}c@{\hspace{.6ex}}|}{\raisebox{-.3ex}{$#1$}}}\raisebox{-.6ex}{$\begin{array}[b]{*{1}c}\cline{1-1}\lr{\overline{5}}\\\cline{1-1}\end{array}$}}$};
\node (node_5) at (8.5bp,447.5bp) [draw,draw=none] {${\def\lr#1{\multicolumn{1}{|@{\hspace{.6ex}}c@{\hspace{.6ex}}|}{\raisebox{-.3ex}{$#1$}}}\raisebox{-.6ex}{$\begin{array}[b]{*{1}c}\cline{1-1}\lr{2}\\\cline{1-1}\end{array}$}}$};
\node (node_7) at (8.5bp,666.5bp) [draw,draw=none] {${\def\lr#1{\multicolumn{1}{|@{\hspace{.6ex}}c@{\hspace{.6ex}}|}{\raisebox{-.3ex}{$#1$}}}\raisebox{-.6ex}{$\begin{array}[b]{*{1}c}\cline{1-1}\lr{5}\\\cline{1-1}\end{array}$}}$};
\draw [spinach,->] (node_0) ..controls (8.5bp,208.04bp) and (8.5bp,189.45bp)  .. (node_8);
\draw (17.0bp,192.0bp) node {$3$};
\draw [red,->] (node_1) ..controls (8.5bp,573.04bp) and (8.5bp,554.45bp)  .. (node_9);
\draw (17.0bp,557.0bp) node {$4$};
\draw [brown,->] (node_2) ..controls (8.5bp,354.04bp) and (8.5bp,335.45bp)  .. (node_3);
\draw (17.0bp,338.0bp) node {$1$};
\draw [black,->] (node_3) ..controls (8.5bp,281.04bp) and (8.5bp,262.45bp)  .. (node_0);
\draw (17.0bp,265.0bp) node {$2$};
\draw [blue,->] (node_4) ..controls (8.5bp,62.042bp) and (8.5bp,43.449bp)  .. (node_6);
\draw (17.0bp,46.0bp) node {$5$};
\draw [black,->] (node_5) ..controls (8.5bp,427.04bp) and (8.5bp,408.45bp)  .. (node_2);
\draw (17.0bp,411.0bp) node {$2$};
\draw [blue,->] (node_7) ..controls (8.5bp,646.04bp) and (8.5bp,627.45bp)  .. (node_1);
\draw (17.0bp,630.0bp) node {$5$};
\draw [red,->] (node_8) ..controls (8.5bp,135.04bp) and (8.5bp,116.45bp)  .. (node_4);
\draw (17.0bp,119.0bp) node {$4$};
\draw [spinach,->] (node_9) ..controls (8.5bp,500.04bp) and (8.5bp,481.45bp)  .. (node_5);
\draw (17.0bp,484.0bp) node {$3$};
\end{tikzpicture}
,\quad
\typed[5]\colon
\begin{tikzpicture}[anchorbase,>=latex,line join=bevel,scale=0.6,every path/.style={very thick}]
\node (node_0) at (28.5bp,228.5bp) [draw,draw=none] {${\def\lr#1{\multicolumn{1}{|@{\hspace{.6ex}}c@{\hspace{.6ex}}|}{\raisebox{-.3ex}{$#1$}}}\raisebox{-.6ex}{$\begin{array}[b]{*{1}c}\cline{1-1}\lr{\overline{4}}\\\cline{1-1}\end{array}$}}$};
\node (node_8) at (28.5bp,155.5bp) [draw,draw=none] {${\def\lr#1{\multicolumn{1}{|@{\hspace{.6ex}}c@{\hspace{.6ex}}|}{\raisebox{-.3ex}{$#1$}}}\raisebox{-.6ex}{$\begin{array}[b]{*{1}c}\cline{1-1}\lr{\overline{3}}\\\cline{1-1}\end{array}$}}$};
\node (node_1) at (28.5bp,520.5bp) [draw,draw=none] {${\def\lr#1{\multicolumn{1}{|@{\hspace{.6ex}}c@{\hspace{.6ex}}|}{\raisebox{-.3ex}{$#1$}}}\raisebox{-.6ex}{$\begin{array}[b]{*{1}c}\cline{1-1}\lr{2}\\\cline{1-1}\end{array}$}}$};
\node (node_9) at (28.5bp,447.5bp) [draw,draw=none] {${\def\lr#1{\multicolumn{1}{|@{\hspace{.6ex}}c@{\hspace{.6ex}}|}{\raisebox{-.3ex}{$#1$}}}\raisebox{-.6ex}{$\begin{array}[b]{*{1}c}\cline{1-1}\lr{3}\\\cline{1-1}\end{array}$}}$};
\node (node_2) at (8.5bp,301.5bp) [draw,draw=none] {${\def\lr#1{\multicolumn{1}{|@{\hspace{.6ex}}c@{\hspace{.6ex}}|}{\raisebox{-.3ex}{$#1$}}}\raisebox{-.6ex}{$\begin{array}[b]{*{1}c}\cline{1-1}\lr{5}\\\cline{1-1}\end{array}$}}$};
\node (node_3) at (48.5bp,301.5bp) [draw,draw=none] {${\def\lr#1{\multicolumn{1}{|@{\hspace{.6ex}}c@{\hspace{.6ex}}|}{\raisebox{-.3ex}{$#1$}}}\raisebox{-.6ex}{$\begin{array}[b]{*{1}c}\cline{1-1}\lr{\overline{5}}\\\cline{1-1}\end{array}$}}$};
\node (node_4) at (28.5bp,82.5bp) [draw,draw=none] {${\def\lr#1{\multicolumn{1}{|@{\hspace{.6ex}}c@{\hspace{.6ex}}|}{\raisebox{-.3ex}{$#1$}}}\raisebox{-.6ex}{$\begin{array}[b]{*{1}c}\cline{1-1}\lr{\overline{2}}\\\cline{1-1}\end{array}$}}$};
\node (node_6) at (28.5bp,9.5bp) [draw,draw=none] {${\def\lr#1{\multicolumn{1}{|@{\hspace{.6ex}}c@{\hspace{.6ex}}|}{\raisebox{-.3ex}{$#1$}}}\raisebox{-.6ex}{$\begin{array}[b]{*{1}c}\cline{1-1}\lr{\overline{1}}\\\cline{1-1}\end{array}$}}$};
\node (node_5) at (28.5bp,374.5bp) [draw,draw=none] {${\def\lr#1{\multicolumn{1}{|@{\hspace{.6ex}}c@{\hspace{.6ex}}|}{\raisebox{-.3ex}{$#1$}}}\raisebox{-.6ex}{$\begin{array}[b]{*{1}c}\cline{1-1}\lr{4}\\\cline{1-1}\end{array}$}}$};
\node (node_7) at (28.5bp,593.5bp) [draw,draw=none] {${\def\lr#1{\multicolumn{1}{|@{\hspace{.6ex}}c@{\hspace{.6ex}}|}{\raisebox{-.3ex}{$#1$}}}\raisebox{-.6ex}{$\begin{array}[b]{*{1}c}\cline{1-1}\lr{1}\\\cline{1-1}\end{array}$}}$};
\draw [spinach,->] (node_0) ..controls (28.5bp,208.04bp) and (28.5bp,189.45bp)  .. (node_8);
\draw (37.0bp,192.0bp) node {$3$};
\draw [red,->] (node_1) ..controls (28.5bp,500.04bp) and (28.5bp,481.45bp)  .. (node_9);
\draw (37.0bp,484.0bp) node {$2$};
\draw [brown,->] (node_2) ..controls (6.188bp,282.54bp) and (5.5353bp,267.87bp)  .. (9.5bp,256.0bp) .. controls (10.612bp,252.67bp) and (12.299bp,249.4bp)  .. (node_0);
\draw (18.0bp,265.0bp) node {$5$};
\draw [black,->] (node_3) ..controls (42.989bp,280.94bp) and (37.677bp,262.08bp)  .. (node_0);
\draw (49.0bp,265.0bp) node {$4$};
\draw [blue,->] (node_4) ..controls (28.5bp,62.042bp) and (28.5bp,43.449bp)  .. (node_6);
\draw (37.0bp,46.0bp) node {$1$};
\draw [black,->] (node_5) ..controls (16.25bp,360.08bp) and (11.723bp,353.66bp)  .. (9.5bp,347.0bp) .. controls (6.8052bp,338.93bp) and (6.2436bp,329.57bp)  .. (node_2);
\draw (18.0bp,338.0bp) node {$4$};
\draw [brown,->] (node_5) ..controls (34.011bp,353.94bp) and (39.323bp,335.08bp)  .. (node_3);
\draw (49.0bp,338.0bp) node {$5$};
\draw [blue,->] (node_7) ..controls (28.5bp,573.04bp) and (28.5bp,554.45bp)  .. (node_1);
\draw (37.0bp,557.0bp) node {$1$};
\draw [red,->] (node_8) ..controls (28.5bp,135.04bp) and (28.5bp,116.45bp)  .. (node_4);
\draw (37.0bp,119.0bp) node {$2$};
\draw [spinach,->] (node_9) ..controls (28.5bp,427.04bp) and (28.5bp,408.45bp)  .. (node_5);
\draw (37.0bp,411.0bp) node {$3$};
\end{tikzpicture}
,\quad\typea[5]\colon
\begin{tikzpicture}[anchorbase,>=latex,line join=bevel,scale=0.55,every path/.style={very thick}]
\node (node_0) at (98.485bp,268.93bp) [draw,draw=none] {${\def\lr#1{\multicolumn{1}{|@{\hspace{.6ex}}c@{\hspace{.6ex}}|}{\raisebox{-.3ex}{$#1$}}}\raisebox{-.6ex}{$\begin{array}[b]{*{1}c}\cline{1-1}\lr{2}\\\cline{1-1}\lr{6}\\\cline{1-1}\end{array}$}}$};
\node (node_12) at (58.485bp,184.38bp) [draw,draw=none] {${\def\lr#1{\multicolumn{1}{|@{\hspace{.6ex}}c@{\hspace{.6ex}}|}{\raisebox{-.3ex}{$#1$}}}\raisebox{-.6ex}{$\begin{array}[b]{*{1}c}\cline{1-1}\lr{3}\\\cline{1-1}\lr{6}\\\cline{1-1}\end{array}$}}$};
\node (node_1) at (38.485bp,691.69bp) [draw,draw=none] {${\def\lr#1{\multicolumn{1}{|@{\hspace{.6ex}}c@{\hspace{.6ex}}|}{\raisebox{-.3ex}{$#1$}}}\raisebox{-.6ex}{$\begin{array}[b]{*{1}c}\cline{1-1}\lr{1}\\\cline{1-1}\lr{2}\\\cline{1-1}\end{array}$}}$};
\node (node_4) at (38.485bp,607.14bp) [draw,draw=none] {${\def\lr#1{\multicolumn{1}{|@{\hspace{.6ex}}c@{\hspace{.6ex}}|}{\raisebox{-.3ex}{$#1$}}}\raisebox{-.6ex}{$\begin{array}[b]{*{1}c}\cline{1-1}\lr{1}\\\cline{1-1}\lr{3}\\\cline{1-1}\end{array}$}}$};
\node (node_2) at (38.485bp,15.276bp) [draw,draw=none] {${\def\lr#1{\multicolumn{1}{|@{\hspace{.6ex}}c@{\hspace{.6ex}}|}{\raisebox{-.3ex}{$#1$}}}\raisebox{-.6ex}{$\begin{array}[b]{*{1}c}\cline{1-1}\lr{5}\\\cline{1-1}\lr{6}\\\cline{1-1}\end{array}$}}$};
\node (node_3) at (8.4845bp,353.48bp) [draw,draw=none] {${\def\lr#1{\multicolumn{1}{|@{\hspace{.6ex}}c@{\hspace{.6ex}}|}{\raisebox{-.3ex}{$#1$}}}\raisebox{-.6ex}{$\begin{array}[b]{*{1}c}\cline{1-1}\lr{3}\\\cline{1-1}\lr{4}\\\cline{1-1}\end{array}$}}$};
\node (node_10) at (38.485bp,268.93bp) [draw,draw=none] {${\def\lr#1{\multicolumn{1}{|@{\hspace{.6ex}}c@{\hspace{.6ex}}|}{\raisebox{-.3ex}{$#1$}}}\raisebox{-.6ex}{$\begin{array}[b]{*{1}c}\cline{1-1}\lr{3}\\\cline{1-1}\lr{5}\\\cline{1-1}\end{array}$}}$};
\node (node_11) at (18.485bp,522.58bp) [draw,draw=none] {${\def\lr#1{\multicolumn{1}{|@{\hspace{.6ex}}c@{\hspace{.6ex}}|}{\raisebox{-.3ex}{$#1$}}}\raisebox{-.6ex}{$\begin{array}[b]{*{1}c}\cline{1-1}\lr{2}\\\cline{1-1}\lr{3}\\\cline{1-1}\end{array}$}}$};
\node (node_13) at (58.485bp,522.58bp) [draw,draw=none] {${\def\lr#1{\multicolumn{1}{|@{\hspace{.6ex}}c@{\hspace{.6ex}}|}{\raisebox{-.3ex}{$#1$}}}\raisebox{-.6ex}{$\begin{array}[b]{*{1}c}\cline{1-1}\lr{1}\\\cline{1-1}\lr{4}\\\cline{1-1}\end{array}$}}$};
\node (node_5) at (128.48bp,353.48bp) [draw,draw=none] {${\def\lr#1{\multicolumn{1}{|@{\hspace{.6ex}}c@{\hspace{.6ex}}|}{\raisebox{-.3ex}{$#1$}}}\raisebox{-.6ex}{$\begin{array}[b]{*{1}c}\cline{1-1}\lr{1}\\\cline{1-1}\lr{6}\\\cline{1-1}\end{array}$}}$};
\node (node_6) at (18.485bp,184.38bp) [draw,draw=none] {${\def\lr#1{\multicolumn{1}{|@{\hspace{.6ex}}c@{\hspace{.6ex}}|}{\raisebox{-.3ex}{$#1$}}}\raisebox{-.6ex}{$\begin{array}[b]{*{1}c}\cline{1-1}\lr{4}\\\cline{1-1}\lr{5}\\\cline{1-1}\end{array}$}}$};
\node (node_9) at (38.485bp,99.827bp) [draw,draw=none] {${\def\lr#1{\multicolumn{1}{|@{\hspace{.6ex}}c@{\hspace{.6ex}}|}{\raisebox{-.3ex}{$#1$}}}\raisebox{-.6ex}{$\begin{array}[b]{*{1}c}\cline{1-1}\lr{4}\\\cline{1-1}\lr{6}\\\cline{1-1}\end{array}$}}$};
\node (node_7) at (98.485bp,438.03bp) [draw,draw=none] {${\def\lr#1{\multicolumn{1}{|@{\hspace{.6ex}}c@{\hspace{.6ex}}|}{\raisebox{-.3ex}{$#1$}}}\raisebox{-.6ex}{$\begin{array}[b]{*{1}c}\cline{1-1}\lr{1}\\\cline{1-1}\lr{5}\\\cline{1-1}\end{array}$}}$};
\node (node_14) at (68.485bp,353.48bp) [draw,draw=none] {${\def\lr#1{\multicolumn{1}{|@{\hspace{.6ex}}c@{\hspace{.6ex}}|}{\raisebox{-.3ex}{$#1$}}}\raisebox{-.6ex}{$\begin{array}[b]{*{1}c}\cline{1-1}\lr{2}\\\cline{1-1}\lr{5}\\\cline{1-1}\end{array}$}}$};
\node (node_8) at (28.485bp,438.03bp) [draw,draw=none] {${\def\lr#1{\multicolumn{1}{|@{\hspace{.6ex}}c@{\hspace{.6ex}}|}{\raisebox{-.3ex}{$#1$}}}\raisebox{-.6ex}{$\begin{array}[b]{*{1}c}\cline{1-1}\lr{2}\\\cline{1-1}\lr{4}\\\cline{1-1}\end{array}$}}$};
\draw [red,->] (node_0) ..controls (85.78bp,241.71bp) and (77.066bp,223.73bp)  .. (node_12);
\draw (89.985bp,226.65bp) node {$2$};
\draw [red,->] (node_1) ..controls (38.485bp,664.7bp) and (38.485bp,647.18bp)  .. (node_4);
\draw (46.985bp,649.41bp) node {$2$};
\draw [black,->] (node_3) ..controls (5.1232bp,327.94bp) and (4.9331bp,313.74bp)  .. (9.4845bp,302.21bp) .. controls (12.243bp,295.22bp) and (17.255bp,288.77bp)  .. (node_10);
\draw (17.985bp,311.21bp) node {$4$};
\draw [blue,->] (node_4) ..controls (26.062bp,589.18bp) and (21.646bp,581.49bp)  .. (19.485bp,573.86bp) .. controls (17.189bp,565.75bp) and (16.482bp,556.61bp)  .. (node_11);
\draw (27.985bp,564.86bp) node {$1$};
\draw [spinach,->] (node_4) ..controls (44.779bp,580.15bp) and (49.025bp,562.63bp)  .. (node_13);
\draw (58.985bp,564.86bp) node {$3$};
\draw [blue,->] (node_5) ..controls (119.0bp,326.38bp) and (112.55bp,308.63bp)  .. (node_0);
\draw (123.98bp,311.21bp) node {$1$};
\draw [brown,->] (node_6) ..controls (16.168bp,159.19bp) and (16.107bp,145.03bp)  .. (19.485bp,133.1bp) .. controls (20.633bp,129.05bp) and (22.418bp,124.98bp)  .. (node_9);
\draw (27.985bp,142.1bp) node {$5$};
\draw [brown,->] (node_7) ..controls (107.97bp,410.93bp) and (114.42bp,393.18bp)  .. (node_5);
\draw (123.98bp,395.76bp) node {$5$};
\draw [blue,->] (node_7) ..controls (88.999bp,410.93bp) and (82.548bp,393.18bp)  .. (node_14);
\draw (93.985bp,395.76bp) node {$1$};
\draw [red,->] (node_8) ..controls (16.062bp,420.08bp) and (11.646bp,412.39bp)  .. (9.4845bp,404.76bp) .. controls (7.1886bp,396.65bp) and (6.4817bp,387.51bp)  .. (node_3);
\draw (17.985bp,395.76bp) node {$2$};
\draw [black,->] (node_8) ..controls (41.189bp,410.81bp) and (49.903bp,392.83bp)  .. (node_14);
\draw (60.985bp,395.76bp) node {$4$};
\draw [black,->] (node_9) ..controls (38.485bp,72.845bp) and (38.485bp,55.321bp)  .. (node_2);
\draw (46.985bp,57.551bp) node {$4$};
\draw [spinach,->] (node_10) ..controls (26.062bp,250.97bp) and (21.646bp,243.29bp)  .. (19.485bp,235.65bp) .. controls (17.189bp,227.55bp) and (16.482bp,218.41bp)  .. (node_6);
\draw (27.985bp,226.65bp) node {$3$};
\draw [brown,->] (node_10) ..controls (44.779bp,241.95bp) and (49.025bp,224.42bp)  .. (node_12);
\draw (58.985bp,226.65bp) node {$5$};
\draw [spinach,->] (node_11) ..controls (17.575bp,497.3bp) and (17.707bp,483.4bp)  .. (19.485bp,471.31bp) .. controls (19.863bp,468.73bp) and (20.383bp,466.08bp)  .. (node_8);
\draw (27.985bp,480.31bp) node {$3$};
\draw [spinach,->] (node_12) ..controls (52.19bp,157.4bp) and (47.944bp,139.87bp)  .. (node_9);
\draw (58.985bp,142.1bp) node {$3$};
\draw [black,->] (node_13) ..controls (71.189bp,495.37bp) and (79.903bp,477.38bp)  .. (node_7);
\draw (89.985bp,480.31bp) node {$4$};
\draw [blue,->] (node_13) ..controls (48.999bp,495.48bp) and (42.548bp,477.73bp)  .. (node_8);
\draw (54.985bp,480.31bp) node {$1$};
\draw [brown,->] (node_14) ..controls (77.97bp,326.38bp) and (84.421bp,308.63bp)  .. (node_0);
\draw (93.985bp,311.21bp) node {$5$};
\draw [red,->] (node_14) ..controls (58.999bp,326.38bp) and (52.548bp,308.63bp)  .. (node_10);
\draw (64.985bp,311.21bp) node {$2$};
\end{tikzpicture}
,\quad
\typeb[2]\colon
\begin{tikzpicture}[anchorbase,>=latex,line join=bevel,scale=0.7,every path/.style={very thick}]
\node (node_0) at (28.485bp,438.03bp) [draw,draw=none] {${\def\lr#1{\multicolumn{1}{|@{\hspace{.6ex}}c@{\hspace{.6ex}}|}{\raisebox{-.3ex}{$#1$}}}\raisebox{-.6ex}{$\begin{array}[b]{*{1}c}\cline{1-1}\lr{1}\\\cline{1-1}\lr{0}\\\cline{1-1}\end{array}$}}$};
\node (node_4) at (8.4845bp,353.48bp) [draw,draw=none] {${\def\lr#1{\multicolumn{1}{|@{\hspace{.6ex}}c@{\hspace{.6ex}}|}{\raisebox{-.3ex}{$#1$}}}\raisebox{-.6ex}{$\begin{array}[b]{*{1}c}\cline{1-1}\lr{1}\\\cline{1-1}\lr{\overline{2}}\\\cline{1-1}\end{array}$}}$};
\node (node_7) at (48.485bp,353.48bp) [draw,draw=none] {${\def\lr#1{\multicolumn{1}{|@{\hspace{.6ex}}c@{\hspace{.6ex}}|}{\raisebox{-.3ex}{$#1$}}}\raisebox{-.6ex}{$\begin{array}[b]{*{1}c}\cline{1-1}\lr{2}\\\cline{1-1}\lr{0}\\\cline{1-1}\end{array}$}}$};
\node (node_1) at (28.485bp,522.58bp) [draw,draw=none] {${\def\lr#1{\multicolumn{1}{|@{\hspace{.6ex}}c@{\hspace{.6ex}}|}{\raisebox{-.3ex}{$#1$}}}\raisebox{-.6ex}{$\begin{array}[b]{*{1}c}\cline{1-1}\lr{1}\\\cline{1-1}\lr{2}\\\cline{1-1}\end{array}$}}$};
\node (node_2) at (28.485bp,15.276bp) [draw,draw=none] {${\def\lr#1{\multicolumn{1}{|@{\hspace{.6ex}}c@{\hspace{.6ex}}|}{\raisebox{-.3ex}{$#1$}}}\raisebox{-.6ex}{$\begin{array}[b]{*{1}c}\cline{1-1}\lr{\overline{2}}\\\cline{1-1}\lr{\overline{1}}\\\cline{1-1}\end{array}$}}$};
\node (node_3) at (48.485bp,268.93bp) [draw,draw=none] {${\def\lr#1{\multicolumn{1}{|@{\hspace{.6ex}}c@{\hspace{.6ex}}|}{\raisebox{-.3ex}{$#1$}}}\raisebox{-.6ex}{$\begin{array}[b]{*{1}c}\cline{1-1}\lr{0}\\\cline{1-1}\lr{0}\\\cline{1-1}\end{array}$}}$};
\node (node_5) at (48.485bp,184.38bp) [draw,draw=none] {${\def\lr#1{\multicolumn{1}{|@{\hspace{.6ex}}c@{\hspace{.6ex}}|}{\raisebox{-.3ex}{$#1$}}}\raisebox{-.6ex}{$\begin{array}[b]{*{1}c}\cline{1-1}\lr{0}\\\cline{1-1}\lr{\overline{2}}\\\cline{1-1}\end{array}$}}$};
\node (node_9) at (8.4845bp,268.93bp) [draw,draw=none] {${\def\lr#1{\multicolumn{1}{|@{\hspace{.6ex}}c@{\hspace{.6ex}}|}{\raisebox{-.3ex}{$#1$}}}\raisebox{-.6ex}{$\begin{array}[b]{*{1}c}\cline{1-1}\lr{2}\\\cline{1-1}\lr{\overline{2}}\\\cline{1-1}\end{array}$}}$};
\node (node_6) at (28.485bp,99.827bp) [draw,draw=none] {${\def\lr#1{\multicolumn{1}{|@{\hspace{.6ex}}c@{\hspace{.6ex}}|}{\raisebox{-.3ex}{$#1$}}}\raisebox{-.6ex}{$\begin{array}[b]{*{1}c}\cline{1-1}\lr{0}\\\cline{1-1}\lr{\overline{1}}\\\cline{1-1}\end{array}$}}$};
\node (node_8) at (8.4845bp,184.38bp) [draw,draw=none] {${\def\lr#1{\multicolumn{1}{|@{\hspace{.6ex}}c@{\hspace{.6ex}}|}{\raisebox{-.3ex}{$#1$}}}\raisebox{-.6ex}{$\begin{array}[b]{*{1}c}\cline{1-1}\lr{2}\\\cline{1-1}\lr{\overline{1}}\\\cline{1-1}\end{array}$}}$};
\draw [red,->] (node_0) ..controls (16.062bp,420.08bp) and (11.646bp,412.39bp)  .. (9.4845bp,404.76bp) .. controls (7.1886bp,396.65bp) and (6.4817bp,387.51bp)  .. (node_4);
\draw (17.985bp,395.76bp) node {$2$};
\draw [blue,->] (node_0) ..controls (34.779bp,411.05bp) and (39.025bp,393.53bp)  .. (node_7);
\draw (48.985bp,395.76bp) node {$1$};
\draw [red,->] (node_1) ..controls (28.485bp,495.6bp) and (28.485bp,478.08bp)  .. (node_0);
\draw (36.985bp,480.31bp) node {$2$};
\draw [red,->] (node_3) ..controls (48.485bp,241.95bp) and (48.485bp,224.42bp)  .. (node_5);
\draw (56.985bp,226.65bp) node {$2$};
\draw [blue,->] (node_4) ..controls (8.4845bp,326.5bp) and (8.4845bp,308.98bp)  .. (node_9);
\draw (16.985bp,311.21bp) node {$1$};
\draw [blue,->] (node_5) ..controls (42.19bp,157.4bp) and (37.944bp,139.87bp)  .. (node_6);
\draw (48.985bp,142.1bp) node {$1$};
\draw [red,->] (node_6) ..controls (28.485bp,72.845bp) and (28.485bp,55.321bp)  .. (node_2);
\draw (36.985bp,57.551bp) node {$2$};
\draw [red,->] (node_7) ..controls (48.485bp,326.5bp) and (48.485bp,308.98bp)  .. (node_3);
\draw (56.985bp,311.21bp) node {$2$};
\draw [red,->] (node_8) ..controls (11.243bp,158.98bp) and (13.341bp,145.07bp)  .. (16.485bp,133.1bp) .. controls (17.151bp,130.57bp) and (17.941bp,127.95bp)  .. (node_6);
\draw (24.985bp,142.1bp) node {$2$};
\draw [blue,->] (node_9) ..controls (8.4845bp,241.95bp) and (8.4845bp,224.42bp)  .. (node_8);
\draw (16.985bp,226.65bp) node {$1$};
\end{tikzpicture}
.
$}
\end{gather}
These also exemplify the tableaux realizations for $k=1$,
which should be compared to \autoref{E:FiniteTypesCrystalList}, and also gives
the cases $\typea[5]$ and $\typeb[2]$ with $\fweight=2$. In type $\typeb$ we write $0=e+1$ for typographical reasons.

We only discuss type $\typeb$, the other cases are similar
and are omitted.
Note that it suffices to concentrate on the column words
obtained by reading a tableau from top to bottom. If $w=w_{1}\dots w_{k}$ is a column word, then the edge $w\xrightarrow{i}w^{\prime}$ is given by replacing $i$ with $i+1$, or $\overline{i+1}$ with $\overline{i}$.

Assume first that we have a square $ij=ji$, with $i<j\neq e$. Locally, the column words vertices in the square must be of the form:
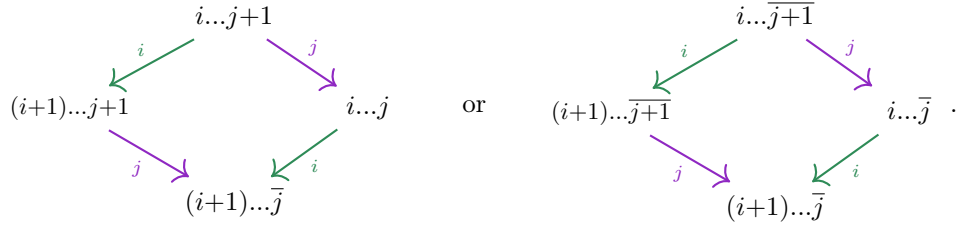
\begin{gather*}
\begin{tikzcd}[column sep=small,ampersand replacement=\&]
\& i\dots j{+}1 \arrow[dl,spinach,"i"',thick]\arrow[dr,orchid,"j",thick] \&
\\
\scalebox{0.9}{$(i{+}1)\dots j{+}1$}\arrow[dr,orchid,"j"',thick] \& \& \phantom{.}i\dots j\phantom{.}\arrow[dl,spinach,"i",thick]
\\
\& (i{+}1)\dots\overline{j} \&
\end{tikzcd}
\qquad\text{or}\qquad
\begin{tikzcd}[column sep=small,ampersand replacement=\&]
\& i\dots\overline{j{+}1} \arrow[dl,spinach,"i"',thick]\arrow[dr,orchid,"j",thick] \&
\\
\scalebox{0.9}{$(i{+}1)\dots\overline{j{+}1}$}\arrow[dr,orchid,"j"',thick] \& \& \phantom{.}i\dots\overline{j}\phantom{.}\arrow[dl,spinach,"i",thick]
\\
\& (i{+}1)\dots\overline{j} \&
\end{tikzcd}
.
\end{gather*}
or an overlined version of the first diagram. In the first case $i+1$ and $j+1$ do not appear in the initial word and in the second case $i+1$ and $\overline{j+1}$ do not appear. If $i$ and $j$ are adjacent, then the first case does not arise and
in order to get from the highest weight vertex $12\dots k$ to the vertex $i\dots\overline{j+1}$ there must be previous edge labeled $j$ since only this edge could have erased
to entry $j$. The edge labeled $j$ must come right before the top vertex $i\dots\overline{j+1}$ in the diagram above since otherwise the order condition would not be satisfied. Thus, we get the desired conditions when $j$ precedes $i$.

The other cases with $j<i\neq e$ are similar.
To conclude the square case, for $i=e-1$ and $j=e$ there is no possible square since $e$ needs to appear in the column word directly before of after $e-1$,
and similarly for $i=e$ and $j=e-1$.

A minimal octagon can
only appear for adjacent $i$ and $j$ since $i$ and $j$ must be applied twice because of the column word combinatorics. Hence,
by \cite[Theorem 1]{St-local-crystal-doubly-laced}, it remains to look at
decagons and tetradecagons as in {\loccit} However, these cannot appear
because they must involve three consecutive applications of, say, $i$ which is
impossible by the column word combinatorics.

\textit{Exceptional types.} As there are only finitely many fundamental crystal graphs of exceptional type, this can be verified directly by inspection.
For type $\typeg[2]$, these crystal graphs are:
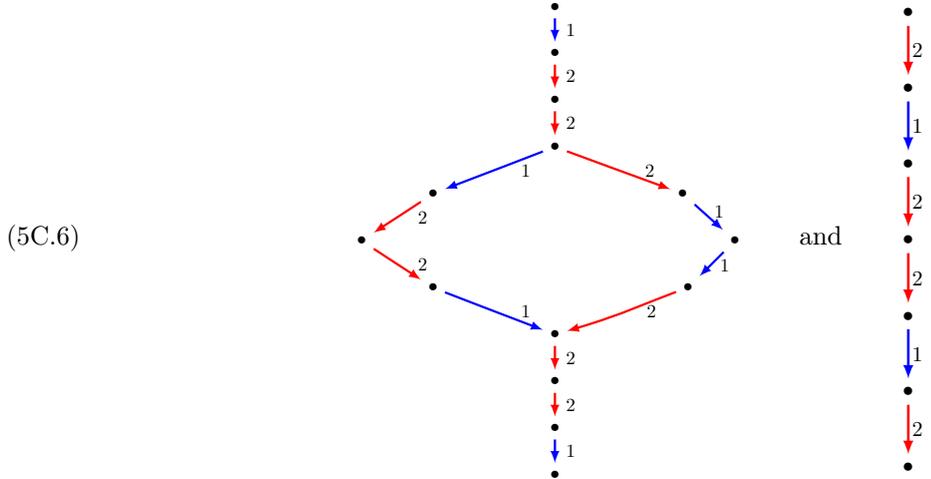
\begin{gather}\label{Eq:FiniteTypesG2}
\scalebox{0.7}{$\begin{tikzpicture}[anchorbase,>=latex,line join=bevel,xscale=1,yscale=0.35,
every path/.style={very thick}]
\node (node_0) at (212.5bp,583.0bp) [draw,draw=none] {$\bullet$};
\node (node_3) at (212.5bp,511.0bp) [draw,draw=none] {$\bullet$};
\node (node_1) at (147.5bp,439.0bp) [draw,draw=none] {$\bullet$};
\node (node_8) at (109.5bp,367.0bp) [draw,draw=none] {$\bullet$};
\node (node_2) at (212.5bp,725.5bp) [draw,draw=none] {$\bullet$};
\node (node_6) at (212.5bp,654.5bp) [draw,draw=none] {$\bullet$};
\node (node_9) at (280.5bp,439.0bp) [draw,draw=none] {$\bullet$};
\node (node_4) at (147.5bp,295.0bp) [draw,draw=none] {$\bullet$};
\node (node_12) at (212.5bp,223.0bp) [draw,draw=none] {$\bullet$};
\node (node_5) at (308.5bp,367.0bp) [draw,draw=none] {$\bullet$};
\node (node_7) at (283.5bp,295.0bp) [draw,draw=none] {$\bullet$};
\node (node_10) at (212.5bp,79.5bp) [draw,draw=none] {$\bullet$};
\node (node_11) at (212.5bp,8.5bp) [draw,draw=none] {$\bullet$};
\node (node_13) at (212.5bp,151.0bp) [draw,draw=none] {$\bullet$};
\draw [red,->] (node_0) ..controls (212.5bp,563.43bp) and (212.5bp,544.61bp)  .. (node_3);
\draw (221.0bp,547.0bp) node {$2$};
\draw [red,->] (node_1) ..controls (137.25bp,419.12bp) and (126.59bp,399.48bp)  .. (node_8);
\draw (142.0bp,403.0bp) node {$2$};
\draw [blue,->] (node_2) ..controls (212.5bp,706.44bp) and (212.5bp,687.5bp)  .. (node_6);
\draw (221.0bp,690.0bp) node {$1$};
\draw [blue,->] (node_3) ..controls (194.67bp,490.8bp) and (175.67bp,470.34bp)  .. (node_1);
\draw (197.0bp,475.0bp) node {$1$};
\draw [red,->] (node_3) ..controls (231.38bp,490.56bp) and (251.89bp,469.46bp)  .. (node_9);
\draw (263.0bp,475.0bp) node {$2$};
\draw [blue,->] (node_4) ..controls (165.33bp,274.8bp) and (184.33bp,254.34bp)  .. (node_12);
\draw (197.0bp,259.0bp) node {$1$};
\draw [blue,->] (node_5) to (node_7);
\draw (303.0bp,329.0bp) node {$1$};
\draw [red,->] (node_6) ..controls (212.5bp,635.42bp) and (212.5bp,616.63bp)  .. (node_0);
\draw (221.0bp,619.0bp) node {$2$};
\draw [red,->] (node_7) ..controls (264.18bp,274.95bp) and (248.68bp,253.75bp)  .. (node_12);
\draw (264.0bp,259.0bp) node {$2$};
\draw [red,->] (node_8) ..controls (119.75bp,347.12bp) and (130.41bp,327.48bp)  .. (node_4);
\draw (142.0bp,331.0bp) node {$2$};
\draw [blue,->] (node_9) to (node_5);
\draw (300.0bp,411.0bp) node {$1$};
\draw [blue,->] (node_10) ..controls (212.5bp,60.442bp) and (212.5bp,41.496bp)  .. (node_11);
\draw (221.0bp,44.0bp) node {$1$};
\draw [red,->] (node_12) ..controls (212.5bp,203.43bp) and (212.5bp,184.61bp)  .. (node_13);
\draw (221.0bp,187.0bp) node {$2$};
\draw [red,->] (node_13) ..controls (212.5bp,131.19bp) and (212.5bp,112.59bp)  .. (node_10);
\draw (221.0bp,115.0bp) node {$2$};
\end{tikzpicture}$}
\qquad\text{and}\qquad
\scalebox{0.8}{$\begin{tikzpicture}[anchorbase,>=latex,line join=bevel,scale=0.5,every path/.style={very thick}]
\node (node_0) at (53.5bp,79.5bp) [draw,draw=none] {$\bullet$};
\node (node_6) at (53.5bp,8.5bp) [draw,draw=none] {$\bullet$};
\node (node_1) at (53.5bp,222.0bp) [draw,draw=none] {$\bullet$};
\node (node_3) at (53.5bp,150.5bp) [draw,draw=none] {$\bullet$};
\node (node_2) at (53.5bp,293.5bp) [draw,draw=none] {$\bullet$};
\node (node_4) at (53.5bp,364.5bp) [draw,draw=none] {$\bullet$};
\node (node_5) at (53.5bp,435.5bp) [draw,draw=none] {$\bullet$};
\draw [red,->] (node_0) ..controls (53.5bp,60.442bp) and (53.5bp,41.496bp)  .. (node_6);
\draw (62.0bp,44.0bp) node {$2$};
\draw [red,->] (node_1) ..controls (53.5bp,202.19bp) and (53.5bp,183.59bp)  .. (node_3);
\draw (62.0bp,186.0bp) node {$2$};
\draw [red,->] (node_2) ..controls (53.5bp,274.42bp) and (53.5bp,255.63bp)  .. (node_1);
\draw (62.0bp,258.0bp) node {$2$};
\draw [blue,->] (node_3) ..controls (53.5bp,131.44bp) and (53.5bp,112.5bp)  .. (node_0);
\draw (62.0bp,115.0bp) node {$1$};
\draw [blue,->] (node_4) ..controls (53.5bp,345.44bp) and (53.5bp,326.5bp)  .. (node_2);
\draw (62.0bp,329.0bp) node {$1$};
\draw [red,->] (node_5) ..controls (53.5bp,416.44bp) and (53.5bp,397.5bp)  .. (node_4);
\draw (62.0bp,400.0bp) node {$2$};
\end{tikzpicture}$}
.
\end{gather}
These are the two relevant fundamental crystals in type $\typeg[2]$.
In type $\typee[i]$ the statements follow from \cite{St-local-crystal-simply-laced} and \cite{St-local-crystal-doubly-laced}.
All other types were checked using the code available at \cite{MaTu-sagemath-finite-type-klrw}.
\end{proof}

As we will see, all cases in \autoref{P:FiniteTypesSSR} correspond
to relations in the wKLRW algebra.

Recall from \autoref{SS:RecollectionKacMoody} that $P^{+}$ is the set of dominant weights for $\quiver$. For the rest of this section, fix a dominant weight $\Lambda\in P^{+}$ and let $\crystalgraph[\Lambda]$ be the crystal graph of highest weight $\Lambda$.  A \emph{rooted path} in $\crystalgraph[\Lambda]$ is a path with source the highest weight vector in $\crystalgraph[\Lambda]$. Let $\Parts{n}{\Lambda}$ be the set of rooted paths in $\crystalgraph[\Lambda]$. If $\bpath\in\Parts{n}{\Lambda}$. The \emph{residue sequence} of $\bpath$ is the sequence $\Res[\bpath]\in\vertices^{n}$ of edges labels for $\bpath$. As before, the path $\bpath$ is uniquely determined by its residue sequence. The \emph{sink} $\omega(\bpath)$ of $\bpath$ is the final vertex in the path $\bpath$.

We want to associate an idempotent diagram to each path $\bpath$ in $\crystalgraph[\Lambda]$. Before we can do this we need to introduce red strings and affine red strings, which we postpone to
\autoref{S:Construction}.

%%%%%%%%%%%%%%%%%%%%%%%%%%%%%%%%%%%%%%%%%

\subsection{Face permutations}\label{SS:FacePermutations}

%%%%%%%%%%%%%%%%%%%%%%%%%%%%%%%%%%%%%%%%%

One of the basic building blocks for our cellular basis are the permutations of the residue sequences of paths in the crystal, which we now introduce.

\begin{Definition}\label{D:FiniteTypesSSPerm}
Let $\lambda$ and $\mu$ be vertices of $\crystalgraph[\Lambda]$. A \emph{basic face permutation} between paths $\bpath,\cpath\in\Parts{n}{\Lambda}$ is the equivalence classes of permutations that act on the residue sequences around a two color face as follows:
\begin{gather*}
\text{square}\colon ij\to ji
,\quad
\text{octagon}\colon ijji\to jiij\text{ or }ij\to ji
,
\\
\text{decagon}\colon
\left\{
\begin{gathered}
23332\to 32233
,
23332\to 32323
,
32233\to 32323
,
\\
23\to 32
,
32\to 23
,
233\to 323
,
323\to 233
,
3223\to 3232
,
3232\to 3223
,
\end{gathered}
\right.
\\
\text{tetradecagon}\colon
\text{all permutations of }
2332323,3223332,3232332,2333223.
\end{gather*}
A \emph{partial face permutation} is a (residue preserving) basic face permutation that is applied only partially around a face starting from the top of the face. A \emph{face permutation} of two paths $\bpath$ and $\cpath$ is a composition of basic and partial face permutations. Let $\Face[\Lambda](\bpath,\cpath)$ be the set of face permutations from $\bpath$ to $\cpath$.
\end{Definition}

For fundamental weight, decagons and tetradecagons appear only in type $\typef[4]$ precisely for the faces listed in \autoref{P:FiniteTypesSSR}, for the residues $2$ and $3$.

\begin{Remark}\label{R:FiniteTypesSSPerm}
A face permutation from $\bpath$ to $\cpath$ is determined by its action on the residue sequences. For squares, there is a unique permutation that sends $ij$ to $ji$, but for each of the other basic face permutations there is more than one permutation in the equivalence class. For example, $(1,3)(2,4)$ and $(1,2)(3,4)$ both send $ijji$ to $jiij$. In \autoref{D:ConstructionPermBasis} below we fix a representative for each equivalence class of face permutations. In what follows we abuse terminology and call a face permutation any permutation in the corresponding equivalence class.
\end{Remark}

\begin{Example}\label{E:FiniteTypesPermutation}
In \autoref{Eq:MainExample}, the path $\qpath[\lambda]$ is related to $\qpath[\mu]$
by using a partial face permutation around the $1221=2112$ octagon. This can pictured as:
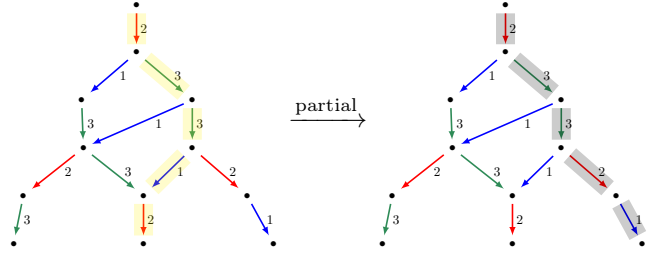
\begin{gather*}
\scalebox{0.5}{$\begin{tikzpicture}[anchorbase,>=latex,line join=bevel,xscale=0.65,yscale=0.5,every path/.style={very thick}]
\node (node_0) at (65.5bp,81.0bp) [draw,draw=none] {$\bullet$};
\node (node_1) at (54.5bp,9.0bp) [draw,draw=none] {$\bullet$};
\node (node_2) at (132.5bp,225.0bp) [draw,draw=none] {$\bullet$};
\node (node_6) at (134.5bp,153.0bp) [draw,draw=none] {$\bullet$};
\node (node_3) at (203.5bp,9.0bp) [draw,draw=none] {$\bullet$};
\node (node_4) at (352.5bp,9.0bp) [draw,draw=none] {$\bullet$};
\node (node_5) at (322.5bp,81.0bp) [draw,draw=none] {$\bullet$};
\node (node_7) at (203.5bp,81.0bp) [draw,draw=none] {$\bullet$};
\node (node_8) at (259.5bp,153.0bp) [draw,draw=none] {$\bullet$};
\node (node_9) at (259.5bp,225.0bp) [draw,draw=none] {$\bullet$};
\node (node_10) at (195.5bp,367.5bp) [draw,draw=none] {$\bullet$};
\node (node_11) at (195.5bp,296.5bp) [draw,draw=none] {$\bullet$};
\draw [spinach,->] (node_0) ..controls (62.582bp,61.429bp) and (59.624bp,42.606bp)  .. (node_1);
\draw (70.0bp,45.0bp) node {$3$};
\draw [spinach,->] (node_2) ..controls (133.02bp,205.68bp) and (133.57bp,186.46bp)  .. (node_6);
\draw (143.0bp,189.0bp) node {$3$};
\draw [blue,->] (node_5) ..controls (330.44bp,61.471bp) and (338.92bp,41.695bp)  .. (node_4);
\draw (349.0bp,45.0bp) node {$1$};
\draw [red,->] (node_6) ..controls (115.47bp,132.7bp) and (95.029bp,111.96bp)  .. (node_0);
\draw (117.0bp,117.0bp) node {$2$};
\draw [spinach,->] (node_6) ..controls (153.76bp,132.46bp) and (174.85bp,111.06bp)  .. (node_7);
\draw (186.0bp,117.0bp) node {$3$};
\draw [red,->] (node_7) ..controls (203.5bp,61.68bp) and (203.5bp,42.459bp)  .. (node_3);
\draw (212.0bp,45.0bp) node {$2$};
\draw [red,->] (node_8) ..controls (276.66bp,132.93bp) and (295.77bp,111.7bp)  .. (node_5);
\draw (306.0bp,117.0bp) node {$2$};
\draw [blue,->] (node_8) ..controls (244.33bp,133.04bp) and (227.58bp,112.1bp)  .. (node_7);
\draw (246.0bp,117.0bp) node {$1$};
\draw [blue,->] (node_9) ..controls (223.89bp,204.06bp) and (183.75bp,181.58bp)  .. (node_6);
\draw (221.0bp,189.0bp) node {$1$};
\draw [spinach,->] (node_9) ..controls (259.5bp,205.2bp) and (259.5bp,185.79bp)  .. (node_8);
\draw (268.0bp,189.0bp) node {$3$};
\draw [red,->] (node_10) ..controls (195.5bp,348.44bp) and (195.5bp,329.5bp)  .. (node_11);
\draw (204.0bp,332.0bp) node {$2$};
\draw [blue,->] (node_11) ..controls (178.34bp,276.57bp) and (159.23bp,255.49bp)  .. (node_2);
\draw (180.0bp,261.0bp) node {$1$};
\draw [spinach,->] (node_11) ..controls (212.71bp,276.81bp) and (231.53bp,256.37bp)  .. (node_9);
\draw (243.0bp,261.0bp) node {$3$};
\draw[yellow,line width=0.5cm,opacity=0.2,rounded corners] (node_10) to (node_11) to (node_9) to (node_8) to (node_7) to (node_3);
\end{tikzpicture}$}
\xrightarrow{\text{partial}}
\scalebox{0.5}{$\begin{tikzpicture}[anchorbase,>=latex,line join=bevel,xscale=0.65,yscale=0.5,every path/.style={very thick}]
\node (node_0) at (65.5bp,81.0bp) [draw,draw=none] {$\bullet$};
\node (node_1) at (54.5bp,9.0bp) [draw,draw=none] {$\bullet$};
\node (node_2) at (132.5bp,225.0bp) [draw,draw=none] {$\bullet$};
\node (node_6) at (134.5bp,153.0bp) [draw,draw=none] {$\bullet$};
\node (node_3) at (203.5bp,9.0bp) [draw,draw=none] {$\bullet$};
\node (node_4) at (352.5bp,9.0bp) [draw,draw=none] {$\bullet$};
\node (node_5) at (322.5bp,81.0bp) [draw,draw=none] {$\bullet$};
\node (node_7) at (203.5bp,81.0bp) [draw,draw=none] {$\bullet$};
\node (node_8) at (259.5bp,153.0bp) [draw,draw=none] {$\bullet$};
\node (node_9) at (259.5bp,225.0bp) [draw,draw=none] {$\bullet$};
\node (node_10) at (195.5bp,367.5bp) [draw,draw=none] {$\bullet$};
\node (node_11) at (195.5bp,296.5bp) [draw,draw=none] {$\bullet$};
\draw [spinach,->] (node_0) ..controls (62.582bp,61.429bp) and (59.624bp,42.606bp)  .. (node_1);
\draw (70.0bp,45.0bp) node {$3$};
\draw [spinach,->] (node_2) ..controls (133.02bp,205.68bp) and (133.57bp,186.46bp)  .. (node_6);
\draw (143.0bp,189.0bp) node {$3$};
\draw [blue,->] (node_5) ..controls (330.44bp,61.471bp) and (338.92bp,41.695bp)  .. (node_4);
\draw (349.0bp,45.0bp) node {$1$};
\draw [red,->] (node_6) ..controls (115.47bp,132.7bp) and (95.029bp,111.96bp)  .. (node_0);
\draw (117.0bp,117.0bp) node {$2$};
\draw [spinach,->] (node_6) ..controls (153.76bp,132.46bp) and (174.85bp,111.06bp)  .. (node_7);
\draw (186.0bp,117.0bp) node {$3$};
\draw [red,->] (node_7) ..controls (203.5bp,61.68bp) and (203.5bp,42.459bp)  .. (node_3);
\draw (212.0bp,45.0bp) node {$2$};
\draw [red,->] (node_8) ..controls (276.66bp,132.93bp) and (295.77bp,111.7bp)  .. (node_5);
\draw (306.0bp,117.0bp) node {$2$};
\draw [blue,->] (node_8) ..controls (244.33bp,133.04bp) and (227.58bp,112.1bp)  .. (node_7);
\draw (246.0bp,117.0bp) node {$1$};
\draw [blue,->] (node_9) ..controls (223.89bp,204.06bp) and (183.75bp,181.58bp)  .. (node_6);
\draw (221.0bp,189.0bp) node {$1$};
\draw [spinach,->] (node_9) ..controls (259.5bp,205.2bp) and (259.5bp,185.79bp)  .. (node_8);
\draw (268.0bp,189.0bp) node {$3$};
\draw [red,->] (node_10) ..controls (195.5bp,348.44bp) and (195.5bp,329.5bp)  .. (node_11);
\draw (204.0bp,332.0bp) node {$2$};
\draw [blue,->] (node_11) ..controls (178.34bp,276.57bp) and (159.23bp,255.49bp)  .. (node_2);
\draw (180.0bp,261.0bp) node {$1$};
\draw [spinach,->] (node_11) ..controls (212.71bp,276.81bp) and (231.53bp,256.37bp)  .. (node_9);
\draw (243.0bp,261.0bp) node {$3$};
\draw[black,line width=0.5cm,opacity=0.2,rounded corners] (node_10) to (node_11) to (node_9) to (node_8) to (node_5) to (node_4);
\end{tikzpicture}$}
\end{gather*}
Applying basic face permutations to go from $\lambda$ to $\nu$ needs three steps:
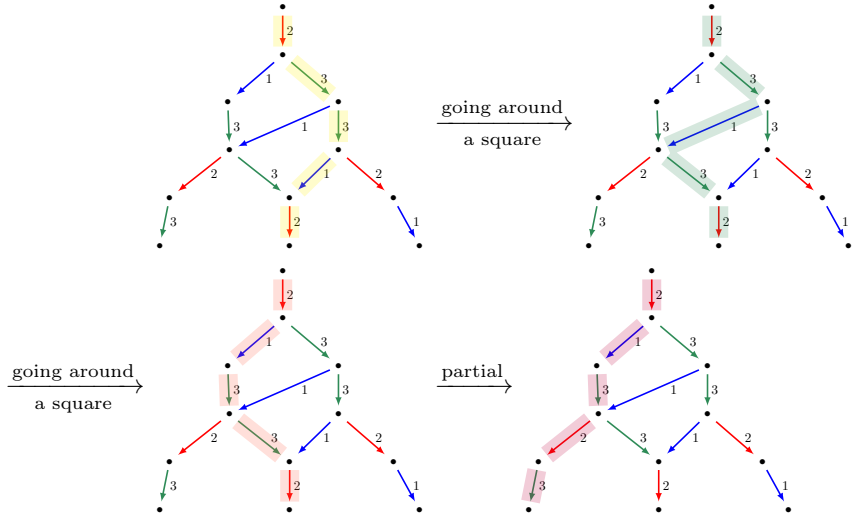
\begin{align*}
\scalebox{0.5}{$\begin{tikzpicture}[anchorbase,>=latex,line join=bevel,xscale=0.65,yscale=0.5,every path/.style={very thick}]
\node (node_0) at (65.5bp,81.0bp) [draw,draw=none] {$\bullet$};
\node (node_1) at (54.5bp,9.0bp) [draw,draw=none] {$\bullet$};
\node (node_2) at (132.5bp,225.0bp) [draw,draw=none] {$\bullet$};
\node (node_6) at (134.5bp,153.0bp) [draw,draw=none] {$\bullet$};
\node (node_3) at (203.5bp,9.0bp) [draw,draw=none] {$\bullet$};
\node (node_4) at (352.5bp,9.0bp) [draw,draw=none] {$\bullet$};
\node (node_5) at (322.5bp,81.0bp) [draw,draw=none] {$\bullet$};
\node (node_7) at (203.5bp,81.0bp) [draw,draw=none] {$\bullet$};
\node (node_8) at (259.5bp,153.0bp) [draw,draw=none] {$\bullet$};
\node (node_9) at (259.5bp,225.0bp) [draw,draw=none] {$\bullet$};
\node (node_10) at (195.5bp,367.5bp) [draw,draw=none] {$\bullet$};
\node (node_11) at (195.5bp,296.5bp) [draw,draw=none] {$\bullet$};
\draw [spinach,->] (node_0) ..controls (62.582bp,61.429bp) and (59.624bp,42.606bp)  .. (node_1);
\draw (70.0bp,45.0bp) node {$3$};
\draw [spinach,->] (node_2) ..controls (133.02bp,205.68bp) and (133.57bp,186.46bp)  .. (node_6);
\draw (143.0bp,189.0bp) node {$3$};
\draw [blue,->] (node_5) ..controls (330.44bp,61.471bp) and (338.92bp,41.695bp)  .. (node_4);
\draw (349.0bp,45.0bp) node {$1$};
\draw [red,->] (node_6) ..controls (115.47bp,132.7bp) and (95.029bp,111.96bp)  .. (node_0);
\draw (117.0bp,117.0bp) node {$2$};
\draw [spinach,->] (node_6) ..controls (153.76bp,132.46bp) and (174.85bp,111.06bp)  .. (node_7);
\draw (186.0bp,117.0bp) node {$3$};
\draw [red,->] (node_7) ..controls (203.5bp,61.68bp) and (203.5bp,42.459bp)  .. (node_3);
\draw (212.0bp,45.0bp) node {$2$};
\draw [red,->] (node_8) ..controls (276.66bp,132.93bp) and (295.77bp,111.7bp)  .. (node_5);
\draw (306.0bp,117.0bp) node {$2$};
\draw [blue,->] (node_8) ..controls (244.33bp,133.04bp) and (227.58bp,112.1bp)  .. (node_7);
\draw (246.0bp,117.0bp) node {$1$};
\draw [blue,->] (node_9) ..controls (223.89bp,204.06bp) and (183.75bp,181.58bp)  .. (node_6);
\draw (221.0bp,189.0bp) node {$1$};
\draw [spinach,->] (node_9) ..controls (259.5bp,205.2bp) and (259.5bp,185.79bp)  .. (node_8);
\draw (268.0bp,189.0bp) node {$3$};
\draw [red,->] (node_10) ..controls (195.5bp,348.44bp) and (195.5bp,329.5bp)  .. (node_11);
\draw (204.0bp,332.0bp) node {$2$};
\draw [blue,->] (node_11) ..controls (178.34bp,276.57bp) and (159.23bp,255.49bp)  .. (node_2);
\draw (180.0bp,261.0bp) node {$1$};
\draw [spinach,->] (node_11) ..controls (212.71bp,276.81bp) and (231.53bp,256.37bp)  .. (node_9);
\draw (243.0bp,261.0bp) node {$3$};
\draw[yellow,line width=0.5cm,opacity=0.2,rounded corners] (node_10) to (node_11) to (node_9) to (node_8) to (node_7) to (node_3);
\end{tikzpicture}$}
&\xrightarrow[\text{a square}]{\text{going around}}
\scalebox{0.5}{$\begin{tikzpicture}[anchorbase,>=latex,line join=bevel,xscale=0.65,yscale=0.5,every path/.style={very thick}]
\node (node_0) at (65.5bp,81.0bp) [draw,draw=none] {$\bullet$};
\node (node_1) at (54.5bp,9.0bp) [draw,draw=none] {$\bullet$};
\node (node_2) at (132.5bp,225.0bp) [draw,draw=none] {$\bullet$};
\node (node_6) at (134.5bp,153.0bp) [draw,draw=none] {$\bullet$};
\node (node_3) at (203.5bp,9.0bp) [draw,draw=none] {$\bullet$};
\node (node_4) at (352.5bp,9.0bp) [draw,draw=none] {$\bullet$};
\node (node_5) at (322.5bp,81.0bp) [draw,draw=none] {$\bullet$};
\node (node_7) at (203.5bp,81.0bp) [draw,draw=none] {$\bullet$};
\node (node_8) at (259.5bp,153.0bp) [draw,draw=none] {$\bullet$};
\node (node_9) at (259.5bp,225.0bp) [draw,draw=none] {$\bullet$};
\node (node_10) at (195.5bp,367.5bp) [draw,draw=none] {$\bullet$};
\node (node_11) at (195.5bp,296.5bp) [draw,draw=none] {$\bullet$};
\draw [spinach,->] (node_0) ..controls (62.582bp,61.429bp) and (59.624bp,42.606bp)  .. (node_1);
\draw (70.0bp,45.0bp) node {$3$};
\draw [spinach,->] (node_2) ..controls (133.02bp,205.68bp) and (133.57bp,186.46bp)  .. (node_6);
\draw (143.0bp,189.0bp) node {$3$};
\draw [blue,->] (node_5) ..controls (330.44bp,61.471bp) and (338.92bp,41.695bp)  .. (node_4);
\draw (349.0bp,45.0bp) node {$1$};
\draw [red,->] (node_6) ..controls (115.47bp,132.7bp) and (95.029bp,111.96bp)  .. (node_0);
\draw (117.0bp,117.0bp) node {$2$};
\draw [spinach,->] (node_6) ..controls (153.76bp,132.46bp) and (174.85bp,111.06bp)  .. (node_7);
\draw (186.0bp,117.0bp) node {$3$};
\draw [red,->] (node_7) ..controls (203.5bp,61.68bp) and (203.5bp,42.459bp)  .. (node_3);
\draw (212.0bp,45.0bp) node {$2$};
\draw [red,->] (node_8) ..controls (276.66bp,132.93bp) and (295.77bp,111.7bp)  .. (node_5);
\draw (306.0bp,117.0bp) node {$2$};
\draw [blue,->] (node_8) ..controls (244.33bp,133.04bp) and (227.58bp,112.1bp)  .. (node_7);
\draw (246.0bp,117.0bp) node {$1$};
\draw [blue,->] (node_9) ..controls (223.89bp,204.06bp) and (183.75bp,181.58bp)  .. (node_6);
\draw (221.0bp,189.0bp) node {$1$};
\draw [spinach,->] (node_9) ..controls (259.5bp,205.2bp) and (259.5bp,185.79bp)  .. (node_8);
\draw (268.0bp,189.0bp) node {$3$};
\draw [red,->] (node_10) ..controls (195.5bp,348.44bp) and (195.5bp,329.5bp)  .. (node_11);
\draw (204.0bp,332.0bp) node {$2$};
\draw [blue,->] (node_11) ..controls (178.34bp,276.57bp) and (159.23bp,255.49bp)  .. (node_2);
\draw (180.0bp,261.0bp) node {$1$};
\draw [spinach,->] (node_11) ..controls (212.71bp,276.81bp) and (231.53bp,256.37bp)  .. (node_9);
\draw (243.0bp,261.0bp) node {$3$};
\draw[spinach,line width=0.5cm,opacity=0.2,rounded corners] (node_10) to (node_11) to (node_9) to (node_6) to (node_7) to (node_3);
\end{tikzpicture}$}
\\
\xrightarrow[\text{a square}]{\text{going around}}
\scalebox{0.5}{$\begin{tikzpicture}[anchorbase,>=latex,line join=bevel,xscale=0.65,yscale=0.5,every path/.style={very thick}]
\node (node_0) at (65.5bp,81.0bp) [draw,draw=none] {$\bullet$};
\node (node_1) at (54.5bp,9.0bp) [draw,draw=none] {$\bullet$};
\node (node_2) at (132.5bp,225.0bp) [draw,draw=none] {$\bullet$};
\node (node_6) at (134.5bp,153.0bp) [draw,draw=none] {$\bullet$};
\node (node_3) at (203.5bp,9.0bp) [draw,draw=none] {$\bullet$};
\node (node_4) at (352.5bp,9.0bp) [draw,draw=none] {$\bullet$};
\node (node_5) at (322.5bp,81.0bp) [draw,draw=none] {$\bullet$};
\node (node_7) at (203.5bp,81.0bp) [draw,draw=none] {$\bullet$};
\node (node_8) at (259.5bp,153.0bp) [draw,draw=none] {$\bullet$};
\node (node_9) at (259.5bp,225.0bp) [draw,draw=none] {$\bullet$};
\node (node_10) at (195.5bp,367.5bp) [draw,draw=none] {$\bullet$};
\node (node_11) at (195.5bp,296.5bp) [draw,draw=none] {$\bullet$};
\draw [spinach,->] (node_0) ..controls (62.582bp,61.429bp) and (59.624bp,42.606bp)  .. (node_1);
\draw (70.0bp,45.0bp) node {$3$};
\draw [spinach,->] (node_2) ..controls (133.02bp,205.68bp) and (133.57bp,186.46bp)  .. (node_6);
\draw (143.0bp,189.0bp) node {$3$};
\draw [blue,->] (node_5) ..controls (330.44bp,61.471bp) and (338.92bp,41.695bp)  .. (node_4);
\draw (349.0bp,45.0bp) node {$1$};
\draw [red,->] (node_6) ..controls (115.47bp,132.7bp) and (95.029bp,111.96bp)  .. (node_0);
\draw (117.0bp,117.0bp) node {$2$};
\draw [spinach,->] (node_6) ..controls (153.76bp,132.46bp) and (174.85bp,111.06bp)  .. (node_7);
\draw (186.0bp,117.0bp) node {$3$};
\draw [red,->] (node_7) ..controls (203.5bp,61.68bp) and (203.5bp,42.459bp)  .. (node_3);
\draw (212.0bp,45.0bp) node {$2$};
\draw [red,->] (node_8) ..controls (276.66bp,132.93bp) and (295.77bp,111.7bp)  .. (node_5);
\draw (306.0bp,117.0bp) node {$2$};
\draw [blue,->] (node_8) ..controls (244.33bp,133.04bp) and (227.58bp,112.1bp)  .. (node_7);
\draw (246.0bp,117.0bp) node {$1$};
\draw [blue,->] (node_9) ..controls (223.89bp,204.06bp) and (183.75bp,181.58bp)  .. (node_6);
\draw (221.0bp,189.0bp) node {$1$};
\draw [spinach,->] (node_9) ..controls (259.5bp,205.2bp) and (259.5bp,185.79bp)  .. (node_8);
\draw (268.0bp,189.0bp) node {$3$};
\draw [red,->] (node_10) ..controls (195.5bp,348.44bp) and (195.5bp,329.5bp)  .. (node_11);
\draw (204.0bp,332.0bp) node {$2$};
\draw [blue,->] (node_11) ..controls (178.34bp,276.57bp) and (159.23bp,255.49bp)  .. (node_2);
\draw (180.0bp,261.0bp) node {$1$};
\draw [spinach,->] (node_11) ..controls (212.71bp,276.81bp) and (231.53bp,256.37bp)  .. (node_9);
\draw (243.0bp,261.0bp) node {$3$};
\draw[tomato,line width=0.5cm,opacity=0.2,rounded corners] (node_10) to (node_11) to (node_2) to (node_6) to (node_7) to (node_3);
\end{tikzpicture}$}
&\xrightarrow{\text{partial}}
\scalebox{0.5}{$\begin{tikzpicture}[anchorbase,>=latex,line join=bevel,xscale=0.65,yscale=0.5,every path/.style={very thick}]
\node (node_0) at (65.5bp,81.0bp) [draw,draw=none] {$\bullet$};
\node (node_1) at (54.5bp,9.0bp) [draw,draw=none] {$\bullet$};
\node (node_2) at (132.5bp,225.0bp) [draw,draw=none] {$\bullet$};
\node (node_6) at (134.5bp,153.0bp) [draw,draw=none] {$\bullet$};
\node (node_3) at (203.5bp,9.0bp) [draw,draw=none] {$\bullet$};
\node (node_4) at (352.5bp,9.0bp) [draw,draw=none] {$\bullet$};
\node (node_5) at (322.5bp,81.0bp) [draw,draw=none] {$\bullet$};
\node (node_7) at (203.5bp,81.0bp) [draw,draw=none] {$\bullet$};
\node (node_8) at (259.5bp,153.0bp) [draw,draw=none] {$\bullet$};
\node (node_9) at (259.5bp,225.0bp) [draw,draw=none] {$\bullet$};
\node (node_10) at (195.5bp,367.5bp) [draw,draw=none] {$\bullet$};
\node (node_11) at (195.5bp,296.5bp) [draw,draw=none] {$\bullet$};
\draw [spinach,->] (node_0) ..controls (62.582bp,61.429bp) and (59.624bp,42.606bp)  .. (node_1);
\draw (70.0bp,45.0bp) node {$3$};
\draw [spinach,->] (node_2) ..controls (133.02bp,205.68bp) and (133.57bp,186.46bp)  .. (node_6);
\draw (143.0bp,189.0bp) node {$3$};
\draw [blue,->] (node_5) ..controls (330.44bp,61.471bp) and (338.92bp,41.695bp)  .. (node_4);
\draw (349.0bp,45.0bp) node {$1$};
\draw [red,->] (node_6) ..controls (115.47bp,132.7bp) and (95.029bp,111.96bp)  .. (node_0);
\draw (117.0bp,117.0bp) node {$2$};
\draw [spinach,->] (node_6) ..controls (153.76bp,132.46bp) and (174.85bp,111.06bp)  .. (node_7);
\draw (186.0bp,117.0bp) node {$3$};
\draw [red,->] (node_7) ..controls (203.5bp,61.68bp) and (203.5bp,42.459bp)  .. (node_3);
\draw (212.0bp,45.0bp) node {$2$};
\draw [red,->] (node_8) ..controls (276.66bp,132.93bp) and (295.77bp,111.7bp)  .. (node_5);
\draw (306.0bp,117.0bp) node {$2$};
\draw [blue,->] (node_8) ..controls (244.33bp,133.04bp) and (227.58bp,112.1bp)  .. (node_7);
\draw (246.0bp,117.0bp) node {$1$};
\draw [blue,->] (node_9) ..controls (223.89bp,204.06bp) and (183.75bp,181.58bp)  .. (node_6);
\draw (221.0bp,189.0bp) node {$1$};
\draw [spinach,->] (node_9) ..controls (259.5bp,205.2bp) and (259.5bp,185.79bp)  .. (node_8);
\draw (268.0bp,189.0bp) node {$3$};
\draw [red,->] (node_10) ..controls (195.5bp,348.44bp) and (195.5bp,329.5bp)  .. (node_11);
\draw (204.0bp,332.0bp) node {$2$};
\draw [blue,->] (node_11) ..controls (178.34bp,276.57bp) and (159.23bp,255.49bp)  .. (node_2);
\draw (180.0bp,261.0bp) node {$1$};
\draw [spinach,->] (node_11) ..controls (212.71bp,276.81bp) and (231.53bp,256.37bp)  .. (node_9);
\draw (243.0bp,261.0bp) node {$3$};
\draw[purple,line width=0.5cm,opacity=0.2,rounded corners] (node_10) to (node_11) to (node_2) to (node_6) to (node_0) to (node_1);
\end{tikzpicture}$}
\end{align*}
So, the face permutation from $\qpath[\lambda]$ to $\qpath[\nu]$ is the composition of two basic face permutations and a partial face permutation.
\end{Example}

\begin{Lemma}\label{L:FiniteTypesPermutations}
Let $\bpath$ and $\cpath$ be two paths in $\crystalgraph[\Lambda]$ ending at the vertex $\omega(\bpath)=\omega(\cpath)$. Then $\Res[\bpath]=w\Res[\cpath]$, for some face permutation $w$.
\end{Lemma}

\begin{proof}
For simply laced quivers this is \cite[Theorem 4.2]{HeLe-crystal-path} and the remaining types follow by \cite[Remark 4.3(2)]{HeLe-crystal-path}.
Note that even though this remark only mentions doubly laced quivers it also applies to $\typeg$ because the conditions required by this remark always hold in finite type.
\end{proof}

Recall from \autoref{SS:RecollectionKacMoody} that $Q^{+}=\bigoplus_{i\in\vertices}\N\alpha_{i}$ is the positive root lattice of $\quiver$.

\begin{Definition}\label{D:FiniteTypesPermutations}
Let $\bpath$ be a rooted path in a crystal graph $\crystal$.
Let $\alpha_{\bpath}=\sum_{k=1}^{n}\alpha_{i_{k}}\in Q^{+}$, where $\Res[\bpath]=(i_{1},\dots,i_{n})$. If $\lambda$ is any vertex of $\crystal$ let $\alpha_{\lambda}=\alpha_{\bpath}$, where $\bpath$ is any path in $\crystal$ with sink $\omega(\bpath)=\lambda$.
Two paths $\bpath$ and $\cpath$ in $\crystal$ are in the same \emph{block} if $\alpha_{\bpath}=\alpha_{\cpath}$. Similarly, two vertices $\lambda$ and $\mu$ of $\crystal$ are in the same block if $\alpha_{\lambda}=\alpha_{\mu}$.
\end{Definition}

If $\lambda$ is a vertex of $\crystal$, then the root $\alpha_{\lambda}=\alpha_{\bpath}$ is independent of the choice of path $\bpath$ with $\omega(\bpath)=\lambda$ by \autoref{L:FiniteTypesPermutations}.

\begin{Proposition}\label{P:FiniteTypesPermutations}
Two paths $\bpath$ and $\cpath$ in the crystal $\crystal$
are in the same block if and only if $\Res[\bpath]=w\Res[\bpath]$, for some face permutation $w$.
\end{Proposition}

\begin{proof}
It is enough to consider the case when $\crystal=\crystalgraph[\Lambda]$, for some dominant weight $\Lambda\in P^{+}$.

\textit{$(\Leftarrow)$.} The condition $\Res[\bpath]=w\Res[\cpath]$ implies that $\alpha_{\bpath}=\alpha_{\cpath}$, so they are in the same block.

\textit{$(\Rightarrow)$.}
Since $\Gamma$ is of finite type, the crystal graph $\crystalgraph[\Lambda]$ has unique lowest weight vertex $\zeta$ of weight $-\fweight$. Let $\bpath_{\zeta}$ and $\cpath_{\zeta}$ be extensions of $\bpath$ and $\cpath$, respectively, to paths to $\zeta$. By \autoref{L:FiniteTypesPermutations}, there exists a face permutation $w_{\zeta}$ such that $\Res[\bpath_{\zeta}]=w_{\zeta}\Res[\cpath_{\zeta}]$. Since $w_\zeta$ is a face permutation, we can restrict $w_{\zeta}$ to $\cpath$ give the result.
\end{proof}

%%%%%%%%%%%%%%%%%%%%%%%%%%%%%%%%%%%%%%%%%

\subsection{Tensor product crystals}\label{SS:TensorProducts}

%%%%%%%%%%%%%%%%%%%%%%%%%%%%%%%%%%%%%%%%%

This section attaches crystal graphs to the tuples $\abfweight=(\Lambda_{\affine{\rho}_1},\dots,\Lambda_{\affine{\rho}_{\hell}})$ and $\bfweight=(\Lambda_{\rho_1},\dots,\Lambda_{\rho_\ell})$ from \autoref{D:hell}. By definition, $\affine{\rho}_i=\rho_i$, for $1\le i\le \ell$.

By \cite{Ka-crystal-bases}, if $\mathcal{G}$ and $\mathcal{H}$ are crystals then their tensor product $\mathcal{G}\otimes\mathcal{H}$ is also a crystal. More explicitly, the vertices of $\mathcal{G}\otimes\mathcal{H}$ are tensor products of the vertices of $\mathcal{G}$ and $\mathcal{H}$ and the edges are given by the following variant of Kashiwara's tensor product rule~\cite[Proposition~6]{Ka-crystal-bases}:
\begin{equation}\label{E:TensorProduct}
e_{i}(a\otimes x) = \begin{cases*}
e_{i}(a)\otimes x & if $\varphi_{i}(a)<\varepsilon_{i}(x)$,\\
a\otimes e_{i}(x) & if $\varphi_{i}(a)\geq\varepsilon_{i}(x)$,
\\
\end{cases*}
\quad\text{and}\quad
f_{i}(a\otimes x) = \begin{cases*}
f_{i}(a)\otimes x & if $\varphi_{i}(a)\geq\varepsilon_{i}(x)$,\\
a\otimes f_{i}(x) & if $\varphi_{i}(a)<\varepsilon_{i}(x)$.
\\
\end{cases*}
\end{equation}

Let $\Crystalgraph=\crystalgraph[{\fweight[\affine{\rho}_{1}]}]\otimes\dots\otimes\crystalgraph[{\fweight[\affine{\rho}_{\hell}]}]$
be the tensor product of the crystal graphs $\crystalgraph[{\fweight[\affine{\rho}_{1}]}],\dots,\crystalgraph[{\fweight[\affine{\rho}_{\hell}]}]$. By classical results, $\Crystalgraph[\abfweight]$ is a direct sum of highest weight crystals. In view of \autoref{E:TensorProduct}, the vertices of $\Crystalgraph[\abfweight]$ are tensor products $\lambda_{1}\otimes\dots\otimes\lambda_{\hell}$, where $\lambda_{k}$ is a vertex of $\crystalgraph[{\fweight[\affine{\rho}_{k}]}]$, and the labeled edges of $\Crystalgraph[\abfweight]$ are of the form $\lambda_{1}\otimes\dots\otimes\lambda_{\hell}\xrightarrow{i}\lambda^{\prime}_{1}\otimes\dots\otimes\lambda^{\prime}_{\hell}$ such that there exists an $m$ such that $\lambda_{k}=\lambda_{k}^{\prime}$ for $k\ne m$ and $\lambda_{m}\xrightarrow{i}\lambda^{\prime}_{m}$ is an edge in $\crystalgraph[{\fweight[\affine{\rho}_{m}]}]$.
Similarly, let $\Crystalgraph[\bfweight]=\crystalgraph[{\fweight[\rho_{1}]}]\otimes\dots\otimes\crystalgraph[{\fweight[\rho_{\ell}]}]$.

Our choice of tensor product rule in \autoref{E:TensorProduct} prefers using edges coming from the left-hand tensor factors before edges in the right-hand tensor factors. In terms of KLRW diagrams, this corresponds to preferring to put solid strings to the left of as many red strings as possible.

\begin{Lemma}\label{L:TensorInjection}
There is an injective map $\tau\colon\Crystalgraph[\bfweight]\to\Crystalgraph$ such that $\tau\bigl(f_{i}(x)\bigr)=f_{i}\bigl(\tau(x)\bigr)$, whenever $f_{i}(x)\neq0$ for $x\in\Crystalgraph[\bfweight]$ and $i\in\vertices$.
\end{Lemma}

\begin{proof}
For $1\leq r\leq\hell$ let $v_{\affine{\rho}_{r}}$ be the highest weight vector of weight $\Lambda_{\affine{\rho}_{r}}$ in $\crystalgraph[{\fweight[\affine{\rho}_{r}]}]$.
The vertices of $\Crystalgraph[\bfweight]$ are of the form $x=x_{1}\otimes\dots\otimes x_{\ell}$, for $x_{r}\in\crystalgraph[{\fweight[\rho_{r}]}]$. Define
$\tau(x)=x_{1}\otimes\dots\otimes x_{\ell}\otimes v_{\affine{\rho}_{\ell+1}}\otimes\dots\otimes
v_{\affine{\rho}_{\hell}}$,
a vertex of $\Crystalgraph$.
By \autoref{E:TensorProduct}, if $f_{i}(x)\neq0$ then $\tau\bigl(f_{i}(x)\bigr)=f_{i}\bigl(\tau(x)\bigr)$, for $i\in\vertices$ since $\varepsilon_{i}(v_{\affine{\rho}_{r}})=0$ for $\ell<r\leq\hell$. This completes the prooof.
\end{proof}

The map $\tau$ of \autoref{L:TensorInjection} is not a crystal embedding because the image of $\tau$ is not necessarily a connected component of $\Crystalgraph$. This will not matter for the results that follow.

Recall that $\crystalgraph[\Lambda_i]$ is the crystal graph of highest weight $\Lambda_i$. Abusing notation, let $\crystalgraph[\abfweight]$ and $\crystalgraph[\bfweight]$ be the crystal graphs of highest weights $\sum_{k=1}^{\hell}\Lambda_{\affine{\rho}_{k}}$ and $\sum_{k=1}^{\ell}\Lambda_{\rho_{k}}$, respectively.

\begin{Lemma}\label{L:Embeddings}
There are injective maps
$\crystalgraph[\bfweight]\hookrightarrow\Crystalgraph[\bfweight]$
and
$\crystalgraph[\abfweight]\hookrightarrow\Crystalgraph[\abfweight]$
that commute with the Kashiwara operators whenever they are nonzero.
\end{Lemma}

\begin{proof}
The existence of a crystal embedding $\tau\colon\crystalgraph[\bfweight]\hookrightarrow\Crystalgraph[\abfweight]$ is a well-known consequence of Kashiwara's tensor product rule~\autoref{E:TensorProduct}. It can be proved by repeating the arguments of \autoref{L:RootedPaths}.
\end{proof}

Recall from \autoref{D:FiniteTypesResidueSequence} that $\Parts{n}{\fweight}$ is the set of rooted paths of length $n$ in $\crystalgraph[\fweight]$. Similarly, let $\Parts{n}{\abfweight}$ and $\Parts{n}{\bfweight}$ be the sets of rooted paths of length $n$ in $\crystalgraph[\abfweight]$ and $\crystalgraph[\bfweight]$, respectively. We consider $\Parts{n}{\abfweight}$ and $\Parts{n}{\bfweight}$ as subcrystals of $\Crystalgraph[\abfweight]$ and $\Crystalgraph[\bfweight]$, respectively.

As the paths in $\Crystalgraph$ are given by applying the Kashiwara operators $f_{i}$, for $i\in\vertices$, the map $\tau$ of \autoref{L:TensorInjection} induces an injective, residue sequence preserving, map from $\Parts{n}{\bfweight}$ to $\Parts{n}{\abfweight}$.
Hereafter, we identify the rooted paths in
$\Parts{n}{\bfweight}$ and
$\Parts{n}{\abfweight}$ with the corresponding paths in $\Crystalgraph$ under the map $\tau$ of \autoref{L:TensorInjection}. Similarly, we identify the crystals
$\crystalgraph[\bfweight]$, $\Crystalgraph[\bfweight]$, and $\crystalgraph[\abfweight]$
 with their images in $\Crystalgraph$. In particular, we always think of the crystal graphs $\crystalgraph[\bfweight]$ and $\crystalgraph[\abfweight]$ as being embedded in $\Crystalgraph$. These embeddings are the crystal graph analog of the difference between $\Lambda=\sum_{k=1}^{\ell}\Lambda_{\rho_{k}}\in P^{+}$ and the ordered sequence $\rho\in\vertices^{\ell}$, which corresponds to the order of the residues on the red strings.

Let $\bpath$ be a path in $\Crystalgraph[\abfweight]$. Recall from the end of \autoref{SS:FacePermutations} that the \emph{sink}
$\omega(\bpath)$ of $\bpath$ is the final vertex
$\lambda_{1}\otimes\dots\otimes\lambda_{\hell}$ in the path $\bpath$. As
described above, each edge in $\bpath$ corresponds to a unique edge in
one of the tensor factors, $\crystalgraph[\Lambda_{\affine{\rho}_{k}}]$.
Let $\bpath_{k}$ be the (connected) path in
$\crystalgraph[\Lambda_{\affine{\rho}_{k}}]$ consisting of all of the
edges that come from $\crystalgraph[\Lambda_{\affine{\rho}_{k}}]$, for
$1\leq k\leq\hell$. We write $\bpath=(\bpath_{1},\dots,\bpath_{\hell})$.
%The \emph{residue sequence} $\Res[\bpath]$ of $\bpath$ is $\Res[\bpath]=\bigl(\Res[\bpath_{1}],\dots,\Res[\bpath_{\ell}],\Res[\bpath_{\ell+1}]\Res[\bpath_{\hell}]\bigr)\in\vertices^{n,\ell}$.

As above, if $\bpath$ is a path in $\Crystalgraph[\bfweight]$ let $\omega(\bpath)$ be its \emph{sink} and write $\bpath=(\bpath_{1},\dots,\bpath_{\ell})$, where $\bpath_{k}$ is the component path in $\crystalgraph[{\fweight[\rho_{k}]}]$. %Let $\Res[\bpath]=\bigl(\Res[\bpath_{1}],\dots,\Res[\bpath_{\ell}]\bigr)\in\vertices^{n,\ell}$ be the \emph{residue sequence} of $\bpath$.

We warn the reader that writing $\bpath=(\bpath_{1},\dots,\bpath_{\hell})$ is a slight abuse of notation because we have not shown that the path $\bpath$ is uniquely determined by the paths $\bpath_{1},\dots,\bpath_{\hell}$. For example, suppose that $\bfweight=(\Lambda_{i},\Lambda_{j})$, where $i,j\in\vertices$ are nonadjacent, and let $\bpath_{i}$ and $\bpath_{j}$ be the rooted paths of length $1$ in $\crystalgraph[\Lambda_{i}]$ and $\crystalgraph[\Lambda_{j}]$, respectively. Then there are, potentially, two paths in $\Crystalgraph[\bfweight]$ such that $\bpath=(\bpath_{1},\bpath_{2})$, which have residue sequences $(i,j)$ and $(j,i)$, respectively.

\begin{Lemma}\label{L:RootedPaths}
Suppose that $\bpath=(\bpath_{1},\dots,\bpath_{\hell})\in\Parts{n}{\abfweight}$.
Then $\bpath_{k}$ is a rooted path in $\crystalgraph[{\fweight[\affine{\rho}_{k}]}]$, for $1\leq k\leq\hell$.
\end{Lemma}

\begin{proof}
Let $\lambda_{k}$ be the source of $\bpath_{k}$, corresponding to a vertex $\blam=\lambda_{1}\otimes\dots\otimes\lambda_{\hell}$ in $\Crystalgraph$. By way of contradiction, suppose that $\lambda_{k}$ is not the highest weight vector in $\crystalgraph[{\fweight[\affine{\rho}_{k}]}]$. Then $\varepsilon_{j}(\lambda_{k})>0$ for some $j\in\vertices$. By \autoref{E:TensorProduct}, we can find $N=n_{1}+\dots+n_{k-1}+1$ such that $e_{j}^N\blam=e_{j}^{n_{1}}\lambda_{1}\otimes\dots\otimes e_{j}^{n_{k-1}}\lambda_{1}\otimes e_{j}\lambda_{k}\otimes\dots\otimes\lambda_{\hell}$ is nonzero, so $\lambda_{k}$ is not the source of $\bpath_{k}$, which is a contradiction. Hence, $\lambda_{k}$ is the highest weight vector in $\crystalgraph[{\fweight[\affine{\rho}_{k}]}]$, so $\bpath_{k}$ is a rooted path.
\end{proof}

To give a partial converse to \autoref{L:RootedPaths}, a \emph{rooted path} in $\Crystalgraph$, or in $\Crystalgraph[\bfweight]$, is a path with sink the highest weight vector of weight $\abfweight$, or $\bfweight$, respectively. If $\cpath$ is a path in a crystal graph, let $|\cpath|$ be its the \emph{length}.

\begin{Lemma}\label{L:TensorFactors}
For $1\leq m\leq\hell$, suppose that $\bpath_{m}$ is a rooted path in $\crystalgraph[{\fweight[\affine{\rho}_{m}]}]$. Then there is a rooted path $\bpath\in\Parts{n}{\abfweight}$ such that $\bpath=(\bpath_{1},\dots,\bpath_{\hell})$. Moreover, if $|\bpath_{m}|=0$ for $\ell<m\leq\hell$ then $\bpath\in\Parts{n}{\bfweight}$.
\end{Lemma}

\begin{proof}
We argue by induction on the total length $|\bpath_{1}|+\dots+|\bpath_{\hell}|$ of the paths. If they all have length zero then there is nothing to prove. Suppose that the total length of these paths is positive and fix $i\in\vertices$ so that at least one of the paths $\bpath_{1},\dots,\bpath_{\hell}$ ends in a path of residue $i$. Let $\lambda_{m}$ be the source of $\bpath_{m}$, for $1\leq m\leq\hell$, and set $\blam=\lambda_{1}\otimes\dots\otimes\lambda_{\hell}\in\Crystalgraph$. Let~$m$ be minimal such that
\begin{gather*}
\varepsilon_{i}(\lambda_{m})+\sum_{1\leq l<m}\bigl(\varepsilon_{i}(\lambda_{l})-\varphi_{i}(\lambda_{l})\bigr)\geq
\varepsilon_{i}(\lambda_{k})+\sum_{1\leq l<k}\bigl(\varepsilon_{i}(\lambda_{l})-\varphi_{i}(\lambda_{l})\bigr)
\end{gather*}
for $k<m$. Allowing for our different conventions for the tensor product of crystals~\autoref{E:TensorProduct}, \cite[Lemma~2.33]{BuSc-crystal-bases} implies that $e_{i}\blam=\lambda_{1}\otimes\dots\otimes e_{i}\lambda_{m}\otimes\dots\otimes\lambda_{\hell}$ in $\Crystalgraph$. Let $\bpath^{\prime}=\bpath_{k}$ if $k\neq m$ and let $\bpath^{\prime}_{m}$ be the path obtained from $\bpath_{m}$ by removing the last edge, which has residue $i$. By induction, there exists a path $\bpath^{\prime}$ in $\Crystalgraph$ such that $\bpath^{\prime}=(\bpath^{\prime}_{1},\dots,\bpath^{\prime}_{\hell})$. By what we have just said, $f_{i}\omega(\bpath^{\prime})=\blam$ in $\Crystalgraph$. Therefore, $\bpath=(\bpath_{1},\dots,\bpath_{\hell})$ where $\bpath$ is obtained from $\bpath^{\prime}$ by adding an edge from $\omega(\bpath^{\prime})$ to $\blam=\omega(\bpath)$. The final claim that $\bpath\in\Parts{n}{\bfweight}$ if $\bpath_{m}$ is empty for $\ell<m\le\hell$ is immediate from the construction.
\end{proof}

We emphasize that \autoref{L:TensorFactors} only claims that there exists at least one path $\bpath$ in $\Parts{n}{\abfweight}$ such that $\bpath=(\bpath_{1},\dots,\bpath_{\hell})$, where $\bpath_{m}$ is a path in $\Parts{n}{\fweight[\affine{\rho}_{m}]}$ for $1\leq m\leq\hell$. As shown in \autoref{Ex:Octogon} below, there may be more than one path $\bpath$ such that $\bpath=(\bpath_{1},\dots,\bpath_{\hell})$.

%%%%%%%%%%%%%%%%%%%%%%%%%%%%%%%%%%%%%%%%%

\section{Sandwich cellular bases}\label{S:Construction}

%%%%%%%%%%%%%%%%%%%%%%%%%%%%%%%%%%%%%%%%%

We are almost ready to describe homogeneous (affine) sandwich cellular bases for
$\WA[n](X)$ and $\WAc[n](X)$.

%%%%%%%%%%%%%%%%%%%%%%%%%%%%%%%%%%%%%%%%%

\subsection{Affine red strings}\label{SS:ConstructionAffine}

%%%%%%%%%%%%%%%%%%%%%%%%%%%%%%%%%%%%%%%%%

\begin{Notation}\label{N:ConstructionAffineHell}
Until \autoref{SS:Affine}, unless otherwise stated, we fix
a finite type Dynkin diagram $\quiver=(\vertices,\edges)$ and set $e=\#\vertices$.
We also fix $n,\ell\in\N$,
the number of solid and red
strings, $\charge\in\R^{\ell}$ with $\kappa_{1}<\dots<\kappa_{\ell}$
and $\brho\in\vertices^{\ell}$, giving the positions and the residues
of the red strings, respectively.

Recall from \autoref{D:hell} that we fixed $\bfweight=(\fweight[{\rho_{1}}],\dots,\fweight[{\rho_{\ell}}])$,
an $\ell$ tuple of fundamental weights.
As described in the next paragraph, we also fix ghost shifts for our choice of $\varepsilon$ in \autoref{Eq:ConstructionEpsilon}.

We will use a ghost shift of $1$
for all edges except for the fishtail edges (the ones pointing into fishtail vertices) in types
$\typed[e>3]$, $\typee[6]$, $\typee[7]$ or $\typee[8]$.
For the fishtail edges we chose a ghost shift of
$1$ for the edge pointing to the right in \autoref{SS:FiniteTypesConventions} and a ghost shift of $1-\varepsilon^{2}$ for the other edge.
The resulting ghost strings for these two edges are very close and in what follows it does not matter which of these strings is on the left because these strings satisfy an honest Reidemeister II relation. We draw these strings as a doubled line and do not distinguish between them. For example, in type $\typee[6]$ we write:
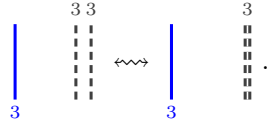
\begin{gather*}
\begin{tikzpicture}[anchorbase,smallnodes,rounded corners]
\draw[ghost](-0.2,0)--++(0,1)node[above,yshift=-1pt]{$3$};
\draw[ghost](0,0)--++(0,1)node[above,yshift=-1pt]{$3$};
\draw[solid](-1,0)node[below]{$3$}--++(0,1);
\end{tikzpicture}
\leftrightsquigarrow
\begin{tikzpicture}[anchorbase,smallnodes,rounded corners]
\draw[dghost](0,0)node[below]{$\phantom{i}$}--++(0,1)node[above,yshift=-1pt]{$3$};
\draw[solid](-1,0)node[below]{$3$}--++(0,1);
\end{tikzpicture}
.
\end{gather*}
For convenience, we only draw one ghost dot on the
doubled ghost string rather than a
dot on each ghost string, as required by the definitions.
\end{Notation}

Recall from \autoref{SS:RecollectionSteady} that steady and unsteady diagrams are nonzero in the infinite dimensional
affine wKLRW algebra $\WA(X)$, while only the steady
diagrams are nonzero in their finite dimension cyclotomic quotient $\WAc(X)$, since this algebra is  defined by quotienting out by the ideal generated by the unsteady diagrams.

In order to distinguish between steady and unsteady diagrams, we now introduce affine red strings. The notation that we use to describe these is taken, {\muta}, from \cite{MaTu-klrw-algebras}
and \cite{MaTu-klrw-algebras-bad}.
For everything affine, or unsteady, we use an underline
as part of the notation.

In \autoref{D:hell} we fixed an affine charge $\affine{\charge}=
(\affine{\kappa}_{1},\dots,\affine{\kappa}_{\hell})\in\Z^{\hell}$
and affine red labels
$\affine{\brho}=(\affine{\rho}_{1},\dots,\affine{\rho}_{\hell})\in \vertices^{\hell}$, which we
use to control the unsteady diagrams in $\WA(X)$. More precisely, as in \autoref{D:hell},
we place \emph{red strings} at positions $\kappa_{1},\dots,\kappa_{\ell}$, with residues $\rho_{1},\dots,\rho_{\ell},$. In addition, we imagine that there are \emph{affine red strings} at positions $\kappa_{\ell+1},\dots,\kappa_{\hell}$, with residues $\rho_{\ell+1},\dots,\rho_{\hell}$, that are  illustrated as:
\begin{gather*}
\text{genuine red string}:
\begin{tikzpicture}[anchorbase,smallnodes,rounded corners]
\draw[redstring] (0,0)node[below]{$i$} to (0,0.5)node[above,yshift=-1pt]{$\phantom{i}$};
\end{tikzpicture}
,\quad
\text{affine red string}:
\begin{tikzpicture}[anchorbase,smallnodes,rounded corners]
\draw[affine] (0,0)node[below]{$i$} to (0,0.5)node[above,yshift=-1pt]{$\phantom{i}$};
\end{tikzpicture}
.
\end{gather*}
The affine red strings are not part of the wKLRW diagrams:
they do not exist and are included only as a visual aid to mark the positions of the unsteady parts of the diagrams. By definition, unsteady strings can be pulled arbitrarily far to the right, but we think of them as being blocked by the affine red strings as in \autoref{E:RecollectionUnsteady}. By \autoref{Eq:RecollectionReidemeisterII}, a dotted $i$-string is not blocked by a red $i$-string. In contrast, we think of an $i$-string that is immediately to the left of an affine $i$-string, with arbitrarily many dots, as being blocked. On the other hand, when two solid $i$-strings are immediately adjacent to an affine $i$ string, as in \autoref{N:RecollectionPartner}.(a), then the left-hand $i$ string can be pulled past the affine $i$-string. Like red strings, an affine $i$-string does not block ghost strings or $j$-strings, for $j\ne i$. That is, we have:
\[
\begin{tikzpicture}[anchorbase,smallnodes,rounded corners]
\draw[solid,dot=0.5](0.5,1)--++(0,-1)node[below]{$i$};
\node at (0.5,-0.5){blocked};
\node at (0.5,1.5){};
\draw[affine](1,0)node[below]{$i$}--++(0,1);
\end{tikzpicture}~,
\hspace*{20mm}
\begin{tikzpicture}[anchorbase,smallnodes,rounded corners]
\draw[solid](0,1)node[above,yshift=-1pt]{$\phantom{i}$}--++(0,-1)node[below]{$i$};
\draw[solid](0.5,1)--++(0,-1)node[below]{$i$};
\draw[affine](1,0)node[below]{$i$}--++(0,1);
\end{tikzpicture}
=
\begin{tikzpicture}[anchorbase,smallnodes,rounded corners]
\draw[solid,dot=0.5](0,1)node[above,yshift=-1pt]{$\phantom{i}$}--++(0,-0.12)--++(1.2,-0.24)--++(0,-0.28)--++(-1.2,-0.24)--++(0,-0.12)node[below]{$i$};
\draw[solid,dot=0.1](0.5,1)--++(0,-1) node[below]{$i$};
\draw[affine](1,0)node[below]{$i$}--++(0,1);
\end{tikzpicture}
-
\begin{tikzpicture}[anchorbase,smallnodes,rounded corners]
\draw[solid,dot=0.5,dot=0.9](0,1)node[above,yshift=-1pt]{$\phantom{i}$}--++(0,-0.12)--++(1.2,-0.24)--++(0,-0.28)--++(-1.2,-0.24)--++(0,-0.12)node[below]{$i$};
\draw[solid](0.5,1)--++(0,-1) node[below]{$i$};
\draw[affine](1,0)node[below]{$i$}--++(0,1);
\end{tikzpicture}.
\]

The number of affine red strings is $\hell=\ell+ne$. This depends on $n$ because, as we will see, this ensures that there are enough affine red strings to block every solid string in any undotted diagram
with residue sequence $\bi\in\vertices^{n}$.

\begin{Example}\label{E:ConstructionAffineRed}
For $n=2$, $\ell=1$ and $e=2$ we have $\hell=1+2\cdot 2=5$. For $\charge=(0)$
and $\brho=(1)$ we have $\affine{\charge}=(0,3,6,9,12)$ and $\affine{\brho}=(1,1,1,2,2)$. Every diagram in this case has one
red string and four affine red strings positioned as
(the number above each string is its position):
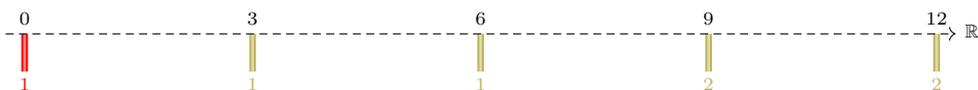
\begin{gather*}
\begin{tikzpicture}[anchorbase,smallnodes,rounded corners]
\draw[redstring] (0,0)node[below]{$1$} to (0,0.5)node[above,yshift=-1pt,black]{0};
\draw[affine] (3,0)node[below]{$1$} to (3,0.5)node[above,yshift=-1pt,black]{3};
\draw[affine] (6,0)node[below]{$1$} to (6,0.5)node[above,yshift=-1pt,black]{6};
\draw[affine] (9,0)node[below]{$2$} to (9,0.5)node[above,yshift=-1pt,black]{9};
\draw[affine] (12,0)node[below]{$2$} to (12,0.5)node[above,yshift=-1pt,black]{12};
\draw[densely dashed,->] (-0.25,0.5) to (12.25,0.5)node[right]{$\R$};
\end{tikzpicture}
.
\end{gather*}
Finally, if $\bfweight=(\fweight[1])$, then
$\abfweight=(\fweight[1],\fweight[1],\fweight[1],\fweight[2],\fweight[2])$. Note that we have enough affine red strings to block all solid strings (for all possible residues of these strings). For example, the relations in \autoref{L:RecollectionMovingStringsDots} give
\begin{gather*}
\begin{tikzpicture}[anchorbase,smallnodes,rounded corners]
\draw[solid](0,0)node[below]{$2$}--++(0,1)node[above,yshift=-1pt]{$\phantom{i}$};
\draw[redstring](0.5,0)node[below]{$1$}--++(0,1);
\end{tikzpicture}
=
\begin{tikzpicture}[anchorbase,smallnodes,rounded corners]
\draw[solid](0,1)node[above,yshift=-1pt]{$\phantom{i}$}--++(0.8,-0.5)--++(-0.8,-0.5)node[below]{$2$};
\draw[redstring](0.5,1)--++(0,-1) node[below]{$1$};
\end{tikzpicture}
,\quad
\begin{tikzpicture}[anchorbase,smallnodes,rounded corners]
\draw[solid](0,0)node[below]{$2$}--++(0,1)node[above,yshift=-1pt]{$\phantom{i}$};
\draw[affine](0.5,0)node[below]{$1$}--++(0,1);
\end{tikzpicture}
=
\begin{tikzpicture}[anchorbase,smallnodes,rounded corners]
\draw[solid](0,1)node[above,yshift=-1pt]{$\phantom{i}$}--++(0.8,-0.5)--++(-0.8,-0.5)node[below]{$2$};
\draw[affine](0.5,1)--++(0,-1) node[below]{$1$};
\end{tikzpicture}
,\quad
\begin{tikzpicture}[anchorbase,smallnodes,rounded corners]
\draw[solid](0,1)node[above,yshift=-1pt]{$\phantom{i}$}--++(0,-1)node[below]{$2$};
\draw[solid](0.5,1)--++(0,-1)node[below]{$2$};
\draw[affine](1,0)node[below]{$2$}--++(0,1);
\end{tikzpicture}
=
\begin{tikzpicture}[anchorbase,smallnodes,rounded corners]
\draw[solid,dot=0.5](0,1)node[above,yshift=-1pt]{$\phantom{i}$}--++(0,-0.12)--++(1.2,-0.24)--++(0,-0.28)--++(-1.2,-0.24)--++(0,-0.12)node[below]{$2$};
\draw[solid,dot=0.1](0.5,1)--++(0,-1) node[below]{$2$};
\draw[affine](1,0)node[below]{$2$}--++(0,1);
\end{tikzpicture}
-
\begin{tikzpicture}[anchorbase,smallnodes,rounded corners]
\draw[solid,dot=0.5,dot=0.9](0,1)node[above,yshift=-1pt]{$\phantom{i}$}--++(0,-0.12)--++(1.2,-0.24)--++(0,-0.28)--++(-1.2,-0.24)--++(0,-0.12)node[below]{$2$};
\draw[solid](0.5,1)--++(0,-1) node[below]{$2$};
\draw[affine](1,0)node[below]{$2$}--++(0,1);
\end{tikzpicture}
.
\end{gather*}
Unlike with red strings, a dotted solid $i$-string does not lose its dot when it is pushed through an affine $i$-string (after all, the affine strings are not really there).
These identities imply that we can pull two solid $2$-strings all the way to the right, so that the right-hand diagram factors through the diagram:
\begin{gather*}
\begin{tikzpicture}[anchorbase,smallnodes,rounded corners]
\draw[solid] (8.8,0)node[below]{$2$} to (8.8,0.5)node[above,yshift=-1pt]{$\phantom{i}$};
\draw[solid] (11.8,0)node[below]{$2$} to (11.8,0.5)node[above,yshift=-1pt]{$\phantom{i}$};
\draw[redstring] (0,0)node[below]{$1$} to (0,0.5)node[above,yshift=-1pt,black]{0};
\draw[affine] (3,0)node[below]{$1$} to (3,0.5)node[above,yshift=-1pt,black]{3};
\draw[affine] (6,0)node[below]{$1$} to (6,0.5)node[above,yshift=-1pt,black]{6};
\draw[affine] (9,0)node[below]{$2$} to (9,0.5)node[above,yshift=-1pt,black]{9};
\draw[affine] (12,0)node[below]{$2$} to (12,0.5)node[above,yshift=-1pt,black]{12};
\draw[densely dashed,->] (-0.25,0.5) to (12.25,0.5)node[right]{$\R$};
\end{tikzpicture}
.
\end{gather*}
Hence, the two affine red $2$-strings are needed to block the two solid $2$-strings.
\end{Example}

As in \autoref{E:ConstructionAffineRed}, in general there are $n$ affine red strings of any possible residue, which are placed far enough to the right so that they do not interfere with the steady strings in the diagram.

%%%%%%%%%%%%%%%%%%%%%%%%%%%%%%%%%%%%%%%%%

\subsection{Blocking strings}\label{SS:ConstructionBlocking}

%%%%%%%%%%%%%%%%%%%%%%%%%%%%%%%%%%%%%%%%%

The string $s$ in the next definition
corresponds to the rightmost string on the left-hand side of the identities in \autoref{L:RecollectionMovingStringsDots}.(b).

\begin{Definition}\label{D:ConstructionBlockingStrings}
For $i\in\vertices$,
assume there is a solid, ghost or red $j$-string $s$ in one of the following situations:
\begin{gather*}
s=
\begin{tikzpicture}[anchorbase,smallnodes,rounded corners]
\draw[redstring] (0,0)node[below]{$j$} to (0,0.5)node[above,yshift=-1pt]{$\phantom{i}$};
\end{tikzpicture}
\text{ and }
i=j
,\quad
s=
\begin{tikzpicture}[anchorbase,smallnodes,rounded corners]
\draw[solid] (0,0)node[below]{$j$} to (0,0.5)node[above,yshift=-1pt]{$\phantom{i}$};
\end{tikzpicture}
\text{ and }
i=j
,\quad
s=
\begin{tikzpicture}[anchorbase,smallnodes,rounded corners]
\draw[ghost] (0,0)node[below]{$\phantom{i}$} to (0,0.5)node[above,yshift=-1pt]{$j$};
\end{tikzpicture}
\text{ and }
i\leftsquigarrow j
,\quad\text{ or }\quad
s=
\begin{tikzpicture}[anchorbase,smallnodes,rounded corners]
\draw[solid] (0,0)node[below]{$j$} to (0,0.5)node[above,yshift=-1pt]{$\phantom{i}$};
\end{tikzpicture}
\text{ and }
i\rightsquigarrow j.
\end{gather*}
In the first three cases from the left, we say $s$
\emph{blocks solid $i$-strings}, and in the rightmost
case $s$ \emph{blocks ghost $i$-strings}. We allow affine
red strings to block strings in the same way that red strings do, as in the leftmost case above.
\end{Definition}

Similarly, solid $i$-strings are blocked by affine $i$-strings. Unlike red strings, which satisfy \autoref{Eq:RecollectionReidemeisterII}, we think of a dotted solid $i$-string as being blocked by an affine $i$-string.

Solid and ghost strings can
block as many strings with different residues as they have neighbors
in $\quiver$. In particular, if $j$ corresponds to a fishtail vertex in $\Gamma$, then it will block $i$-strings corresponding to its two neighbors.

\begin{Notation}\label{N:ConstructionBlockingStrings}
There will always be two cases depending on whether the string $s$ in \autoref{D:ConstructionBlockingStrings} blocks solid or ghost $i$-strings.
All notions involving these two cases run in parallel
and we say $s$ \emph{blocks $i$-strings} to cover both cases, allowing us to treat both cases {\muta}.
\end{Notation}

Given a positioning $\boldsymbol{x}$ and $\boldsymbol{j}\in\vertices^{n}$, the
\emph{idempotent diagram} $\idem[{(\boldsymbol{x},\boldsymbol{j})}]$ is the diagram with vertical solid $j_{k}$-strings at position $x_{k}$, for $1\leq k\leq n$, together with their ghost strings. The idempotent diagrams described in \autoref{E:MainExampleTheBeastItself} are all of this form.

\begin{Definition}\label{D:ConstructionParking}
Fix an idempotent diagram $\idem[{(\boldsymbol{x},\boldsymbol{j})}]$.
Let $s$ be the
$j_{k}$-string of $\idem[{(\boldsymbol{x},\boldsymbol{j})}]$ positioned
at $x_{j}\in\R$ and suppose that $s$ blocks $i$-strings. Then
the coordinates $x_{k}-m\varepsilon$
for $m\in\set{1,\dots,n}$ are the \emph{$i$-parking positions},
or \emph{parking positions} if $i$ is not specified.
The position $x_{k}-m\varepsilon$ is \emph{occupied} if $x_{k}-m\varepsilon\in\set{x_{1},\dots,x_{n}}$ and otherwise it is \emph{free}.
\end{Definition}

\begin{Example}\label{E:ConstructionParking}
Let $i\rightsquigarrow j$. The illustration
\begin{gather*}
\begin{tikzpicture}[anchorbase,smallnodes,rounded corners]
\draw[ghost] (-0.3,0)to (-0.3,0.5)node[above,yshift=-1pt]{$k$};
\draw[spinach,densely dotted] (-0.1,0)to (-0.1,0.5);
\draw[spinach,densely dotted] (-0.2,0)to (-0.2,0.5);
\draw[spinach,densely dotted] (-0.4,0)to (-0.4,0.5);
\draw[spinach,densely dotted] (-0.5,0)to (-0.5,0.5);
\draw[spinach,densely dotted] (-0.6,0)to (-0.6,0.5);
\draw[spinach,densely dotted] (-0.7,0)to (-0.7,0.5);
\draw[spinach,densely dotted] (-0.8,0)to (-0.8,0.5);
\draw[spinach,densely dotted] (-0.9,0)to (-0.9,0.5);
\draw[solid] (0,0)node[below]{$j$} to (0,0.5)node[above,yshift=-1pt]{$\phantom{i}$};
\draw[redstring] (-1,0)node[below]{$\rho$}to (-1,0.5);
\end{tikzpicture}
\end{gather*}
shows a solid $j$-string that blocks $i$-strings.
The first few free $i$-parking positions are illustrated by
dotted and colored strings. The third
parking position $x_{1}-3\varepsilon$ is occupied by a
ghost $k$-string, and the tenth position is occupied by a
red $\rho$-string. All other positions are free.
\end{Example}

\begin{Definition}\label{D:ConstructionParkingTwo}
If a $j$-string blocks $i$-strings, then an $i$-parking position for a $j$-string is \emph{admissible}
if it is the rightmost free parking position such that if an $i$-string is placed in this position, then the $i$ and $j$ strings can be pulled arbitrarily close to one another, in some neighborhood containing no other strings, using only isotopies and honest Reidemeister II relations.

If an $i$-string is in an admissible parking position for
a $j$-string, then the $i$-string is \emph{blocked}
by the $j$-string.
More generally, a $i$-string is blocked by a $j$-string if its solid or ghost string is blocked by the $j$-string.
We extend being blocked transitively, so that if an $i_{k+1}$-string is blocked by an $i_{k}$-string, for $1\leq k<m$, then the $i_{m}$-string is blocked by the $i_{1}$-string.
\end{Definition}

Note that if an $i$-string is in an admissible parking position for $j$-string, then the $i$-string is blocked by the $j$-string.
In particular, an honest Reidemeister II relation cannot be applied to pull the $i$-string past the $j$-string. Importantly, if the $j$-string is in position $x_k$, then the admissible parking position is not necessarily at $x_k-\varepsilon$ because, for example, the associated ghost and solid strings might be blocked and so prevent the $i$-string from reaching this position.

\begin{Example}\label{E:LeptinCheck}
An $i$-string being blocked by a $j$-string is not the same as the $i$-string being close to the $j$-string. For example, if $\Gamma$ is a quiver of type $\typed[e]$, then the both the ghost $(e-1)$ and $e$-strings are blocked by the solid $(e-2)$-string in the diagram:
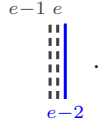
\begin{gather*}
\begin{tikzpicture}[anchorbase,smallnodes,rounded corners]
\draw[ghost] (0.8,0)--(0.8,1)node[above,shift={(-0.3,-1pt)}]{$e{-}1$};
\draw[ghost] (0.9,0)-- (0.9,1)node[above,yshift=-1pt]{$e$};
\draw[solid] (1,0)node[below]{$e{-}2$}--(1,1);
\end{tikzpicture}
.
\end{gather*}
Note that the ghost $(e-1)$ and $e$-strings do not block each other.
\end{Example}

Note that if a string is placed in an admissible $i$-parking position, then it can be blocked by a
$j$-string that is placed in an admissible $j$-parking position. If the $j$ string is not an (affine) red string, then, in turn, it can be blocked by a $k$-string that sits in an admissible parking position. As there are only finitely many parking positions, continuing in this way we can associate an (affine) red string to every admissible parking position, which we call its \emph{anchor}.

\begin{Lemma}\label{L:ConstructionParkingRed}
Every admissible parking position is associated to
a unique red string, its anchor.
\end{Lemma}

\begin{proof}
By \autoref{Eq:ConstructionEpsilon}, $\varepsilon$ is very small with respect to the charge $\affine\charge$, so the admissible parking position belongs to $\kappa_{m}-\set{1,\dots,n}\varepsilon+\Z$ for a unique $m$. By the paragraph above, the (affine) red string at position $\kappa_{m}$ is the anchor associated with this parking position.
\end{proof}

We use \autoref{L:ConstructionParkingRed} to associate
(affine) red strings as anchors for parking positions. Anchors block strings in the following sense: if an $i$-string is placed in an admissible parking position, then this $i$-string will be blocked by its anchor if we inductively add strings using the process described before \autoref{L:ConstructionParkingRed}.

As in \autoref{SS:RecollectionDiagrams}, the steady strings of the diagrams are blocked by the (genuine) red strings and the unsteady strings are blocked by the affine red strings. We will see that a string in an admissible parking possible is steady if and only if it is anchored by a red string and it is unsteady if and only it is anchored by an affine red string.
If $D$ is a diagram, then its \emph{steady part} consists of those strings that are anchored on the red strings, corresponding to $(\fweight[\rho_{1}],\dots,\fweight[\rho_{\ell}])$, and its \emph{unsteady part} consists of the strings anchored on the affine red strings, corresponding to $(\fweight[\rho_{\ell+1}],\dots,\fweight[\rho_{\hell}])$.

%%%%%%%%%%%%%%%%%%%%%%%%%%%%%%%%%%%%%%%%%

\subsection{Idempotent diagrams}\label{SS:ConstructionIdempotents}

%%%%%%%%%%%%%%%%%%%%%%%%%%%%%%%%%%%%%%%%%

We now describe how to construct idempotent diagrams
$\idem[\bpath]$ such that every string in $\idem[\bpath]$ is in an admissible parking position. To define the indexing set for these diagrams
if $\bi=(i_{1},\dots,i_{a})\in\vertices^{a}$ and $\bj=(j_{1},\dots,j_{b})\in\vertices^{b}$ let $\bi\bj=(i_{1},\dots,i_{a},j_{1},\dots,j_{b})\in\vertices^{a+b}$ be the concatenation of these two sequences. Define
\begin{gather*}
\vertices^{n,\ell}=\set[\big]{\bpath=(\bpath_{1},\dots,\bpath_{\ell})|
\bpath_{1}\dots\bpath_{\ell}\in\vertices^{n}}.
\end{gather*}
We will eventually see that the elements of $\vertices^{n,\ell}$ are the residue sequences of paths in an $\hell$-tuple of crystal graphs. As before, we call the elements of $\vertices^{n,\ell}$ \emph{residue sequences} and $\bpath_{1},\dots,\bpath_{\ell}$ are the \emph{components} of $\bpath$. For convenience, if $\bpath\in\vertices^{n,\ell}$ we write
$\bpath=(i_{11}\dots i_{1n_{1}},\dots,i_{\ell1}\dots i_{\ell n_{\ell}})$ if
$\bpath=\bigl((i_{1,1},\dots,i_{1,n_{1}}),\dots,(i_{\ell,1}\dots i_{\ell,n_{\ell}})\bigr)$. If $1\leq a\leq n$, then the $a$th residue in $\bpath$ is~$i_{k,m}$, where $a=k+\sum_{l<m}n_{l}$.

\begin{Definition}\label{D:ConstructionIdempotents}
Let $\bpath\in\vertices^{n,\ell}$. We inductively construct the idempotent diagram $\idem[{\qpath[]}]$ as follows.
\begin{enumerate}

\item Let $\idem[\bpath]^{0}$ be the diagram with no solid or ghost strings
and $\ell$ red strings, of residues $\rho_{m}$, placed at $\kappa_{m}$
in the steady part for $m\in\set{1,\dots,\ell}$. We imagine that there are $\hell-\ell$ affine red string of residues $\rho_{m}$ at positions $\kappa_{m}$, in the unsteady part, for $m\in\set{\ell+1,\dots,\hell}$.

\item Assume that $\idem[\bpath]^{a-1}$ is constructed
for $a\in\set{1,\dots,n}$, and suppose that $i_{k,m}$ is the $a$th residue in $\bpath$.
Then $\idem[\bpath]^{a}$ is the diagram obtained from $\idem[\bpath]^{a-1}$
by adding a solid $i_{k,m}$-string, and its ghosts, in the leftmost admissible parking position for an $i_{k,m}$-string in component $l$, where $l\leq\hell$ is minimal such that $l\geq m$.

\end{enumerate}
We set $\idem[\bpath]=\idem[\bpath]^{n}$. Let $\bx_{\bpath}=(x_{1}<\dots<x_{n})$ be the coordinates of the solid strings of $\idem[\bpath]$, read from left to right.
\end{Definition}

It is not clear that this construction is well-defined because, a priori, the diagram $\idem[\bpath]^{a-1}$ may not have an admissible parking position for an $i_{k,m}$-string. We prove this below. We also note that, in general, $\bx_{\bpath}\notin X$, so $\idem[\bpath]$ is not an element of $\WA(X)$. Nonetheless, we will show that every element of $\WA(X)$ \emph{factors through} these idempotents in the sense that every element of $\WA(X)$ is a linear combination of diagrams of the form $D\idem[\bpath]E$, for some diagrams $D$ and $E$ and $\bpath\in\vertices^{n,\ell}$.

We will eventually see that it is enough to consider idempotent diagrams associated to paths in crystal graphs.
If $\bpath=(\bpath_{1},\dots,\bpath_{\hell})$ is a path in the crystal graph $\crystalgraph[\abfweight]$, let $\idem[\bpath]$ be the idempotent diagram for the residue sequence $\Res[\bpath]=\big(\Res[\bpath_{1}],\dots,\Res[\bpath_{\ell}]\dots\Res[\bpath_{\hell}]\big)\in\vertices^{n,\ell}$. Similarly, if $\bpath=(\bpath_{1},\dots,\bpath_{\ell})$ is a path in $\crystalgraph[\bfweight]$ define the idempotent diagram $\idem[\bpath]$ using the residue sequence $\Res[\bpath]=(\Res[\bpath_{1}],\dots,\Res[\bpath_{\ell}])\in\vertices^{n,\ell}$.

\begin{Remark}\label{R:ConstructionIdempotents}
The conditions in \autoref{D:ConstructionIdempotents}(b) are not as convoluted as they appear. One should think of placing an $i_{k,m}$-string on the left of component $m$, and then pulling it to the right until it is blocked. The definition chooses the leftmost admissible parking position for the solid $i_{k,m}$-string in the $l$th component where~$l$ is minimal such that $m\leq l\leq\hell$.
\end{Remark}

\begin{Example}\label{E:ConstructionExampleOne}
In this example, let $\ell=1$, $\kappa=(0)$, $\brho=(1)$, and assume that ghost shifts are $1$. Below, we describe a path $\bpath$ in the crystal graph by giving its residue sequence $\Res[\bpath]$.
\begin{enumerate}

\item The main example \autoref{Eq:MainExample} shows how to construct the idempotent diagram $\idem[\mu]$ by adding one string at a time.

\item For the crystal graphs
in \autoref{Eq:ConstructionMainCrystals},
\begin{gather*}
\typea[4],\qpath[]=1234\colon
\begin{tikzpicture}[anchorbase,smallnodes,yscale=0.5]
\draw[ghost](1,0)--++(0,1)node[above,yshift=-1pt]{$1$};
\draw[ghost](1.9,0)--++(0,1)node[above,yshift=-1pt]{$2$};
\draw[ghost](2.8,0)--++(0,1)node[above,yshift=-1pt]{$3$};
\draw[ghost](3.7,0)--++(0,1)node[above,yshift=-1pt]{$4$};
\draw[solid](0,0)node[below]{$1$}--++(0,1);
\draw[solid](0.9,0)node[below]{$2$}--++(0,1);
\draw[solid](1.8,0)node[below]{$3$}--++(0,1);
\draw[solid](2.7,0)node[below]{$4$}--++(0,1);
\draw[redstring](0.15,0)node[below]{$1$}--++(0,1);
\end{tikzpicture}
,\quad
\typeb[4],\qpath[]=12344321\colon
\begin{tikzpicture}[anchorbase,smallnodes,yscale=0.5]
\draw[ghost](0.3,0)--++(0,1)node[above,yshift=-1pt]{$1$};
\draw[ghost](1,0)--++(0,1)node[above,yshift=-1pt]{$1$};
\draw[ghost](1.4,0)--++(0,1)node[above,yshift=-1pt]{$2$};
\draw[ghost](1.9,0)--++(0,1)node[above,yshift=-1pt]{$2$};
\draw[ghost](2.5,0)--++(0,1)node[above,yshift=-1pt]{$3$};
\draw[ghost](2.8,0)--++(0,1)node[above,yshift=-1pt]{$3$};
\draw[solid](-0.7,0)node[below]{$1$}--++(0,1);
\draw[solid](0,0)node[below]{$1$}--++(0,1);
\draw[solid](0.4,0)node[below]{$2$}--++(0,1);
\draw[solid](0.9,0)node[below]{$2$}--++(0,1);
\draw[solid](1.5,0)node[below]{$3$}--++(0,1);
\draw[solid](1.8,0)node[below]{$3$}--++(0,1);
\draw[solid](2.6,0)node[below,xshift=-0.03cm]{$4$}--++(0,1);
\draw[solid](2.7,0)node[below]{$4$}--++(0,1);
\draw[redstring](0.15,0)node[below]{$1$}--++(0,1);
\end{tikzpicture}
,
\\
\typec[4],\qpath[]=4321234\colon
\begin{tikzpicture}[anchorbase,smallnodes,yscale=0.5]
\draw[ghost](-2.3,0)--++(0,1)node[above,yshift=-1pt]{$1$};
\draw[ghost](-1.4,0)--++(0,1)node[above,yshift=-1pt]{$2$};
\draw[ghost](-1.2,0)--++(0,1)node[above,yshift=-1pt]{$2$};
\draw[ghost](-0.3,0)--++(0,1)node[above,yshift=-1pt]{$3$};
\draw[ghost](-0.1,0)--++(0,1)node[above,yshift=-1pt]{$3$};
\draw[solid](-3.3,0)node[below]{$1$}--++(0,1);
\draw[solid](-2.4,0)node[below]{$2$}--++(0,1);
\draw[solid](-2.2,0)node[below]{$2$}--++(0,1);
\draw[solid](-1.3,0)node[below]{$3$}--++(0,1);
\draw[solid](-1.1,0)node[below]{$3$}--++(0,1);
\draw[solid](-0.4,0)node[below]{$4$}--++(0,1);
\draw[solid](0,0)node[below]{$4$}--++(0,1);
\draw[redstring](0.15,0)node[below]{$4$}--++(0,1);
\end{tikzpicture}
,\quad
\typed[4],
\left\{
\begin{aligned}
\qpath[]=123421\colon&
\begin{tikzpicture}[anchorbase,smallnodes,yscale=0.5]
\draw[ghost](0.5,0)--++(0,1)node[above,yshift=-1pt]{$1$};
\draw[ghost](1,0)--++(0,1)node[above,yshift=-1pt]{$1$};
\draw[dghost](1.6,0)--++(0,1)node[above,yshift=-1pt]{$2$};
\draw[dghost](1.9,0)--++(0,1)node[above,yshift=-1pt]{$2$};
\draw[solid](-0.5,0)node[below]{$1$}--++(0,1);
\draw[solid](0,0)node[below]{$1$}--++(0,1);
\draw[solid](0.6,0)node[below]{$2$}--++(0,1);
\draw[solid](0.9,0)node[below]{$2$}--++(0,1);
\draw[solid](1.7,0)node[below,xshift=-0.03cm]{$4$}--++(0,1);
\draw[solid](1.8,0)node[below]{$3$}--++(0,1);
\draw[redstring](0.15,0)node[below]{$1$}--++(0,1);
\end{tikzpicture}
,
\\
\qpath[]=124321\colon&
\begin{tikzpicture}[anchorbase,smallnodes,yscale=0.5]
\draw[ghost](0.5,0)--++(0,1)node[above,yshift=-1pt]{$1$};
\draw[ghost](1,0)--++(0,1)node[above,yshift=-1pt]{$1$};
\draw[dghost](1.6,0)--++(0,1)node[above,yshift=-1pt]{$2$};
\draw[dghost](1.9,0)--++(0,1)node[above,yshift=-1pt]{$2$};
\draw[solid](-0.5,0)node[below]{$1$}--++(0,1);
\draw[solid](0,0)node[below]{$1$}--++(0,1);
\draw[solid](0.6,0)node[below]{$2$}--++(0,1);
\draw[solid](0.9,0)node[below]{$2$}--++(0,1);
\draw[solid](1.7,0)node[below,xshift=-0.03cm]{$3$}--++(0,1);
\draw[solid](1.8,0)node[below]{$4$}--++(0,1);
\draw[redstring](0.15,0)node[below]{$1$}--++(0,1);
\end{tikzpicture}
,
\end{aligned}
\right.
\end{gather*}
are the associated idempotent diagrams.
The two idempotent diagrams in type $\typed[4]$ are almost exactly the same except that the solid $3$ and $4$-strings are swapped. These two residue sequences correspond to different paths to the same vertex in the crystal graph $\crystalgraph[{\fweight[1]}]$, see \autoref{Eq:ConstructionMainCrystals} for this graph.

\item In type $\typeg[2]$ and
$\bpath=212212$ the idempotent diagram and crystal graph $\crystalgraph[{\fweight[2]}]$ are:
\begin{gather*}
\typeg[2],\bpath=212212\colon
\begin{tikzpicture}[anchorbase,smallnodes,yscale=0.5,every path/.style={very thick}]
\draw[ghost](-0.4,0)--++(0,1)node[above,yshift=-1pt]{$1$};
\draw[ghost](-0.1,0)--++(0,1)node[above,yshift=-1pt]{$1$};
\draw[solid](-0.5,0)node[below]{$2$}--++(0,1);
\draw[solid](-0.3,0)node[below]{$2$}--++(0,1);
\draw[solid](-0.2,0)node[below]{$2$}--++(0,1);
\draw[solid](-1.6,0)node[below]{$1$}--++(0,1);
\draw[solid](-1.1,0)node[below]{$1$}--++(0,1);
\draw[solid](0,0)node[below]{$2$}--++(0,1);
\draw[redstring](0.15,0)node[below]{$2$}--++(0,1);
\end{tikzpicture}
,\quad
\scalebox{0.75}{$
\begin{tikzpicture}[anchorbase,>=latex,line join=bevel,scale=0.35,xscale=1,every path/.style={very thick}]
\node (node_0) at (53.5bp,79.5bp) [draw,draw=none] {$\bullet$};
\node (node_6) at (53.5bp,8.5bp) [draw,draw=none] {$\bullet$};
\node (node_1) at (53.5bp,222.0bp) [draw,draw=none] {$\bullet$};
\node (node_3) at (53.5bp,150.5bp) [draw,draw=none] {$\bullet$};
\node (node_2) at (53.5bp,293.5bp) [draw,draw=none] {$\bullet$};
\node (node_4) at (53.5bp,364.5bp) [draw,draw=none] {$\bullet$};
\node (node_5) at (53.5bp,435.5bp) [draw,draw=none] {$\bullet$};
\draw [red,->] (node_0) ..controls (53.5bp,60.442bp) and (53.5bp,41.496bp)  .. (node_6);
\draw (62.0bp,44.0bp) node {\;$2$};
\draw [red,->] (node_1) ..controls (53.5bp,202.19bp) and (53.5bp,183.59bp)  .. (node_3);
\draw (62.0bp,186.0bp) node {\;$2$};
\draw [red,->] (node_2) ..controls (53.5bp,274.42bp) and (53.5bp,255.63bp)  .. (node_1);
\draw (62.0bp,258.0bp) node {\;$2$};
\draw [blue,->] (node_3) ..controls (53.5bp,131.44bp) and (53.5bp,112.5bp)  .. (node_0);
\draw (62.0bp,115.0bp) node {\;$1$};
\draw [blue,->] (node_4) ..controls (53.5bp,345.44bp) and (53.5bp,326.5bp)  .. (node_2);
\draw (62.0bp,329.0bp) node {\;$1$};
\draw [red,->] (node_5) ..controls (53.5bp,416.44bp) and (53.5bp,397.5bp)  .. (node_4);
\draw (62.0bp,400.0bp) node {\;$2$};
\end{tikzpicture}$}
.
\end{gather*}
\end{enumerate}
The reader should compare the examples in (b)
with the crystal graphs in \autoref{E:FiniteTypesCrystalList}.
\end{Example}

\begin{Example}\label{E:ConstructionExampleAnchor}
Suppose we are in type $\typea[e]$, and set $n=2$, $\ell=2$, for $e>2$, $\charge=(0,3)$, and $\brho=(1,2)$. Consider the residue sequences $(1,2)$ and $(12,\emptyset)$ in $\vertices^{2,2}$. The associated idempotent diagrams are:
\begin{gather*}
\idem[(1,2)]=
\begin{tikzpicture}[anchorbase,smallnodes,yscale=0.5]
%\draw[spinach,densely dotted] (0.9,0)--++(0,1);
\draw[ghost](1,0)--++(0,1)node[above,yshift=-1pt]{$1$};
\draw[ghost](4,0)--++(0,1)node[above,yshift=-1pt]{$2$};
\draw[solid](0,0)node[below]{$1$}--++(0,1);
\draw[solid](3,0)node[below]{$2$}--++(0,1);
\draw[redstring](0.15,0)node[below]{$1$}--++(0,1);
\draw[redstring](3.15,0)node[below]{$2$}--++(0,1);
\end{tikzpicture}
\qquad\text{and}\qquad
\idem[(12,\emptyset)]=
\begin{tikzpicture}[anchorbase,smallnodes,yscale=0.5]
%\draw[spinach,densely dotted] (0.9,0)--++(0,1);
\draw[ghost](1,0)--++(0,1)node[above,yshift=-1pt]{$1$};
\draw[ghost](1.85,0)--++(0,1)node[above,yshift=-1pt]{$2$};
\draw[solid](0,0)node[below]{$1$}--++(0,1);
\draw[solid](0.85,0)node[below]{$2$}--++(0,1);
\draw[redstring](0.15,0)node[below]{$1$}--++(0,1);
\draw[redstring](3.15,0)node[below]{$2$}--++(0,1);
\end{tikzpicture}
.
\end{gather*}
Note that in $\idem[(1,2)]$ the solid $2$-string is anchored on the red $2$-string even though
there is an admissible $2$-parking position to the left of the ghost
$1$-string.

In contrast, when placing strings in the unsteady part of the diagram we will always require that the leftmost admissible parking positions in the steady and unsteady part of the diagram are occupied first. In particular, applying our definitions gives the diagram
\begin{gather*}
\idem[(\emptyset,12)]=
\begin{tikzpicture}[anchorbase,smallnodes,yscale=0.5]
\draw[redstring](-5.85,0)node[below]{$1$}--++(0,1);
\draw[redstring](-2.85,0)node[below]{$2$}--++(0,1);
\draw[ghost](1,0)--++(0,1)node[above,yshift=-1pt]{$1$};
\draw[ghost](-2,0)--++(0,1)node[above,yshift=-1pt]{$2$};
\draw[solid](0,0)node[below]{$1$}--++(0,1);
\draw[solid](-3,0)node[below]{$2$}--++(0,1);
\draw[affine](0.15,0)node[below]{$1$}--++(0,1);
\draw[affine](3.15,0)node[below]{$2$}--++(0,1);
\end{tikzpicture}
.
\end{gather*}
This distinction between placing strings in the steady and unsteady parts of the diagram is necessary because the affine red strings are a notational sleight of hand, which requires slightly different combinatorics, because we want every string to be blocked by either a red string or an affine red string. In this example, the solid $1$-string is unsteady and is blocked by the first affine $1$-string, whereas the solid $2$-string is blocked by the red $2$-string.
\end{Example}

\begin{Example}\label{E:ConstructionExampleAffine}
In the setting of \autoref{E:ConstructionAffineRed}, consider a special case of relation \autoref{L:RecollectionMovingStringsDots}:
\begin{gather}\label{Eq:ConstructionEquator}
\begin{tikzpicture}[anchorbase,smallnodes,rounded corners]
\draw[solid](0,1)node[above,yshift=-1pt]{$\phantom{i}$}--++(0,-1)node[below]{$1$};
\draw[solid](0.5,1)--++(0,-1)node[below]{$1$};
\draw[redstring](0.9,0)node[below]{$1$}--++(0,1);
\end{tikzpicture}
=
\begin{tikzpicture}[anchorbase,smallnodes,rounded corners]
\draw[solid,dot](0,1)node[above,yshift=-1pt]{$\phantom{i}$}--++(0.8,-0.5)--++(-0.8,-0.5)node[below]{$1$};
\draw[solid,dot=0.1](0.5,1)--++(0,-1) node[below]{$1$};
\draw[redstring](0.9,0)node[below]{$1$}--++(0,1);
\end{tikzpicture}
-
\begin{tikzpicture}[anchorbase,smallnodes,rounded corners]
\draw[solid,dot,dot=0.9](0,1)node[above,yshift=-1pt]{$\phantom{i}$}--++(0.8,-0.5)--++(-0.8,-0.5)node[below]{$1$};
\draw[solid](0.5,1)--++(0,-1) node[below]{$1$};
\draw[redstring](0.9,0)node[below]{$1$}--++(0,1);
\end{tikzpicture}
=
\begin{tikzpicture}[anchorbase,smallnodes,rounded corners]
\draw[solid](0,1)node[above,yshift=-1pt]{$\phantom{i}$}--++(1.1,-0.5)--++(-1.1,-0.5)node[below]{$1$};
\draw[solid,dot=0.1](0.5,1)--++(0,-1) node[below]{$1$};
\draw[redstring](0.9,0)node[below]{$1$}--++(0,1);
\end{tikzpicture}
-
\begin{tikzpicture}[anchorbase,smallnodes,rounded corners]
\draw[solid,dot=0.9](0,1)node[above,yshift=-1pt]{$\phantom{i}$}--++(1.1,-0.5)--++(-1.1,-0.5)node[below]{$1$};
\draw[solid](0.5,1)--++(0,-1) node[below]{$1$};
\draw[redstring](0.9,0)node[below]{$1$}--++(0,1);
\end{tikzpicture}
,
\end{gather}
where we have ignored the ghost strings.
The left-hand side of this relation has a diagram with two consecutive solid $1$-strings. The crystal graph $\crystalgraph[{\fweight[1]}]$ does not contain a path with residue sequence $11$. On the other hand, if $\bpath=(11)\in\vertices^{2,1}$, then \autoref{D:ConstructionIdempotents} gives:
\begin{gather*}
\idem[\bpath]=
\begin{tikzpicture}[anchorbase,smallnodes,rounded corners]
\draw[solid](-0.2,0.5)--++(0,-0.5) node[below]{$1$};
\draw[solid](1.8,0.5)--++(0,-0.5) node[below]{$1$};
\draw[redstring] (0,0)node[below]{$1$} to (0,0.5)node[above,yshift=-1pt]{$\phantom{i}$};
\draw[affine] (2,0)node[below]{$1$} to (2,0.5);
\draw[affine] (4,0)node[below]{$1$} to (4,0.5);
\draw[affine] (6,0)node[below]{$2$} to (6,0.5);
\draw[affine] (8,0)node[below]{$2$} to (8,0.5);
\end{tikzpicture}
.
\end{gather*}
Note that this diagram is unsteady. We think of $\idem[\bpath]$ as being obtained from the relation above by pulling the solid $1$-string through the red string until it is close to the affine red string.
Observe that $\idem[\bpath]$ appears as the equator of the diagrams
on the right-hand side of \autoref{Eq:ConstructionEquator} above. We emphasize that the affine strings are not really part of the diagram $\idem[\bpath]$, however, the first affine $1$-string serves as a good visual aid because it appears to block the rightmost solid $1$-string, which could otherwise be pulled arbitrarily far to the right. In general, the affine strings give parking positions for the unsteady strings.
\end{Example}

\begin{Lemma}\label{L:ConstructionIdempotents}
Let $\bpath\in\vertices^{n,\ell}$.
The construction of $\idem[\bpath]$ in \autoref{D:ConstructionIdempotents} is well-defined.
\end{Lemma}

\begin{proof}
The definition of the affine charge in
\autoref{SS:ConstructionAffine} ensures that there are
$n$ affine red strings for each $i\in\vertices$.
The inductive construction of \autoref{D:ConstructionIdempotents} places the $i_{k,m}$-string from component $m$ to the left of the $m$th red string and then pulls it to the right into component~$l$, where $l\geq m$. In particular, unsteady strings are placed, from left to right, into components $\ell+1,\dots,\hell$.
As each $\bpath\in\vertices^{n,\ell}$ has $n$ solid strings, this implies the claim because there are enough admissible parking positions anchored on the affine red strings so that every string in $\bpath$ becomes anchored on either a red string or an affine string when it is pulled to the right.
\end{proof}

Let $\bot$ be the equivalence relation on $\vertices$ that is trivial except that $(e-1)\bot e$ in $\typed[e>3]$ and $4\bot 5$ in types $\typee[e]$, for $e\in\set{6,7,8}$, where we use the labeling of the Dynkin diagrams given in \autoref{SS:FiniteTypesConventions}. Note that $\bot$ is only concerned with fishtail vertices.

\begin{Lemma}\label{L:ProofsAnchor}
Let $s$ be a solid string in $\idem[\bpath]$, for $\bpath\in\vertices^{n,\ell}$. Then the coordinate of $s$ is uniquely determined by its components and its residue, up to $\bot$-equivalence.
\end{Lemma}

\begin{proof}
This follows directly from the definitions.
That is, the component can be read off from the
position of the anchor as
in \autoref{L:ConstructionParkingRed}, which is
given by the $\kappa_{m}$, and the residue
can be read off from the parking positions, up to $\bot$-equivalence.
\end{proof}

\iffalse
We can finally complete the definition of the wKLRW
algebras by defining the sets $X$ and $X$ of endpoints for the wKLRW diagrams.

\begin{Definition}\label{D:ConstructionEndpoints}
Let $\bpath\in\vertices^{n,\ell}$. Then $\idem[\bpath]$ is a \emph{red idempotent} if all of the strings in $\idem[\bpath]$ are anchored on red strings. Set
$X=\set{\bx_{\bpath}|\bpath\in\vertices^{n,\ell}}$ and
$X=\set{\bx_{\bpath}|\idem[\bpath]\text{ is a red idempotent}}$.
\end{Definition}

By definition, if $\bpath\in\vertices^{n,\ell}$, then $\idem[\bpath]$ is an idempotent in $\aWA$. In particular, as the examples above show, $\idem[\bpath]$ is not necessarily a red idempotent if $\bpath\in\vertices^{n,\ell}$. As \autoref{E:ConstructionExampleAffine} shows, different elements of $\vertices^{n,\ell}$ can give rise to the same idempotent in $\WA(X)$ or in $\aWA$.

Let $X^{\text{all}}\subset\R^{n}$ be the set of all allowed endpoints of diagrams satisfying the constraints of \autoref{SS:RecollectionDiagrams}.(b)(iii).
Then $\aWA$ is the idempotent subalgebra of the affine wKLRW algebra $\WA(X^{\text{all}})$, with truncation given by the sum of the idempotent diagrams $\idem[\bpath]$ for $\bpath\in\vertices^{n,\ell}$. The corresponding cyclotomic quotient $\WAc(X)$ can also be viewed as an idempotent subalgebra of $\WAc(X^{\text{all}})$.
\fi

%%%%%%%%%%%%%%%%%%%%%%%%%%%%%%%%%%%%%%%%%

\subsection{Plactic crystals}\label{SS:Plactic}

%%%%%%%%%%%%%%%%%%%%%%%%%%%%%%%%%%%%%%%%%

This section introduces a condition on idempotent diagrams that can be checked on the crystal graph $\crystalgraph[\Lambda]$, for $\Lambda\in P^{+}$. We use will this condition to control the detour permutations in $\WA(X)$.

A \emph{weighted plactic monoid move} is a multilocal change of strings of the form
\begin{gather}\label{Eq:ConstructionPlactic}
i\not\!\rightsquigarrow j\colon
\begin{tikzpicture}[anchorbase,smallnodes,rounded corners]
\draw[ghost](1,1)node[above,yshift=-1pt]{$i$}--++(0,-1)node[below]{$\phantom{i}$};
\draw[solid,smallnodes,rounded corners](1.5,1)--++(0,-1)node[below]{$j$};
\end{tikzpicture}
\to
\begin{tikzpicture}[anchorbase,smallnodes,rounded corners]
\draw[ghost](1.5,1)node[above,yshift=-1pt]{$i$}--++(0,-1)node[below]{$\phantom{i}$};
\draw[solid,smallnodes,rounded corners](1,1)--++(0,-1)node[below]{$j$};
\end{tikzpicture}
,\quad
i\rightarrow j\colon
\begin{tikzpicture}[anchorbase,smallnodes,rounded corners]
\draw[ghost](1.5,1)node[above,yshift=-1pt]{$i$}--++(0,-1)node[below]{$\phantom{i}$};
\draw[ghost](2,1)node[above,yshift=-1pt]{$i$}--++(0,-1)node[below]{$\phantom{i}$};
\draw[solid,smallnodes,rounded corners](2.5,1)--++(0,-1)node[below]{$j$};
\end{tikzpicture}
\to
\begin{tikzpicture}[anchorbase,smallnodes,rounded corners]
\draw[ghost](1,1)node[above,yshift=-1pt]{$i$}--++(0,-1)node[below]{$\phantom{i}$};
\draw[ghost](2,1)node[above,yshift=-1pt]{$i$}--++(0,-1)node[below]{$\phantom{i}$};
\draw[solid,smallnodes,rounded corners](1.5,1)--++(0,-1)node[below]{$j$};
\end{tikzpicture}
,\quad
i\Rightarrow j\colon
\begin{tikzpicture}[anchorbase,smallnodes,rounded corners]
\draw[ghost](1,1)node[above,yshift=-1pt]{$i$}--++(0,-1)node[below]{$\phantom{i}$};
\draw[ghost](1.5,1)node[above,yshift=-1pt]{$i$}--++(0,-1)node[below]{$\phantom{i}$};
\draw[ghost](2,1)node[above,yshift=-1pt]{$i$}--++(0,-1)node[below]{$\phantom{i}$};
\draw[solid,smallnodes,rounded corners](2.5,1)--++(0,-1)node[below]{$j$};
\end{tikzpicture}
\to
\begin{tikzpicture}[anchorbase,smallnodes,rounded corners]
\draw[ghost](0.5,1)node[above,yshift=-1pt]{$i$}--++(0,-1)node[below]{$\phantom{i}$};
\draw[ghost](1,1)node[above,yshift=-1pt]{$i$}--++(0,-1)node[below]{$\phantom{i}$};
\draw[ghost](2,1)node[above,yshift=-1pt]{$i$}--++(0,-1)node[below]{$\phantom{i}$};
\draw[solid,smallnodes,rounded corners](1.5,1)--++(0,-1)node[below]{$j$};
\end{tikzpicture}
,
\end{gather}
and a similar relation in case $i\Rrightarrow j$ (we do not need this one explicitly), or any of their partner moves. Two diagrams are \emph{plactic equivalent} if they related by a sequence of weighted plactic monoid moves.

\begin{Remark}\label{R:ConstructionPlactic}
This nomenclature comes from a special case of the \emph{Knuth relations} in the plactic monoid (see, for example, \cite[Chapter 8]{BuSc-crystal-bases}),
which are the relations $iji=iij$ and $iji=jii$.
\end{Remark}

The \emph{plactic moves on residues sequences} are:
\begin{gather}\label{Eq:ConstructionPlacticTwo}
i\not\!\rightsquigarrow j\colon
ij\to ji
,\quad
i\rightarrow j\colon
iij\to iji
,\quad
i\Rightarrow j\colon
iiij\to iiji
,\quad
i\Rrightarrow j\colon
iiiij\to iiiji.
\end{gather}
In practise, the weighted plactic monoid moves on idempotent diagrams can be detected by looking at the residue sequences:

\begin{Lemma}\label{L:ConstructionPlacticMoves}
If two residue sequences are related by plactic moves, then
the corresponding idempotent diagrams are related by weighted plactic moves.
\end{Lemma}

\begin{proof}
By comparing \autoref{Eq:ConstructionPlactic} with
\autoref{Eq:ConstructionPlacticTwo}, this follows by
the construction of the idempotent diagram $\idem[\bpath]$ from its residue sequences, for $\bpath\in\vertices^{n,\ell}$.
\end{proof}

Recall the definition of adjacent squares $ij=ji$, for $i\rightsquigarrow j$, in the crystal graph $\crystalgraph[{\fweight[k]}]$ from \autoref{SS:FiniteTypesCrystalsSSR}. By \autoref{D:FiniteTypesSSR}, the source of a square is its unique vertex $\blam$ of minimal distance $d(\blam)$ from $\fweight[k]$.

\begin{Definition}\label{D:ConstructionPlactic}
A crystal graph $\crystal$ is \emph{plactic} if whenever $\lambda$ is the source of an adjacent square $ij=ji$ in $\crystal$ then one of the following holds:
\begin{enumerate}

\item There exists a path $\bpath$ to $\blam$ that ends in an edge labeled $i$ or $j$.

\item There exists a path $\bpath$ to $\blam$ and an idempotent diagram $\idem[(\bx,\bi)]$ with $\bi=\Res[\bpath]$
that ends with an edge labeled $i$ or $j$ such that $\idem[\bpath]$ and $\idem[(\bx,\bi)]$ are plactic equivalent.

\end{enumerate}
\end{Definition}

In particular, crystal graphs without adjacent squares are plactic, which includes all fundamental crystals for quivers of type $\typea[e\ge1]$ by \autoref{P:FiniteTypesSSR}.

\begin{Example}\label{E:ConstructionPlactic}
Let us list the three fundamental crystal graphs in type $\typec[3]$:
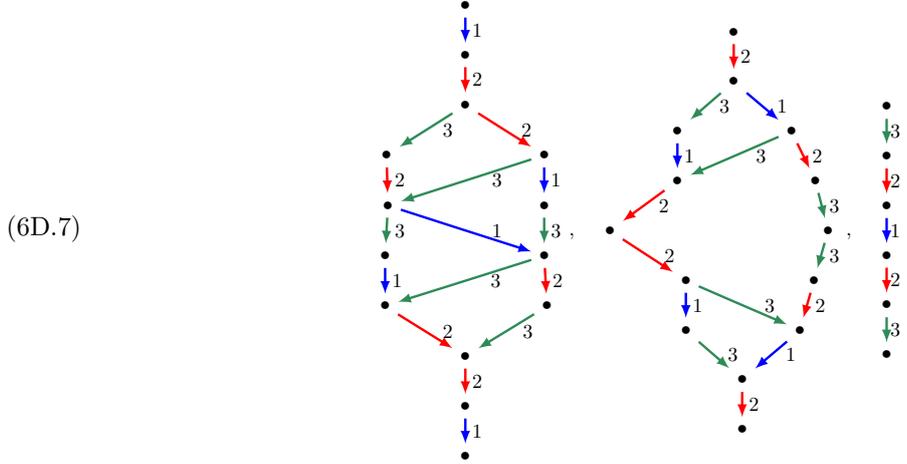
\begin{gather}\label{Eq:ConstructionTypeC3}
\scalebox{0.75}{$\begin{tikzpicture}[anchorbase,>=latex,line join=bevel,xscale=0.7,yscale=0.35,every path/.style={very thick}]
\node (node_0) at (61.5bp,367.0bp) [draw,draw=none] {$\bullet$};
\node (node_3) at (59.5bp,295.0bp) [draw,draw=none] {$\bullet$};
\node (node_7) at (172.5bp,295.0bp) [draw,draw=none] {$\bullet$};
\node (node_1) at (116.5bp,653.5bp) [draw,draw=none] {$\bullet$};
\node (node_4) at (116.5bp,582.5bp) [draw,draw=none] {$\bullet$};
\node (node_2) at (172.5bp,367.0bp) [draw,draw=none] {$\bullet$};
\node (node_8) at (59.5bp,223.0bp) [draw,draw=none] {$\bullet$};
\node (node_6) at (116.5bp,511.0bp) [draw,draw=none] {$\bullet$};
\node (node_5) at (60.5bp,439.0bp) [draw,draw=none] {$\bullet$};
\node (node_13) at (172.5bp,439.0bp) [draw,draw=none] {$\bullet$};
\node (node_12) at (174.5bp,223.0bp) [draw,draw=none] {$\bullet$};
\node (node_10) at (116.5bp,151.0bp) [draw,draw=none] {$\bullet$};
\node (node_9) at (116.5bp,79.5bp) [draw,draw=none] {$\bullet$};
\node (node_11) at (116.5bp,8.5bp) [draw,draw=none] {$\bullet$};
\draw [spinach,->] (node_0) ..controls (60.963bp,347.2bp) and (60.408bp,327.79bp)  .. (node_3);
\draw (70.0bp,331.0bp) node {$3$};
\draw [blue,->] (node_0) ..controls (92.952bp,346.17bp) and (128.14bp,323.98bp)  .. (node_7);
\draw (139.0bp,331.0bp) node {$1$};
\draw [blue,->] (node_1) ..controls (116.5bp,634.44bp) and (116.5bp,615.5bp)  .. (node_4);
\draw (125.0bp,618.0bp) node {$1$};
\draw [spinach,->] (node_2) ..controls (172.5bp,347.68bp) and (172.5bp,328.46bp)  .. (node_7);
\draw (181.0bp,331.0bp) node {$3$};
\draw [blue,->] (node_3) ..controls (59.5bp,275.57bp) and (59.5bp,256.05bp)  .. (node_8);
\draw (68.0bp,259.0bp) node {$1$};
\draw [red,->] (node_4) ..controls (116.5bp,563.42bp) and (116.5bp,544.63bp)  .. (node_6);
\draw (125.0bp,547.0bp) node {$2$};
\draw [red,->] (node_5) ..controls (60.765bp,419.43bp) and (61.034bp,400.61bp)  .. (node_0);
\draw (70.0bp,403.0bp) node {$2$};
\draw [spinach,->] (node_6) ..controls (101.23bp,490.91bp) and (85.078bp,470.72bp)  .. (node_5);
\draw (104.0bp,475.0bp) node {$3$};
\draw [red,->] (node_6) ..controls (131.96bp,490.67bp) and (148.62bp,469.85bp)  .. (node_13);
\draw (160.0bp,475.0bp) node {$2$};
\draw [spinach,->] (node_7) ..controls (140.18bp,273.98bp) and (103.51bp,251.26bp)  .. (node_8);
\draw (138.0bp,259.0bp) node {$3$};
\draw [red,->] (node_7) ..controls (173.03bp,275.43bp) and (173.57bp,256.61bp)  .. (node_12);
\draw (182.0bp,259.0bp) node {$2$};
\draw [red,->] (node_8) ..controls (74.843bp,203.16bp) and (91.636bp,182.53bp)  .. (node_10);
\draw (104.0bp,187.0bp) node {$2$};
\draw [blue,->] (node_9) ..controls (116.5bp,60.442bp) and (116.5bp,41.496bp)  .. (node_11);
\draw (125.0bp,44.0bp) node {$1$};
\draw [red,->] (node_10) ..controls (116.5bp,131.19bp) and (116.5bp,112.59bp)  .. (node_9);
\draw (125.0bp,115.0bp) node {$2$};
\draw [spinach,->] (node_12) ..controls (158.68bp,202.91bp) and (141.96bp,182.72bp)  .. (node_10);
\draw (161.0bp,187.0bp) node {$3$};
\draw [spinach,->] (node_13) ..controls (141.54bp,418.48bp) and (105.85bp,395.97bp)  .. (node_0);
\draw (139.0bp,403.0bp) node {$3$};
\draw [blue,->] (node_13) ..controls (172.5bp,419.57bp) and (172.5bp,400.05bp)  .. (node_2);
\draw (181.0bp,403.0bp) node {$1$};
\end{tikzpicture}
,\qquad
\begin{tikzpicture}[anchorbase,>=latex,line join=bevel,xscale=0.7,yscale=0.35,every path/.style={very thick}]
\node (node_0) at (238.5bp,293.0bp) [draw,draw=none] {$\bullet$};
\node (node_4) at (228.5bp,221.5bp) [draw,draw=none] {$\bullet$};
\node (node_1) at (177.5bp,79.5bp) [draw,draw=none] {$\bullet$};
\node (node_13) at (177.5bp,8.5bp) [draw,draw=none] {$\bullet$};
\node (node_2) at (131.5bp,435.5bp) [draw,draw=none] {$\bullet$};
\node (node_5) at (131.5bp,364.5bp) [draw,draw=none] {$\bullet$};
\node (node_3) at (137.5bp,150.5bp) [draw,draw=none] {$\bullet$};
\node (node_11) at (218.5bp,150.5bp) [draw,draw=none] {$\bullet$};
\node (node_6) at (83.5bp,293.0bp) [draw,draw=none] {$\bullet$};
\node (node_10) at (137.5bp,221.5bp) [draw,draw=none] {$\bullet$};
\node (node_7) at (171.5bp,506.5bp) [draw,draw=none] {$\bullet$};
\node (node_12) at (212.5bp,435.5bp) [draw,draw=none] {$\bullet$};
\node (node_8) at (229.5bp,364.5bp) [draw,draw=none] {$\bullet$};
\node (node_9) at (171.5bp,577.5bp) [draw,draw=none] {$\bullet$};
\draw [spinach,->] (node_0) ..controls (235.78bp,273.09bp) and (233.06bp,254.21bp)  .. (node_4);
\draw (243.0bp,257.0bp) node {$3$};
\draw [red,->] (node_1) ..controls (177.5bp,60.442bp) and (177.5bp,41.496bp)  .. (node_13);
\draw (186.0bp,44.0bp) node {$2$};
\draw [blue,->] (node_2) ..controls (131.5bp,416.44bp) and (131.5bp,397.5bp)  .. (node_5);
\draw (140.0bp,400.0bp) node {$1$};
\draw [spinach,->] (node_3) ..controls (148.15bp,131.13bp) and (159.61bp,111.36bp)  .. (node_1);
\draw (171.0bp,115.0bp) node {$3$};
\draw [red,->] (node_4) ..controls (225.88bp,202.44bp) and (223.14bp,183.5bp)  .. (node_11);
\draw (233.0bp,186.0bp) node {$2$};
\draw [red,->] (node_5) ..controls (118.73bp,345.01bp) and (105.0bp,325.13bp)  .. (node_6);
\draw (122.0bp,329.0bp) node {$2$};
\draw [red,->] (node_6) ..controls (98.524bp,272.66bp) and (114.03bp,252.71bp)  .. (node_10);
\draw (126.0bp,257.0bp) node {$2$};
\draw [spinach,->] (node_7) ..controls (160.85bp,487.13bp) and (149.39bp,467.36bp)  .. (node_2);
\draw (165.0bp,471.0bp) node {$3$};
\draw [blue,->] (node_7) ..controls (182.41bp,487.13bp) and (194.16bp,467.36bp)  .. (node_12);
\draw (206.0bp,471.0bp) node {$1$};
\draw [spinach,->] (node_8) ..controls (231.84bp,345.42bp) and (234.27bp,326.63bp)  .. (node_0);
\draw (243.0bp,329.0bp) node {$3$};
\draw [red,->] (node_9) ..controls (171.5bp,558.44bp) and (171.5bp,539.5bp)  .. (node_7);
\draw (180.0bp,542.0bp) node {$2$};
\draw [blue,->] (node_10) ..controls (137.5bp,202.44bp) and (137.5bp,183.5bp)  .. (node_3);
\draw (146.0bp,186.0bp) node {$1$};
\draw [spinach,->] (node_10) ..controls (159.67bp,201.62bp) and (184.53bp,180.44bp)  .. (node_11);
\draw (197.0bp,186.0bp) node {$3$};
\draw [blue,->] (node_11) ..controls (207.59bp,131.13bp) and (195.84bp,111.36bp)  .. (node_1);
\draw (212.0bp,115.0bp) node {$1$};
\draw [spinach,->] (node_12) ..controls (190.33bp,415.62bp) and (165.47bp,394.44bp)  .. (node_5);
\draw (191.0bp,400.0bp) node {$3$};
\draw [red,->] (node_12) ..controls (216.95bp,416.44bp) and (221.62bp,397.5bp)  .. (node_8);
\draw (231.0bp,400.0bp) node {$2$};
\end{tikzpicture}
,\qquad
\begin{tikzpicture}[anchorbase,>=latex,line join=bevel,scale=0.35,every path/.style={very thick}]
\node (node_0) at (27.5bp,79.5bp) [draw,draw=none] {$\bullet$};
\node (node_1) at (27.5bp,8.5bp) [draw,draw=none] {$\bullet$};
\node (node_2) at (27.5bp,150.5bp) [draw,draw=none] {$\bullet$};
\node (node_3) at (27.5bp,292.5bp) [draw,draw=none] {$\bullet$};
\node (node_4) at (27.5bp,221.5bp) [draw,draw=none] {$\bullet$};
\node (node_5) at (27.5bp,363.5bp) [draw,draw=none] {$\bullet$};
\draw [spinach,->] (node_0) ..controls (27.5bp,60.442bp) and (27.5bp,41.496bp)  .. (node_1);
\draw (36.0bp,44.0bp) node {\;$3$};
\draw [red,->] (node_2) ..controls (27.5bp,131.44bp) and (27.5bp,112.5bp)  .. (node_0);
\draw (36.0bp,115.0bp) node {\;$2$};
\draw [red,->] (node_3) ..controls (27.5bp,273.44bp) and (27.5bp,254.5bp)  .. (node_4);
\draw (36.0bp,257.0bp) node {\;$2$};
\draw [blue,->] (node_4) ..controls (27.5bp,202.44bp) and (27.5bp,183.5bp)  .. (node_2);
\draw (36.0bp,186.0bp) node {\;$1$};
\draw [spinach,->] (node_5) ..controls (27.5bp,344.44bp) and (27.5bp,325.5bp)  .. (node_3);
\draw (36.0bp,328.0bp) node {\;$3$};
\end{tikzpicture}$}
.
\end{gather}
All of these are plactic: for the two graphs on the right this is clear, because they do not have any adjacent squares.
As predicted by \autoref{P:FiniteTypesSSR},
the left-hand graph contains two adjacent squares $23=32$ one of which has an incoming edge labeled $2$ and the other has an incoming edge labeled by $3$. In both cases, these squares are preceded by an edge that also appears in the square itself. Hence, all of these crystals satisfy \autoref{D:ConstructionPlactic}.(a) and so are plactic

An example of a crystal that satisfies \autoref{D:ConstructionPlactic}.(b) is the crystal graph $\crystalgraph[{\fweight[2]}]$ of type $\typef[4]$, which
contains the path $231233241432$ to the source $\lambda$ of a $34=43$ square. (With $1274$ vertices, this crystal graph is too large to display but can be viewed using SageMath~\cite{sage}.) No path to $\lambda$ ends in an edge labeled $3$ or $4$, however, one can check that the following
is a sequence of plactic moves on this residue sequence:
\begin{gather*}
2312332\colorbox{spinach!50}{41}432
\xrightarrow{41\to 14}
23123321\colorbox{tomato!50}{443}2
\xrightarrow{443\to 434}
2312332143\colorbox{orchid!50}{42}
\xrightarrow{42\to 24}
231233214324
,
\end{gather*}
where we use \autoref{L:ConstructionPlacticMoves}. The plactic moves allow us to move the residue $4$ to the end of the sequence.
\end{Example}

\begin{Lemma}\label{L:ConstructionPlactic}
We have the following.
\begin{enumerate}

\item All fundamental crystal graphs of classical type and of type $\typeg$ are plactic.

\item All fundamental crystal graphs of types $\typee[6]$ and
$\typef$ are plactic.

\item In type $\typee[7]$ the fundamental crystal graphs $\crystalgraph[{\fweight[k]}]$
for $k\in\set{1,2,3,6,7}$ are plactic.

\item In type $\typee[8]$ the fundamental crystal graphs
$\crystalgraph[{\fweight[k]}]$ for $k\in\set{1,8}$ are plactic.

\end{enumerate}
\end{Lemma}

\begin{proof}
Part (a) follows from \autoref{P:FiniteTypesSSR} for classical types and from \autoref{Eq:FiniteTypesG2} for type $\typeg[2]$. Parts (b)--(d) were verified with the code from \cite{MaTu-sagemath-finite-type-klrw}, which uses \autoref{L:ConstructionPlacticMoves}.
\end{proof}

For types $\typee[7]$ and $\typee[8]$ we are only claiming that the fundamental crystals listed above are plactic, and not that the missing fundamental crystals are not plactic. It turns out that knowing that the crystals listed are plactic is enough to show that all crystals are plactic (and computationally, the  remaining cases take a long time to check).

\begin{Lemma}\label{L:PlacticTensor}
Suppose that $\mathcal{G}$ and $\mathcal{H}$ are plactic crystal graphs. Then $\mathcal{G}\otimes\mathcal{H}$ is plactic.
\end{Lemma}

\begin{proof}
The vertices of the crystal graph $\mathcal{G}\otimes\mathcal{H}$ are tensor products of the vertices of $\mathcal{G}$ and $\mathcal{H}$ and the Kashiwara operator $f_{i}$, for $i\in\vertices$, acts via \autoref{E:TensorProduct}.
Therefore, any $ij$-square in $\mathcal{G}\otimes\mathcal{H}$ either comes from an $ij$-square in $\mathcal{G}$ or $\mathcal{H}$ or it comes from the tensor product of subgraphs:
\begin{gather*}
\begin{tikzpicture}[anchorbase]
\node (a) at (0,0){$a$};
\node (b) at (1,0){$b$};
\draw[spinach,directed=0.55,line width=1.2](a)--node[above]{$i$}(b);
\end{tikzpicture}
\quad\otimes\quad
\begin{tikzpicture}[anchorbase]
\node (x) at (0,0){$x$};
\node (y) at (1,0){$y$};
\draw[orchid,directed=0.55,line width=1.2](x)--node[above,xshift=0.02cm]{$j$}(y);
\end{tikzpicture}
\quad=\quad
\begin{tikzpicture}[anchorbase]
\node (ax) at (0,2){$a\otimes x$};
\node (bx) at (-1,1){$b\otimes x$};
\node (ay) at (1,1){$a\otimes y$};
\node (by) at (0,0){$b\otimes y$};
\draw[spinach,directed=0.55,line width=1.2](ax)--node[above,xshift=-0.2cm]{$i$}(bx);
\draw[orchid,directed=0.55,line width=1.2](ax)--node[above,xshift=0.2cm]{$j$}(ay);
\draw[spinach,directed=0.55,line width=1.2,orchid](bx)--node[below]{$j$}(by);
\draw[spinach,directed=0.55,line width=1.2](ay)--node[below]{$i$}(by);
\end{tikzpicture}
\end{gather*}
By \autoref{D:ConstructionPlactic}(a), to prove that this square is plactic it is enough to show that there is an edge labeled by $i$ or $j$ coming into the vertex $a\otimes x$. That is, it is enough to show that $\varepsilon_{i}(a\otimes x)>0$ or $\varepsilon_{j}(a\otimes x)>0$. By \autoref{E:TensorProduct}, if $\varepsilon_{i}(a)>0$ then there is an $i$-edge coming into $a\otimes x$. Similarly, if $\varepsilon_{i}(x)>0$ there is a $j$-edge coming into $a\otimes x$. Hence, the $ij$-square is plactic in both of these cases, so it remains to consider the case when $\varepsilon_{j}(a)=0=\varepsilon_{i}(x)$.

We first assume that the quiver $\quiver$ is simply laced. By the Stembridge axioms~(see \cite[Theorem~2.4]{St-local-crystal-simply-laced} or \cite[Proposition~4.5]{BuSc-crystal-bases}), we are in one of the following four mutually exclusive cases:
\begin{itemize}
\item $\varphi_{j}(a)=\varphi_{j}(b)$ and $\varepsilon_{j}(a)=\varepsilon_{j}(b)+1$, and
$\varphi_{i}(x)=\varphi_{i}(y)$ and $\varepsilon_{i}(x)=\varepsilon_{i}(y)+1$,
\item $\varphi_{j}(a)=\varphi_{j}(b)$ and $\varepsilon_{j}(a)=\varepsilon_{j}(b)+1$, and
$\varphi_{i}(x)=\varphi_{i}(y)+1$ and $\varepsilon_{i}(x)=\varepsilon_{i}(y)$,
\item $\varphi_{j}(a)=\varphi_{j}(b)+1$ and $\varepsilon_{j}(a)=\varepsilon_{j}(b)$, and
$\varphi_{i}(x)=\varphi_{i}(y)$ and $\varepsilon_{i}(x)=\varepsilon_{i}(y)+1$,
\item $\varphi_{j}(a)=\varphi_{j}(b)+1$ and $\varepsilon_{j}(a)=\varepsilon_{j}(b)$, and
$\varphi_{i}(x)=\varphi_{i}(y)+1$ and $\varepsilon_{i}(x)=\varepsilon_{i}(y)$.
\end{itemize}
Since $\varepsilon_{j}(a)=0=\varepsilon_{i}(x)$, the Stembridge axioms imply that  $\varphi_{j}(a)=\varphi_{j}(b)+1$ and $\varphi_{i}(x)=\varphi_{i}(y)+1$. However, this gives a contradiction because \autoref{E:TensorProduct} now implies that $\mathcal{G}\otimes\mathcal{H}$ contains the edge $a\otimes x\xrightarrow{j}f_{j}a\otimes x$ and so, in particular, the $ij$-square above cannot not appear in the crystal $\mathcal{G}\otimes\mathcal{H}$.

Now suppose that $\quiver$ is doubly laced. By assumption, $\varepsilon_{j}(a)=0=\varepsilon_{i}(x)$ in the displayed $ij$-square above. If $i$ and $j$ are both simply laced edges then the argument of the last paragraph using the Stembridge axioms shows that this $ij$-square cannot appear in $\mathcal{G}\otimes\mathcal{H}$. If one of $i$ or $j$ is a doubly laced edge then Sternberg~\cite[Theorem~1]{St-local-crystal-doubly-laced} shows that $a_{ij}=0$ since $\varepsilon_{j}(a)=0=\varepsilon_{i}(x)$. Hence, $\mathcal{G}\otimes\mathcal{H}$ is plactic.

Finally, it remains to consider the quiver of type $\typeg$. This follows easily from the definitions because crystal graphs for $\typeg$ only have two colors.
\end{proof}

\begin{Proposition}\label{P:PlacticCrystals}
Let $\Lambda\in P^{+}$ be a dominant weight. Then $\crystalgraph[\Lambda]$ is plactic.
\end{Proposition}

\begin{proof}
Write $\Lambda=\sum_{i\in\vertices}l_{i}\Lambda_{i}$.
By Kashiwara's tensor product theorem \cite{Ka-crystal-bases}, $\crystalgraph[\Lambda]$ is a connected subgraph of the tensor product crystal $\bigotimes_{i\in\vertices}\crystalgraph[\Lambda_{i}]^{\otimes l_{i}}$. Therefore, by \autoref{L:PlacticTensor}, it is enough to show that the fundamental crystals are plactic. By \cite[Section 5]{BuSc-crystal-bases}, all of the fundamental crystals appear in tensor products of the fundamental crystals considered in \autoref{L:ConstructionPlactic}, so another application of \autoref{L:PlacticTensor} completes the proof.
\end{proof}

Applying \autoref{L:PlacticTensor} to \autoref{P:PlacticCrystals} we obtain:

\begin{Corollary}
The two crystal graphs $\Crystalgraph[\bfweight]$ and $\Crystalgraph$ are plactic.
\end{Corollary}

%%%%%%%%%%%%%%%%%%%%%%%%%%%%%%%%%%%%%%%%%

\subsection{Dots on idempotents}\label{SS:ConstructionDots}

%%%%%%%%%%%%%%%%%%%%%%%%%%%%%%%%%%%%%%%%%

Following the recipe of \autoref{SS:MainExampleTwo}, the next step is to put dots on the path idempotents $\idem[\bpath]$, for $\bpath\in\vertices^{n,\ell}$. This is necessary because
flanking two adjacent solid $i$-strings with
crossings annihilates them by
\autoref{Eq:RecollectionReidemeisterII}. The next definition allows us to avoid this problem.

Recall from \autoref{SS:RecollectionDegree} that $(a_{ij})_{i,j=1}^{e}$ is the Cartan matrix of $\quiver$.

\begin{Definition}\label{D:ConstructionRepeated}
Let $\bpath\in\vertices^{n,\ell}$ and $i\in\vertices$.
Write $\Res[\bpath]=(i_{1},\dots,i_{n})$. Fix $1\leq m\leq n$. The \emph{dotted $i_{m}$-subsequence} of $\Res[\bpath]$ is the maximal subsequence $(i_{p},i_{p+1},\dots,i_{m})$ such that $1\leq p\leq m$, $i_{p}$ is in the same component of $\bpath$ as $i_{m}$, and either $a_{i_{r}i_{m}}=0$, or $i_{r}=i_{m}$ and the $i_{r}$-string has the same anchor as the $i_{m}$-string in $\idem[\bpath]$, for $p\leq r<m$. Define $d_{m}(\bpath)=\#\set{p\leq r<m|i_{r}=i_{m}}$.
\end{Definition}

Note that every residue $i=i_{m}$ of $\bpath$ is contained in a unique dotted $i_{m}$-subsequence, so $d_{m}(\bpath)$ is defined for $1\leq m\leq n$. In particular, $d_{m}(\bpath)=0$ if $i_{m}$ appears only once in its dotted $i_{m}$-sequence. The condition $a_{i_{r}i_{m}}=0$ implies that close $i_{r}$ and $i_{m}$-strings satisfy an honest Reidemeister II relation. The condition on anchors implies that if $i_{r}=i_{m}$ then these strings are close in $\idem[\bpath]$.

\begin{Remark}\label{R:ConstructionRepeated}
In many examples we consider below, such as in \autoref{E:ConstructionExampleTwo}, all of the dotted $i$-subsequences are of the form $i\dots i$. For these examples, $d_{m}(\bpath)+1$ is just the length of this $i$-subsequence. The point of \autoref{D:ConstructionRepeated} is that, if $a_{ij}=0$ and $i_{m}=i$, then $d_{m}(\bpath)+1$ is equal to the number of times that $i$ appears in the sequences $ji$, $iji$, $ijiji$, \dots.
\end{Remark}

\begin{Notation}\label{N:ConstructionRepeated}
Let $\bpath\in\vertices^{n,\ell}$.
From now on, the \emph{$m$th solid string of $\dotidem[\bpath]$}
is the solid string that corresponds to the $m$th entry
of $\Res[\bpath]$.
% Its \textbf{component}, which is steady or unsteady, is the component of $\Res[\bpath]$ that contains this string.
Note, in particular, that this changes the convention that we use for the (linear combination of) diagrams $f(\bu)D$ as in \autoref{SS:RecollectionQPoly}, for $\bu\in\N[\bu]$. This convention will also be used when we associate permutation diagrams $D_{w}$ to permutations $w\in\sym$ in \autoref{SS:ConstructionPermutation}.
\end{Notation}

\begin{Example}\label{E:ConstructionExamplekthString}
In the second displayed equation of \autoref{E:ConstructionExampleTwo} below, the third solid string of the (dotted) idempotent diagram in type $\typeb[4]$ with residue sequence $12344321$ is the rightmost solid $3$-string.
\end{Example}

\begin{Definition}\label{D:ConstructionDots}
Let $\bpath\in\vertices^{n,\ell}$.
The \emph{dotted idempotent} $\dotidem[\bpath]$ is the diagram
$\dotidem[\bpath]=y_{1}^{d_{1}(\bpath)}\dots y_{n}^{d_{n}(\bpath)}\idem[\bpath]$.
\end{Definition}

\autoref{D:ConstructionDots} will remind some readers of the idempotents in the nil-Hecke algebra. As the following example shows, the dotted idempotent for a residue sequence containing a substring of the form $ii$ has only one dot on the second $i$-string, while the substring $iii$ obtains three dots, and $iiii$ gets six dots {\etc}

Locally, the dot placement from \autoref{D:ConstructionDots} takes the form:
\begin{gather*}
\colorbox{tomato!50}{i}\colon
\begin{tikzpicture}[anchorbase,smallnodes,yscale=0.5]
\draw[solid](0,0)node[below]{$i$}--++(0,1)node[above,yshift=-1pt]{\phantom{$1$}};
\end{tikzpicture}
,\quad
\colorbox{tomato!50}{ii}\colon
\begin{tikzpicture}[anchorbase,smallnodes,yscale=0.5]
\draw[solid](0,0)node[below]{$i$}--++(0,1)node[above,yshift=-1pt]{\phantom{$1$}};
\draw[solid,dot](0.5,0)node[below]{$i$}--++(0,1);
\end{tikzpicture}
,\quad
\colorbox{tomato!50}{iii}\colon
\begin{tikzpicture}[anchorbase,smallnodes,yscale=0.5]
\draw[solid](0,0)node[below]{$i$}--++(0,1)node[above,yshift=-1pt]{\phantom{$1$}};
\draw[solid,dot](0.5,0)node[below]{$i$}--++(0,1);
\draw[solid,dot=0.3,dot=0.7](1,0)node[below]{$i$}--++(0,1);
\end{tikzpicture}
,\quad
\colorbox{tomato!50}{iiii}\colon
\begin{tikzpicture}[anchorbase,smallnodes,yscale=0.5]
\draw[solid](0,0)node[below]{$i$}--++(0,1)node[above,yshift=-1pt]{\phantom{$1$}};
\draw[solid,dot](0.5,0)node[below]{$i$}--++(0,1);
\draw[solid,dot=0.3,dot=0.7](1,0)node[below]{$i$}--++(0,1);
\draw[solid,dot=0.16,dot=0.5,dot=0.84](1.5,0)node[below]{$i$}--++(0,1);
\end{tikzpicture}
,
\end{gather*}
where we assume that the illustrated number of
solid $i$-strings
is maximal.

\begin{Example}\label{E:ConstructionExampleTwo}
Continuing \autoref{E:ConstructionExampleOne}:
\begin{gather*}
\left.
\begin{gathered}
\typeb[4],
\\[-0.15cm]
\Res[\bpath]=123\colorbox{tomato!50}{44}321
\end{gathered}\right\}
\colon
\dotidem[\bpath]=
\begin{tikzpicture}[anchorbase,smallnodes,yscale=0.5]
\draw[ghost](0.3,0)--++(0,1)node[above,yshift=-1pt]{$1$};
\draw[ghost](1,0)--++(0,1)node[above,yshift=-1pt]{$1$};
\draw[ghost](1.4,0)--++(0,1)node[above,yshift=-1pt]{$2$};
\draw[ghost](1.9,0)--++(0,1)node[above,yshift=-1pt]{$2$};
\draw[ghost](2.5,0)--++(0,1)node[above,yshift=-1pt]{$3$};
\draw[ghost](2.8,0)--++(0,1)node[above,yshift=-1pt]{$3$};
\draw[solid](-0.7,0)node[below]{$1$}--++(0,1);
\draw[solid](0,0)node[below]{$1$}--++(0,1);
\draw[solid](0.4,0)node[below]{$2$}--++(0,1);
\draw[solid](0.9,0)node[below]{$2$}--++(0,1);
\draw[solid](1.5,0)node[below]{$3$}--++(0,1);
\draw[solid](1.8,0)node[below]{$3$}--++(0,1);
\draw[solid](2.6,0)node[below,xshift=-0.03cm]{$4$}--++(0,1);
\draw[solid,dot](2.7,0)node[below]{$4$}--++(0,1);
\draw[redstring](0.15,0)node[below]{$1$}--++(0,1);
\end{tikzpicture}
,\quad
\left.
\begin{gathered}
\typeg[2],
\\[-0.15cm]
\Res[\bpath]=21\colorbox{tomato!50}{22}12
\end{gathered}\right\}
\colon
\dotidem[\bpath]=
\begin{tikzpicture}[anchorbase,smallnodes,yscale=0.5]
\draw[ghost](-0.4,0)--++(0,1)node[above,yshift=-1pt]{$1$};
\draw[ghost](-0.1,0)--++(0,1)node[above,yshift=-1pt]{$1$};
\draw[solid](-1.4,0)node[below]{$1$}--++(0,1);
\draw[solid](-1.1,0)node[below]{$1$}--++(0,1);
\draw[solid](-0.5,0)node[below]{$2$}--++(0,1);
\draw[solid](-0.3,0)node[below,xshift=-0.03cm]{$2$}--++(0,1);
\draw[solid,dot](-0.2,0)node[below]{$2$}--++(0,1);
\draw[solid](0,0)node[below]{$2$}--++(0,1);
\draw[redstring](0.15,0)node[below]{$2$}--++(0,1);
\end{tikzpicture}
\end{gather*}
as the dotted idempotent diagrams. The other idempotents
in \autoref{E:ConstructionExampleOne} do not get extra dots.
For shorter residue sequences we have
\begin{gather*}
\left.
\begin{gathered}
\typeb[4],
\\[-0.15cm]
\Res[\bpath]=1234
\end{gathered}\right\}
\colon
\dotidem[\bpath]=
\begin{tikzpicture}[anchorbase,smallnodes,yscale=0.5]
\draw[ghost](1,0)--++(0,1)node[above,yshift=-1pt]{$1$};
\draw[ghost](1.9,0)--++(0,1)node[above,yshift=-1pt]{$2$};
\draw[ghost](2.8,0)--++(0,1)node[above,yshift=-1pt]{$3$};
\draw[solid](0,0)node[below]{$1$}--++(0,1);
\draw[solid](0.9,0)node[below]{$2$}--++(0,1);
\draw[solid](1.8,0)node[below]{$3$}--++(0,1);
\draw[solid](2.7,0)node[below]{$4$}--++(0,1);
\draw[redstring](0.15,0)node[below]{$1$}--++(0,1);
\end{tikzpicture}
,\quad
\left.
\begin{gathered}
\typeg[2],
\\[-0.15cm]
\Res[\bpath]=212
\end{gathered}\right\}
\colon
\dotidem[\bpath]=
\begin{tikzpicture}[anchorbase,smallnodes,yscale=0.5]
\draw[ghost](-0.1,0)--++(0,1)node[above,yshift=-1pt]{$1$};
\draw[solid](-1.1,0)node[below]{$1$}--++(0,1);
\draw[solid](-0.2,0)node[below]{$2$}--++(0,1);
\draw[solid](0,0)node[below]{$2$}--++(0,1);
\draw[redstring](0.15,0)node[below]{$2$}--++(0,1);
\end{tikzpicture}.
\end{gather*}
The extra dots only appear when, up to commuting strings, there are close strings of the same residue.
\end{Example}

%%%%%%%%%%%%%%%%%%%%%%%%%%%%%%%%%%%%%%%%%

\subsection{Sandwiched dots}\label{SS:ConstructionSandwichDots}

%%%%%%%%%%%%%%%%%%%%%%%%%%%%%%%%%%%%%%%%%

Recall from \autoref{D:ConstructionRepeated} that $d_{m}(\bpath)+1$ is the length of a dotted $i$-sequence in $\Res[\bpath]$. The following definition should be compared to
\autoref{L:RecollectionMovingStringsDots}.

\begin{Definition}\label{D:ConstructionSanddots}
Let $\bpath\in\vertices^{n,\ell}$ and $1\leq m\leq n$. Assume that the $m$th string of $\dotidem[\bpath]$ is an $i$-string. Set
$c_{m}(\bpath)=d_{i}-d_{m}(\bpath)-1$.
Define the sets of \emph{finite dots} and \emph{affine dots} to be
\begin{align*}
\Finch[\bpath]&=\set[\big]{\bfi=(f_{1},\dots,f_{n})\in\N^{n}|0\leq f_{m}\leq c_{m}(\bpath)},
\\
\Affch[\bpath]&=\set[\bigg]{\ba=(a_{1},\dots,a_{n})\in\N^{n}|\text{\begin{tabular}{l} $a_{m}=0$ whenever the $m$th string is not immediately\\ to the left of an affine string of the same residue\end{tabular}}}.
\end{align*}
The set of \emph{sandwiched dots} is $\Sandch[\bpath]=\Affch[\bpath]\cup\Finch[\bpath]$. We set $\sandwich{\bpath}{\ba}{\bfi}=\daffine\zeetwo\dotidem[\bpath]$
for $\bpath\in\vertices^{n,\ell}$, where $\ba\in\Affch[\bpath]$ and $\bfi\in\Finch[\bpath]$. When $\ba=(0,\dots,0)$ we write
$\fsandwich{\bpath}{\bfi}=\sandwich{\bpath}{(0,\dots,0)}{\bfi}$.
\end{Definition}

Again, these definitions are not as complicated as they appear because they can be checked (multi)locally in the diagrams. It follows from the proof of \autoref{C:ConstructionRepeated} below, that $c_{m}(\bpath)$ is always non-negative. When computing $c_{m}(\bpath)$ note that solid strings are only blocked by ghost strings, and ghost strings are only blocked by solid strings.

\begin{Example}\label{Ex:D5}
Let $\Gamma$ be the quiver of type $\typed[5]$ and consider the path $\Res[\bpath]=234532112$ in the fundamental crystal $\crystalgraph[{\fweight[2]}]$. Then:
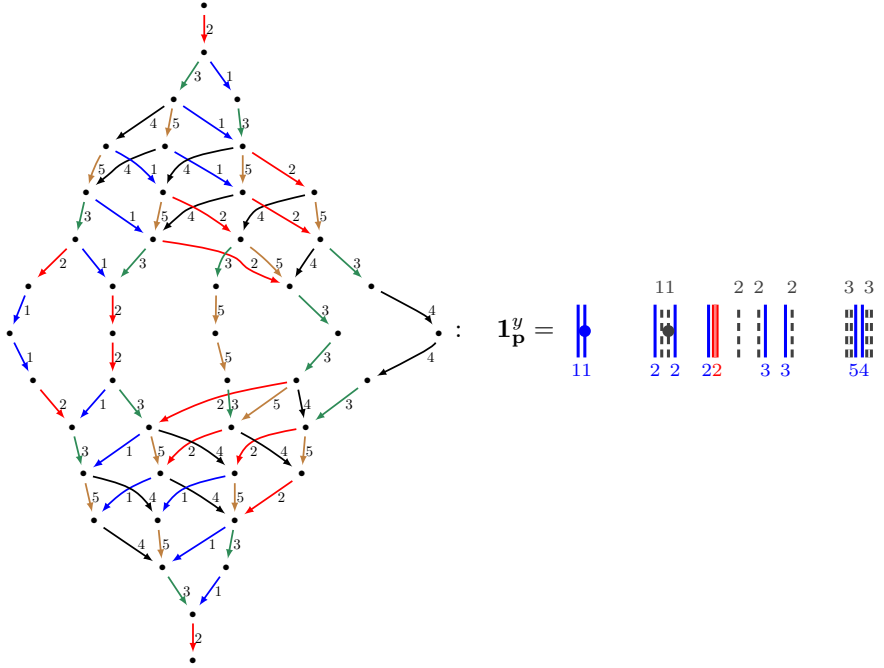
\begin{gather*}
\scalebox{0.55}{$\begin{tikzpicture}[anchorbase,>=latex,line join=bevel,scale=0.45,every path/.style={very thick}]
\node (node_0) at (485.5bp,434.5bp) [draw,draw=none] {$\bullet$};
\node (node_12) at (499.5bp,363.5bp) [draw,draw=none] {$\bullet$};
\node (node_22) at (387.5bp,363.5bp) [draw,draw=none] {$\bullet$};
\node (node_44) at (263.5bp,363.5bp) [draw,draw=none] {$\bullet$};
\node (node_1) at (379.5bp,150.5bp) [draw,draw=none] {$\bullet$};
\node (node_15) at (329.5bp,79.5bp) [draw,draw=none] {$\bullet$};
\node (node_2) at (152.5bp,648.5bp) [draw,draw=none] {$\bullet$};
\node (node_40) at (81.5bp,577.5bp) [draw,draw=none] {$\bullet$};
\node (node_41) at (208.5bp,577.5bp) [draw,draw=none] {$\bullet$};
\node (node_3) at (208.5bp,506.0bp) [draw,draw=none] {$\bullet$};
\node (node_27) at (208.5bp,434.5bp) [draw,draw=none] {$\bullet$};
\node (node_4) at (180.5bp,221.5bp) [draw,draw=none] {$\bullet$};
\node (node_38) at (283.5bp,150.5bp) [draw,draw=none] {$\bullet$};
\node (node_5) at (381.5bp,434.5bp) [draw,draw=none] {$\bullet$};
\node (node_6) at (512.5bp,719.5bp) [draw,draw=none] {$\bullet$};
\node (node_25) at (401.5bp,648.5bp) [draw,draw=none] {$\bullet$};
\node (node_30) at (521.5bp,648.5bp) [draw,draw=none] {$\bullet$};
\node (node_7) at (269.5bp,648.5bp) [draw,draw=none] {$\bullet$};
\node (node_21) at (475.5bp,577.5bp) [draw,draw=none] {$\bullet$};
\node (node_8) at (364.5bp,577.5bp) [draw,draw=none] {$\bullet$};
\node (node_33) at (363.5bp,506.0bp) [draw,draw=none] {$\bullet$};
\node (node_9) at (592.5bp,434.5bp) [draw,draw=none] {$\bullet$};
\node (node_10) at (598.5bp,577.5bp) [draw,draw=none] {$\bullet$};
\node (node_17) at (700bp,506.0bp) [draw,draw=none] {$\bullet$};
\node (node_11) at (280.5bp,292.5bp) [draw,draw=none] {$\bullet$};
\node (node_37) at (392.5bp,221.5bp) [draw,draw=none] {$\bullet$};
\node (node_24) at (392.5bp,292.5bp) [draw,draw=none] {$\bullet$};
\node (node_34) at (493.5bp,292.5bp) [draw,draw=none] {$\bullet$};
\node (node_13) at (287.5bp,790.5bp) [draw,draw=none] {$\bullet$};
\node (node_20) at (404.5bp,719.5bp) [draw,draw=none] {$\bullet$};
\node (node_42) at (168.5bp,719.5bp) [draw,draw=none] {$\bullet$};
\node (node_14) at (88.5bp,434.5bp) [draw,draw=none] {$\bullet$};
\node (node_16) at (147.5bp,363.5bp) [draw,draw=none] {$\bullet$};
\node (node_32) at (329.5bp,8.5bp) [draw,draw=none] {$\bullet$};
\node (node_28) at (164.5bp,292.5bp) [draw,draw=none] {$\bullet$};
\node (node_18) at (276.5bp,221.5bp) [draw,draw=none] {$\bullet$};
\node (node_19) at (346.5bp,1003.5bp) [draw,draw=none] {$\bullet$};
\node (node_26) at (346.5bp,932.5bp) [draw,draw=none] {$\bullet$};
\node (node_39) at (548.5bp,506.0bp) [draw,draw=none] {$\bullet$};
\node (node_23) at (396.5bp,861.5bp) [draw,draw=none] {$\bullet$};
\node (node_43) at (404.5bp,790.5bp) [draw,draw=none] {$\bullet$};
\node (node_35) at (300.5bp,861.5bp) [draw,draw=none] {$\bullet$};
\node (node_29) at (198.5bp,790.5bp) [draw,draw=none] {$\bullet$};
\node (node_31) at (284.5bp,719.5bp) [draw,draw=none] {$\bullet$};
\node (node_36) at (53.5bp,506.0bp) [draw,draw=none] {$\bullet$};
\draw [black,->] (node_0) ..controls (489.16bp,415.44bp) and (493.01bp,396.5bp)  .. (node_12);
\draw (502.0bp,399.0bp) node {$4$};
\draw [brown,->] (node_0) ..controls (458.23bp,414.3bp) and (426.89bp,392.23bp)  .. (node_22);
\draw (456.0bp,399.0bp) node {$5$};
\draw [red,->] (node_0) ..controls (414.59bp,421.76bp) and (386.81bp,415.76bp)  .. (362.5bp,408.0bp) .. controls (336.03bp,399.55bp) and (307.08bp,386.33bp)  .. (node_44);
\draw (371.0bp,399.0bp) node {$2$};
\draw [blue,->] (node_1) ..controls (366.12bp,131.03bp) and (351.59bp,110.98bp)  .. (node_15);
\draw (368.0bp,115.0bp) node {$1$};
\draw [red,->] (node_2) ..controls (133.18bp,628.72bp) and (111.68bp,607.83bp)  .. (node_40);
\draw (134.0bp,613.0bp) node {$2$};
\draw [blue,->] (node_2) ..controls (167.57bp,628.93bp) and (184.07bp,608.6bp)  .. (node_41);
\draw (196.0bp,613.0bp) node {$1$};
\draw [red,->] (node_3) ..controls (208.5bp,486.19bp) and (208.5bp,467.59bp)  .. (node_27);
\draw (217.0bp,470.0bp) node {$2$};
\draw [black,->] (node_4) ..controls (209.24bp,201.25bp) and (242.4bp,179.03bp)  .. (node_38);
\draw (252.0bp,186.0bp) node {$4$};
\draw [spinach,->] (node_5) ..controls (383.07bp,415.44bp) and (384.72bp,396.5bp)  .. (node_22);
\draw (393.0bp,399.0bp) node {$3$};
\draw [black,->] (node_6) ..controls (454.73bp,706.23bp) and (431.61bp,699.71bp)  .. (423.5bp,693.0bp) .. controls (415.34bp,686.26bp) and (409.9bp,675.72bp)  .. (node_25);
\draw (432.0bp,684.0bp) node {$4$};
\draw [brown,->] (node_6) ..controls (514.86bp,700.44bp) and (517.33bp,681.5bp)  .. (node_30);
\draw (526.0bp,684.0bp) node {$5$};
\draw [red,->] (node_7) ..controls (352.05bp,636.48bp) and (373.25bp,630.73bp)  .. (391.5bp,622.0bp) .. controls (402.9bp,616.55bp) and (402.65bp,610.47bp)  .. (413.5bp,604.0bp) .. controls (422.72bp,598.5bp) and (433.35bp,593.65bp)  .. (node_21);
\draw (422.0bp,613.0bp) node {$2$};
\draw [spinach,->] (node_7) ..controls (252.99bp,628.82bp) and (234.77bp,608.21bp)  .. (node_41);
\draw (254.0bp,613.0bp) node {$3$};
\draw [brown,->] (node_8) ..controls (364.24bp,558.42bp) and (363.97bp,539.63bp)  .. (node_33);
\draw (372.0bp,542.0bp) node {$5$};
\draw [spinach,->] (node_9) ..controls (566.69bp,414.35bp) and (537.14bp,392.43bp)  .. (node_12);
\draw (565.0bp,399.0bp) node {$3$};
\draw [black,,->] (node_10) ..controls (636.07bp,557.16bp) and (680.41bp,534.33bp)  .. (node_17);
\draw (690.0bp,542.0bp) node {$4$};
\draw [blue,->] (node_11) ..controls (249.29bp,279.5bp) and (236.18bp,273.32bp)  .. (225.5bp,266.0bp) .. controls (213.7bp,257.91bp) and (202.13bp,246.63bp)  .. (node_4);
\draw (234.0bp,257.0bp) node {$1$};
\draw [black,->] (node_11) ..controls (305.97bp,279.11bp) and (318.66bp,272.53bp)  .. (329.5bp,266.0bp) .. controls (344.89bp,256.73bp) and (361.6bp,245.12bp)  .. (node_37);
\draw (364.0bp,257.0bp) node {$4$};
\draw [red,->] (node_12) ..controls (438.09bp,352.55bp) and (423.58bp,346.57bp)  .. (412.5bp,337.0bp) .. controls (404.72bp,330.28bp) and (399.79bp,319.9bp)  .. (node_24);
\draw (421.0bp,328.0bp) node {$2$};
\draw [brown,->] (node_12) ..controls (497.93bp,344.44bp) and (496.28bp,325.5bp)  .. (node_34);
\draw (505.0bp,328.0bp) node {$5$};
\draw [blue,->] (node_13) ..controls (310.27bp,776.78bp) and (322.12bp,770.08bp)  .. (332.5bp,764.0bp) .. controls (349.77bp,753.89bp) and (369.17bp,742.15bp)  .. (node_20);
\draw (370.0bp,755.0bp) node {$1$};
\draw [black,->] (node_13) ..controls (250.62bp,777.62bp) and (235.26bp,771.47bp)  .. (222.5bp,764.0bp) .. controls (211.6bp,757.62bp) and (210.33bp,753.94bp)  .. (200.5bp,746.0bp) .. controls (195.8bp,742.21bp) and (190.7bp,738.12bp)  .. (node_42);
\draw (231.0bp,755.0bp) node {$4$};
\draw [red,->] (node_14) ..controls (104.38bp,414.93bp) and (121.76bp,394.6bp)  .. (node_16);
\draw (134.0bp,399.0bp) node {$2$};
\draw [red,->] (node_15) ..controls (329.5bp,60.442bp) and (329.5bp,41.496bp)  .. (node_32);
\draw (338.0bp,44.0bp) node {$2$};
\draw [spinach,->] (node_16) ..controls (151.95bp,344.44bp) and (156.62bp,325.5bp)  .. (node_28);
\draw (167.0bp,328.0bp) node {$3$};
\draw [black,->] (node_17) ..controls (692.38bp,484.73bp) and (645.64bp,461.69bp)  .. (node_9);
\draw (688.0bp,470.0bp) node {$4$};
\draw [brown,->] (node_18) ..controls (278.33bp,202.44bp) and (280.25bp,183.5bp)  .. (node_38);
\draw (289.0bp,186.0bp) node {$5$};
\draw [red,->] (node_19) ..controls (346.5bp,984.44bp) and (346.5bp,965.5bp)  .. (node_26);
\draw (355.0bp,968.0bp) node {$2$};
\draw [black,->] (node_20) ..controls (350.81bp,706.98bp) and (333.2bp,701.17bp)  .. (318.5bp,693.0bp) .. controls (314.71bp,690.9bp) and (298.07bp,675.88bp)  .. (node_7);
\draw (327.0bp,684.0bp) node {$4$};
\draw [red,->] (node_20) ..controls (437.32bp,699.14bp) and (475.49bp,676.64bp)  .. (node_30);
\draw (485.0bp,684.0bp) node {$2$};
\draw [spinach,->] (node_21) ..controls (495.24bp,557.7bp) and (517.0bp,536.99bp)  .. (node_39);
\draw (529.0bp,542.0bp) node {$3$};
\draw [red,->] (node_22) ..controls (345.63bp,350.84bp) and (331.16bp,344.97bp)  .. (319.5bp,337.0bp) .. controls (308.7bp,329.62bp) and (298.88bp,318.63bp)  .. (node_11);
\draw (328.0bp,328.0bp) node {$2$};
\draw [black,->] (node_22) ..controls (411.42bp,350.05bp) and (423.33bp,343.47bp)  .. (433.5bp,337.0bp) .. controls (448.03bp,327.75bp) and (463.77bp,316.28bp)  .. (node_34);
\draw (467.0bp,328.0bp) node {$4$};
\draw [spinach,->] (node_23) ..controls (398.59bp,842.44bp) and (400.79bp,823.5bp)  .. (node_43);
\draw (409.0bp,826.0bp) node {$3$};
\draw [blue,->] (node_24) ..controls (335.24bp,280.15bp) and (320.44bp,274.45bp)  .. (308.5bp,266.0bp) .. controls (298.72bp,259.08bp) and (290.73bp,248.23bp)  .. (node_18);
\draw (317.0bp,257.0bp) node {$1$};
\draw [brown,->] (node_24) ..controls (392.5bp,273.44bp) and (392.5bp,254.5bp)  .. (node_37);
\draw (401.0bp,257.0bp) node {$5$};
\draw [spinach,->] (node_25) ..controls (385.14bp,635.34bp) and (378.54bp,629.05bp)  .. (374.5bp,622.0bp) .. controls (370.0bp,614.14bp) and (367.52bp,604.39bp)  .. (node_8);
\draw (383.0bp,613.0bp) node {$3$};
\draw [brown,->] (node_25) ..controls (418.87bp,635.0bp) and (427.42bp,628.42bp)  .. (434.5bp,622.0bp) .. controls (444.19bp,613.22bp) and (454.27bp,602.54bp)  .. (node_21);
\draw (460.0bp,613.0bp) node {$5$};
\draw [blue,->] (node_26) ..controls (359.88bp,913.03bp) and (374.41bp,892.98bp)  .. (node_23);
\draw (385.0bp,897.0bp) node {$1$};
\draw [spinach,->] (node_26) ..controls (334.19bp,913.03bp) and (320.82bp,892.98bp)  .. (node_35);
\draw (337.0bp,897.0bp) node {$3$};
\draw [blue,->] (node_27) ..controls (191.99bp,414.82bp) and (173.77bp,394.21bp)  .. (node_16);
\draw (194.0bp,399.0bp) node {$1$};
\draw [spinach,->] (node_27) ..controls (223.3bp,414.93bp) and (239.51bp,394.6bp)  .. (node_44);
\draw (251.0bp,399.0bp) node {$3$};
\draw [brown,->] (node_28) ..controls (168.69bp,273.44bp) and (173.08bp,254.5bp)  .. (node_4);
\draw (183.0bp,257.0bp) node {$5$};
\draw [black,->] (node_28) ..controls (215.8bp,280.09bp) and (230.6bp,274.34bp)  .. (242.5bp,266.0bp) .. controls (252.63bp,258.89bp) and (261.2bp,247.89bp)  .. (node_18);
\draw (269.0bp,257.0bp) node {$4$};
\draw [blue,->] (node_29) ..controls (223.56bp,777.27bp) and (234.6bp,770.98bp)  .. (243.5bp,764.0bp) .. controls (254.02bp,755.75bp) and (264.33bp,744.78bp)  .. (node_31);
\draw (271.0bp,755.0bp) node {$1$};
\draw [brown,->] (node_29) ..controls (190.24bp,776.74bp) and (186.34bp,770.16bp)  .. (183.5bp,764.0bp) .. controls (179.64bp,755.62bp) and (176.21bp,745.95bp)  .. (node_42);
\draw (192.0bp,755.0bp) node {$5$};
\draw [spinach,->] (node_30) ..controls (542.57bp,628.62bp) and (566.2bp,607.44bp)  .. (node_10);
\draw (577.0bp,613.0bp) node {$3$};
\draw [black,->] (node_30) ..controls (511.53bp,634.61bp) and (506.6bp,628.02bp)  .. (502.5bp,622.0bp) .. controls (496.46bp,613.13bp) and (490.11bp,602.95bp)  .. (node_21);
\draw (511.0bp,613.0bp) node {$4$};
\draw [brown,->] (node_31) ..controls (280.57bp,700.44bp) and (276.46bp,681.5bp)  .. (node_7);
\draw (287.0bp,684.0bp) node {$5$};
\draw [red,->] (node_31) ..controls (313.86bp,706.2bp) and (327.74bp,699.78bp)  .. (339.5bp,693.0bp) .. controls (354.93bp,684.1bp) and (371.39bp,672.44bp)  .. (node_25);
\draw (375.0bp,684.0bp) node {$2$};
\draw [brown,->] (node_33) ..controls (368.4bp,486.09bp) and (373.29bp,467.21bp)  .. (node_5);
\draw (383.0bp,470.0bp) node {$5$};
\draw [red,->] (node_34) ..controls (465.39bp,272.3bp) and (433.09bp,250.23bp)  .. (node_37);
\draw (463.0bp,257.0bp) node {$2$};
\draw [brown,->] (node_35) ..controls (297.1bp,842.44bp) and (293.53bp,823.5bp)  .. (node_13);
\draw (304.0bp,826.0bp) node {$5$};
\draw [black,->] (node_35) ..controls (272.12bp,841.3bp) and (239.49bp,819.23bp)  .. (node_29);
\draw (270.0bp,826.0bp) node {$4$};
\draw [blue,->] (node_35) ..controls (329.52bp,841.25bp) and (363.0bp,819.03bp)  .. (node_43);
\draw (373.0bp,826.0bp) node {$1$};
\draw [blue,->] (node_36) ..controls (63.131bp,485.88bp) and (72.907bp,466.46bp)  .. (node_14);
\draw (84.0bp,470.0bp) node {$1$};
\draw [spinach,->] (node_37) ..controls (389.1bp,202.44bp) and (385.53bp,183.5bp)  .. (node_1);
\draw (396.0bp,186.0bp) node {$3$};
\draw [blue,->] (node_37) ..controls (362.0bp,201.2bp) and (326.68bp,178.83bp)  .. (node_38);
\draw (359.0bp,186.0bp) node {$1$};
\draw [spinach,->] (node_38) ..controls (295.81bp,131.03bp) and (309.18bp,110.98bp)  .. (node_15);
\draw (320.0bp,115.0bp) node {$3$};
\draw [spinach,->] (node_39) ..controls (530.78bp,485.45bp) and (512.19bp,464.94bp)  .. (node_0);
\draw (532.0bp,470.0bp) node {$3$};
\draw [blue,->] (node_40) ..controls (74.176bp,558.32bp) and (66.497bp,539.26bp)  .. (node_36);
\draw (80.0bp,542.0bp) node {$1$};
\draw [red,->] (node_41) ..controls (208.5bp,558.42bp) and (208.5bp,539.63bp)  .. (node_3);
\draw (217.0bp,542.0bp) node {$2$};
\draw [spinach,->] (node_42) ..controls (164.31bp,700.44bp) and (159.92bp,681.5bp)  .. (node_2);
\draw (171.0bp,684.0bp) node {$3$};
\draw [blue,->] (node_42) ..controls (196.61bp,699.3bp) and (228.91bp,677.23bp)  .. (node_7);
\draw (239.0bp,684.0bp) node {$1$};
\draw [red,->] (node_43) ..controls (434.72bp,770.2bp) and (469.72bp,747.83bp)  .. (node_6);
\draw (479.0bp,755.0bp) node {$2$};
\draw [brown,->] (node_43) ..controls (404.5bp,771.44bp) and (404.5bp,752.5bp)  .. (node_20);
\draw (413.0bp,755.0bp) node {$5$};
\draw [black,->] (node_43) ..controls (336.98bp,778.19bp) and (322.9bp,772.52bp)  .. (311.5bp,764.0bp) .. controls (302.5bp,757.28bp) and (295.8bp,746.58bp)  .. (node_31);
\draw (320.0bp,755.0bp) node {$4$};
\draw [brown,->] (node_44) ..controls (267.95bp,344.44bp) and (272.62bp,325.5bp)  .. (node_11);
\draw (282.0bp,328.0bp) node {$5$};
\draw [black,->] (node_44) ..controls (305.71bp,350.68bp) and (322.47bp,344.63bp)  .. (336.5bp,337.0bp) .. controls (351.18bp,329.02bp) and (366.03bp,317.27bp)  .. (node_24);
\draw (370.0bp,328.0bp) node {$4$};
\draw [blue,->] (node_44) ..controls (235.95bp,343.3bp) and (204.29bp,321.23bp)  .. (node_28);
\draw (234.0bp,328.0bp) node {$1$};
\end{tikzpicture}$}
:\quad
\dotidem[\bpath]=
\begin{tikzpicture}[anchorbase,smallnodes,xscale=1.1,yscale=0.7]
\draw[redstring](0,0)node[below]{$\,2$}--+(0,1);
\draw[solid](-0.08,0)node[below]{$2\,$}--+(0,1);
\draw[ghost](0.92,0)--+(0,1)node[above]{$2$};
\draw[solid](0.84,0)node[below]{$3$}--+(0,1);
\draw[ghost,double](1.84,0)--+(0,1)node[above]{$3$};
\draw[solid](1.76,0)node[below]{$\,4$}--+(0,1);
\draw[solid](1.68,0)node[below]{$5\,$}--+(0,1);
\draw[solid](0.60,0)node[below]{$3$}--+(0,1);
\draw[ghost,double](1.60,0)--+(0,1)node[above]{$3$};
\draw[solid](-0.48,0)node[below]{$2$}--+(0,1);
\draw[ghost](0.52,0)--+(0,1)node[above]{$2$};
\draw[solid,dot=0.5](-1.56,0)node[below]{$\,1$}--+(0,1);
\draw[ghost,dot=0.5](-0.56,0)--+(0,1)node[above]{$\,1$};
\draw[solid](-1.64,0)node[below]{$1\,$}--+(0,1);
\draw[ghost](-0.64,0)--+(0,1)node[above]{$1\,$};
\draw[solid](-0.72,0)node[below]{$2$}--+(0,1);
\draw[ghost](0.28,0)--+(0,1)node[above]{$2$};
\end{tikzpicture}
.
\end{gather*}
In this case, $\Finch[\bpath]=\set{(0,0,0,0,0,0,0,0,0),(0,0,0,0,0,0,1,0,0),(0,0,0,0,0,0,0,0,1),(0,0,0,0,0,0,1,0,1)}$ and $\Affch[\bpath]=\set{(0,0,0,0,0,0,0,0,0)}$.
\end{Example}

\begin{Definition}\label{D:ConstructionSandwichAlgebra}
In the setting of \autoref{D:ConstructionSanddots}, the
\emph{(affine) sandwiched algebra} $\sand[\bpath]$ at $\bpath\in\vertices^{n,\ell}$ is
the algebra $\WA(X)$ generated by
$\asandbasis[\bpath]=\set[\big]{\sandwich{\bpath}{\ba}{\bfi}
|\ba\in\Affch[\bpath],\bfi\in\Finch[\bpath]}$.

The \emph{(finite) sandwiched algebra} $\sand[\bpath]$ at $\bpath\in\vertices^{n,\ell}$ is
the algebra generated by
$\sandbasis[\bpath]=\set[\big]{\fsandwich{\bpath}{\bfi}|\bfi\in\Finch[\bpath]}$ for $\bpath\in\vertices^{n,\ell}$.
\end{Definition}

As we will see, the sandwich bases of $\aWA$ and $\WAc(X)$ that we construct in \autoref{T:ConstructionMain} below factor through the sandwiched algebras. This will imply that the diagrams in $\asandbasis[\bpath]$ and $\sandbasis[\bpath]$ are linearly independent. Later, we will fix preferred paths $\ppath$ for a vertex $\blam$ in $\crystalgraph$ we will write $\Res[\blam]=\Res[\ppath]$, $\dotidem[\blam]=\dotidem[\ppath]$, $\sand[\blam]=\sand[\ppath]$ etc.
We give some examples.

\begin{Example}\label{E:ConstructionSandwichpart}
The diagrams $\dotidem[\blam]$, $\dotidem[\mu]$ and $\dotidem[\nu]$ in \autoref{E:MainExampleTheBeastItself} have five strings, all of which
are in the steady part of these diagrams. We have $c_{1}(\nu)=c_{2}(\nu)=c_{4}(\nu)=0$,
while $c_{3}(\nu)=c_{5}(\nu)=1$, which we will see implies that the sandwiched algebra $\sand[\nu]$ is nontrivial, agreeing with \autoref{E:MainExampleTheBeastItself}(c). The sandwich
algebras $\sand[\lambda]$ and $\sand[\mu]$ are both trivial.
\end{Example}

\begin{Example}\label{E:ConstructionSandwichpartTwo}
Returning to \autoref{E:ConstructionExampleTwo},
for type $\typeb[4]$ and $\Res[\bpath]=12344321$ we have
$c_{k}(\bpath)=0$ for all $k\in\set{1,\dots,8}$, so $\sand[\bpath]=\ring$ is trivial. In contrast, again in type $\typeb[4]$, but with $\Res[\bpath]=1234$, we have
$c_{k}(\bpath)=0$, for $k\in\set{1,2,3}$, and
$c_{4}(\bpath)=1$, so $\sand[\bpath]=\ring[y_4]/(y_{4}^{2})$. Similarly, in the type $\typeg[2]$ example,
we have a trivial sandwich algebra for $\Res[\bpath]=212212$. On the other hand, $c_{3}(\bpath)=2$ for $\Res[\bpath]=212$, since $d_2=3$ because (the partner relation of) \autoref{L:RecollectionReidemeisterIII} applies,
so $\sand[\bpath]=\ring[y_{3}]/(y_{3}^{2})$.
\end{Example}

\begin{Example}\label{E:ConstructionSandwichAlgebra}
Continuing \autoref{E:ConstructionExampleAffine}, the sandwiched algebra
is spanned by
\begin{gather*}
\sandbasis[\blam]=
\set[\bigg]{
\begin{tikzpicture}[anchorbase,smallnodes,rounded corners]
\draw[solid](-0.2,0.5)--++(0,-0.5) node[below]{$1$};
\draw[solid,dot](1.8,0.5)--++(0,-0.5) node[below]{$1$};
\draw[redstring] (0,0)node[below]{$1$} to (0,0.5)node[above,yshift=-1pt]{$\phantom{i}$};
\draw[affine] (2,0)node[below]{$1$} to (2,0.5);
\draw[affine] (4,0)node[below]{$1$} to (4,0.5);
\draw[affine] (6,0)node[below]{$2$} to (6,0.5);
\draw[affine] (8,0)node[below]{$2$} to (8,0.5);
\node at (1.65,0.25) {$a$};
\end{tikzpicture}
|a\in\N}
\end{gather*}
where $a$ is the number of (affine) dots. Hence, $\sand[\blam]\cong\ring[y]$ is a polynomial ring.
\end{Example}

%%%%%%%%%%%%%%%%%%%%%%%%%%%%%%%%%%%%%%%%%

\subsection{Preferred paths and permutation diagrams}\label{SS:ConstructionPermutation}

%%%%%%%%%%%%%%%%%%%%%%%%%%%%%%%%%%%%%%%%%

We are now ready to define our sandwich cellular bases. We start with defining the sets that label the cell modules in our basis.

Recall from \autoref{D:ConstructionIdempotents} that each path $\bpath\in\Parts{n}{\abfweight}$, via its residue sequence $\Res[\bpath]$, determines an idempotent diagram $\idem[\bpath]$.
Recall from \autoref{D:ConstructionIdempotents} that $\bx_{\bpath}=(x_{1}<\dots<x_{n})$ records the coordinates of the solid strings in $\idem[\bpath]$, for $\bpath\in\Parts{n}{\abfweight}$.

%We now define the total order for our sandwich cellular bases.

%\begin{Definition}\label{D:ConstructionChoice}
%Let $\bpath,\cpath\in\Parts{n}{\abfweight}$.
%Then $\bpath$ \emph{dominates} $\cpath$, written $\bpath\gsandorder\cpath$ %if
%$\bx_{\bpath}\geq\bx_{\cpath}$, where $\geq$ is the lexicographic order on %$\R^{n}$ where coordinates are read right to left.
%\end{Definition}

Let $w\in\sym$. As in \cite[Definition 3B.1]{MaTu-klrw-algebras}, given $\bx,\by\in X$ and $\bi\in\vertices^{n}$ let $D_{\bx,\by,\bi}(w)$ be the \emph{permutation diagram} with residue sequence $\bi$ and with a minimal number of crossings such that has a solid $i_{k}$-string connecting $(x_{k},1)$ to $(y_{w(k)},0)$, for $1\leq k\leq n$, together with the corresponding ghost and red strings. The diagram $D_{\bx,\by,\bi}(w)$ depends on the choice of crossings but by \cite[Lemma 3B.3]{MaTu-klrw-algebras} different choices give diagrams that are equal modulo ``bigger'' terms. To simplify notation, we often write $D_{w}$ instead of $D_{\bx,\by,\bi}(w)$ when $\bx$, $\by$ and $\bi$ are clear.

\begin{Definition}\label{D:ConstructionPermutationDiagrams}
Let $\bpath\in\Parts{n}{\abfweight}$. Let $E_{\bpath}$ be the unique diagram such that $E_{\bpath}=\idem[\bpath]E_{\bpath}\idem[{\bx,\Res[\bpath]}]$ and the $r$th string in $\idem[{\bx,\Res[\bpath]}]$ is joined to the $r$th string in $\idem[\bpath]$, for $1\le r\le n$. As paths in $\crystalgraph[\abfweight]$ are uniquely determined by the residue sequences, set $E_{\Res[\bpath]}=E_{\bpath}$.
\end{Definition}

\begin{Example}
    Take $\bpath=\ppath$ in \autoref{E:MainExampleTheBeastItself}. Then there are three paths $\cpath$ in the crystal graph $\crystalgraph[\Lambda_{2}]$ with sink $\lambda$, which have residue sequences $21332$, $23132$ and $23312$.
    \begin{gather*}
        E_{21332} =
        \begin{tikzpicture}[anchorbase,smallnodes,xscale=1.5,yscale=0.8]
            \foreach \r/\x/\c in {2/-0.1/1,1/-1.2/2,3/0.8/3,3/0.7/4,2/-0.4/5} {
               \draw[solid](\x,2)--+(0,-0.5)--({-\c*1.2},0)node[below]{$\r$};
               \ifnum\r<3
                 \draw[ghost]({\x+1},2)node[above]{$\r$}--+(0,-0.4)--({1-\c*1.2},0);
               \fi
            }
            \draw[redstring](0,0)node[below]{$2$}--++(0,2);
        \end{tikzpicture},
        \\
        E_{23132} =
        \begin{tikzpicture}[anchorbase,smallnodes,xscale=1.5,yscale=0.8]
            \foreach \r/\x/\c in {2/-0.1/1,1/-1.2/3,3/0.8/2,3/0.7/4,2/-0.4/5} {
               \draw[solid](\x,2)--+(0,-0.5)--({-1.2*\c},0)node[below]{$\r$};
               \ifnum\r<3
                 \draw[ghost]({\x+1},2)node[above]{$\r$}--+(0,-0.4)--({1-\c*1.2},0);
               \fi
            }
            \draw[redstring](0,0)node[below]{$2$}--++(0,2);
        \end{tikzpicture},
        \\
        E_{23312} =
        \begin{tikzpicture}[anchorbase,smallnodes,xscale=1.5,yscale=0.8]
            \foreach \r/\x/\c in {2/-0.1/1,1/-1.2/4,3/0.8/3,3/0.7/2,2/-0.4/5} {
               \draw[solid](\x,2)--+(0,-0.5)--({-\c*1.2},0)node[below]{$\r$};
               \ifnum\r<3
                 \draw[ghost]({\x+1},2)node[above]{$\r$}--+(0,-0.4)--({1-\c*1.2},0);
               \fi
            }
            \draw[redstring](0,0)node[below]{$2$}--++(0,2);
        \end{tikzpicture}.
    \end{gather*}
\end{Example}

For $\bpath,\cpath\in\Parts{n}{\abfweight}$ let $\idem[\bpath]$ and $\idem[\cpath]$ be the idempotent diagrams constructed in \autoref{D:ConstructionIdempotents}. If $w\in\sym$ is a permutation such that $\Res[\bpath]=w\,\Res[\cpath]$, then let $D_{\bpath,\cpath}(w)=D_{\bx_{\bpath},\bx_{\cpath}, \Res[\bpath]}(w)$ be the permutation diagram that connects the $k$th solid string of $\idem[\bpath]$ to the $w(k)$th solid string of $\idem[\cpath]$, for $1\leq k\leq n$. Then $D_{\bpath,\cpath}(w)=\idem[\bpath]D_{\bpath,\cpath}(w)\idem[\cpath]$. Note that the coordinates of the strings at the top and bottom of $D_{\bpath,\cpath}(w)$ are given by the parking functions for $\bpath$ and $\cpath$, respectively.

Let $\bi,\bj\in\vertices^{n}$ and write $\bi=(i_{1},\dots,i_{n})$. Let $\sim$ be the
equivalence relation on $\vertices^{n}$ that is the transitive closure of the relation $\bi\sim\bj$ if $\bi=(k,k+1)\,\bj$ for some $k$ such that $a_{i_{k},i_{k+1}}=0$. Write $\bpath\sim\cpath$ if $\bpath,\cpath\in\Parts{n}{\abfweight}$ and $\Res[\bpath]\sim\Res[\cpath]$.

\begin{Lemma}\label{L:ConstructionPermBasisAdSquares}
Suppose that $\bpath,\cpath\in\Parts{n}{\abfweight}$ and that $w\in\sym$ is a path permutation from $\bpath$ to $\cpath$ such that $\Res[\bpath]\sim w\,\Res[\cpath]$.
Then $D_{w^{-1}}D_{w}\idem[\bpath]=\idem[\bpath]=\idem[\bpath]D_{w^{-1}}D_{w}$. Consequently, if  $\omega(\bpath)=\omega(\cpath)$ then $\idem[\bpath]=\idem[\cpath]$.
\end{Lemma}

\begin{proof}
By \autoref{P:FiniteTypesPermutations}, $w$ is a face permutation, so it is enough to consider the case when $w$ is a partial face permutation. By assumption, if $i$ and $j$ are the two colors appearing in $w$, then $a_{ij}=0$. In particular, if all of the strings in $\idem[\bpath]$ corresponding to the edges permuted by $w$ are in the same component, then the result follows by the honest Reidemeister II relations,
as in \autoref{Eq:RecollectionReidemeisterII}. If the strings are in different components, then the result still holds because the number of $i$ and $j$'s appearing the two arms of the face is the same, since $\Res[\bpath]\sim w\,\Res[\cpath]$. Hence, by \autoref{D:ConstructionIdempotents}, the edges allowed by the tensor product rule \autoref{E:TensorProduct} occur between the same vertices, the corresponding strings are placed in the same components. Therefore, the strings in $\idem[\bpath]$ and $\idem[\cpath]$ again differ by honest Reidemeister II relations. (In particular, this shows that the idempotent diagrams $\idem[\bpath]$ and $\idem[\cpath]$ coincide unless $i$ and $j$ are nonadjacent fishtail vertices.) Hence, $D_{w^{-1}}D_{w}\idem[\bpath]=\idem[\bpath]=\idem[\bpath]D_{w^{-1}}D_{w}$.

For the second statement, if $\omega(\bpath)=\omega(\cpath)$ then $\Res[\bpath]=\Res[\cpath]$. Therefore, the construction of the idempotent diagrams from \autoref{D:ConstructionIdempotents} ensures that $\idem[\cpath]=D_{w}\idem[\bpath]D_{w^{-1}}$, so the result follows from the first claim.
\end{proof}

% Recall from \autoref{D:FiniteTypesSSPerm} that if $\bpath$ and $\cpath$ are paths in $\crystalgraph[{\fweight[k]}]$, then $\Face[k](\bpath,\cpath)$ is the set of face permutations between these paths.
%
% \begin{Definition}\label{D:ConstructionPermutationDiagrams}
% Let $\bpath=(\bpath_{1},\dots,\bpath_{\hell})$ and $\cpath=(\cpath_{1},\dots,\cpath_{\hell})$ be paths in $\Crystalgraph[\abfweight]$. A \emph{face permutation} from $\bpath$ to $\cpath$ is a composition of face permutations in the disjoint union $\bigsqcup_{k=1}^{\hell}\Face[\affine{\rho}_{k}](\bpath_{k},\cpath_{k})$ that do not involve basic face permutations around nonadjacent squares.
% Let $\hFace(\bpath,\cpath)$ be the set of face permutations from $\bpath$ to $\cpath$. If $\bpath,\cpath\in\Crystalgraph[\bfweight]$ set $\Face(\bpath,\cpath)=\hFace(\bpath,\cpath)$.
% \end{Definition}

We can now define the indexing sets for our cellular bases.

\begin{Definition}\label{D:ConstructionChoice}
\leavevmode
\begin{enumerate}

%\item  Let $j\in\vertices$ and let $\bpath$ be a rooted path in $\crystalgraph[{\fweight[i]}]$. Then $j$ is an \emph{good residue} of $\bpath$ if $\bpath\xrightarrow{j}\lambda$ is a rooted path in $\crystalgraph[{\fweight[i]}]$, for some $\lambda$.

%\item Let $\bpath\in\Parts{n}{\abfweight}$ have residue sequence $\Res[\bpath]=(i_{1},\dots,i_{n})$. The path $\bpath$ is \emph{good} if either $\bpath$ is the empty path, or if $m$ is maximal such that $\bpath_{m}$ is nonempty, then either $1\leq m\leq\ell$ or $\ell<m\leq\hell$ and $i_{n}$ is not a good residue of $\bpath_{k}$ for $1\leq k<m$.

%\item Two %good
%paths $\bpath,\cpath\in\crystalgraph[\abfweight]$ with the same sink are \emph{equivalent} if $\Res[\bpath]=w\,\Res[\cpath]$, where $w\in\sym$ is a composition of basic face permutations for nonadjacent squares. Let $\fParts{n}{\fweight[]}$ and $\fParts{n}{\abfweight}$ be sets of equivalence class representatives for the %good
%paths in $\crystalgraph[\bfweight]$ and $\crystalgraph[\abfweight]$, respectively, chosen so that $\fParts{n}{\bfweight}$ is a subset of $\fParts{n}{\abfweight}$ (under the map $\tau$ of \autoref{L:TensorInjection}).

\item Let $\CrystalVertices=\set{\omega(\bpath)|\bpath\in\Parts{n}{\bfweight}}$ and $\aCrystalVertices=\set{\omega(\bpath)|\bpath\in\Parts{n}{\abfweight}}$ be the sets of sinks of the %good
paths in $\Parts{n}{\bfweight}$ and $\Parts{n}{\abfweight}$, respectively. We think of $\CrystalVertices$ as being contained in $\aCrystalVertices$ via \autoref{L:TensorInjection}.

\item The \emph{(diagrammatic) dominance order}, on the set of idempotent diagrams is the order $\sandorder$ where $\idem[(x,\bi)]$ and $\idem[(y,\bj)]$ if there exists a $k$ such that $x_k<y_k$ and $x_l=y_l$ for $k< l\leq n$. We use $\lsandorder$, $\gsandorder$ {\etc} with the evident meaning. If $\bpath,\cpath\in\Parts{n}{\abfweight}$ write $\bpath\lsandorder\cpath$ if $\idem[\bpath]\lsandorder\idem[\cpath]$.

\item  For each $\blam\in\aCrystalVertices$, let $\ppath\in\Parts{n}{\abfweight}$ be the unique path in $\Parts{n}{\abfweight}$ such that $\omega(\ppath)=\blam$ and $\ppath\lsandorder\cpath$ whenever $\omega(\cpath)=\blam$. Set $\idem[\blam]=\idem[\ppath]$ and $\dotidem[\blam]=\dotidem[\ppath]$.

\item If $\blam,\bmu\in\aCrystalVertices$ write $\blam\lsandorder\bmu$ if $\ppath\lsandorder\ppath[\bmu]$. Write $\blam\sandorder\bmu$ if $\blam\lsandorder\bmu$ and $\blam\neq\bmu$.

\item If $\blam,\bmu\in\aCrystalVertices$, then $\blam\sim\bmu$ if $\bmu=w\blam$, for some face permutation $w\in\sym$, and $\bx_{\ppath[\blam]}=\bx_{\ppath[\bmu]}$. Let $\daCrystalVertices[]=\aCrystalVertices/\sim$ and $\dCrystalVertices=\CrystalVertices/\sim$ be the corresponding equivalence classes of vertices.
\end{enumerate}
\end{Definition}

We also consider $\lsandorder$ as a total order on $\daCrystalVertices$ and $\dCrystalVertices$. We abuse notation and identify paths in $\aCrystalVertices$ and $\CrystalVertices$ with their equivalence class representative in $\daCrystalVertices$ and $\dCrystalVertices$, respectively.

\begin{Example}\label{E:BasesDominance}
In \autoref{Eq:MainExampleOrder} we have $\nu\sandorder\mu\sandorder\lambda$.
\end{Example}

By convention, $\fParts{n}{\bfweight}\subset\fParts{n}{\abfweight}$ and $\CrystalVertices\subset\aCrystalVertices$.
The embeddings
$\crystalgraph[\bfweight]\hookrightarrow\Crystalgraph[\bfweight]$
and
$\crystalgraph[\abfweight]\hookrightarrow\Crystalgraph[\abfweight]$
of \autoref{L:Embeddings} play an important role because they allow us to attach a residue sequence $\Res[\bpath]\in\vertices^{n,\ell}$ to each rooted path in these crystals.
In particular, part (c) defines diagrams $\idem[\blam]$ and $\dotidem[\blam]$ for a choice of preferred $\blam\in\CrystalVertices$.
% The results that follow are independent of the (arbitrary) choices of preferred paths in \autoref{D:ConstructionChoice}.(c) by \autoref{P:ProofsDetour} below.

%Good paths are analogs of good residue sequences for symmetric groups. By definition, every path in $\Parts{n}{\bfweight}$ is good, so $\CrystalVertices$ is the set of sinks of the paths in $\Parts{n}{\bfweight}$, so \autoref{D:ConstructionChoice} is trivial in this case. The main complication in this definition is that good conditions in \autoref{D:ConstructionChoice} are used to choose distinguished paths in the crystal graph $\Crystalgraph[\abfweight]$, which will be used to define the cellular basis of the infinite dimensional wKLRW algebras $\aWA$.

%The extra complications in \autoref{D:ConstructionChoice} appear to have nothing to do with the algebras $\WA(X)$. It turns out that we need the unsteady diagrams in $\aWA$ to prove that the steady diagrams in $\WA(X)$ actually are steady.

\begin{Remark}\label{R:ConstructionChoice}
We will see in \autoref{L:ProofsIdempotentOne}.(b) below that if $\blam,\bmu\in\aCrystalVertices$, then $\idem[\blam]\neq\idem[\bmu]$. In particular, if $\blam\sim\bmu$, then $\bx_{\blam}=\bx_{\bmu}$ but the idempotent diagrams are not equal because they have different residue sequences.
%Moreover, by \autoref{C:ProofsIdempotentTwo} below, $\CrystalVertices$ and $\aCrystalVertices$ are totally ordered by $\sandorder$ from \autoref{D:ConstructionChoice}. \DT{Careful: not true}
\end{Remark}

\begin{Example}\label{E:ConstructionChoice}
We explain how \autoref{D:ConstructionChoice} works in \autoref{E:MainExampleTheBeastItself}. By \autoref{Eq:MainExample}, we can take $\fParts{5}{\bfweight}$ to be $\set{\ppath[\lambda],\ppath[\mu], \ppath[\nu]}$. In $\crystalgraph[{\fweight[2]}]$, there is one path to $\mu$, three paths
to $\lambda$ and one path $\mu$. The three paths to $\lambda$, and the two paths to $\mu$, differ by face permutations around nonadjacent squares, so there is only one equivalence class of paths to each vertex. In contrast, the paths $\ppath[\lambda^{\prime}]$ and $\qpath[\lambda^{\prime}]$ in \autoref{E:MainExampleTheBeastItself}.(f) are not equivalent, so they give different elements of $\fParts{7}{\bfweight}$. There are nontrivial face permutations between these paths.
% By \autoref{D:ConstructionChoice}, all paths in the crystal graph $\crystalgraph[{\fweight[2]}]$ are good.

Paths that are in $\crystalgraph[\bfweight]$ and not in $\crystalgraph[\abfweight]$ only appear when we consider paths for which the idempotent diagrams have strings that are close to the affine red strings. These paths will be used to construct the unsteady basis elements of $\aWA$. For example, the idempotent diagram at the end of \autoref{E:ConstructionAffineRed} corresponds to the path in $\crystalgraph[\abfweight]$ with residue sequence $(2,2)$. Similarly, the idempotent diagram in \autoref{E:ConstructionExampleAffine} corresponds to the path in $\crystalgraph[\abfweight]$ with residue sequence $(1,1)$. Neither of these paths belongs to $\crystalgraph[\bfweight]$: even though the path $\bpath$ with residue sequence $(1,1)$ starts in~$\crystalgraph[\bfweight]$ its sink $\omega(\bpath)$ is not in $\crystalgraph[\bfweight]$.
\end{Example}

\begin{Example}\label{E:TypeDNotGood}
  Let $\Gamma$ be the quiver of type $\typed[4]$ and let $\bfweight=\Lambda_2$ and consider $\WAc[5](x)$. The paths of length five in the crystal graph are given by:
  \[
  \begin{tikzpicture}[>=latex,line join=bevel,scale=0.5]
  \node (node_0) at (287.49bp,716.23bp) [draw,draw=none] {$\bullet$};
  \node (node_1) at (287.49bp,645.62bp) [draw,draw=none] {$\bullet$};
  \node (node_2) at (194.49bp,575.02bp) [draw,draw=none] {$\bullet$};
  \node (node_3) at (376.49bp,575.02bp) [draw,draw=none] {$\bullet$};
  \node (node_4) at (287.49bp,575.02bp) [draw,draw=none] {$\bullet$};
  \node (node_5) at (287.49bp,504.42bp) [draw,draw=none] {$\bullet$};
  \node (node_6) at (167.49bp,504.42bp) [draw,draw=none] {$\bullet$};
  \node (node_7) at (403.49bp,504.42bp) [draw,draw=none] {$\bullet$};
  \node (node_8) at (230.49bp,433.81bp) [draw,draw=none] {$\bullet$};
  \node (node_9) at (352.49bp,433.81bp) [draw,draw=none] {$\bullet$};
  \node (node_10) at (64.491bp,433.81bp) [draw,draw=none] {$\bullet$};
  \node (node_11) at (504.49bp,433.81bp) [draw,draw=none] {$\bullet$};
  \node (node_12) at (178.49bp,362.27bp) [draw,draw=none] {$\bullet$};
  \node (node_13) at (363.49bp,362.27bp) [draw,draw=none] {$\bullet$};
  \node (node_14) at (53.491bp,362.27bp) [draw,draw=none] {$\bullet$};
  \node (node_15) at (548.49bp,362.27bp) [draw,draw=none] {$\bullet$};
  \draw [red,->] (node_0) ..controls (287.49bp,698.08bp) and (287.49bp,679.29bp)  .. (node_1);
  \definecolor{strokecol}{rgb}{0.0,0.0,0.0};
  \pgfsetstrokecolor{strokecol}
  \draw (295.99bp,680.93bp) node {$2$};
  \draw [blue,->] (node_1) ..controls (262.71bp,626.35bp) and (232.9bp,604.36bp)  .. (node_2);
  \draw (260.56bp,610.32bp) node {$1$};
  \draw [green,->] (node_1) ..controls (311.14bp,626.4bp) and (339.47bp,604.56bp)  .. (node_3);
  \draw (351.08bp,610.32bp) node {$3$};
  \draw [black,->] (node_1) ..controls (287.49bp,627.48bp) and (287.49bp,608.69bp)  .. (node_4);
  \draw (295.99bp,610.32bp) node {$4$};
  \draw [green,->] (node_2) ..controls (208.91bp,557.8bp) and (224.2bp,541.77bp)  .. (239.49bp,530.72bp) .. controls (246.16bp,525.9bp) and (253.87bp,521.43bp)  .. (node_5);
  \draw (247.99bp,539.72bp) node {$3$};
  \draw [black,->] (node_2) ..controls (187.66bp,556.67bp) and (180.04bp,537.31bp)  .. (node_6);
  \draw (192.7bp,539.72bp) node {$4$};
  \draw [blue,->] (node_3) ..controls (365.3bp,557.47bp) and (353.1bp,541.24bp)  .. (339.49bp,530.72bp) .. controls (332.8bp,525.54bp) and (324.87bp,521.03bp)  .. (node_5);
  \draw (365.82bp,539.72bp) node {$1$};
  \draw [black,->] (node_3) ..controls (383.32bp,556.67bp) and (390.94bp,537.31bp)  .. (node_7);
  \draw (401.7bp,539.72bp) node {$4$};
  \draw [blue,->] (node_4) ..controls (279.46bp,557.04bp) and (269.8bp,540.08bp)  .. (256.49bp,530.72bp) .. controls (247.33bp,524.27bp) and (236.55bp,519.46bp)  .. (node_6);
  \draw (281.9bp,539.72bp) node {$1$};
  \draw [green,->] (node_4) ..controls (295.56bp,557.09bp) and (305.24bp,540.17bp)  .. (318.49bp,530.72bp) .. controls (327.42bp,524.35bp) and (337.97bp,519.57bp)  .. (node_7);
  \draw (326.99bp,539.72bp) node {$3$};
  \draw [red,->] (node_5) ..controls (272.73bp,485.65bp) and (255.69bp,465.14bp)  .. (node_8);
  \draw (274.27bp,469.12bp) node {$2$};
  \draw [black,->] (node_5) ..controls (304.52bp,485.45bp) and (324.5bp,464.35bp)  .. (node_9);
  \draw (336.22bp,469.12bp) node {$4$};
  \draw [green,->] (node_6) ..controls (187.15bp,486.72bp) and (209.5bp,469.43bp)  .. (231.49bp,460.12bp) .. controls (249.17bp,452.63bp) and (269.26bp,447.29bp)  .. (node_9);
  \draw (239.99bp,469.12bp) node {$3$};
  \draw [red,->] (node_6) ..controls (139.89bp,485.04bp) and (106.43bp,462.75bp)  .. (node_10);
  \draw (136.75bp,469.12bp) node {$2$};
  \draw [blue,->] (node_7) ..controls (393.15bp,491.45bp) and (387.32bp,484.48bp)  .. (382.49bp,478.12bp) .. controls (375.95bp,469.51bp) and (369.06bp,459.61bp)  .. (node_9);
  \draw (390.99bp,469.12bp) node {$1$};
  \draw [red,->] (node_7) ..controls (430.55bp,485.04bp) and (463.37bp,462.75bp)  .. (node_11);
  \draw (474.51bp,469.12bp) node {$2$};
  \draw [black,->] (node_8) ..controls (217.19bp,415.03bp) and (202.13bp,394.88bp)  .. (node_12);
  \draw (219.18bp,398.51bp) node {$4$};
  \draw [red,->] (node_9) ..controls (355.22bp,415.54bp) and (358.19bp,396.79bp)  .. (node_13);
  \draw (367.8bp,398.51bp) node {$2$};
  \draw [green,->] (node_10) ..controls (61.759bp,415.54bp) and (58.793bp,396.79bp)  .. (node_14);
  \draw (68.799bp,398.51bp) node {$3$};
  \draw [blue,->] (node_11) ..controls (515.68bp,415.13bp) and (528.25bp,395.26bp)  .. (node_15);
  \draw (540.23bp,398.51bp) node {$1$};
\end{tikzpicture}
.
  \]
  Applying \autoref{D:ConstructionIdempotents}, the idempotent diagrams in $\WAc[5](X)$ are:
  \begin{align*}
  \idem[21423]&=
\begin{tikzpicture}[anchorbase,smallnodes,xscale=1.1,yscale=0.7]
\draw[redstring](0,0)node[below]{$\,2$}--+(0,1);
\draw[solid](-0.1,0)node[below]{$2\,$}--+(0,1);
\draw[ghost,double](0.9,0)--+(0,1)node[above]{$2$};
\draw[ghost](-0.2,0)--+(0,1)node[above]{$\,1$};
\draw[solid](-1.2,0)node[below]{$1\,$}--+(0,1);
\draw[solid](0.8,0)node[below]{$4\,$}--+(0,1);
\draw[solid](-0.3,0)node[below]{$2\,$}--+(0,1);
\draw[ghost,double](0.7,0)--+(0,1)node[above]{$2$};
\draw[solid](0.6,0)node[below]{$3\,$}--+(0,1);
\end{tikzpicture},
  &
  \idem[23142]&=
\begin{tikzpicture}[anchorbase,smallnodes,xscale=1.1,yscale=0.7]
\draw[use as bounding box,draw=none] (-1.7,-0.2) rectangle (0.9,1.2);
\draw[redstring](0,0)node[below]{$\,2$}--+(0,1);
\draw[solid](-0.1,0)node[below]{$2\,$}--+(0,1);
\draw[ghost,double](0.9,0)--+(0,1)node[above]{$2$};
\draw[solid](-1.2,0)node[below]{$1\,$}--+(0,1);
\draw[ghost](-0.2,0)--+(0,1)node[above]{$\,1$};
\draw[solid](0.8,0)node[below]{$3\,$}--+(0,1);
\draw[solid](0.7,0)node[below]{$4\,$}--+(0,1);
\draw[solid](-0.4,0)node[below]{$2\,$}--+(0,1);
\draw[ghost,double](0.6,0)--+(0,1)node[above]{$2$};
\end{tikzpicture},
  \\
    \idem[21324]&=
\begin{tikzpicture}[anchorbase,smallnodes,xscale=1.1,yscale=0.7]
\draw[redstring](0,0)node[below]{$\,2$}--+(0,1);
\draw[solid](-0.1,0)node[below]{$2\,$}--+(0,1);
\draw[ghost,double](0.9,0)--+(0,1)node[above]{$2$};
\draw[solid](-1.2,0)node[below]{$1\,$}--+(0,1);
\draw[ghost](-0.2,0)--+(0,1)node[above]{$\,1$};
\draw[solid](0.8,0)node[below]{$3\,$}--+(0,1);
\draw[solid](-0.3,0)node[below]{$2\,$}--+(0,1);
\draw[ghost,double](0.7,0)--+(0,1)node[above]{$2$};
\draw[solid](0.6,0)node[below]{$4\,$}--+(0,1);
\end{tikzpicture},
  &
  \idem[23421]&=
\begin{tikzpicture}[anchorbase,smallnodes,xscale=1.1,yscale=0.7]
\draw[redstring](0,0)node[below]{$\,2$}--+(0,1);
\draw[solid](-0.1,0)node[below]{$2\,$}--+(0,1);
\draw[ghost,double](0.9,0)--+(0,1)node[above]{$2$};
\draw[solid](0.8,0)node[below]{$3\,$}--+(0,1);
\draw[solid](0.7,0)node[below]{$4\,$}--+(0,1);
\draw[solid](-0.4,0)node[below]{$2\,$}--+(0,1);
\draw[ghost,double](0.6,0)--+(0,1)node[above]{$2$};
\draw[solid](-1.5,0)node[below]{$1\,$}--+(0,1);
\draw[ghost](-0.5,0)--+(0,1)node[above]{$\,1$};
\end{tikzpicture}.
\end{align*}
  The idempotents for each vertex are independent of the choice of path because all paths to the same vertex differ by non-adjacent squares.
  In particular, note that the $\bx$-coordinates for $\idem[21423]$ and $\idem[21324]$  coincide, so these two idempotents are equivalent under \autoref{D:ConstructionChoice}.(e).
  Note that the idempotent diagrams $\idem[21423]$ and $\idem[21324]$ are not equal because they have different residue sequences. The equivalence of these two idempotents is symptomatic of~$\WAc[5](X)$ not being cellular.
\end{Example}

The next definition associates diagrams in the wKLRW algebra to the face permutations of \autoref{D:ConstructionPermutationDiagrams}. Recall from \autoref{D:FiniteTypesSSPerm} that a two color face permutation is an equivalence class of permutations that act in the same way on the residue sequence of a given path. The next definition implicitly chooses an equivalence class representative for each face permutation. We also need to consider related $1$-color relations, which correspond to the nil-Hecke case.

\begin{Definition}\label{D:ConstructionPermBasis}
Let $w$ be a path permutation from $\bpath=(\bpath_{1},\dots,\bpath_{\hell})$ to $\cpath=(\cpath_{1},\dots,\cpath_{\hell})$, for $\bpath,\cpath\in\Crystalgraph$. Given a \emph{basic diagrammatic path permutation} $w$, the
corresponding diagram is given by:

\begin{enumerate}

\item Putting the idempotent diagrams $\idem[\bpath]$ and $\idem[\cpath]$ along the top and bottom, respectively, of the permutation diagram.

\item Connecting all of the vertices of the same residue using a choice of representative of the longest word, for the solid strings and their ghosts, as in the permutation diagrams below. All other strings, including (affine) red strings, are vertical.

\item  To each partial face permutation we associate the subdiagram of the corresponding basic diagrammatic face permutation.
\end{enumerate}

We make the diagrams for basic face permutations more explicit.  The diagrams in brackets, on the right-hand side below, show how to place these strings separately for the $i$ and $j$ strings in the $2$-color relations, in accordance with condition (b) above.

\begin{enumerate}[label=(\roman*)]

\item $1$-color relations: to a subsequence of consecutive solid $i$-strings, we attach the permutation for (a choice of representative of) the corresponding longest word. Ignoring all other strings:
\begin{gather*}
\begin{tikzpicture}[anchorbase,smallnodes,rounded corners]
\draw[solid,black](0,0)node[below]{$i$}--++(0.5,0.5) node[above,yshift=-1pt]{$i$};
\draw[solid,black](0.5,0)node[below]{$i$}--++(-0.5,0.5)node[above,yshift=-1pt]{$i$};
\end{tikzpicture}
,\quad
\begin{tikzpicture}[anchorbase,smallnodes,rounded corners]
\draw[solid,black](0,0)node[below]{$i$}--++(1,1) node[above,yshift=-1pt]{$i$};
\draw[solid,black](0.5,0)node[below]{$i$}--++(-0.5,0.5)--++(0.5,0.5)node[above,yshift=-1pt]{$i$};
\draw[solid,black](1,0)node[below]{$i$}--++(-1,1)node[above,yshift=-1pt]{$i$};
\end{tikzpicture}
,\quad\text{\etc}
\end{gather*}

\item $2$-color relations:
\begin{gather*}
\text{adjacent square}\colon
\begin{tikzpicture}[anchorbase,smallnodes]
\draw[solid,black](0,0)node[below]{$i$}--++(0.5,0.5) node[above,yshift=-1pt]{$i$};
\draw[solid,spinach](0.5,0)node[below]{$j$}--++(-0.5,0.5)node[above,yshift=-1pt]{$j$};
\end{tikzpicture}
,\quad
\text{octagon}\colon
\begin{tikzpicture}[anchorbase,smallnodes]
\draw[solid,black](0,0)node[below]{$i$}--++(1,0.5) node[above,yshift=-1pt]{$i$};
\draw[solid,black](1.5,0)node[below]{$i$}--++(-1,0.5)node[above,yshift=-1pt]{$i$};
\draw[solid,spinach](0.5,0)node[below]{$j$}--++(1,0.5)node[above,yshift=-1pt]{$j$};
\draw[solid,spinach](1,0)node[below]{$j$}--++(-1,0.5)node[above,yshift=-1pt]{$j$};
\end{tikzpicture}\quad
\leftrightsquigarrow\quad
\left(
\begin{tikzpicture}[anchorbase,smallnodes]
\draw[solid,black](0,0)node[below]{$i$}--++(0.5,0.5) node[above,yshift=-1pt]{$i$};
\draw[solid,black](0.5,0)node[below]{$i$}--++(-0.5,0.5)node[above,yshift=-1pt]{$i$};
\end{tikzpicture}
\text{ and }
\begin{tikzpicture}[anchorbase,smallnodes]
\draw[solid,spinach](0,0)node[below]{$j$}--++(0.5,0.5) node[above,yshift=-1pt]{$j$};
\draw[solid,spinach](0.5,0)node[below]{$j$}--++(-0.5,0.5)node[above,yshift=-1pt]{$j$};
\end{tikzpicture}
\right)
,
\\
\text{decagon}\colon\hspace{-0.15cm}
\begin{tikzpicture}[anchorbase,smallnodes,rounded corners]
\draw[solid,black](0,0)node[below]{$2$}--++(0,1)--++(1,0.5) node[above,yshift=-1pt]{$2$};
\draw[solid,black](2,0)node[below]{$2$}--++(0,1)--++(-1.5,0.5)node[above,yshift=-1pt]{$2$};
\draw[solid,spinach](0.5,0)node[below]{$3$}--++(1,1)--++(0.5,0.5)node[above,yshift=-1pt]{$3$};
\draw[solid,spinach](1,0)node[below]{$3$}--++(-0.5,0.5)--++(0.5,0.5)--++(0.5,0.5)node[above,yshift=-1pt]{$3$};
\draw[solid,spinach](1.5,0)node[below]{$3$}--++(-1,1)--++(-0.5,0.5)node[above,yshift=-1pt]{$3$};
\end{tikzpicture}
,\;
\begin{tikzpicture}[anchorbase,smallnodes,rounded corners]
\draw[solid,black](0,0)node[below]{$2$}--++(0,1)--++(0.5,0.5) node[above,yshift=-1pt]{$2$};
\draw[solid,black](2,0)node[below]{$2$}--++(0,1)--++(-0.5,0.5)node[above,yshift=-1pt]{$2$};
\draw[solid,spinach](0.5,0)node[below]{$3$}--++(1,1)--++(0.5,0.5)node[above,yshift=-1pt]{$3$};
\draw[solid,spinach](1,0)node[below]{$3$}--++(-0.5,0.5)--++(0.5,0.5)--++(0,0.5)node[above,yshift=-1pt]{$3$};
\draw[solid,spinach](1.5,0)node[below]{$3$}--++(-1,1)--++(-0.5,0.5)node[above,yshift=-1pt]{$3$};
\end{tikzpicture}
,\;
\begin{tikzpicture}[anchorbase,smallnodes,rounded corners,yscale=-1]
\draw[solid,black](0,0)node[above,yshift=-1pt]{$2$}--++(0,1)--++(1,0.5) node[below]{$2$};
\draw[solid,black](2,0)node[above,yshift=-1pt]{$2$}--++(0,1)--++(-1.5,0.5)node[below]{$2$};
\draw[solid,spinach](0.5,0)node[above,yshift=-1pt]{$3$}--++(1,1)--++(0.5,0.5)node[below]{$3$};
\draw[solid,spinach](1,0)node[above,yshift=-1pt]{$3$}--++(-0.5,0.5)--++(0.5,0.5)--++(0.5,0.5)node[below]{$3$};
\draw[solid,spinach](1.5,0)node[above,yshift=-1pt]{$3$}--++(-1,1)--++(-0.5,0.5)node[below]{$3$};
\end{tikzpicture}
,\;
\begin{tikzpicture}[anchorbase,smallnodes,rounded corners,yscale=-1]
\draw[solid,black](0,0)node[above,yshift=-1pt]{$2$}--++(0,1)--++(0.5,0.5) node[below]{$2$};
\draw[solid,black](2,0)node[above,yshift=-1pt]{$2$}--++(0,1)--++(-0.5,0.5)node[below]{$2$};
\draw[solid,spinach](0.5,0)node[above,yshift=-1pt]{$3$}--++(1,1)--++(0.5,0.5)node[below]{$3$};
\draw[solid,spinach](1,0)node[above,yshift=-1pt]{$3$}--++(-0.5,0.5)--++(0.5,0.5)--++(0,0.5)node[below]{$3$};
\draw[solid,spinach](1.5,0)node[above,yshift=-1pt]{$3$}--++(-1,1)--++(-0.5,0.5)node[below]{$3$};
\end{tikzpicture}\quad
\leftrightsquigarrow\quad
\left(
\begin{tikzpicture}[anchorbase,smallnodes]
\draw[solid,black](0,0)node[below]{$2$}--++(0.5,0.5) node[above,yshift=-1pt]{$2$};
\draw[solid,black](0.5,0)node[below]{$2$}--++(-0.5,0.5)node[above,yshift=-1pt]{$2$};
\end{tikzpicture}
\text{ and }
\begin{tikzpicture}[anchorbase,smallnodes,rounded corners]
\draw[solid,spinach](0,0)node[below]{$3$}--++(1,1) node[above,yshift=-1pt]{$3$};
\draw[solid,spinach](0.5,0)node[below]{$3$}--++(-0.5,0.5)--++(0.5,0.5)node[above,yshift=-1pt]{$3$};
\draw[solid,spinach](1,0)node[below]{$3$}--++(-1,1)node[above,yshift=-1pt]{$3$};
\end{tikzpicture}
\right)
,
\\
\text{tetradecagon}\colon
\begin{array}{l}
2332323
\\[-0.07cm]
\ =3223332
\\[-0.07cm]
\ =3232332
\\[-0.07cm]
\ =2333223
\end{array}\quad
\leftrightsquigarrow\quad
\left(
\begin{tikzpicture}[anchorbase,smallnodes,rounded corners]
\draw[solid,black](0,0)node[below]{$2$}--++(1,1) node[above,yshift=-1pt]{$2$};
\draw[solid,black](0.5,0)node[below]{$2$}--++(-0.5,0.5)--++(0.5,0.5)node[above,yshift=-1pt]{$2$};
\draw[solid,black](1,0)node[below]{$2$}--++(-1,1)node[above,yshift=-1pt]{$2$};
\end{tikzpicture}
\text{ and }
\begin{tikzpicture}[anchorbase,smallnodes]
\draw[solid,spinach](0,0)node[below]{$3$}--++(1,0.5) node[above,yshift=-1pt]{$3$};
\draw[solid,spinach](1.5,0)node[below]{$3$}--++(-1,0.5)node[above,yshift=-1pt]{$3$};
\draw[solid,spinach](0.5,0)node[below]{$3$}--++(1,0.5)node[above,yshift=-1pt]{$3$};
\draw[solid,spinach](1,0)node[below]{$3$}--++(-1,0.5)node[above,yshift=-1pt]{$3$};
\end{tikzpicture}
\right)
.
\end{gather*}

\item Multi-color relations: As in (i) and (ii), at the top and bottom of the diagram, any subsequences of repeated residues are connected by the permutation diagram on the longest word for this subsequence.
\end{enumerate}

All of these diagrams are classical permutation diagrams with strings labeled by their residues. The colors in the diagrams above are only used to highlight the different residues. Although we have omitted them, there will typically be (affine) red strings in these diagrams.

Finally, a \emph{diagrammatic face permutation} between
two paths $\bpath,\cpath\in\Parts{n}{\abfweight}$ is a composition of basic and partial diagrammatic face permutations. Recall from \autoref{D:FiniteTypesSSPerm}, that $\Face(\bpath,\cpath)$ is the set of face permutations from $\bpath$ to $\cpath$. If $\bT\in\Face(\bpath,\cpath)$ let $D_{\bT}$ be the corresponding diagrammatic face permutation.
By definition, if $\bT\in\Face(\bpath,\cpath)$, then $D_{\bT}=\idem[\bpath]D_{\bT}\idem[\cpath]$. %Set $D_{\bT}=D_{\bT}E_{\cpath}$.
\end{Definition}

The definition of a diagrammatic face permutation in \autoref{D:ConstructionPermBasis} involves a choice because the braid relations do not hold in $\aWA$. The bases that we define below depend upon these choices, but our results hold for all choices. This is explained by \cite[Lemma 3B.3]{MaTu-klrw-algebras}, which says that $\WA(X)$ is a filtered algebra and any two choices for $D_{\bT}$ are equal in the associated graded algebra.

The map $\sigma\mapsto\sigma^{-1}$ gives a bijection $\Face(\blam,\bmu)\xrightarrow{\sim}\Face(\bmu,\blam)$. This is problematic for us because we want to associate each face permutation to a unique cell module. The next definition corrects for this overcounting.

Recall that $\blam=\omega(\bpath)\in\aCrystalVertices$ is the sink of the path $\bpath$ and that $\ppath$ is the preferred path for $\blam$.

\begin{Definition}\label{D:ConstructionPermBasisTwo}
Let $\blam\in\aCrystalVertices$ and $\bmu\in\CrystalVertices$.
The sets of \emph{face permutations} for $\blam$ %and $\bmu$
are
\begin{gather*}
\Face(\blam)=\bigcup_{\bpath\in\fParts{n}{\abfweight},\omega(\bpath)\gsandorder\blam}
\Face(\ppath,\bpath).
\quad\text{and}\quad
\fFace(\bmu)=
\bigcup_{\bpath\in\fParts{n}{\bfweight},\omega(\bpath)\gsandorder\bmu}
\fFace(\ppath[\bmu],\bpath)
\end{gather*}
If $\bnu\in\daCrystalVertices$, set $\Face(\bnu)=\bigcup_{\blam\sim\bmu}\Face(\blam)$. Similarly, if $\bnu\in\dCrystalVertices$ set $\fFace(\bnu)=\bigcup_{\blam\sim\bmu}\fFace(\blam)$.
\end{Definition}

By construction, if $\mu\in\fParts{n}{\bfweight}$, then
$\fFace(\bmu)\subseteq\Face(\bmu)$
Since $\CrystalVertices\subset\aCrystalVertices$.

The face permutations from $\blam$ to $\bnu$ are analogs of semistandard $\blam$-tableaux of type $\bnu$ in the classical types. As we already saw in \autoref{E:ConstructionChoice}, if $\bpath\in\fParts{n}{\abfweight}$ has sink $\blam$, then $\Face(\blam)$ contains all face permutations from $\ppath$ to $\bpath$.

As in \autoref{D:FiniteTypesDetour}, if $\blam\in\aCrystalVertices$
% and $\bmu\in\CrystalVertices$
then
the set of \emph{detour permutations} in $\Face(\blam)$ is
\begin{gather}\label{E:DetourPermutation}
\Detour(\blam) = \bigcup_{\bpath\in\fParts{n}{\abfweight},\omega(\bpath)=\blam}
\Face(\ppath,\bpath).
%\qquad\text{and}\qquad
%\bigcup_{\bpath\in\fParts{n}{\abfweight},\omega(\bpath)=\bmu} %\Face(\ppath[\bmu],\bpath)
\end{gather}
%are the sets of \emph{detour permutations} in $\hFace(\blam)$ and $\Face(\bmu)$, respectively.
By definition, if $\bmu\in\CrystalVertices$, then $\Detour(\bmu)\subseteq\fFace(\bmu)$.

\begin{Example}\label{E:ConstructionPermutationTwo}
Using the notation of \autoref{E:MainExampleTheBeastItself}.(d),
the face permutation $(2,3)(4,5)\in\fFace(\mu)$ and $(2,5,4)\in\fFace(\nu)$, but neither of these permutations belong to $\fFace(\lambda)$ since $\nu\lsandorder\mu\lsandorder\lambda$.
\end{Example}

%%%%%%%%%%%%%%%%%%%%%%%%%%%%%%%%%%%%%%%%%

\subsection{The bases}\label{SS:ConstructionBases}

%%%%%%%%%%%%%%%%%%%%%%%%%%%%%%%%%%%%%%%%%

We use the usual
\emph{diagrammatic antiinvolution} $(\placeholder)^{\star}$ given
by (the $\ring$-linear extension of) reflecting diagrams
bottom-to-top as illustrated by:
\begin{gather*}
\left(
\begin{tikzpicture}[anchorbase,smallnodes,rounded corners]
\draw[ghost](1.5,1)node[above,yshift=-1pt]{$j$}--++(1,-1);
\draw[ghost](2.5,1)node[above,yshift=-1pt]{$i$}--++(-1,-1);
\draw[ghost](4,1)node[above,yshift=-1pt]{$k$}--++(0.5,-1);
\draw[solid](1,1)node[above,yshift=-1pt]{$\phantom{i}$}--++(1,-1)node[below]{$j$};
\draw[solid](2,1)--++(-1,-1)node[below]{$i$};
\draw[solid](3,1)--++(0.5,-1)node[below]{$k$};
\draw[redstring](1.25,1)--++(0,-1)node[below]{$\rho$};
\end{tikzpicture}
\right)^{\star}
=
\begin{tikzpicture}[anchorbase,smallnodes,rounded corners]
\draw[ghost](1.5,1)node[above,yshift=-1pt]{$i$}--++(1,-1);
\draw[ghost](2.5,1)node[above,yshift=-1pt]{$j$}--++(-1,-1);
\draw[ghost](4,0)--++(0.5,1)node[above,yshift=-1pt]{$k$};
\draw[solid](1,1)node[above,yshift=-1pt]{$\phantom{i}$}--++(1,-1)node[below]{$i$};
\draw[solid](2,1)--++(-1,-1)node[below]{$j$};
\draw[solid](3,0)node[below]{$k$}--++(0.5,1);
\draw[redstring](1.25,1)--++(0,-1)node[below]{$\rho$};
\end{tikzpicture}
.
\end{gather*}
Recall that we often omit the word involutive from a sandwich cell datum, but the cell data below
are involutive using the diagrammatic antiinvolution $\star$.

Extending \autoref{SS:ConstructionSandwichDots}, if $\bsig\in\aCrystalVertices$ set
\begin{gather*}
\Affch[\bsig]=\Affch[{\ppath[{\bsig}]}],\quad
\Finch[\bsig]=\Finch[{\ppath[\bsig]}],\quad
c_{m}(\bsig)=c_{m}(\ppath[\bsig])\quad\text{and}\quad
d_{m}(\bsig)=d_{m}(\ppath[\bsig]),
\end{gather*}
for $1\leq m\leq n$.
If $\ba\in\Affch[\bsig]$ and $\bfi\in\Finch[\bsig]$ set
$\sandwich{\bsig}{\ba}{\bfi}=\daffine\zeetwo\dotidem[\bsig]$.

Now comes our main definition, which should be compared with \autoref{Eq:MainExampleSandwich}.

\begin{Definition}\label{R:ConstructionTheBasis}
Let $D_{\bS\bT}^{\ba,\bfi}=(D_{\bS}E_{\bnu})^{\star}\sandwich{\bsig}{\ba}{\bfi}\dotidem[\bsig]D_{\bT}E_{\bmu}$, for $\bsig\in\aCrystalVertices$, $\bS\in\Face(\ppath[\bsig],\bnu)$ and
$\bT\in\Face(\ppath[\bsig],\bmu)$. Define
\begin{gather}\label{Eq:ConstructionTheBasisAffine}
\BX=
\set[\big]{D_{\bS\bT}^{\ba,\bfi}|\blam\in\daCrystalVertices,\bS,\bT\in\Face(\bsig),\ba\in\Affch[\bsig],\bfi\in\Finch[\bsig],\text{ for }\bsig\in\blam}.
\end{gather}
Similarly, let $D_{\bS\bT}^{\bfi}=D_{\bS\bT}^{(0,\dots,0),\bfi}$, and define
\begin{gather}\label{Eq:ConstructionTheBasisFinite}
\BXc=
\set[\big]{D_{\bS\bT}^{\bfi}|\blam\in\dCrystalVertices,\bS,\bT\in\fFace(\bsig),\bfi\in\Finch[\bsig],\text{ for }\bsig\in\blam}.
\end{gather}
\end{Definition}

Our choice for the sandwich cell datum $\affine{\mathscr{C}}$ for $\WA[n](X)$ consists of:
\begin{itemize}

\item The middle set $\daCrystalVertices=(\daCrystalVertices,\lsandorder)$.
\item The positions of the strings at the top and bottom of the diagrams are given by the set $X$, from \autoref{Eq:X}. In the sandwich cellular basis, the top and bottom and top sets are equal to $\bigcup_{\blam\in\daCrystalVertices}\Face(\blam)$ composed with the diagrams $E_{\blam}=\sum_{\bsig\in\blam}E_{\bsig}$.

\item The bases of the sandwiched algebras are $\sandbasis[\blam]=\bigcup_{\bsig\in\blam}\set[\big]{\sandwich{\bsig}{\ba}{\bfi}|\ba\in\Affch[\bsig],\bfi\in\Finch[\bsig]}$, for $\blam\in\daCrystalVertices$, and the sandwiched algebras $\sand[\blam]=\bigoplus_{\bsig\in\blam}\sand[{\ppath[\bsig]}]$ are the subalgebras of the wKLRW algebra generated by these bases.

\item Our basis is the set $\BX[]=\BX$ from \autoref{Eq:ConstructionTheBasisAffine}, viewed as a map.

\item The degree function is $\bS\mapsto\deg D_{\bS}$.
\item The antiinvolution is the diagrammatic antiinvolution $(\placeholder)^{\star}$.

\end{itemize}

For $\WAc[n](X)$, the sandwich cell datum of our choice
$\mathscr{C}$ is essentially the same, except that $\daCrystalVertices$
is replaced with $\dCrystalVertices$, $\Face(\bsig)$ is replaced by $\fFace(\bsig)$, and we take the set in \autoref{Eq:ConstructionTheBasisFinite} as our basis.

\begin{Theorem}\label{T:ConstructionMain}
Suppose that $\quiver$ is a quiver of finite type.
\begin{enumerate}

\item The datum $\affine{\mathscr{C}}$ is an involutive graded affine sandwich cell datum for $\WA[n](X)$. In particular, \autoref{Eq:ConstructionTheBasisAffine} is a homogeneous affine sandwich cellular basis of $\WA[n](X)$.

\item The datum $\mathscr{C}$ is an involutive graded sandwich cell datum for $\WAc[n](X)$. In particular, \autoref{Eq:ConstructionTheBasisFinite} is a homogeneous sandwich cellular basis of $\WAc[n](X)$.

\end{enumerate}
\end{Theorem}

We prove this result in \autoref{S:Proofs}.

Let $\TA[n]$ and $\TAc[n]$, respectively, be
the KLR algebras and their cyclotomic
quotients associated to the fixed Kac--Moody datum.
(These are the algebras defined by Khovanov--Lauda, Rouquier and
Webster in {\eg} \cite{KhLa-cat-quantum-sln-first} and  \cite{Ro-2-kac-moody}.)
By \cite[Section 3F]{MaTu-klrw-algebras}, the algebras $\WA(X)$ and $\WAc(X)$ are isomorphic to the KLR  algebras $\TA$ and $\TAc$, respectively; see \autoref{P:ProofsKLRIsomorphism} below for a more direct argument. Hence, \autoref{T:ConstructionMain} implies:

\begin{Corollary}\label{C:ConstructionMainTwo}
Suppose that $\quiver$ is a quiver of finite type. Then $\TA[n]$ is a homogeneous affine sandwich cellular algebra and $\TAc[n]$ is a sandwich cellular algebra.
%\begin{enumerate}
%\item The set $\EX=\set[\big]{E^{\ba,\bfi}_{\bS\bT}|\blam\in\aCrystalVertices,\bS,\bT\in\hSFace(\blam),\ba\in\Affch,\bfi\inJ\Finch}$ is a homogeneous affine sandwich cellular basis of $\TA[n]$.
%\item The set $\EXc=\set[\big]{E^{\bfi}_{\bS\bT}|\blam\in\CrystalVertices,\bS,\bT\in\SFace(\blam),\bfi\in\Finch}$ is a homogeneous sandwich cellular basis of $\TAc[n]$.\qed
%\end{enumerate}
\end{Corollary}

%Note the surprising fact that there is a bijection $\hSFace(\blam)\xrightarrow{\sim}\hFace(\blam)$, for $\blam\in\CrystalVertices$, which also follows from an appropriate interpretation of \cite[Proposition 2.14]{We-weighted-klr}.

%%%%%%%%%%%%%%%%%%%%%%%%%%%%%%%%%%%%%%%%%

\subsection{Cellularity for general dominant weights}\label{SS:ConstructionGeneral}

%%%%%%%%%%%%%%%%%%%%%%%%%%%%%%%%%%%%%%%%%

The definition of the wKLRW algebra labels the red strings by $i\in\vertices$, which corresponds to labeling them with the corresponding fundamental weights $\fweight[i]$.
Using more general relations, the definition of wKLRW algebras
can be generalized to allow arbitrary dominant weights
as labels of the red strings,
see \cite[Definition 4.5]{We-knot-invariants}. Although we will not define these algebras explicitly, our results imply that these algebras are always cellular.
Given an $\ell$ tuple $\bbeta=(\beta_{1},\dots,\beta_{\ell})$
of dominant weights, let $\gWA[n](X)$ be the corresponding generalized wKLRW algebra and let $\gWAc[n](X)$ be its cyclotomic quotient.
As we now sketch, the cellularity results in this paper apply to $\gWA[n](X)$
an $\gWAc[n](X)$ as well.

\begin{Theorem}\label{T:ConstructionGeneralMain}
The algebra $\gWA[n](X)$
is graded affine sandwich cellular and
$\WAc[n](X)$ is graded sandwich cellular.
\end{Theorem}

\begin{proof}
The generalized wKLRW algebra $\gWA[n](X)$ is the idempotent truncation of a wKLRW algebra that has red strings labeled by
fundamental weights. More precisely, if $\beta=a_{1}\fweight[1]+\dots+a_{e}\fweight[e]$, for $a_{i}\in\N$, then set
\begin{gather*}
\Lambda=(\underbrace{\Lambda_{1},\dots,\Lambda_{1}}_{a_{1}},\dots,\underbrace{\Lambda_{e},\dots,\Lambda_{e}}_{a_{e}})
\end{gather*}
and let $\WA(X)$ be the corresponding wKLRW algebra where the red strings are placed very close together. For example:
\begin{gather*}
\beta=3\fweight[1]+\fweight[3]\colon
\begin{tikzpicture}[anchorbase,smallnodes,rounded corners]
\draw[redstring](0,0)node[below]{$1$}--++(0,1)node[above,yshift=-1pt]{$\phantom{i}$};
\draw[redstring](0.2,0)node[below]{$1$}--++(0,1)node[above,yshift=-1pt]{$\phantom{i}$};
\draw[redstring](0.4,0)node[below]{$1$}--++(0,1)node[above,yshift=-1pt]{$\phantom{i}$};
\draw[redstring](0.6,0)node[below]{$3$}--++(0,1)node[above,yshift=-1pt]{$\phantom{i}$};
\end{tikzpicture}
\xrightarrow[\text{truncation}]{\text{idempotent}}
\begin{tikzpicture}[anchorbase,smallnodes,rounded corners]
\draw[redstring](0,0)node[below]{$\beta$}--++(0,1)node[above,yshift=-1pt]{$\phantom{i}$};
\end{tikzpicture}
.
\end{gather*}
Then $\gWA(X)$ is an idempotent truncation of $\WA(X)$ by \cite[Proposition~3F.3]{MaTu-klrw-algebras}. Hence, the result follows by applying \autoref{T:ConstructionMain} and
\autoref{T:SandwichMain}(a).
\end{proof}

%%%%%%%%%%%%%%%%%%%%%%%%%%%%%%%%%%%%%%%%%

\section{Proof of sandwich cellularity}\label{S:Proofs}

%%%%%%%%%%%%%%%%%%%%%%%%%%%%%%%%%%%%%%%%%

This section proves \autoref{T:ConstructionMain}, our sandwich cellular basis theorem for $\aWA$ and $\WAc(X)$.

%%%%%%%%%%%%%%%%%%%%%%%%%%%%%%%%%%%%%%%%%

\subsection{Ordering in finite types}\label{SS:ProofsOne}

%%%%%%%%%%%%%%%%%%%%%%%%%%%%%%%%%%%%%%%%%

We first need a crucial lemma about the
paths in a fundamental crystal graph in finite type. Let $i\in\vertices$ and suppose that $\bpath$ is a rooted path in $\crystalgraph$. A permutation in the symmetric group $\sym$ is \emph{$\bpath$-nonadjacent} if it only swaps nonadjacent residues in $\bpath$.

\begin{Definition}\label{D:ProofsGraphs}
A crystal $\crystal$, is \emph{$n$-path ordered} if whenever $\bpath$ and $\cpath$ are paths in $\crystal$ such $\Res[\cpath]=w\Res[\bpath]$ for a $\bpath$-nonadjacent permutation $w$, then $\omega(\bpath)=\omega(\cpath)$. The crystal graph $\crystal$ is \emph{path ordered} if it is $n$-path ordered for all $n\in\N$.
\end{Definition}

That is, $\crystal$ is path ordered if the reside sequences of paths in $\crystalgraph[\Lambda]$ to distinct vertices are not related by permutations that swap only nonadjacent residues.
Consequently, by \autoref{P:FiniteTypesSSR} and \autoref{L:ConstructionPermBasisAdSquares}, if $\crystal$ is path ordered, then any two paths in $\crystal$ that differ by nonadjacent permutations
end at the same vertex and  give the same idempotent diagram in $\aWA$.

For this paper, being path ordered is one of the most fundamental properties of a crystal graph of finite type because it is crucial to proving that the wKLRW algebras are sandwich cellular algebras. The proof of the next result is easy but only because it invokes Lusztig's quantum PBW theorem, which is only valid in finite type. Before we give the proof, we recall some of these notions in an example.

\begin{Example}\label{E:ProofsPBW}
In finite type every crystal has a PBW realization that is
obtained from the $\mathcal{B}_{\infty}$ crystal
of PBW monomials that models Verma
modules; see
{\eg} \cite[Chapter 12]{BuSc-crystal-bases}.
For example, for type $\typeb[3]$ choose $w_{0}=r_{1}r_{2}r_{3}r_{2}r_{1}r_{2}r_{3}r_{2}r_{3}$, where
$\set{r_{i}|i\in\vertices}$ is the set of Coxeter generators of the Weyl group
of type $\typeb[3]$. Then
the crystal graph from \autoref{Eq:MainExample} in the PBW notation is
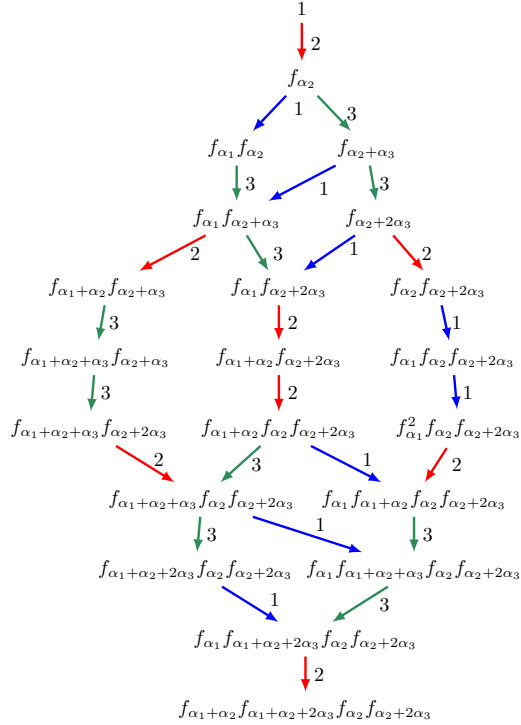
\begin{gather*}
\scalebox{0.75}{$\begin{tikzpicture}[>=latex,line join=bevel,xscale=-0.6,yscale=0.5,every path/.style={very thick},anchorbase]
\node (node_0) at (210.5bp,289.0bp) [draw,draw=none] {$f_{\alpha_{1} + \alpha_{2}} f_{\alpha_{2}} f_{\alpha_{2} + 2 \alpha_{3}}$};
\node (node_14) at (274.5bp,218.0bp) [draw,draw=none] {$f_{\alpha_{1} + \alpha_{2} + \alpha_{3}} f_{\alpha_{2}} f_{\alpha_{2} + 2 \alpha_{3}}$};
\node (node_20) at (98.5bp,218.0bp) [draw,draw=none] {$f_{\alpha_{1}} f_{\alpha_{1} + \alpha_{2}} f_{\alpha_{2}} f_{\alpha_{2} + 2 \alpha_{3}}$};
\node (node_1) at (61.5bp,289.0bp) [draw,draw=none] {$f_{\alpha_{1}}^{2} f_{\alpha_{2}} f_{\alpha_{2} + 2 \alpha_{3}}$};
\node (node_2) at (245.5bp,500.0bp) [draw,draw=none] {$f_{\alpha_{1}} f_{\alpha_{2} + \alpha_{3}}$};
\node (node_5) at (210.5bp,430.0bp) [draw,draw=none] {$f_{\alpha_{1}} f_{\alpha_{2} + 2 \alpha_{3}}$};
\node (node_18) at (352.5bp,430.0bp) [draw,draw=none] {$f_{\alpha_{1} + \alpha_{2}} f_{\alpha_{2} + \alpha_{3}}$};
\node (node_3) at (362.5bp,360.0bp) [draw,draw=none] {$f_{\alpha_{1} + \alpha_{2} + \alpha_{3}} f_{\alpha_{2} + \alpha_{3}}$};
\node (node_7) at (367.5bp,289.0bp) [draw,draw=none] {$f_{\alpha_{1} + \alpha_{2} + \alpha_{3}} f_{\alpha_{2} + 2 \alpha_{3}}$};
\node (node_4) at (98.5bp,148.0bp) [draw,draw=none] {$f_{\alpha_{1}} f_{\alpha_{1} + \alpha_{2} + \alpha_{3}} f_{\alpha_{2}} f_{\alpha_{2} + 2 \alpha_{3}}$};
\node (node_17) at (188.5bp,78.0bp) [draw,draw=none] {$f_{\alpha_{1}} f_{\alpha_{1} + \alpha_{2} + 2 \alpha_{3}} f_{\alpha_{2}} f_{\alpha_{2} + 2 \alpha_{3}}$};
\node (node_11) at (210.5bp,360.0bp) [draw,draw=none] {$f_{\alpha_{1} + \alpha_{2}} f_{\alpha_{2} + 2 \alpha_{3}}$};
\node (node_6) at (191.5bp,710.0bp) [draw,draw=none] {$1$};
\node (node_12) at (191.5bp,640.0bp) [draw,draw=none] {$f_{\alpha_{2}}$};
\node (node_8) at (78.5bp,430.0bp) [draw,draw=none] {$f_{\alpha_{2}} f_{\alpha_{2} + 2 \alpha_{3}}$};
\node (node_16) at (66.5bp,360.0bp) [draw,draw=none] {$f_{\alpha_{1}} f_{\alpha_{2}} f_{\alpha_{2} + 2 \alpha_{3}}$};
\node (node_9) at (137.5bp,570.0bp) [draw,draw=none] {$f_{\alpha_{2} + \alpha_{3}}$};
\node (node_13) at (127.5bp,500.0bp) [draw,draw=none] {$f_{\alpha_{2} + 2 \alpha_{3}}$};
\node (node_10) at (188.5bp,8.0bp) [draw,draw=none] {$f_{\alpha_{1} + \alpha_{2}} f_{\alpha_{1} + \alpha_{2} + 2 \alpha_{3}} f_{\alpha_{2}} f_{\alpha_{2} + 2 \alpha_{3}}$};
\node (node_19) at (245.5bp,570.0bp) [draw,draw=none] {$f_{\alpha_{1}} f_{\alpha_{2}}$};
\node (node_15) at (279.5bp,148.0bp) [draw,draw=none] {$f_{\alpha_{1} + \alpha_{2} + 2 \alpha_{3}} f_{\alpha_{2}} f_{\alpha_{2} + 2 \alpha_{3}}$};
\draw [spinach,->] (node_0) ..controls (227.36bp,269.82bp) and (247.32bp,248.31bp)  .. (node_14);
\draw ($(259.0bp,253.0bp)+(-30bp,-0bp)$) node {$3$};
\draw [blue,->] (node_0) ..controls (180.07bp,269.25bp) and (142.44bp,246.07bp)  .. (node_20);
\draw ($(177.0bp,253.0bp)+(-40bp,5bp)$) node {$1$};
\draw [red,->] (node_1) ..controls (71.743bp,268.9bp) and (82.236bp,249.33bp)  .. (node_20);
\draw ($(93.0bp,253.0bp)+(-30bp,-0bp)$) node {$2$};
\draw [spinach,->] (node_2) ..controls (236.49bp,481.49bp) and (226.17bp,461.45bp)  .. (node_5);
\draw ($(241.0bp,465.0bp)+(-30bp,-0bp)$) node {$3$};
\draw [red,->] (node_2) ..controls (274.49bp,480.58bp) and (310.21bp,457.88bp)  .. (node_18);
\draw ($(320.0bp,465.0bp)+(-40bp,-0bp)$) node {$2$};
\draw [spinach,->] (node_3) ..controls (363.77bp,341.44bp) and (365.2bp,321.69bp)  .. (node_7);
\draw ($(374.0bp,325.0bp)+(-20bp,-0bp)$) node {$3$};
\draw [spinach,->] (node_4) ..controls (122.68bp,128.73bp) and (152.12bp,106.48bp)  .. (node_17);
\draw ($(163.0bp,113.0bp)+(-40bp,-0bp)$) node {$3$};
\draw [red,->] (node_5) ..controls (210.5bp,411.8bp) and (210.5bp,392.61bp)  .. (node_11);
\draw ($(219.0bp,395.0bp)+(-20bp,-0bp)$) node {$2$};
\draw [red,->] (node_6) ..controls (191.5bp,691.8bp) and (191.5bp,672.61bp)  .. (node_12);
\draw ($(200.0bp,675.0bp)+(-20bp,-0bp)$) node {$2$};
\draw [red,->] (node_7) ..controls (342.44bp,269.41bp) and (311.81bp,246.68bp)  .. (node_14);
\draw ($(340.0bp,253.0bp)+(-30bp,5bp)$) node {$2$};
\draw [blue,->] (node_8) ..controls (75.447bp,411.7bp) and (72.011bp,392.23bp)  .. (node_16);
\draw ($(83.0bp,395.0bp)+(-20bp,-0bp)$) node {$1$};
\draw [blue,->] (node_9) ..controls (166.84bp,550.53bp) and (203.13bp,527.68bp)  .. (node_2);
\draw ($(213.0bp,535.0bp)+(-40bp,-5bp)$) node {$1$};
\draw [spinach,->] (node_9) ..controls (134.96bp,551.7bp) and (132.09bp,532.23bp)  .. (node_13);
\draw ($(143.0bp,535.0bp)+(-20bp,-0bp)$) node {$3$};
\draw [red,->] (node_11) ..controls (210.5bp,341.44bp) and (210.5bp,321.69bp)  .. (node_0);
\draw ($(219.0bp,325.0bp)+(-20bp,-0bp)$) node {$2$};
\draw [spinach,->] (node_12) ..controls (177.36bp,621.19bp) and (160.75bp,600.28bp)  .. (node_9);
\draw ($(180.0bp,605.0bp)+(-30bp,-0bp)$) node {$3$};
\draw [blue,->] (node_12) ..controls (205.64bp,621.19bp) and (222.25bp,600.28bp)  .. (node_19);
\draw ($(234.0bp,605.0bp)+(-40bp,5bp)$) node {$1$};
\draw [blue,->] (node_13) ..controls (149.61bp,480.88bp) and (176.22bp,459.09bp)  .. (node_5);
\draw ($(188.0bp,465.0bp)+(-40bp,5bp)$) node {$1$};
\draw [red,->] (node_13) ..controls (114.74bp,481.29bp) and (99.881bp,460.67bp)  .. (node_8);
\draw ($(118.0bp,465.0bp)+(-30bp,-0bp)$) node {$2$};
\draw [blue,->] (node_14) ..controls (225.5bp,198.07bp) and (162.77bp,173.83bp)  .. (node_4);
\draw ($(216.0bp,183.0bp)+(-40bp,10bp)$) node {$1$};
\draw [spinach,->] (node_14) ..controls (275.76bp,199.8bp) and (277.18bp,180.61bp)  .. (node_15);
\draw ($(286.0bp,183.0bp)+(-20bp,-0bp)$) node {$3$};
\draw [blue,->] (node_15) ..controls (255.05bp,128.73bp) and (225.28bp,106.48bp)  .. (node_17);
\draw ($(254.0bp,113.0bp)+(-40bp,5bp)$) node {$1$};
\draw [blue,->] (node_16) ..controls (65.243bp,341.66bp) and (63.856bp,322.51bp)  .. (node_1);
\draw ($(74.0bp,325.0bp)+(-20bp,-0bp)$) node {$1$};
\draw [red,->] (node_17) ..controls (188.5bp,59.799bp) and (188.5bp,40.613bp)  .. (node_10);
\draw ($(197.0bp,43.0bp)+(-20bp,-0bp)$) node {$2$};
\draw [spinach,->] (node_18) ..controls (355.04bp,411.7bp) and (357.91bp,392.23bp)  .. (node_3);
\draw ($(367.0bp,395.0bp)+(-20bp,-0bp)$) node {$3$};
\draw [spinach,->] (node_19) ..controls (245.5bp,551.8bp) and (245.5bp,532.61bp)  .. (node_2);
\draw ($(254.0bp,535.0bp)+(-20bp,-0bp)$) node {$3$};
\draw [spinach,->] (node_20) ..controls (98.5bp,199.8bp) and (98.5bp,180.61bp)  .. (node_4);
\draw ($(107.0bp,183.0bp)+(-20bp,-0bp)$) node {$3$};
\end{tikzpicture}$}
.
\end{gather*}
In this example the three paths $\lambda$, $\mu$ and $\nu$ from \autoref{S:MainExample}
end at the vertices labeled
$f_{\sroot[1]+\sroot[2]}f_{\sroot[2]+2\sroot[3]}$,
$f_{\sroot[1]}f_{\sroot[2]}f_{\sroot[2]+2\sroot[3]}$ and
$f_{\sroot[1]+\sroot[2]+\sroot[3]}f_{\sroot[2]+\sroot[3]}$, respectively.
The generating sequences of $f$ Kashiwara operators are given by the residue sequence of the corresponding path. In general, in the PBW realization of the crystal the vertices are labeled by products $f_{\beta_{1}}\dots f_{\beta_{k}}$ of Kashiwara operators, where $\beta_{1},\dots,\beta_{k}$ are certain sequences of positive roots ordered with respect to a convex order that is determined by a choice of reduced expression for $w_{0}$, the longest element of $W$.
\end{Example}

\begin{Lemma}\label{L:ProofsGraphs}
Suppose that $\Gamma$ is a quiver of finite type, $\Lambda\in P^{+}$ and let $i\in\vertices$. Then $\crystalgraph[\Lambda]$ is path ordered.
\end{Lemma}

\begin{proof}
As in \autoref{E:ProofsPBW},
each vertex $\lambda\in\crystalgraph[\Lambda]$
corresponds to a PBW element in $L(\Lambda)$ by applying the
PBW basis given at the vertex to a highest weight vector. Moreover, every
rooted path $\bpath$ to $\lambda$ corresponds to a sequence $f_{\beta_{1}}\dots f_{\beta_{k}}$ of Kashiwara operators, which gives the associated PBW basis element up to lower order terms.
This depends on a choice of reduced expression for the longest word in the Coxeter group, which exists because we are in finite type, but any choice works in the same way in this proof.

The statement now follows from the PBW theorem as all different PBW elements built from the same sequence of Kashiwara operators are related by at least one nontrivial Serre relation, which involves adjacent residues in $\quiver$. That is, as is explained in many places such as \cite[Lemma 3.3(i)]{Ti-elementary-canonical-basis} or \cite[Chapter 8]{Ja-lectures-qgroups}, if $\bpath$ and $\cpath$ are two paths in $\crystalgraph[\Lambda]$ such that their residue sequences $\Res[\bpath]$ and $\Res[\cpath]$ only differ by successive swaps of $f_{\beta_{j}}$ and $f_{\beta_{k}}$, for nonadjacent $j,k\in\vertices$, then $\bpath$ and $\cpath$ determine the same canonical basis element and, hence, the same vertex of $\crystalgraph[\Lambda]$.
\end{proof}

\begin{Corollary}
The crystals $\crystalgraph[\abfweight]$, $\Crystalgraph[\bfweight]$ and $\Crystalgraph$ are all path ordered.
\end{Corollary}

\begin{proof}
The crystal graph $\crystalgraph[\abfweight]$ is path ordered by \autoref{L:ProofsGraphs}.
As remarked already, $\Crystalgraph[\bfweight]$ and $\Crystalgraph$ are both direct sums of highest weight crystals. Hence, the result again follows from \autoref{L:ProofsGraphs}.
\end{proof}

%%%%%%%%%%%%%%%%%%%%%%%%%%%%%%%%%%%%%%%%%

\subsection{Steady diagrams in level one}\label{SS:ProofsGrothendieckRings}

%%%%%%%%%%%%%%%%%%%%%%%%%%%%%%%%%%%%%%%%%

We now use standard arguments to connect the wKLRW algebras with the quantum groups. In this section we assume that $\ell=1$ and that $\bfweight=(\fweight[j])$, for some $j\in\vertices$.

Recall from \autoref{SS:ConstructionBases} that $\TA$ and $\TAc$ are the infinite and finite dimensional KLR algebras attached to the quiver $\quiver$.

\begin{Proposition}\label{P:ProofsKLRIsomorphism}
As graded algebras,
$\WA[n](X)\cong\TA[n]$ and
$\WAc[n](X)\cong\TAc[n]$.
\end{Proposition}

\begin{proof}
The proof is similar to \cite[Proposition~3F.1]{MaTu-klrw-algebras}, so we only sketch the details. The algebra $\TA$ is generated by elements $\psi_{r}$, $Y_{s}$ and $e_{\bi}$, where $1\leq r<n$, $1\leq s\leq n$ and $\bi\in\vertices^{n}$, so it is enough to define a map $f\colon\TA\to\aWA$ on the generators of $\TA$, show that this map respects the relations of $\TA$ and then define an inverse map. Write $X=\set{\bx}$ and define
\begin{gather}\label{E:KLREmbedding}
f(e_{\bi})=\idem[\bx,\bi],\quad
f(\psi_{r}e_{\bi})=D_{s_{r}}\idem[\bx,\bi],\quad
f(Y_{s}e_{\bi})=y_{s}\idem[\bx,\bi],
\end{gather}
for $\bi\in\vertices^{n}$ and all admissible $r,s$ as above. By definition, the images of the generators of $\TA[n]$ are diagrams in the subalgebra $\aWA$, which is spanned by diagrams that have all of their solid and ghost strings, except for the rightmost ghost string, to the left of the red strings. As discussed in \cite[Proposition~3F.1]{MaTu-klrw-algebras}, this follows using the arguments of \cite[Proposition 6.19]{Bo-many-cellular-structures}. Strictly speaking Bowman is working in type $\aonetypea$ but the argument applies, essentially without change, to quivers of all types. Once it is known that $f$ is an isomorphism it is straightforward to check that it induces an isomorphism between the cyclotomic quotients.
\end{proof}

Let $q$ be an indeterminate over $\Q$. For $a\in\Z$ and $i\in\vertices$ define the \emph{quantum integer} $[a]_{i}=(q_{i}^{a}-q_{i}^{-a})/(q_{i}-q_{i}^{-1})$, where $q_{i}=q^{d_{i}}$. If $a\in\N$ set $[a]_{i}!=[a]_{i}[a-1]_{i}\dots[1]_{i}$. For $a\in\Z$ and $b\in\N$ set $\qbinom[i]{a}{b}=\frac{[a]_{i}[a-1]_{i}\dots[a-b+1]_{i}}{[b]_{i}[b-1]_{i}\dots[1]_{i}}$. By convention, empty products are~$1$, so $[0]_{i}=0=\qbinom[i]{a}{0}$.

Recall from \autoref{SS:RecollectionKacMoody} that $(a_{ij})_{i,j=1}^{e}$ is the Cartan matrix of $\quiver$. Moreover, let
$Q^{+}$ be the positive root lattice associated to $\quiver$.
Following Lusztig \cite[Section 1]{Lu-intro-quantum-groups}, define $\Luf=\bigoplus_{\alpha\in Q^{+}}\Luf_{\alpha}$ to be the $\Q(q)$-algebra generated by $\set{\theta_{i}|i\in\vertices}$ with relations
\begin{gather*}
\sum_{k=0}^{1-a_{ij}}\qbinom[i]{1-a_{ij}}{k}
\theta_{i}^{1-a_{ij}-k}\theta_{j}\theta_{i}^{k}=0,
\qquad\text{for }i,j\in\vertices.
\end{gather*}
Define the divided powers $\theta_{i}^{(k)}=\theta_{i}^{k}/[k]_{i}!$, for $i\in\vertices$ and $k\in\Z_{>0}$.

Using the nil-Hecke algebra, Khovanov--Lauda define a projective $\TA[n]$-module $P_{\bi}$, which is a summand of $\TA[n]e_{\bi}$ for each $\bi\in\vertices^{+}$. Abusing notation, we consider $P_{\bi}$ as a $\aWA$-module using the isomorphism of \autoref{P:ProofsKLRIsomorphism}.

If $A$ is a graded algebra let $[\Proj A]$
be the Grothendieck group of finitely generated graded projective $A$-modules. If $P$ is a projective $A$-module, let $[P]$ be its image in $[\Proj A]$. By extension of scalars, consider $[\Proj A]$ as a $\Q(q)$-module, where $\vpar$ acts as the grading shift functor. We are interested in the cases when $A$ is a wKLRW algebra $\WA(X)$, a KLR algebra $\TA$ or a cyclotomic quotient $\WAc(X)$ and $\TAc$ of one of these algebras.

Let $L(\Lambda_{j})$ be the Weyl $\Luf$-module of highest
weight $\Lambda_{j}$.

Using what are by now standard results for KLR algebras, we can use \autoref{P:ProofsKLRIsomorphism} to give an $\Luf$-action on
$[\Proj\WA[\oplus]]=\bigoplus_{n\geq 0}[\Proj\WA(X)]$ and
$[\Proj\WAc[\oplus]]=\bigoplus_{n\geq 0}[\Proj\WAc(X)]$.
Similarly, set
$[\Proj\TA[\oplus]]=\bigoplus_{n\geq 0}[\Proj\WA(X)]$ and
$[\Proj\TAc[\oplus]]=\bigoplus_{n\geq 0}[\Proj\WAc(X)]$.

\begin{Lemma}\label{L:ProofsKLRIsomorphismTwo}
Suppose that a steady diagram $D$ in $\WAc(X)$ factors through  an idempotent diagram $\idem[\bx,\bi]$. Then $D$ factors through an idempotent diagram that has all of its strings to the left of the red $j$-string.
\end{Lemma}

\begin{proof}
Suppose that $\idem[\bx,\bi]$ has a solid $j$-string $s$ to the right of the red $j$-string. If $s$ is steady then there is a sequence of strings $s_{1},\dots,s_{k}=s$ such that $s_{r+1}$ is blocked by $s_{r}$, for $1\leq r<k$, and $s_{1}$ is blocked by the red $j$-string. Let $s_{r}$ be an $i_{r}$-string. Since $s_{1}$ is blocked by the red string $i_{1}=j$, and $i_{r}\rightsquigarrow i_{r+1}$ since $s_{r+1}$ is blocked by $s_{r}$. Therefore, $j=i_{1}\rightsquigarrow i_{2}\rightsquigarrow\dots\rightsquigarrow i_{k}=j$. However, this is impossible because Dynkin diagrams of finite type do not contain cycles.

Therefore, all of the solid $j$-strings in a steady idempotent diagram are to the left of the red $j$-string. Hence, applying the relations, all of the strings in $\idem[\bx,\bi]$ can be pulled to the left of the red $j$-string, proving the lemma.
\end{proof}

\begin{Proposition}\label{P:ProofsKLRIsomorphismTwo}
Then
there is an injective algebra homomorphism
$\pi\colon\Luf\to[\Proj\WA[\oplus](X)]$ of $\Q(q)$-algebras such that $\pi(\theta_{\bi})=[P_{\bi}]$, for $\bi\in\vertices^{n}$. Moreover, $\Luf$ acts on $[\Proj\WAc(X)]$ and $L(\Lambda_{j})\cong[\Proj\WAc[\oplus]]$ as $\Luf$-modules.
\end{Proposition}

\begin{proof}
By \cite[Theorem 8]{KhLa-cat-quantum-sln-second}, we have $\Luf\cong[{\Proj\TA[\oplus]}]$ as algebras, where the isomorphism sends $\theta_{i}$ to $P_{i}$, for $i\in\vertices$. Hence, the first isomorphism follows in view of \autoref{P:ProofsKLRIsomorphism}. By \cite[Theorem~7.8]{LaVa-crystals-cat} or \cite[Theorem~6.2]{KaKa-categorification-via-klr}, we have $L(\Lambda_{j})\cong[{\Proj\TAc[\oplus]}]$ as $\Luf$-modules. Applying \autoref{P:ProofsKLRIsomorphism} again, the highest weight module $L(\Lambda_{j})$ is isomorphic to $[\Proj\WAc[\oplus]]$. More explicitly
by \autoref{L:ProofsKLRIsomorphismTwo}, the solid strings in a steady idempotent diagram can always be pulled to the left of the red $j$-string. Therefore,  every projective $\WAc(X)$-module can be identified with a projective $\TAc$-module under the isomorphism of \autoref{P:ProofsKLRIsomorphism}.
\end{proof}

%We continue to assume $\bfweight=(\fweight[j])$, for $j\in\vertices$.
\autoref{SS:ConstructionIdempotents} associates an idempotent diagram $\idem[\bpath]$ to each path $\bpath$ in the crystal graph $\crystalgraph[{\fweight[j]}]$. This section shows how face permutations in the crystal are related to permutation diagrams in the wKLRW algebras.

\begin{Lemma}\label{L:ProofsSteady}
Let $\bi\in\vertices^{n}$ and $\bx\in X$. Then $\idem[\bx,\bi]$ is steady if and only if $\bi=\Res[\bpath]$ where $\bpath$ is a rooted path in $\crystalgraph[{\fweight[j]}]$.
\end{Lemma}

\begin{proof}
By \autoref{L:ProofsKLRIsomorphismTwo}, we can assume that all solid strings are to the left of the red $j$-string.
We argue by induction on $n$. If $n=1$, then $\idem[\bx,\bi]$ is steady if and only if $\bi=(j)$ and $x_{1}=\kappa_{1}-\varepsilon$. As the first edge in the crystal graph $\crystalgraph[{\fweight[j]}]$ is labeled~$j$, the result holds in this case. Now assume that $n>1$ and let $\by=(x_{1},\dots,x_{n-1})$ and $\bj=(i_{1},\dots,i_{n-1})$. By assumption, $\idem[\by,\bj]$ is steady so, by induction, $\bj$ is a path in~$\crystalgraph[{\fweight[j]}]$. By the proof of \autoref{P:ProofsKLRIsomorphismTwo}, adding an $i_{n}$-string to the left of the diagram $\idem[\by,\bj]$ is equivalent to tensoring with~$P_{i_{n}}$, which gives an action of~$\theta_{i_{n}}$ on $\idem[\by,\bj]$. By \cite[Theorem~7.5]{LaVa-crystals-cat}, the crystal graph $\crystalgraph[{\fweight[j]}]$ categorifies the simple $\TAc$-modules, for $n\geq0$. Therefore, $\theta_{i_{n}}$ is nonzero on $\idem[\by,\bj]$ if and only if $\bi$ is the residue sequence of a path in $\crystalgraph[{\fweight[j]}]$, giving the result.
\end{proof}

As a special case of face permutations, \autoref{D:FiniteTypesDetour} defines detour permutations to be the face permutations in $\Face(\ppath,\ppath)$, for $\blam\in\CrystalVertices$.
To understand the detour permutations we use the plactic relations of \autoref{L:RecollectionReidemeisterIII} and the plactic monoid moves of \autoref{SS:Plactic}.

\begin{Proposition}\label{P:ProofsDetour}
Let $\bT\in\Face(\lambda)$ be a detour permutation, for $\lambda\in\CrystalVertices$. Then $D_{\bT}$ is invertible up to more dominant diagrams.
That is, there exists a diagram $E_{\bT}$ such that
$E_{\bT}D_{\bT}=\idem[\lambda]+F$, where $F$ is linear combination of diagrams that dominate $\idem[\lambda]$.
\end{Proposition}

\begin{proof}
Let $i,j\in\vertices$ with $i\neq j$. It is enough to consider the case when $\bT$ is a basic face permutation between two paths $\bpath$ and $\cpath$ with sink $\lambda$. To prove the lemma it is enough to show that the idempotent diagrams $\idem[\bpath]$ and $\idem[\cpath]$ factor through each other, up to more dominant diagrams. By \autoref{D:FiniteTypesSSPerm} and \autoref{P:FiniteTypesSSR}, the basic face permutations are of the form
\begin{gather*}
ij\leftrightsquigarrow ji\text{ with $i,j$ nonadjacent},
\\
i\dots ij\leftrightsquigarrow i\dots ji\text{ with $i,j$ simply laced adjacent}
,\quad
ijji\leftrightsquigarrow jiij\text{ with $i,j$ doubly laced adjacent},
\end{gather*}
except in type $\typef[4]$ where there are two additional cases, which we consider at the end of the proof.
We consider each of the basic face permutations in turn.

\Case{$ij\leftrightsquigarrow ji$, with $i$ and $j$ nonadjacent.} The claim is easy to verify in this case. On the one hand,
if $i,j$ are
not fishtail vertices, then $\idem[\bpath]=\idem[\cpath]$ on the nose
by our string placement strategy in \autoref{D:ConstructionIdempotents}. On the other hand, if $i,j$ are
fishtail vertices, then
\begin{gather*}
\begin{tikzpicture}[anchorbase,smallnodes,rounded corners]
\draw[solid](-0.5,1)--++(0,-1)node[below]{$i$};
\draw[solid](0,1)node[above,yshift=-1pt]{$\phantom{i}$}--++(0,-1)node[below]{$j$};
\end{tikzpicture}
=
\begin{tikzpicture}[anchorbase,smallnodes,rounded corners]
\draw[solid](-0.5,1)--++(0.5,-0.5)--++(-0.5,-0.5)node[below]{$i$};
\draw[solid](0,1)node[above,yshift=-1pt]{$\phantom{i}$}--++(-0.5,-0.5)--++(0.5,-0.5)node[below]{$j$};
\end{tikzpicture}
\end{gather*}
and its partner relation hold
by the honest Reidemeister II relation. Hence, this local
picture verifies the claim.

\Case{$ij\leftrightsquigarrow ji$, with $i$ and $j$ adjacent.} This is the only place where we use \autoref{P:PlacticCrystals}, which says that the crystal $\crystalgraph[{\fweight[j]}]$ is plactic. First assume that the local situation
is $iij\leftrightsquigarrow iji$, which corresponds to part~(a) of \autoref{D:ConstructionPlactic}. Assume first that $i\rightarrow j$,
we can use \autoref{L:RecollectionReidemeisterIII} which reads
\begin{gather*}
\begin{tikzpicture}[anchorbase,smallnodes,rounded corners]
\draw[ghost,spinach](1,0)--++(0,1)node[above,yshift=-1pt]{$i$};
\draw[ghost](1.5,0)--++(0,1)node[above,yshift=-1pt]{$i$};
\draw[solid](1.25,0)node[below]{$j$}--++(0,1);
\end{tikzpicture}
=
-
\begin{tikzpicture}[anchorbase,smallnodes,rounded corners]
\draw[ghost](1,0)--++(0.75,0.5)--++(-0.75,0.5)node[above,yshift=-1pt]{$i$};
\draw[ghost,dot](1.5,0)--++(0,1)node[above,yshift=-1pt]{$i$};
\draw[solid](1.25,0)node[below]{$j$}--++(0,1);
\end{tikzpicture}
\underbrace{-\begin{tikzpicture}[anchorbase,smallnodes,rounded corners]
\draw[ghost,dot](1,0)--++(0.75,0.5)--++(-0.75,0.5)node[above,yshift=-1pt]{$i$};
\draw[ghost](1.5,0)--++(0,1)node[above,yshift=-1pt]{$i$};
\draw[solid](1.25,0)node[below]{$j$}--++(0.75,0.5)--++(-0.75,0.5);
\end{tikzpicture}}_{\text{more dominant}}
.
\end{gather*}
Note that this relates $iji$ and $iij$ modulo $jii$, and $jii$ is always more dominant. A similar argument works for its partner relation and for the case when $i\leftarrow j$. This verifies the claim when $iij\leftrightsquigarrow iji$.

Now consider the case when $ij\leftrightsquigarrow ji$ is not preceded by $i$ or $j$, which corresponds to \autoref{D:ConstructionPlactic}(b). By \autoref{L:ConstructionPlacticMoves}, we can use weighted plactic monoid moves to reduce to the case when $ij\leftrightsquigarrow ji$ is not preceded by $i$ or $j$.
Thus, it remains to show that the weighted plactic moves are realized
by invertible operations up to more dominant diagrams. For the first case
we can simply use an honest Reidemeister II relation \autoref{L:RecollectionMovingStringsDots}.(b), {\eg}
\begin{gather*}
\begin{tikzpicture}[anchorbase,smallnodes,rounded corners]
\draw[ghost](1,1)node[above,yshift=-1pt]{$i$}--++(0,-1)node[below]{$\phantom{i}$};
\draw[solid,smallnodes,rounded corners](1.5,1)--++(0,-1)node[below]{$j$};
\end{tikzpicture}
=
\begin{tikzpicture}[anchorbase,smallnodes,rounded corners]
\draw[ghost](1,1)node[above,yshift=-1pt]{$i$}--++(0.5,-0.5)--++(-0.5,-0.5)node[below]{$\phantom{i}$};
\draw[solid,smallnodes,rounded corners](1.5,1)--++(-0.5,-0.5)--++(0.5,-0.5)node[below]{$j$};
\end{tikzpicture}
.
\end{gather*}
At the equator of these diagrams is the sequence $ji$, showing that the corresponding weighted plactic move is realized
by an invertible operation. Similarly, for the other case applying
a Reidemeister III relation \autoref{Eq:RecollectionReidemeisterIII}, and looking at the equators proves the claim.

\Case{$ijji\leftrightsquigarrow jiij$ and $i\rightarrow j$.}
By \autoref{L:RecollectionMovingStringsDots} and \autoref{L:RecollectionMovingStringsDots},
\begin{gather*}
\begin{tikzpicture}[anchorbase,smallnodes,rounded corners]
\draw[ghost](1,0)--++(0,1)node[above,yshift=-1pt]{$j$};
\draw[ghost](1.75,0)--++(0,1)node[above,yshift=-1pt]{$j$};
\draw[solid](1.25,0)node[below]{$i$}--++(0,1);
\draw[solid](1.5,0)node[below]{$i$}--++(0,1);
\end{tikzpicture}
=
-
\begin{tikzpicture}[anchorbase,smallnodes,rounded corners]
\draw[ghost](1,0)--++(0,1)node[above,yshift=-1pt]{$j$};
\draw[ghost](1.75,0)--++(0,1)node[above,yshift=-1pt]{$j$};
\draw[solid,dot](1.25,0)node[below]{$i$}--++(0.25,0.5)--++(-0.25,0.5);
\draw[solid,dot=0.925](1.5,0)node[below]{$i$}--++(-0.25,0.5)--++(0.25,0.5);
\end{tikzpicture}
-
\begin{tikzpicture}[anchorbase,smallnodes,rounded corners]
\draw[ghost](1,0)--++(0,1)node[above,yshift=-1pt]{$j$};
\draw[ghost](1.75,0)--++(0,1)node[above,yshift=-1pt]{$j$};
\draw[solid,dot=0.075,dot](1.25,0)node[below]{$i$}--++(0.25,0.5)--++(-0.25,0.5);
\draw[solid](1.5,0)node[below]{$i$}--++(-0.25,0.5)--++(0.25,0.5);
\end{tikzpicture}
=
-
\begin{tikzpicture}[anchorbase,smallnodes,rounded corners]
\draw[ghost](1,0)--++(0,1)node[above,yshift=-1pt]{$j$};
\draw[ghost](1.75,0)--++(-0.25,0.5)--++(0.25,0.5)node[above,yshift=-1pt]{$j$};
\draw[solid](1.25,0)node[below]{$i$}--++(0.5,0.5)--++(-0.5,0.5);
\draw[solid,dot=0.925](1.5,0)node[below]{$i$}--++(-0.25,0.5)--++(0.25,0.5);
\end{tikzpicture}
-
\begin{tikzpicture}[anchorbase,smallnodes,rounded corners]
\draw[ghost](1,0)--++(0,1)node[above,yshift=-1pt]{$j$};
\draw[ghost](1.75,0)--++(-0.25,0.5)--++(0.25,0.5)node[above,yshift=-1pt]{$j$};
\draw[solid,dot=0.075](1.25,0)node[below]{$i$}--++(0.5,0.5)--++(-0.5,0.5);
\draw[solid](1.5,0)node[below]{$i$}--++(-0.25,0.5)--++(0.25,0.5);
\end{tikzpicture}
\underbrace{-
\begin{tikzpicture}[anchorbase,smallnodes,rounded corners]
\draw[ghost](1,0)--++(0,1)node[above,yshift=-1pt]{$j$};
\draw[ghost,dot](1.75,0)--++(0,1)node[above,yshift=-1pt]{$j$};
\draw[solid](1.25,0)node[below]{$i$}--++(0.25,0.5)--++(-0.25,0.5);
\draw[solid,dot=0.925](1.5,0)node[below]{$i$}--++(-0.25,0.5)--++(0.25,0.5);
\end{tikzpicture}
-
\begin{tikzpicture}[anchorbase,smallnodes,rounded corners]
\draw[ghost](1,0)--++(0,1)node[above,yshift=-1pt]{$j$};
\draw[ghost,dot](1.75,0)--++(0,1)node[above,yshift=-1pt]{$j$};
\draw[solid,dot=0.075](1.25,0)node[below]{$i$}--++(0.25,0.5)--++(-0.25,0.5);
\draw[solid](1.5,0)node[below]{$i$}--++(-0.25,0.5)--++(0.25,0.5);
\end{tikzpicture}}_{\text{more dominant}}
,
\end{gather*}
or its partner relation. Hence, we have related $jiij$ and $jiji$ up to higher order terms.  We can now use the plactic relations
from \autoref{L:RecollectionReidemeisterIII} to relate $jiji$ to $ijji$ and $jjii$,  with $jjii$ being a more dominant term. The argument when $j\rightarrow i$ is identical.

\Case{$ijji\leftrightsquigarrow jiij$ and $i\Rightarrow j$.}
If $i$ and $j$ are doubly laced adjacent then essentially the same argument
applies except that we make use of the fact that
the rightmost of the two middle strings
carries a dot by \autoref{SS:ConstructionDots}. More precisely, as above, we can
relate $jiij$ to $jiji$ up to more dominant terms since
\autoref{L:RecollectionMovingStringsDots}.(b) gives
\begin{gather*}
\begin{tikzpicture}[anchorbase,smallnodes,rounded corners]
\draw[ghost,dot,spinach](0,1)node[above,yshift=-1pt]{$i$}--++(0,-1);
\draw[solid](0.5,1)--++(0,-1)node[below]{$j$};
\end{tikzpicture}
=
\begin{tikzpicture}[anchorbase,smallnodes,rounded corners]
\draw[ghost](0,1)node[above,yshift=-1pt]{$i$}--++(0.8,-0.5)--++(-0.8,-0.5);
\draw[solid](0.5,1)--++(0,-1) node[below]{$j$};
\end{tikzpicture}
\underbrace{+
\begin{tikzpicture}[anchorbase,smallnodes,rounded corners]
\draw[ghost](0,1)node[above,yshift=-1pt]{$i$}--++(0,-1);
\draw[solid,dot=0.4,dot=0.6](0.5,1)--++(0,-1)node[below]{$j$};
\end{tikzpicture}}_{\text{more dominant}}
,
\end{gather*}
and the same diagrammatic argument used above works (the more dominant terms
will have two dots, but this does not affect the argument).
Note that in the diagram the strings with residues $jij$ are such that the $i$-string has a dot, so using \autoref{L:RecollectionMovingStringsDots}.(b)
we can relate this to $ijj$ up to higher order terms. Hence, we have
related $jiij$ to $ijji$ up to higher order terms.

\Case{the decagon in type $\typef[4]$.} We consider only
the case when the basic diagrammatic face permutation is:
\begin{gather*}
\begin{tikzpicture}[anchorbase,smallnodes,rounded corners]
\draw[solid,black](0,0)node[below]{$2$}--++(0,1)--++(0.5,0.5) node[above,yshift=-1pt]{$2$};
\draw[solid,black](2,0)node[below]{$2$}--++(0,1)--++(-0.5,0.5)node[above,yshift=-1pt]{$2$};
\draw[solid,spinach](0.5,0)--++(1,1)--++(0.5,0.5)node[above,yshift=-1pt]{$3$};
\draw[solid,spinach](1,0)--++(-0.5,0.5)--++(0.5,0.5)--++(0,0.5)node[above,yshift=-1pt]{$3$};
\draw[solid,spinach](1.5,0)--++(-1,1)--++(-0.5,0.5)node[above,yshift=-1pt]{$3$};
\end{tikzpicture}
.
\end{gather*}
From \autoref{L:RecollectionMovingStringsDots}.(a) and the Reidemeister II relation we get
\begin{align*}
\begin{tikzpicture}[anchorbase,smallnodes,rounded corners]
\draw[ghost](0,1)node[above,yshift=-1pt]{$3$}--++(0,-1);
\draw[ghost](0.5,1)node[above,yshift=-1pt]{$3$}--++(0,-1);
\draw[ghost](1,1)node[above,yshift=-1pt]{$3$}--++(0,-1);
\draw[solid](-0.5,1)--++(0,-1)node[below]{$2$};
\draw[solid](1.5,1)--++(0,-1)node[below]{$2$};
\end{tikzpicture}
&=
\begin{tikzpicture}[anchorbase,smallnodes,rounded corners]
\draw[ghost,dot=0.45,dot=0.55](0,1)node[above,yshift=-1pt]{$3$}--++(0,-0.2)--++(1.8,-0.15)--++(0,-0.3)--++(-1.8,-0.15)--++(0,-0.2);
\draw[ghost,dot=0.1](0.5,1)node[above,yshift=-1pt]{$3$}--++(0,-1);
\draw[ghost,dot=0.1](1,1)node[above,yshift=-1pt]{$3$}--++(0,-1);
\draw[solid](-0.5,1)--++(0,-1)node[below]{$2$};
\draw[solid](1.5,1)--++(0,-1)node[below]{$2$};
\end{tikzpicture}
-
\begin{tikzpicture}[anchorbase,smallnodes,rounded corners]
\draw[ghost,dot=0.45,dot=0.55,dot=0.775](0,1)node[above,yshift=-1pt]{$3$}--++(0,-0.2)--++(1.8,-0.15)--++(0,-0.3)--++(-1.8,-0.15)--++(0,-0.2);
\draw[ghost,dot=0.1](0.5,1)node[above,yshift=-1pt]{$3$}--++(0,-1);
\draw[ghost](1,1)node[above,yshift=-1pt]{$3$}--++(0,-1);
\draw[solid](-0.5,1)--++(0,-1)node[below]{$2$};
\draw[solid](1.5,1)--++(0,-1)node[below]{$2$};
\end{tikzpicture}
-
\begin{tikzpicture}[anchorbase,smallnodes,rounded corners]
\draw[ghost,dot=0.45,dot=0.55,dot=0.9](0,1)node[above,yshift=-1pt]{$3$}--++(0,-0.2)--++(1.8,-0.15)--++(0,-0.3)--++(-1.8,-0.15)--++(0,-0.2);
\draw[ghost](0.5,1)node[above,yshift=-1pt]{$3$}--++(0,-1);
\draw[ghost,dot=0.1](1,1)node[above,yshift=-1pt]{$3$}--++(0,-1);
\draw[solid](-0.5,1)--++(0,-1)node[below]{$2$};
\draw[solid](1.5,1)--++(0,-1)node[below]{$2$};
\end{tikzpicture}
+
\begin{tikzpicture}[anchorbase,smallnodes,rounded corners]
\draw[ghost,dot=0.45,dot=0.55,dot=0.775,dot=0.9](0,1)node[above,yshift=-1pt]{$3$}--++(0,-0.2)--++(1.8,-0.15)--++(0,-0.3)--++(-1.8,-0.15)--++(0,-0.2);
\draw[ghost](0.5,1)node[above,yshift=-1pt]{$3$}--++(0,-1);
\draw[ghost](1,1)node[above,yshift=-1pt]{$3$}--++(0,-1);
\draw[solid](-0.5,1)--++(0,-1)node[below]{$2$};
\draw[solid](1.5,1)--++(0,-1)node[below]{$2$};
\end{tikzpicture}
\\
&=
\begin{tikzpicture}[anchorbase,smallnodes,rounded corners]
\draw[ghost,dot=0.45,dot=0.55](0,1)node[above,yshift=-1pt]{$3$}--++(0,-0.2)--++(1.8,-0.15)--++(0,-0.3)--++(-1.8,-0.15)--++(0,-0.2);
\draw[ghost,dot=0.1](0.5,1)node[above,yshift=-1pt]{$3$}--++(0,-1);
\draw[ghost,dot=0.1](1,1)node[above,yshift=-1pt]{$3$}--++(0,-1);
\draw[solid](-0.5,1)--++(0,-1)node[below]{$2$};
\draw[solid](1.5,1)--++(0,-1)node[below]{$2$};
\end{tikzpicture}
+
\text{more dominant terms,}
.
\end{align*}
Now \autoref{L:RecollectionMovingStringsDots}.(a) implies that
\begin{gather*}
\begin{tikzpicture}[anchorbase,smallnodes,rounded corners]
\draw[ghost,dot=0.45,dot=0.55](0,1)node[above,yshift=-1pt]{$3$}--++(0,-0.2)--++(1.8,-0.15)--++(0,-0.3)--++(-1.8,-0.15)--++(0,-0.2);
\draw[ghost,dot=0.1](0.5,1)node[above,yshift=-1pt]{$3$}--++(0,-1);
\draw[ghost,dot=0.1](1,1)node[above,yshift=-1pt]{$3$}--++(0,-1);
\draw[solid](-0.5,1)--++(0,-1)node[below]{$2$};
\draw[solid](1.5,1)--++(0,-1)node[below]{$2$};
\end{tikzpicture}
+
\text{more dominant terms}
=
\begin{tikzpicture}[anchorbase,smallnodes,rounded corners]
\draw[ghost,dot=0.45,dot=0.55](0,1)node[above,yshift=-1pt]{$3$}--++(0,-0.2)--++(1.8,-0.15)--++(0,-0.3)--++(-1.8,-0.15)--++(0,-0.2);
\draw[ghost,dot=0.1](0.5,1)node[above,yshift=-1pt]{$3$}--++(0,-1);
\draw[ghost](1,1)node[above,yshift=-1pt]{$3$}--++(0,-0.35)--++(-1.8,-0.05)--++(0,-0.2)--++(+1.8,-0.05)--++(0,-0.35);
\draw[solid](-0.5,1)--++(0,-1)node[below]{$2$};
\draw[solid](1.5,1)--++(0,-1)node[below]{$2$};
\end{tikzpicture}
+
\text{more dominant terms},
\end{gather*}
which is what we needed to establish, up to symmetry. The argument when $j\Rightarrow i$ is identical.

\Case{The tetradecagon in type $\typef[4]$.} Similar to
the previous case, we use \autoref{L:RecollectionMovingStringsDots}.(a) and the Reidemeister II relation on close strings of the same residue. The details are omitted.
\end{proof}

%%%%%%%%%%%%%%%%%%%%%%%%%%%%%%%%%%%%%%%%%

\subsection{Steady diagrams and crystals}\label{SS:ProofsSteadyDiagrams}

%%%%%%%%%%%%%%%%%%%%%%%%%%%%%%%%%%%%%%%%%

We now return to the case where $\bfweight$ is of arbitrary level $\ell\ge1$.
Recall from \autoref{SS:TensorProducts} that we write each rooted path $\bpath\in\Parts{n}{\abfweight}$ as a tuple $\bpath=(\bpath_{1},\dots,\bpath_{\hell})$, where $\bpath_{m}$ is a rooted path in $\crystalgraph[\Lambda_{\affine{\rho}_{m}}]$, for $1\le m\le\hell$.

\begin{Lemma}\label{C:PathAnchors}
Suppose that $\bpath=(\bpath_{1},\dots,\bpath_{\hell})\in\Parts{n}{\abfweight}$. Then the solid strings in $\idem[\bpath]$ corresponding to $\bpath_{m}$ are anchored around the $m$th (affine) red string, for $1\leq m\leq\hell$. Moreover, $1_{\bpath}$ is steady if and only if $\bpath_{m}$ is empty for $\ell<m\leq\hell$.
\end{Lemma}

\begin{proof}
By \autoref{L:RootedPaths}, $\bpath_{m}$ is a rooted path in $\crystalgraph[{\fweight[\rho_{m}]}]$. By \autoref{P:ProofsKLRIsomorphism}, $\WAc\cong\TAc$ is isomorphic to the wKLRW algebra obtained by putting $n$ solid strings to the left of the $m$th red string using an appropriate modification of \autoref{E:KLREmbedding}. In view of \autoref{L:ConstructionParkingRed}, we can apply \autoref{L:ProofsSteady} locally to the strings anchored on each (affine) red string.
Hence, the strings in $\idem[\bpath]$ that correspond to the edges of the path $\bpath_{m}$ in $\crystalgraph[\affine{\rho}_{m}]$ are anchored on the $m$th (affine) red string by \autoref{L:ProofsSteady}. Finally, to show that $1_{\bpath}$ is steady if and only if $\bpath_{m}$ is empty for $\ell<m\leq\hell$ first note that $1_{\bpath}$ is not steady if $\bpath_{m}$ is not empty for any $m>\ell$. Conversely, if $\bpath_{m}$ is empty for $\ell<m\leq\hell$ then all of the strings anchored on the $k$th red string, for $1\le k\le\ell$, are steady by \autoref{L:ProofsSteady}, so no string in $\idem[\bpath]$ can be pulled arbitrarily far to the right. Hence, $\idem[\bpath]$ is steady.
\end{proof}

A diagram $D$ \emph{factors through} a diagram $B$ if $D=D^{\prime}BD^{\prime\prime}$ for some diagrams $D^{\prime}$ and $D^{\prime\prime}$. The diagram $D$ \emph{dominates} $\idem[\blam]$ if $D$ factors through an idempotent diagram $\idem[\bmu]$ such that $\blam\lsandorder\bmu$. We do not require that $B$ belongs to the wKLRW algebra $\WAc(X)$, only that it is a wKLRW diagram.

\begin{Proposition}\label{P:PathIdempotents}
Let $\bi\in\vertices^{n}$. Then $\idem[\bx,\bi]$ factors through $\idem[\bpath]$ for some path $\bpath$ in $\crystalgraph[\abfweight]$ such that $\idem[\bx,\bi]\lsandorder\idem[\bpath]$. Moreover, $\idem[\bx,\bi]$ is steady if and only if $\bpath\in\crystalgraph[\bfweight]$.
\end{Proposition}

\begin{proof}
By pulling strings to the right, we can assume that every string is anchored on a red string or an affine red string. Therefore, $\idem[\bx,\bi]$ factors through a diagram $\idem[\bj]$, where $\bj=(\bj_{1},\dots,\bj_{\hell})\in\vertices^{n,\hell}$ and $\bj_{m}$ correspond to the strings anchored on the (affine) red $\affine{\rho}_{m}$-string with coordinate $\affine{\kappa}_{m}$. By construction, $\idem[\bx,\bi]\lsandorder\idem[\bj]$. By \autoref{L:ProofsSteady}, $\bj_{m}=\Res[\bpath_{m}]$ is the residue sequence of a rooted path $\bpath_{m}$ in $\crystalgraph[{\fweight[\rho_{m}]}]$. Therefore, by \autoref{L:TensorFactors}, we can find a rooted path $\bpath\in\crystalgraph[\bfweight]$ with $\bpath=(\bpath_{1},\dots,\bpath_{\hell})$. Finally, $\idem[\bx,\bi]$ is steady if and only if $\idem[\bpath]$ is steady, which is if and only if $\bpath_{m}$ is empty for $\ell<m\leq\hell$ by \autoref{C:PathAnchors}. Hence, $\idem[\bx,\bi]$ is steady if and only if $\bpath\in\crystalgraph[\bfweight]$ by \autoref{L:TensorFactors}.
\end{proof}

\begin{Remark}
  In fact, the argument of \autoref{P:PathIdempotents} shows that if $\Gamma$ is a symmetrizable quiver, then an idempotent diagram $\idem[\bx,\bi]$ in the wKLRW algebra $\WA(X)$ is steady if and only if $\bi$ is the residue sequence of a rooted path in the quiver $\crystalgraph[\bfweight]$.
\end{Remark}

Let $\bpath\in\vertices^{n,\ell}$. Recall from \autoref{D:ConstructionRepeated} that $d_{m}(\bpath)+1$ is the length of a dotted $i$-sequence in $\Res[\bpath]$ and $c_{m}(\bpath)=d_{i}-d_{m}-1$ is the maximal number of finite sandwich dots by \autoref{D:ConstructionSanddots}.

\begin{Corollary}\label{C:ConstructionRepeated}
Let $\bpath$ be a rooted path in the crystal graph $\crystalgraph[\bfweight]$. Then
\begin{gather*}
d_{m}(\bpath)\leq
\begin{cases*}
2 & in types $\typef[4]$ and $\typeg[2]$,
\\
1 & in types $\typeb[e>1]$, $\typec[e>2]$, $\typed[e>3]$ and $\typee[k]$,
\\
0 & in type $\typea[e]$.
\end{cases*}
\end{gather*}
In particular, $\dotidem[\bpath]=\idem[\bpath]$ in type $\typea[e]$.
\end{Corollary}

\begin{proof}
Since every string in $\idem[\bpath]$ is anchored on some (affine) red string, it is enough to consider rooted paths in $\crystalgraph$, for $i\in\vertices$.
By definition, $d_{m}(\bpath)$ is determined by the crystal graph of $\crystalgraph[{\fweight[i]}]$.
The Serre relations show that $F_{j}^{1-a_{jk}}$ annihilates $L(\fweight[i])$, for $j,k\in\vertices$. Therefore, no path in $\crystalgraph[{\fweight[i]}]$ has a dotted $j$-subsequence of length $1-a_{jk}$, so if $i_{m}=j$, then $d_{m}(\bpath)\leq\max\set{{|a_{jk}|}|k\in\vertices}$. This gives all of the bounds above except for types $\typea[e]$, $\typeb[e>1]$, $\typec[e>2]$ and $\typeg[2]$, where we can do better. In the classical cases, these bounds can be verified using tableaux combinatorics, as in the proof of \autoref{P:FiniteTypesSSR}. For type $\typeg[2]$, the bound follows directly from \eqref{Eq:FiniteTypesG2}.
\end{proof}

%\begin{Lemma}\label{L:ProofsIdeal}
%Let $\blam\in\aCrystalVertices$. Then $\WA[n](X)^{\rsandorder\blam}$ and $%%%\WA[n](X)^{\gsandorder\blam}$ are two-sided ideals of $\aWA$.
%\end{Lemma}

%\begin{proof}
%Let $D\in\WA[n](X)^{\gsandorder\blam}$. By definition, $EDF\in\WA[n](X)^{\gsandorder\blam}$
%for all $E,F\in\WA[n](X)$. To complete the proof it is enough to check that $\WA[n](X)^{\gsandorder\blam}$ is invariant under the relations of $\WA[n](X)$, which follows by inspecting the relations in \autoref{SS:RecollectionRelations} because in every relation the endpoints of the diagrams are unchanged. The same argument applies to $\WA[n](X)^{\rsandorder\blam}$, so the statement follows.
%\end{proof}

Recall that \autoref{SS:ConstructionPermutation} associates a permutation diagram $D_{\bT}$ to each face permutation $\bT\in\Face(\bpath,\cpath)$, for $\bpath,\cpath\in\Parts{n}{\abfweight}$.
By \autoref{D:FiniteTypesSSPerm}, a face permutation $\bT\in\Face(\bpath,\cpath)$ is an equivalence class of permutations, and \autoref{D:ConstructionPermBasis} defines the diagrammatic face permutation $D_{\bT}$ for~$\bT$. %Recall from \autoref{D:ConstructionPermBasis} that each permutation $w\in\sym$ defines a diagram $D_{w}=\idem[\bpath]D_{w}\idem[\lambda]$.

\begin{Lemma}\label{L:ProofsSpanningPermutations}
Let $\bT\in\Face(\bpath,\cpath)$, for $\bpath,\cpath\in\Parts{n}{\abfweight}$. If $w\in\bT$, then
$D_{w}$ can be obtained from $D_{\bT}$, up to more dominant terms, by applying crossings and dots.
\end{Lemma}

\begin{proof}
As all of the relations in \autoref{D:hell} are bilocal, we can ignore the red strings in the relations below. By definition, the basic diagrammatic face permutations permute the minimal number of strings with different residues and the maximal number of strings with the same residue. Thus, if the permutation diagram $D_{w}$ has more $i$-$j$ crossings, then we can apply these crossings to $D_{\bT}$ to give all of the $i$-$j$ crossings in $D_{w}$. On the other hand, we can remove $i$-$i$ crossings from $D_{\bT}$ by successively applying
relation \autoref{Eq:RecollectionDotCrossing}, which gives
\begin{gather}\label{Eq:ProofsSpanningPermutations}
\begin{tikzpicture}[anchorbase,smallnodes,rounded corners]
\draw[solid](0.5,0.5)node[above,yshift=-1pt]{$\phantom{i}$}--(0,0) node[below]{$i$};
\draw[solid,dot=0.25](0,0.5)--(0.5,0) node[below]{$i$};
\end{tikzpicture}
=
\begin{tikzpicture}[anchorbase,smallnodes,rounded corners]
\draw[solid](0,0.5)node[above,yshift=-1pt]{$\phantom{i}$}--(0,0) node[below]{$i$};
\draw[solid](0.5,0.5)--(0.5,0) node[below]{$i$};
\end{tikzpicture}
+
\underbrace{\begin{tikzpicture}[anchorbase,smallnodes,rounded corners]
\draw[solid](0.5,0.5)node[above,yshift=-1pt]{$\phantom{i}$}--(0,0) node[below]{$i$};
\draw[solid,dot=0.75](0,0.5)--(0.5,0) node[below]{$i$};
\end{tikzpicture}}_{\text{more dominant}}
\quad\text{and}\quad
\begin{tikzpicture}[anchorbase,smallnodes,rounded corners,xscale=-1]
\draw[solid](0.5,0.5)node[above,yshift=-1pt]{$\phantom{i}$}--(0,0) node[below]{$i$};
\draw[solid,dot=0.75](0,0.5)--(0.5,0) node[below]{$i$};
\end{tikzpicture}
=
\begin{tikzpicture}[anchorbase,smallnodes,rounded corners]
\draw[solid](0,0.5)node[above,yshift=-1pt]{$\phantom{i}$}--(0,0) node[below]{$i$};
\draw[solid](0.5,0.5)--(0.5,0) node[below]{$i$};
\end{tikzpicture}
+
\underbrace{\begin{tikzpicture}[anchorbase,smallnodes,rounded corners,xscale=-1]
\draw[solid](0.5,0.5)node[above,yshift=-1pt]{$\phantom{i}$}--(0,0) node[below]{$i$};
\draw[solid,dot=0.25](0,0.5)--(0.5,0) node[below]{$i$};
\end{tikzpicture}}_{\text{more dominant}}
.
\end{gather}
That the rightmost diagrams in these equations are more dominant follows by
pulling strings and jumping dots to the right, as in \autoref{L:RecollectionMovingStringsDots}.
For example, in the simplified notation as in \autoref{D:ConstructionPermBasis}:
\begin{gather*}
\begin{tikzpicture}[anchorbase,smallnodes]
\draw[solid,black](0,0)node[below]{$i$}--++(1,0.5) node[above,yshift=-1pt]{$i$};
\draw[solid,black,dot=0.9](1.5,0)node[below]{$i$}--++(-1,0.5)node[above,yshift=-1pt]{$i$};
\draw[solid,spinach,dot=0.1](0.5,0)node[below]{$j$}--++(1,0.5)node[above,yshift=-1pt]{$j$};
\draw[solid,spinach](1,0)node[below]{$j$}--++(-1,0.5)node[above,yshift=-1pt]{$j$};
\end{tikzpicture}
=
\begin{tikzpicture}[anchorbase,smallnodes]
\draw[solid,black](0,0)node[below]{$i$}--++(0.5,0.5) node[above,yshift=-1pt]{$i$};
\draw[solid,black](1.5,0)node[below]{$i$}--++(-0.5,0.5)node[above,yshift=-1pt]{$i$};
\draw[solid,spinach](0.5,0)node[below]{$j$}--++(-0.5,0.5)node[above,yshift=-1pt]{$j$};
\draw[solid,spinach](1,0)node[below]{$j$}--++(0.5,0.5)node[above,yshift=-1pt]{$j$};
\end{tikzpicture}
+
\underbrace{\begin{tikzpicture}[anchorbase,smallnodes]
\draw[solid,black](0,0)node[below]{$i$}--++(0.5,0.5) node[above,yshift=-1pt]{$i$};
\draw[solid,black](1.5,0)node[below]{$i$}--++(-0.5,0.5)node[above,yshift=-1pt]{$i$};
\draw[solid,spinach,dot=0.9](0.5,0)node[below]{$j$}--++(1,0.5)node[above,yshift=-1pt]{$j$};
\draw[solid,spinach](1,0)node[below]{$j$}--++(-1,0.5)node[above,yshift=-1pt]{$j$};
\end{tikzpicture}
+
\begin{tikzpicture}[anchorbase,smallnodes]
\draw[solid,black](0,0)node[below]{$i$}--++(1,0.5) node[above,yshift=-1pt]{$i$};
\draw[solid,black,dot=0.1](1.5,0)node[below]{$i$}--++(-1,0.5)node[above,yshift=-1pt]{$i$};
\draw[solid,spinach](0.5,0)node[below]{$j$}--++(-0.5,0.5)node[above,yshift=-1pt]{$j$};
\draw[solid,spinach](1,0)node[below]{$j$}--++(0.5,0.5)node[above,yshift=-1pt]{$j$};
\end{tikzpicture}
+
\begin{tikzpicture}[anchorbase,smallnodes]
\draw[solid,black](0,0)node[below]{$i$}--++(1,0.5) node[above,yshift=-1pt]{$i$};
\draw[solid,black,dot=0.1](1.5,0)node[below]{$i$}--++(-1,0.5)node[above,yshift=-1pt]{$i$};
\draw[solid,spinach,dot=0.9](0.5,0)node[below]{$j$}--++(1,0.5)node[above,yshift=-1pt]{$j$};
\draw[solid,spinach](1,0)node[below]{$j$}--++(-1,0.5)node[above,yshift=-1pt]{$j$};
\end{tikzpicture}}_{\text{more dominant}}
.
\end{gather*}
So, we get the diagram with only two crossings from the left-hand one by applying dots.
\end{proof}
\begin{Proposition}\label{P:ProofsPermutationsSteady}
Let $D$ be a steady permutation diagram in $\WA[n](X)$. Then $D$ is a linear combination of diagrams of the form $E_{\bpath}^{\star}D_{\bT}E_{\cpath}$, for $\bpath,\cpath\in\Parts{n}{\abfweight}$ and face permutations $\bT\in\Face(\bpath,\cpath)$.
\end{Proposition}

\begin{proof}
Any diagram can be factored through an idempotent diagram by making a horizontal cut through the diagram. By \autoref{P:PathIdempotents}, we can pull the strings of $D$ to the right until each string becomes blocked. This gives paths $\bpath,\cpath\in\Parts{n}{\abfweight}$ such that $D$ is a linear combination of diagrams of the form $E_{\bpath}^{\star}D^{\prime}E_{\cpath}$, for some diagram $D^{\prime}$. In particular, the residues along the top and bottom of a nonzero diagram uniquely determine paths $\bpath$ and $\cpath$ in $\Parts{n}{\abfweight}$ by \autoref{P:PathIdempotents}. Each diagram $D^{\prime}$ determines a permutation $w=w_{D^{\prime}}\in\sym$ such that $\Res[\bpath]=w\Res[\cpath]$. By \autoref{L:FiniteTypesPermutations} and \autoref{P:FiniteTypesPermutations}, $w=\bT$ is a face permutation in $\Face(\bpath,\cpath)$, so~$D^{\prime}=D_{\bT}$ in view of \autoref{L:ProofsSpanningPermutations}. (As noted after \autoref{D:ConstructionPermBasis}, the definition of diagrammatic face permutations involves several choices, but these give the same diagrams modulo more dominant terms.)
%As the diagrammatic face permutation $D_{\bT}$ involves a choice up to more dominant terms, $D^{\prime}\equiv D_{\bT}\pmod{\WA[n](X)^{\gsandorder\blam}}$ giving the result.
\end{proof}

To summarize the results so far:
\begin{enumerate}[label=$\triangleright$]

\item In all finite types the path idempotents $\idem$
are ordered by how far strings are to the right.

\item The end points of steady diagrams correspond to paths in crystals.

\item Up to more dominant terms, the steady permutation diagrams are compositions
of diagrammatic face permutations, where diagrammatic face permutations correspond to
face permutations in the crystals.

\item In level one, the weighted plactic monoid moves were used to show that steady permutation
diagrams for detour permutations are the same up to more dominant terms.

\end{enumerate}

\begin{Example}\label{Ex:Octogon}
Let $\Gamma$ be a quiver of type $\typea[5]$ and set $\Lambda=\Lambda_{1}+\Lambda_{2}$, $\brho=(1,2)$, and $\charge=(0,3)$. Then $\bpath=(1,2,2,1)$ and $\cpath=(2,1,1,2)$ are both paths in $\crystalgraph[\Lambda]$. Using the notation of \autoref{SS:TensorProducts}, $\bpath=(12,21)$ and $\cpath=(12,21)$. Hence, $\idem[\bpath]=\idem[\cpath]$. There are two face permutations in $\fFace(\bpath,\cpath)$. Namely, the trivial permutation $\idem[\bpath]$ and the permutation
\begin{gather*}
D_{\bpath,\cpath}\big((1,2)(3,4)\big)=
\begin{tikzpicture}[anchorbase,smallnodes,yscale=1,xscale=2]
\draw[redstring](0,-2)node[below]{$1$}--(0,0);
\draw[redstring](3.1,-2)node[below]{$2$}--(3.1,0);
\foreach \so/\i/\j in {-1.3/2/1, -0.1/1/2, } {
\draw[solid](\so,0)--++(3.1,-2)node[below]{$\j$};
\draw[ghost]({\so+1},0)node[above]{$\j$}--++(3.1,-2);
\draw[solid](\so,-2)node[below]{$\i$}--++(3.1,2);
\draw[ghost]({\so+1},-2)--++(3.1,2)node[above]{$\i$};
}
\end{tikzpicture}
\end{gather*}
The face permutation $D_{\bpath,\cpath}\big((1,2)(3,4)\big)$ has degree zero and is invertible.
\end{Example}

\begin{Example}\label{Ex:PathSwapping}
Extending \autoref{E:MainExampleTheBeastItself}, let $\Gamma$ be the quiver of type $\typeb[3]$ and let $\bfweight=(\fweight[2],\fweight[2])$. Let $n=10$, $\ell=2$, $\charge=(0,3)$ and $\brho=(2,2)$. Let $w=(1,6)(2,7)(3,8)(4,9)(5,10)$, Then
\begin{gather*}
D_{(\ppath[\mu],\ppath[\mu]),(\ppath[\mu],\ppath[\mu])}(w) =
\begin{tikzpicture}[anchorbase,smallnodes,yscale=1,xscale=2]
%\foreach \shift in {0,3.1} {
%\draw[ghost]({\shift-0.5},0)--++(0,1)node[above,yshift=-1pt]{$1$};
%\draw[ghost]({\shift+0.6},0)--++(0,1)node[above,yshift=-1pt]{$2$};
%\draw[ghost]({\shift+0.9},0)--++(0,1)node[above,yshift=-1pt]{$2$};
%\draw[solid]({\shift-1.5},0)--++(0,1);
%\draw[solid]({\shift-0.4},0)--++(0,1);
%\draw[solid]({\shift-0.1},0)--++(0,1);
%\draw[solid]({\shift+0.7},0)--++(0,1);
%\draw[solid]({\shift+0.8},0)--++(0,1);
%\draw[redstring](\shift,-2)node[below]{$2$}--(\shift,1);
%}
\draw[redstring](0,-2)node[below]{$2$}--(0,0);
\draw[redstring](3.1,-2)node[below]{$2$}--(3.1,0);
\foreach \go/\i in {-0.5/1, 0.6/2, 0.9/2 }{
\draw[ghost](\go,0)node[above]{$\i$}--++(3.1,-2);
\draw[ghost](\go,-2)--++(3.1,2)node[above]{$\i$};
}
\foreach \so/\i in {-1.5/1, -0.4/2, -0.1/2\;, 0.7/3\;, 0.8/3} {
\draw[solid](\so,0)--++(3.1,-2)node[below]{$\i$};
\draw[solid](\so,-2)node[below]{$\i$}--++(3.1,2);
}
\draw[decoration={brace, mirror}, decorate]
(-1.5,-2.35)--node[below]{$\idem[\mu]$}(0.9,-2.35);
\draw[decoration={brace, mirror}, decorate]
(1.6,-2.35)--node[below]{$\idem[\mu]$}(4.0,-2.35);
\end{tikzpicture}
\end{gather*}
Looking at the equator of this diagram, and comparing with the crystal graph in \autoref{E:MainExampleTheBeastItself}, this diagram is unsteady in view of \autoref{P:PathIdempotents}. In particular, this diagram is not the diagram of a face path permutation by \autoref{L:FiniteTypesPermutations} and \autoref{P:ProofsPermutationsSteady}.
\end{Example}

%%%%%%%%%%%%%%%%%%%%%%%%%%%%%%%%%%%%%%%%%

\subsection{Combinatorics of diagrams in finite types}\label{SS:ProofsTwo}

%%%%%%%%%%%%%%%%%%%%%%%%%%%%%%%%%%%%%%%%%

As in \cite{MaTu-klrw-algebras} and \cite{MaTu-klrw-algebras-bad}, in \autoref{T:ConstructionGeneralMain} the basis theorem for the infinite dimensional wKLRW algebra $\aWA$ implies the result for the cyclotomic quotient $\WAc(X)$, so we focus on
the infinite dimensional algebra in this section. In this setting
we have the following \emph{standard basis} and \emph{polynomial
module}.

Recall from \autoref{D:ConstructionPermBasis} that each permutation $w\in\sym$ defines a diagram $D_{w}\idem[\bx,\bi]$, for $x\in X$ and $\bi\in\vertices^{n}$.

\begin{Proposition}\label{P:ProofWABasis}
The algebra $\WA[n](X)$ is free as an $\ring$-module with homogeneous basis
\begin{gather}\label{Eq:ProofWABasis}
\WABasis=
\set{D_{w}y_{1}^{a_{1}}\dots y_{n}^{a_{n}}\idem[{\bx,\bi}]|a_{1},\dots,a_{n}\in\N,w\in\sym,\bx\in X,\bi\in\vertices^{n}}.
\end{gather}
\end{Proposition}

\begin{proof}
See \cite[Proposition 3B.12]{MaTu-klrw-algebras}.
\end{proof}

Similarly to \autoref{Eq:MainExampleSandwich},
the mnemonic for \autoref{Eq:ProofWABasis} is
\begin{gather*}
D_{w}\underbrace{y_{1}^{a_{1}}\dots y_{n}^{a_{n}}\idem[{\bx,\bi}]}_{=d}
\quad\leftrightsquigarrow\quad
\begin{tikzpicture}[anchorbase,scale=1]
\draw[line width=0.75,color=black,fill=cream] (0,1) to (0.25,0.5) to (0.75,0.5) to (1,1) to (0,1);
\node at (0.5,0.75){\scalebox{0.85}{$D_{w}$}};
\draw[line width=0.75,color=black,fill=cream] (0.25,0) to (0.25,0.5) to (0.75,0.5) to (0.75,0) to (0.25,0);
\node at (0.5,0.25){$d$};
\end{tikzpicture}
,
\quad\text{where}\quad
\begin{aligned}
\begin{tikzpicture}[anchorbase,scale=1]
\draw[line width=0.75,color=black,fill=cream] (0,1) to (0.25,0.5) to (0.75,0.5) to (1,1) to (0,1);
\node at (0.5,0.75){\scalebox{0.85}{$D_{w}$}};
\end{tikzpicture}
&\text{ a permutation diagram,}
\\
\begin{tikzpicture}[anchorbase,scale=1]
\draw[white,ultra thin] (0,0) to (1,0);
\draw[line width=0.75,color=black,fill=cream] (0.25,0) to (0.25,0.5) to (0.75,0.5) to (0.75,0) to (0.25,0);
\node at (0.5,0.25){$d$};
\end{tikzpicture}
&\text{ dots on an idempotent}.
\end{aligned}
\end{gather*}

Importantly, the algebra $\aWA$ has a faithful polynomial representation. To define it, let $y_{1},\dots,y_{n}$ be indeterminates over $R$ and set
\begin{gather*}
P_{n}(X)=
\bigoplus_{\bx\in X,\bi\in\vertices^{n}}
R[y_{1},\dots,y_{n}]\idem[{\bx,\bi}].
\end{gather*}
The symmetric group $\sym$ acts on $P_{n}(X)$ by place permutations.
Let $\partial_{r,s}=\frac{(r,s)-1}{y_{s}-y_{r}}$
be the \emph{Demazure operator}.
By \cite[Lemma 3B.11]{MaTu-klrw-algebras}, $\aWA$ acts on $P_{n}(X)$ by:
\begin{gather}\label{E:polyAction}
\idem[{\by,\bj}]\cdot f(y)\idem[{\bx,\bi}]=\delta_{\bx\by}\delta_{\bi\bj}f(y)\idem[{\bx,\bi}]
,\quad
\begin{tikzpicture}[anchorbase,smallnodes,rounded corners]
\draw[solid,dot] (0,-0.2)node[below]{$i_{r}$} to (0,0.5)node[above,yshift=-1pt]{$\phantom{i}$};
\end{tikzpicture}
\mapsto y_{r}
,\quad
\begin{tikzpicture}[anchorbase,smallnodes,rounded corners]
\draw[ghost,dot] (0,-0.2)node[below]{$\phantom{i}$} to (0,0.5)node[above,yshift=-1pt]{$i_{r}$};
\end{tikzpicture}
\mapsto 1,
\end{gather}
and all crossings act as the identity except that
\begin{gather*}
\begin{tikzpicture}[anchorbase,smallnodes,rounded corners]
\draw[solid] (0,0)node[below]{$i_{r}$} to (0.5,0.5);
\draw[solid] (0.5,0)node[below]{$i_{s}$} to (0,0.5)node[above,yshift=-1pt]{$\phantom{i}$};
\end{tikzpicture}
\mapsto
\begin{cases*}
\partial_{r,s} & if $i_{r}=i_{s}$,
\\
(r,s) & if $i_{r}\neq i_{s}$,
\end{cases*}
\;
\begin{tikzpicture}[anchorbase,smallnodes,rounded corners]
\draw[solid] (0,0)node[below]{$i_{r}$} to (0.5,0.5);
\draw[redstring] (0.5,0)node[below]{$i_{s}$} to (0,0.5)node[above,yshift=-1pt]{$\phantom{i}$};
\end{tikzpicture}
\mapsto
\begin{cases*}
x_{r} & if $i_{r}=i_{s}$,
\\
1 & if $i_{r}\neq i_{s}$,
\end{cases*}
\;
\begin{tikzpicture}[anchorbase,smallnodes,rounded corners]
\draw[ghost] (0,0) to (0.5,0.5)node[above,yshift=-1pt]{$i_{r}$};
\draw[solid] (0.5,0)node[below]{$i_{s}$} to (0,0.5)node[above,yshift=-1pt]{$\phantom{i}$};
\end{tikzpicture}
\mapsto
\begin{cases}
Q_{i_{r},i_{s}}(y_{r},y_{s})
&\text{if $i\rightsquigarrow j$},
\\
1&\text{otherwise}.
\end{cases}
\end{gather*}

\begin{Proposition}\label{P:ProofsFaithful}
The algebra $\aWA$ acts faithfully on $P_{n}(X)$.
\end{Proposition}

\begin{proof}
See \cite[Corollary 3B.13]{MaTu-klrw-algebras}.
\end{proof}

Recall from \autoref{D:ConstructionChoice} that $\fParts{n}{\abfweight}$ is a fixed choice
of paths in $\crystalgraph[\crystal]$ of length of $n$.

\begin{Lemma}\label{L:ProofsIdempotentOne}
Let $\bpath,\cpath\in\fParts{n}{\abfweight}$. Then
$\idem[\bpath]=\idem[\cpath]$ if and only if
$\bpath=\cpath$.
\end{Lemma}

\begin{proof}
One direction is clear, so let us assume that $\bpath\neq\cpath$.
We will show that $\idem[\bpath]\neq\idem[\cpath]$.

\Case{$\Res[\bpath]=\Res[\cpath]$.} This case
does not arise because paths in irreducible crystals are uniquely
determined by the sequence of their labels.

\Case{$\Res[\bpath]\neq w\,\Res[\cpath]$ for any $w\in\sym$.}
In this case the idempotent diagrams $\idem[\bpath]$ and $\idem[\cpath]$
have different residues of their solid strings. Hence, \autoref{P:ProofWABasis} ensures that $\idem[\bpath]\neq\idem[\cpath]$.

\Case{$\Res[\bpath]=w\,\Res[\cpath]$ for some $1\ne w\in\sym$.}
The permutation $w$ necessarily swaps solid strings of different residues $i\neq j$ since otherwise $\Res[\bpath]=\Res[\cpath]$.  By \autoref{L:TensorFactors}, we can write $\bpath=(\bpath_{1},\dots,\bpath_{\hell})$ and $\cpath=(\cpath_{1},\dots,\cpath_{\hell})$. By \autoref{L:ProofsAnchor}, we can assume that $w$ permutes the residues in just one component, so that $\bpath_{m}=w\cpath_{m}$, for some $m$. The crystal graph $\crystalgraph[{\fweight[\affine{\rho}_{m}]}]$ is path ordered by \autoref{L:ProofsGraphs}, so $\omega(\bpath_{m})=\omega(\cpath_{m})$. Hence, $w$ is a face permutation by \autoref{L:FiniteTypesPermutations} and we can assume that $w$ only permutes two non-adjacent residues $i$ and~$j$. As $i$ and $j$ are not adjacent in the quiver $\quiver$, then $\idem[\bpath]=\idem[\cpath]$ by the construction of \autoref{D:ConstructionIdempotents}.
%If $i$ and $j$ are adjacent, then \autoref{P:ProofWABasis} implies that $\idem[\bpath]\neq\idem[\cpath]$.
%\DT{Careful: proof needs revision}
\end{proof}

%Recall, from after \autoref{D:ConstructionChoice}, that $\gsandorder$ is an order on $\Parts{n}{\abfweight}$. By \autoref{L:ProofsIdempotentOne} we get:

%\begin{Corollary}\label{C:ProofsIdempotentTwo}
%The sets $\fParts{n}{\abfweight}$ and $\aCrystalVertices$
%are totally ordered by $\gsandorder$.
%\DT{Wrong!}
%\end{Corollary}

%\begin{Remark}\label{R:ProofsAnchor}
%In \autoref{L:ProofsAnchor} the vertices of degree $3$ in $\quiver$ require extra care because,
%\textit{a priori}, there could exist inequivalent idempotent diagrams corresponding to different %vertices of the crystal, see {\eg} \autoref{R:ConsequencesGraphs} below for an example of this in %affine types. By the miracle in \autoref{C:ProofsIdempotentTwo}, this does not happen in finite type.
%\end{Remark}

A \emph{straight line diagram} is a diagram without crossings. For example, every idempotent diagram is a straight line diagram. As for idempotent diagrams, a diagram $D$ \emph{factors through} a straight line diagram $S$ if $D=ESF$, for some $E,F\in\WA[n](X)$. Here, $E$ and $F$ are allowed to be linear combinations of diagrams.

\begin{Definition}\label{D:RecollectionFactors}
Let $S$ be a straight line diagram. The \emph{left justification} of $S$ is an idempotent diagram $L(S)$ such that $S$ factors through $L(S)$, and $L(S)$ has its solid strings in admissible parking positions that are as far to the left as possible.
\end{Definition}

By successively pulling strings as far to the left as possible, it follows that every straight line diagram has a unique left justification; {\cf} \cite[Lemma 5D.9]{MaTu-klrw-algebras}.

\begin{Example}\label{E:ProofsLeftJustification}
The notions from \autoref{D:RecollectionFactors} can be illustrated by
\begin{gather*}
S=
\begin{tikzpicture}[anchorbase,smallnodes,rounded corners]
\draw[ghost](0,0)node[below]{$\phantom{i}$}--++(-0.2,1)node[above,yshift=-1pt]{$i_{1}$};
\draw[ghost](1.2,0)node[below]{$\phantom{i}$}--++(-0.3,1)node[above,yshift=-1pt]{$i_{2}$};
\draw[ghost](2.85,0)node[below]{$\phantom{i}$}--++(0.2,1)node[above,yshift=-1pt]{$i_{3}$};
\draw[ghost](4.25,0)node[below]{$\phantom{i}$}--++(0.1,1)node[above,yshift=-1pt]{$i_{4}$};
\draw[solid](-1,0)node[below]{$i_{1}$}--++(-0.2,1) node[above,yshift=-1pt]{$\phantom{i_{1}}$};
\draw[solid](0.2,0)node[below]{$i_{2}$}--++(-0.3,1);
\draw[solid](1.85,0)node[below]{$i_{3}$}--++(0.2,1);
\draw[solid](3.25,0)node[below]{$i_{4}$}--++(0.1,1);
\draw[redstring](-0.5,0)node[below]{$\rho_{1}$}--++(-0.1,1);
\draw[redstring](1.5,0)node[below]{$\rho_{2}$}--++(-0.1,1);
\draw[redstring](2.25,0)node[below]{$\rho_{3}$}--++(0.2,1);
\draw[very thick,mygray,densely dotted] (-1.25,0.5) to (4.5,0.5);
\end{tikzpicture}
,\quad
L(S)=
\begin{tikzpicture}[anchorbase,smallnodes,rounded corners]
\draw[ghost](0,0)node[below]{$\phantom{i}$}--++(0,1)node[above,yshift=-1pt]{$i_{1}$};
\draw[ghost](1.1,0)node[below]{$\phantom{i}$}--++(0,1)node[above,yshift=-1pt]{$i_{2}$};
\draw[ghost](2.3,0)node[below]{$\phantom{i}$}--++(0,1)node[above,yshift=-1pt]{$i_{3}$};
\draw[ghost](3.4,0)node[below]{$\phantom{i}$}--++(0,1)node[above,yshift=-1pt]{$i_{4}$};
\draw[solid](-1.0,0)node[below]{$i_{1}$}--++(0,1) node[above,yshift=-1pt]{$\phantom{i_{1}}$};
\draw[solid](0.1,0)node[below]{$i_{2}$}--++(0,1) node[above,yshift=-1pt]{$\phantom{i_{1}}$};
\draw[solid](1.3,0)node[below,xshift=0.075cm]{$i_{3}$}--++(0,1) node[above,yshift=-1pt]{$\phantom{i_{1}}$};
\draw[solid](2.4,0)node[below]{$i_{4}$}--++(0,1) node[above,yshift=-1pt]{$\phantom{i_{1}}$};
\draw[redstring](-0.9,0)--++(0,1)node[above]{$\rho_{1}$};
\draw[redstring](1.2,0)node[below,xshift=-0.075cm]{$\rho_{2}$}--++(0,1);
\draw[redstring](1.4,0)--++(0,1)node[above]{$\rho_{3}$};
\end{tikzpicture}
.
\end{gather*}
The left justification $L(S)$ can be thought
of as being in a collar neighborhood of the points on the dashed line through the middle of $S$, so $S$ factors through $L(S)$. (Recall that a \emph{collar neighborhood} of a set is the Cartesian product of the set with a half-open interval.)
\end{Example}

Once our main theorem is established, the general theory of
sandwich cellular algebras defines the two-sided ideals $\WA[n](X)^{\gsandorder\blam}$ for $\blam\in\daCrystalVertices$. However, we need to define these ideals now even though we do not yet have a
sandwich cellular basis. For $\blam\in\daCrystalVertices$ let $\WA[n](X)^{\gsandorder\blam}$ be the $\ring$-submodule of $\aWA$ spanned by all diagrams~$D$ such that $D$ factors through a straight line diagram $S$ with $L(S)\gsandorder\idem[\blam]$. Define $\WA[n](X)^{\rsandorder\blam}$ in the same way.

The next two results imply that $\WA[n](X)^{\gsandorder\blam}$ and $\WA[n](X)^{\rsandorder\blam}$ are two-sided ideals of $\WA[n](X)$.

\begin{Lemma}\label{L:ProofsVerticalDominance}
Suppose that $D\in\WA[n](X)$ and that $D$ factors through the idempotent diagram $S$. Then there exists $\blam\in\daCrystalVertices$ such that $D$ factors through $\dotidem$ and $L(S)\gsandorder\dotidem[\blam]$ modulo $\WA[n](X)^{\gsandorder\blam}$.
\end{Lemma}

\begin{proof}
Without loss of generality we can assume that $D=S$. Using the results in \autoref{SS:RecollectionMoveRight} to pull strings to the right, if necessary, we can replace $D$ with a linear combination of more $\gsandorder$-dominant diagrams in which every string is anchored on some (affine) red string. With respect to $\gsandorder$, let $T=L(T)$ be the minimal straight line diagram such that one of the summands of $D$ factors through $T$. By construction, $T\gsandorder L(S)$ and every string in $T$ is anchored on an (affine) red string. By \autoref{L:ProofsSteady}, the strings anchored around the red $\rho_{i}$-string
correspond to the path $\bpath$ in $\crystalgraph[{\fweight[i]}]$, so $T$ naturally corresponds to an $\hell$-tuple $(\bpath_{1},\dots,\bpath_{\hell})$, where $\bpath_{k}\in\crystalgraph[{\fweight[\affine{\rho}_{k}]}]$. Note that this argument also applies to the affine red strings because, by definition, the strings anchored around the affine red strings are blocked by the affine red strings. By construction each component of $T$ is anchored around an (affine) red string so, up to conjugation by straight line diagrams, $T=\idem[\bpath]$ by \autoref{L:ProofsSteady}, for some $\bpath\in\fParts{n}{\abfweight}$. Let $\blam=\omega(\bpath)$. In view of \autoref{L:ProofsSpanningPermutations}, we can assume that $\bpath=\ppath$. By \autoref{L:RecollectionMovingStringsDots}.(a), the dots on dotted idempotents are invertible, modulo higher terms, so $D$ factors through  $\dotidem[\blam]$ as claimed.
\end{proof}

From here on we can use almost the same arguments as in \cite[Sections 6D and 7E]{MaTu-klrw-algebras}.

For the next lemma, recall from \autoref{D:ConstructionSanddots} that $c_{k}(\blam)$
is the maximum number of dots allowed on the $k$th solid string in the sandwich algebra $\sand[\blam]$.
The next lemma says that the idempotent diagram $\dotidem[\blam]$ is a minimal diagram in $\WA[n](X)^{\gsandorder\blam}$,
with respect to the $\sandorder$-order, and adding dots or
crossings gives a linear combination of more dominant diagrams.

\begin{Lemma}\label{L:ProofsGettingBigger}
\leavevmode
\begin{enumerate}

\item Suppose that $\blam\in\daCrystalVertices$ and $1\leq k\leq n$ such that the $k$th string is not immediately adjacent to an affine string of the same residue. Then $y_{k}y_{k}^{c_{k}(\blam)}\dotidem[\blam]\in\WA[n](X)^{\rsandorder\blam}$.

\item Suppose that $\blam\in\daCrystalVertices$ and $w\in\Face(\ppath,\bpath)$ for $\idem[\bpath]\rsandorder\dotidem[\blam]$. Then $D_{w}\dotidem[\blam],\dotidem[\blam]D_{w}^{\star}\in\WA[n](X)^{\rsandorder\blam}$.
\end{enumerate}
\end{Lemma}

\begin{proof}
\textit{(a).} By the construction of $\dotidem[\blam]$ in \autoref{Eq:ConstructionTheBasisAffine}, putting a dot on the $k$th solid string allows us to pull the string to the right by
\autoref{L:RecollectionMovingStringsDots}, applied to the $k$th solid string or its ghosts. Now apply \autoref{L:ProofsVerticalDominance}.

\textit{(b).} This follows directly from \autoref{P:ProofsPermutationsSteady} and
\autoref{L:ProofsVerticalDominance}.
\end{proof}

\begin{Lemma}\label{L:ProofsKLRWLinearIndependent}
The diagrams in \upshape{\autoref{Eq:ConstructionTheBasisAffine}} form
an $\ring$-basis of $\WA[n](X)$.
\end{Lemma}

\begin{proof}
\textit{Spanning.} By \autoref{P:ProofWABasis},
$\WA[n](X)$ is spanned by the diagrams
$D_{w}^{\star}y^{\bb}\idem[\blam]D_{v}$, for $w,v\in\sym$ and $\bb\in\N$. \autoref{L:ProofsGettingBigger}(a) implies that we can restrict
the dots $y^{\bb}$ to $\bb=\ba+\bfi$, where $\ba\in\Affch$ and $\bfi\in\Finch$ as in
\autoref{Eq:ConstructionTheBasisAffine}. By \autoref{P:ProofsPermutationsSteady},
we can assume that $D_{w}$ and $D_{v}$ are diagrammatic face permutations. By \autoref{L:ProofsGettingBigger}(b), we can further assume that  $w,v\in\bigcup_{\bnu\lsandorder\blam}\Face(\blam,\bnu)$. Hence, the diagrams in \autoref{Eq:ConstructionTheBasisAffine} span $\WA[n](X)$.

\textit{Linear independence.} Using \autoref{P:ProofsFaithful},
it is enough to show that the images of the elements from \autoref{Eq:ConstructionTheBasisAffine} in the endomorphism algebra of the polynomial module $P_{n}(X)$
are linearly independent. Note that it suffices to check
linear independence of the
elements from \autoref{Eq:ConstructionTheBasisAffine}
in the associated graded algebra of $\WA[n](X)$, see
\cite[Section 3B]{MaTu-klrw-algebras} for details.
To this end, for $\bS,\bT\in\Face(\blam)$ let $\idem[\bS]$ and
$\idem[\bT]$ be the idempotents whose endpoints
are determined by the top of $D_{\bS}^{\star}\idem[\blam]$
and the bottom of $\idem[\blam]D_{\bT}$, {\ie}
\begin{gather*}
D_{\bS}^{\star}\idem[\blam]
\leftrightsquigarrow
\begin{tikzpicture}[anchorbase,scale=1]
\draw[white] (0,0) to (0,-0.25);
\draw[line width=0.75,color=black,fill=cream] (0,1) to (0.25,0.5) to (0.75,0.5) to (1,1) to (0,1);
\node at (0.5,0.75){\scalebox{0.75}{\,$D_{\bS}^{\star}$}};
\draw[line width=0.75,color=black,fill=cream] (0.25,0) to (0.25,0.5) to (0.75,0.5) to (0.75,0) to (0.25,0);
\node at (0.5,0.25){$\idem[\blam]$};
\draw[->] (1.3,1)node[right]{\text{endpoints}} to (1.1,1);
\end{tikzpicture}
,\quad
\idem[\blam]D_{\bT}
\leftrightsquigarrow
\begin{tikzpicture}[anchorbase,yscale=-1]
\draw[white] (0,0) to (0,-0.25);
\draw[line width=0.75,color=black,fill=cream] (0,1) to (0.25,0.5) to (0.75,0.5) to (1,1) to (0,1);
\node at (0.5,0.75){\scalebox{0.85}{$D_{\bT}$}};
\draw[line width=0.75,color=black,fill=cream] (0.25,0) to (0.25,0.5) to (0.75,0.5) to (0.75,0) to (0.25,0);
\node at (0.5,0.25){$\idem[\blam]$};
\draw[->] (1.3,1)node[right]{\text{endpoints}} to (1.1,1);
\end{tikzpicture}
.
\end{gather*}
The action of these elements on the polynomial module is given by \autoref{E:polyAction}. Hence, in the associated graded, $D_{\bT}$ acts on $\idem[\bT]P_{n}(X)$ by sending $\idem[\bT]f(y_{1},\dots,y_{n})$ to $\idem[\blam]f(y_{w_{\bT}(1)},\dots,y_{w_{\bT}(n)})$, and
$D_{\bS}^{\star}$ acts similarly except that the permutation is
$w_{\bS}^{-1}$. Thus, in the associated graded
module, the element $D_{\bS\bT}^{\ba,\bfi}$ from
in \autoref{Eq:ConstructionTheBasisAffine} acts as the endomorphism
\begin{gather*}
\idem[\bS]f(y_{1},\dots,y_{n})\mapsto
\idem[\bT]y^{\ba}y^{\bfi}f(y_{w_{\bS}^{-1}w_{\bT}(1)},\dots,y_{w_{\bS}^{-1}w_{\bT}(n)}).
\end{gather*}
Linear independence follows.
\end{proof}

\begin{proof}[Proof of \autoref{T:ConstructionMain}]
\textit{(a).} By the results and definitions above, the only statement that is not immediate
is (AC$_{3}$). In (AC$_{3}$) it is sufficient to let $x$ be
a crossing or a dot since, up to isotopy, these generate
$\WA[n](X)$. The $x$ will be on top of the diagram
$D_{\bS\bT}^{\ba,\bfi}$.

If $x$ is a crossing, then (AC$_{3}$) is immediate from \autoref{P:ProofsPermutationsSteady} and \autoref{L:ProofsGettingBigger}(b).

Now suppose that $x$ is a dot. This follows almost immediately
from \autoref{L:ProofsGettingBigger}(a) with the caveat that
we first need to pull the dot $x$ to the equator. Let $i,j,k\in\vertices$ with
$i\neq j$, and recall from
\autoref{SS:RecollectionRelations} that we have
\begin{gather*}
\begin{tikzpicture}[anchorbase,smallnodes,rounded corners]
\draw[solid](0.5,0.5)node[above,yshift=-1pt]{$\phantom{i}$}--(0,0) node[below]{$i$};
\draw[solid,dot=0.25](0,0.5)--(0.5,0) node[below]{$i$};
\end{tikzpicture}
=
\begin{tikzpicture}[anchorbase,smallnodes,rounded corners]
\draw[solid](0.5,0.5)node[above,yshift=-1pt]{$\phantom{i}$}--(0,0) node[below]{$i$};
\draw[solid,dot=0.75](0,0.5)--(0.5,0) node[below]{$i$};
\end{tikzpicture}
+
\underbrace{\begin{tikzpicture}[anchorbase,smallnodes,rounded corners]
\draw[solid](0,0.5)node[above,yshift=-1pt]{$\phantom{i}$}--(0,0) node[below]{$i$};
\draw[solid](0.5,0.5)--(0.5,0) node[below]{$i$};
\end{tikzpicture}}_{\text{error term}}
,
\quad
\begin{tikzpicture}[anchorbase,smallnodes,rounded corners]
\draw[solid](0.5,0.5)node[above,yshift=-1pt]{$\phantom{i}$}--(0,0) node[below]{$j$};
\draw[solid,dot=0.25](0,0.5)--(0.5,0) node[below]{$i$};
\end{tikzpicture}
=
\begin{tikzpicture}[anchorbase,smallnodes,rounded corners]
\draw[solid](0.5,0.5)node[above,yshift=-1pt]{$\phantom{i}$}--(0,0) node[below]{$j$};
\draw[solid,dot=0.75](0,0.5)--(0.5,0) node[below]{$i$};
\end{tikzpicture}
,
\quad
\begin{tikzpicture}[anchorbase,smallnodes,rounded corners]
\draw[ghost](0.5,0.5)node[above,yshift=-1pt]{$k$}--(0,0) node[below]{$\phantom{i}$};
\draw[solid,dot=0.25](0,0.5)--(0.5,0) node[below]{$i$};
\end{tikzpicture}
=
\begin{tikzpicture}[anchorbase,smallnodes,rounded corners]
\draw[ghost](0.5,0.5)node[above,yshift=-1pt]{$k$}--(0,0) node[below]{$\phantom{i}$};
\draw[solid,dot=0.75](0,0.5)--(0.5,0) node[below]{$i$};
\end{tikzpicture}
\quad
\begin{tikzpicture}[anchorbase,smallnodes,rounded corners]
\draw[solid,dot=0.25](0,0.5)--(0.5,0) node[below]{$i$};
\draw[redstring](0.5,0.5)node[above,yshift=-1pt]{$\phantom{i}$}--(0,0) node[below]{$k$};
\end{tikzpicture}
=
\begin{tikzpicture}[anchorbase,smallnodes,rounded corners]
\draw[solid,dot=0.75](0,0.5)--(0.5,0) node[below]{$i$};
\draw[redstring](0.5,0.5)node[above,yshift=-1pt]{$\phantom{i}$}--(0,0) node[below]{$k$};
\end{tikzpicture}
,
\end{gather*}
and all the partner relations. In other words,
we can slide dots freely to the equator up to additional error terms
for $(i,i)$-crossings. Let us first ignore these error
terms, and first consider the case when $x$ reaches the equator, and say
it does so on the $k$th solid string. If $k$ is in the unsteady part of the diagram, or if $k$ is in the steady part of the diagram and number of dots on the equator with $x$ is strictly
smaller than $c_{k}(\blam)$, then we just park the dot where it is
to obtain another basis element. In all other cases
\autoref{L:ProofsGettingBigger}(a) applies, giving (AC$_{3}$).
Finally, for the error terms we use \autoref{L:RecollectionMovingStringsDots}
to pull strings to the right, again showing that (AC$_{3}$) holds.

\textit{(b).} By \autoref{P:PathIdempotents}, a diagram in \autoref{Eq:ConstructionTheBasisAffine} is steady if and only if is in \autoref{Eq:ConstructionTheBasisFinite}. Hence, the sandwich cellularity result for $\WAc[n](X)$ follows immediately from \autoref{T:ConstructionMain}.(a), the corresponding result for the algebra $\WA[n](X)$.
\end{proof}

%%%%%%%%%%%%%%%%%%%%%%%%%%%%%%%%%%%%%%%%%

\section{Some consequences of sandwich cellularity}\label{S:Consequences}

%%%%%%%%%%%%%%%%%%%%%%%%%%%%%%%%%%%%%%%%%

This section proves results for $\WA[n](X)$ and $\WAc[n](X)$ that follow from \autoref{T:ConstructionMain}. Some results are restricted to level one.
We remind the reader that we are using the labeling of the Dynkin diagrams given in \autoref{SS:FiniteTypesConventions}.

%%%%%%%%%%%%%%%%%%%%%%%%%%%%%%%%%%%%%%%%%

\subsection{Maximal number of strings}\label{SS:ConsequencesStrings}

%%%%%%%%%%%%%%%%%%%%%%%%%%%%%%%%%%%%%%%%%

We now focus on $\WAc[n](X)$.
We start by observing that $\WAc[n](X)$ is zero for most $n$.

\begin{Definition}\label{D:ConsequencesStrings}
Attach the following numbers $n(i)$, for $i\in\vertices$ to the
fundamental weights in finite types:
\begin{center}
\rowcolors{2}{mygray!25}{}
\begin{tabular}{c|c}
\rowcolor{mygray!75}
Type & $n(i)$ for $i\in\vertices$
\\
$\typea[e]$ &  $i(e-i+1)$
\\
$\typeb[e>1]$ & $(e+i)(e-i+1)$ for $i\neq e$, $\tfrac{e(e+1)}{2}$ otherwise
\\
$\typec[e>2]$ & $(e+i-1)(e-i+1)$
\\
$\typed[e>3]$ & $2ie-i(i+1)$ for $i\notin\set{e-1,e}$, $\frac{e(e-1)}{2}$ otherwise
\\
$\typee[6]$ & $16,30,42,22,30,16$
\\
$\typee[7]$ & $34,66,96,51,75,52,27$
\\
$\typee[8]$ & $92,182,270,136,220,168,114,58$
\\
$\typef[4]$ & $22,42,30,16$
\\
$\typeg[2]$ & $10,6$
\end{tabular}
,
\end{center}
where we list $n(i)$ starting from $i=1$ for the exceptional types.
\end{Definition}

\begin{Proposition}\label{P:ConsequencesStrings}
Let $\bfweight=(\fweight[\rho_{1}],\dots,\fweight[\rho_{\ell}])$.
Then $\WAc[n](X)\not\cong 0$ if and only if $1\leq n\leq n(\rho_{1})+\dots+n(\rho_{\ell})$.
\end{Proposition}

\begin{proof}
The numbers in \autoref{D:ConsequencesStrings} appear as
\begin{gather*}
\fweight[i]=\sum_{j=1}^{e}m_{ij}\cdot\sroot[j],
\quad
n(i)=2\sum_{j=1}^{e}m_{ij},
\end{gather*}
see {\eg} \cite[Plates I--IX]{Bo-chapters-4-6}, with the slight caveat that the numbers for some types need to be permuted because of our differences in labeling in \autoref{SS:FiniteTypesConventions}.
Since $\fweight[i]$ and $-\fweight[i]$
correspond to the highest and lowest weight vectors in the
representation for $\crystalgraph[{\fweight[i]}]$, it follows
that the paths in $\crystalgraph[{\fweight[i]}]$ have length at most $n(i)-1$. Therefore, $\WAc(X)\cong 0$ by \autoref{T:ConstructionMain} if and only if $n>n(\rho_{1})+\dots+n(\rho_{\ell})$ or $n<1$.
\end{proof}

In particular, \autoref{P:ConsequencesStrings} implies
that for each fixed rank $e$
there are only finitely many $n$ such that $\WAc[n](X)$ is nonzero.

\begin{Example}\label{E:ConsequencesStrings}
Let $\ell=1$ and assume that the rank of $\quiver$ is strictly less than nine. Then the maximal number of solid strings in a nonzero diagram is $270$, which happens in type $\typee[8]$.
\end{Example}

%%%%%%%%%%%%%%%%%%%%%%%%%%%%%%%%%%%%%%%%%

\subsection{Classification of simples}\label{SS:ConsequencesSimples}

%%%%%%%%%%%%%%%%%%%%%%%%%%%%%%%%%%%%%%%%%

\begin{Remark}\label{R:ConstructionMainTwo}
By comparing \autoref{T:ConstructionMain} and
\autoref{C:ConstructionMainTwo}, all
the results below have obvious analogs
for $\TA[n]$ and $\TAc[n]$. We omit details for the sake of brevity.
\end{Remark}

Let $\lambda\in\aCrystalVertices$ and let $a(\blam)$ be the number of unsteady strings in $\idem[\blam]$. Equivalently, $a(\blam)$ is the number of strings in $\idem[\blam]$, or the number of positions in $\Affch$, that can carry arbitrarily many dots.
In the following result, simple modules are counted up to homogeneous isomorphism.

\begin{Proposition}\label{P:ConsequencesSimples}
Let $\ring$ be a field.
\begin{enumerate}

\item The set of apexes of $\WA[n](X)$ and $\WAc[n](X)$ are in 1:1-correspondence
with $\aCrystalVertices$ and $\CrystalVertices$, respectively.
In other words, $\Pcal^{\neq 0}=\Pcal$ in the notation
of \autoref{T:SandwichMain}(b).

\item %There is a 1:1-correspondence between simple $\WA[n](X)$-modules with apex $\blam\in\aCrystalVertices$ and $R^{a(\blam)}$. There is a unique $\WAc[n](X)$-module with apex $\blam\in\CrystalVertices$.
There is a 1:1-correspondence
between simple $\WA[n](X)$-modules with apex $\blam\in\daCrystalVertices$ and $\bigoplus_{\bsig\in\blam}R^{a(\bsig)}$. For every $\blam\in\dCrystalVertices$, there
are $\#\blam$ isomorphisms classes of simple $\WAc[n](X)$-modules with apex $\blam$.

\item There is exactly one graded simple $\WA[n](X)$, and one graded simple $\WAc[n](X)$-module, for each vertex in~$\aCrystalVertices$, or $\CrystalVertices$, respectively.

\end{enumerate}
\end{Proposition}

\begin{proof}
\textit{(a).} Suppose that a simple $\WA(X)$-module has apex $\blam\in\daCrystalVertices$c. Then  $\sand[\blam]=\bigoplus_{\bsig\in\blam}\sand[{\ppath[\bsig]}]$ is a direct sum of polynomial rings. Hence, by \autoref{T:SandwichMain}(b), the number of simple $\sand[\blam]$-modules is equal $\#\blam$, which is the size of the equivalence class of~$\blam$. Hence, it is enough to show that if $\bsig\in\aCrystalVertices$ then there is a simple $\WA(X)$-module with apex $\bsig$. Similarly, we need to show that if $\bsig\in\aCrystalVertices$ then there is a simple $\WAc(X)$-module with apex $\bsig$.

Assume first that we are not in type $\typeg[2]$ and fix $\bsig\in\aCrystalVertices$ or
$\bsig\in\CrystalVertices$. If there are
no repeated residues, then $\dotidem[\bsig]=\idem[\bsig]$
is the idempotent associated with $\bsig$. In the case of
repeated residues, first recall that $\bsig$ has at most two consecutive residues
by \autoref{C:ConstructionRepeated}, and secondly note that the relation
\begin{gather*}
\bigg(
\begin{tikzpicture}[anchorbase,smallnodes,rounded corners]
\draw[solid](0.5,0.5)node[above,yshift=-1pt]{$\phantom{i}$}--(0,0) node[below]{$i$};
\draw[solid,dot=0.75](0,0.5)--(0.5,0) node[below]{$i$};
\end{tikzpicture}
\bigg)^{2}
=
-
\begin{tikzpicture}[anchorbase,smallnodes,rounded corners]
\draw[solid](0.5,0.5)node[above,yshift=-1pt]{$\phantom{i}$}--(0,0) node[below]{$i$};
\draw[solid,dot=0.75](0,0.5)--(0.5,0) node[below]{$i$};
\end{tikzpicture}
\end{gather*}
and its partner relation follow from the dot sliding relation \autoref{Eq:RecollectionDotCrossing}. Thus, the crossing with a dot
is a pseudo idempotent. By \autoref{D:ConstructionPermBasis} we always
flank $\dotidem[\bsig]$ with crossings for repeated residues.
Since all $\bsig\in\aCrystalVertices$ and all $\bsig\in\CrystalVertices$ have an associated
idempotent by the above discussion, the statement is a direct consequence of the basis theorems
\autoref{T:ConstructionMain} and \autoref{C:ConstructionMainTwo}, combined with the classification of simple modules in sandwich cellular algebras from  \autoref{T:SandwichMain}(b) and \cite[Theorem 2.16]{TuVa-handlebody}.

In type $\typeg[2]$ we can use essentially the same argument, with the only difference being that three repeated residues are possible by \autoref{C:ConstructionRepeated}.
To deal with these one checks that
\begin{gather*}
\begin{tikzpicture}[anchorbase,smallnodes,rounded corners]
\draw[solid](0,0)node[below]{$i$} to (1,0.5)node[above,yshift=-1pt]{$\phantom{i}$};
\draw[solid,dot=0.15](0.5,0)node[below]{$i$} to (0,0.5);
\draw[solid,dot=0.15,dot=0.45](1,0)node[below]{$i$} to (0.5,0.5);
\end{tikzpicture}
\end{gather*}
is a pseudo idempotent, after which the rest of the argument works {\ver}.

\textit{(b)+(c).} By \autoref{T:ConstructionMain} and
\autoref{C:ConstructionMainTwo}, the sandwiched algebras are
\begin{align*}
\sand[\blam]&\cong \bigoplus_{\bsig\in\blam}R[X_{1},\dots,X_{a(\bsig)}]\otimes R[Y_{1},\dots,Y_{f(\bsig)}]/(Y_{1}^{z_{1}},\dots,Y_{f(\bsig)}^{z_{f(\bsig)}})
,\\
\fsand[\blam]&\cong\bigoplus_{\bsig\in\blam} R[Y_{1},\dots,Y_{f(\bsig)}]/(Y_{1}^{z_{1}},\dots,Y_{f(\bsig)}^{z_{f(\bsig)}}),
\end{align*}
where the exponents
$z_{i}$ are either two or three, depending on the type. Here, $f(\blam)$ is the number of dots on the steady strings, or the possible nonzero
positions in $\Finch$. Thus, both claims follow
from \autoref{T:SandwichMain}(b).
\end{proof}

For the next statement, recall the notions
of \emph{quasi-hereditary} from \cite{ClPaSc-h-weight-qh} and \emph{affine quasi-hereditary} as in \cite{Kl-affine-quasi-hereditary}.

\begin{Proposition}\label{P:ConsequencesQuasiHereditary}
Let $\ring$ be a field and suppose that $\CrystalVertices=\dCrystalVertices$. The algebras $\WA[n](X)$ and $\WAc[n](X)$ are affine quasi
hereditary algebras. If $\quiver$ is of type $\typea[e]$, then $\WAc[n](X)$ is quasi-hereditary.
\end{Proposition}

\begin{proof}
That $\WA[n](X)$ is affine quasi
hereditary follows by combining \autoref{T:ConstructionMain} and \autoref{P:ConsequencesSimples}, and similarly for $\WAc[n](X)$. For the second claim, the assumption that $\quiver$ is of type $\typea[e]$ implies that the sandwiched
algebras are isomorphic to $\ring$, see the proof of \autoref{P:ConsequencesSimples} and combine this with \autoref{C:ConstructionRepeated}.
Thus, $\WAc[n](X)$ is actually an honest cellular algebra. Hence,
the result follows by \autoref{C:ConstructionMainTwo}, \autoref{P:ConsequencesSimples} and \cite[Lemma 2.1]{KoXi-cellular-quasi-hereditary}.
\end{proof}

\begin{Remark}\label{R:ConsequencesQuasiHereditary}
Maintain the assumptions of \autoref{P:ConsequencesQuasiHereditary}.

\begin{enumerate}

\item Since the sandwiched algebras are often nontrivial, the algebra $\WAc[n](X)$ is in general not quasi-hereditary.
%This follows from \autoref{T:SandwichMain}(b), which adds
%``dummy cells'' with no associated simple $\WAc[n](X)$-modules for
%the nontrivial sandwiched algebras.

\item By \autoref{P:ConsequencesQuasiHereditary}, $\WA[n](X)$ is affine quasi-hereditary. This generalizes \cite[Main theorem]{KlLo-klr-affine-cellular-finite-type}, which proves that the infinite dimensional KLR algebras of finite type are affine quasi-hereditary algebras.

\item In affine type $A$, \cite[Corollary 2.26]{We-rouquier-dia-algebra} and
\cite[Corollary 6.24]{Bo-many-cellular-structures} prove that $\WAc(X)$ is quasi-heredity. Using a similar argument to that above, \cite[Corollary 5C.2]{MaTu-klrw-algebras} gives another proof that $\WAc(X)$ is quasi-hereditary in this case.
\end{enumerate}
\end{Remark}

\begin{Proposition}\label{P:ConsequencesMorita}
Let $\ring$ be a field. Then $\WAc[n](X)$ and $\TAc$ are graded Morita equivalent.
\end{Proposition}

\begin{proof}
By \cite[Theorem 7.5]{LaVa-crystals-cat} the simple $\TAc$-modules are indexed by the crystal graph $\crystalgraph[\Lambda]$. Therefore, by \autoref{P:ConsequencesSimples}, the isomorphism of \autoref{P:ProofsKLRIsomorphism} gives a truncation functor that maps simple modules to simple modules, which makes it an equivalence.
\end{proof}

%%%%%%%%%%%%%%%%%%%%%%%%%%%%%%%%%%%%%%%%%

\subsection{Semisimplicity}\label{SS:ConsequencesSemisimple}

%%%%%%%%%%%%%%%%%%%%%%%%%%%%%%%%%%%%%%%%%

We will now show that $\WAc[n](X)$ is rarely semisimple outside of $\typea[e]$.

\begin{Definition}\label{D:ConsequencesEasyCrystals}
A crystal graph $\crystalgraph[{\fweight[i]}]$ is \emph{entirely semisimple} if $\Lambda_{i}$ is in the following list:
\begin{center}
\rowcolors{2}{mygray!25}{}
\begin{tabular}{c|c}
\rowcolor{mygray!75}
Type & $\fweight[i]$
\\
$\typea[e]$ &  $\fweight[i]$, $i\in\set{1,\dots,e}$
\\
$\typeb[e>1]$ & $\fweight[1]$, $\fweight[e]$
\\
$\typec[e>2]$ & $\fweight[e]$
\\
$\typed[e>3]$ & $\fweight[1]$, $\fweight[{e-1}]$, $\fweight[e]$
\\
$\typee[6]$ & $\fweight[1]$, $\fweight[6]$
\\
$\typee[7]$ & $\fweight[7]$
\\
$\typee[8]$ & None
\\
$\typef[4]$ & None
\\
$\typeg[2]$ & None
\end{tabular}
.
\end{center}
\end{Definition}

With the exception of $\fweight[e]$ in type $\typeb[e>1]$, the entirely semisimple crystal graphs are the crystal graphs for minuscule weights ($\fweight[1]$ is not minuscule in type $\typeb[e>1]$).

\begin{Proposition}\label{P:ConsequencesEasyCrystals}
Let $\ring$ be a field.
\begin{enumerate}

\item Let $\ell=1$ (so that we have a level one crystal graph
$\crystalgraph[{\fweight[i]}]$). The algebra $\WAc[n](X)$ is semisimple for all $n\geq 0$
if and only if $\crystalgraph[{\fweight[i]}]$ is an entirely semisimple crystal graph.

\item In general, the algebra $\WAc[n](X)$ is semisimple for all $n\geq 0$ if and only if the components of $\bfweight$ are pairwise distinct entirely semisimple crystal graphs.

\end{enumerate}
\end{Proposition}

\begin{proof}
We will use \autoref{T:SandwichMain}.(c) silently below.

\textit{(a).}
By the proof of \autoref{P:ConsequencesSimples} and
\autoref{P:ConsequencesQuasiHereditary} it follows
that $\WAc[n](X)$ is semisimple if and only if the
all sandwich algebras are trivial and the paths in $\crystalgraph[{\fweight[i]}]$ of length $n$ are in different blocks, in the sense of \autoref{D:FiniteTypesPermutations}.
For the duration of this proof, we call the first condition the sandwich property and the second one the block property.

Suppose that $\crystalgraph[{\fweight[i]}]$ is entirely semisimple.
In \autoref{D:ConstructionChoice} we fixed a choice of paths, say $\bpath_{1},\dots,\bpath_{k}$, for the vertices
of $\crystalgraph[{\fweight[i]}]$.
For entirely semisimple crystal graphs, the paths $\Res[\bpath_{1}],\dots,\Res[\bpath_{k}]$ are never related by face permutations by the explicit tableau realization of these crystal graphs outlined in the proof of \autoref{P:FiniteTypesSSR} and described in detail, for example, in \cite[Section 5.4]{BuSc-crystal-bases}. For $\typea[e]$ this is immediate from \autoref{P:FiniteTypesSSR}.(a), which says that the crystal graph $\crystalgraph[{\fweight[i]}]$ only contains nonadjacent squares, whereas for types $\typeb[e>1]$, $\typec[c>2]$ and $\typed[e>3]$ only one box tableaux appear in the entirely semisimple crystal graphs, so they also contain only nonadjacent squares. The type $\typee[6]$ and $\typee[7]$ crystal graphs can be inspected, for example using SageMath~\cite{sage} or \cite[Section 5.6]{BuSc-crystal-bases}. Hence, the block property holds for all entirely semisimple crystal graphs.
Moreover, it also follows from these descriptions of the entirely semisimple crystal graphs, and our construction of the bases of $\WAc(X)$ in \autoref{SS:ConstructionBases}, that the sandwiched algebras are trivial in these cases. Therefore, if $\crystalgraph[{\fweight[i]}]$ is semisimple, then $\WAc(X)$ is semisimple for $n\geq 0$.

Conversely, suppose that $\crystalgraph[{\fweight[i]}]$ is not entirely semisimple.
In classical types, except for $\fweight[1]$ in $\typec[3]$, the tableau model for the crystals shows that $\crystalgraph[{\fweight[i]}]$ contains
a $ijji=jiij$ octagon, so the block property fails since
we have a partial face permutation to different vertices in the octagon.
In types $\typef[4]$ and the $\typee[k]$ the block property fails by observation (we used SageMath~\cite{sage}).
For type $\typeg[2]$ we can use \autoref{Eq:FiniteTypesG2}: in $\crystalgraph[{\fweight[1]}]$
the block property fails while in $\crystalgraph[{\fweight[2]}]$ the sandwich property does not hold. The remaining case is $\fweight[1]$ for type $\typec[3]$,
which has the crystal graph:
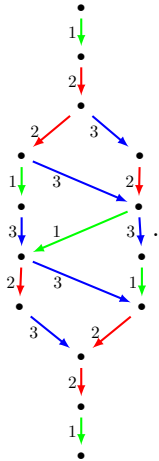
\begin{gather*}
\scalebox{0.65}{$\begin{tikzpicture}[anchorbase,>=latex,line join=bevel,xscale=-0.6,yscale=0.4,every path/.style={very thick}]
\node (node_0) at (63.657bp,295.7bp) [draw,draw=none] {$\bullet$};
\node (node_10) at (63.657bp,223.1bp) [draw,draw=none] {$\bullet$};
\node (node_1) at (180.66bp,223.1bp) [draw,draw=none] {$\bullet$};
\node (node_5) at (121.66bp,150.51bp) [draw,draw=none] {$\bullet$};
\node (node_2) at (65.657bp,440.89bp) [draw,draw=none] {$\bullet$};
\node (node_8) at (66.657bp,368.29bp) [draw,draw=none] {$\bullet$};
\node (node_3) at (178.66bp,440.89bp) [draw,draw=none] {$\bullet$};
\node (node_9) at (178.66bp,368.29bp) [draw,draw=none] {$\bullet$};
\node (node_4) at (178.66bp,295.7bp) [draw,draw=none] {$\bullet$};
\node (node_7) at (121.66bp,78.906bp) [draw,draw=none] {$\bullet$};
\node (node_6) at (121.66bp,8.3018bp) [draw,draw=none] {$\bullet$};
\node (node_11) at (121.66bp,655.69bp) [draw,draw=none] {$\bullet$};
\node (node_12) at (121.66bp,585.09bp) [draw,draw=none] {$\bullet$};
\node (node_13) at (121.66bp,513.49bp) [draw,draw=none] {$\bullet$};
\draw [green,->] (node_0) ..controls (63.657bp,276.73bp) and (63.657bp,256.51bp)  .. (node_10);
\draw (72.157bp,259.4bp) node {$1$};
\draw [blue,->] (node_1) ..controls (164.46bp,202.73bp) and (147.18bp,182.05bp)  .. (node_5);
\draw (167.16bp,186.8bp) node {$3$};
\draw [red,->] (node_2) ..controls (65.922bp,421.15bp) and (66.191bp,402.15bp)  .. (node_8);
\draw (75.157bp,404.59bp) node {$2$};
\draw [blue,->] (node_3) ..controls (148.7bp,421.01bp) and (112.48bp,398.18bp)  .. (node_8);
\draw (145.16bp,404.59bp) node {$3$};
\draw [green,->] (node_3) ..controls (178.66bp,421.92bp) and (178.66bp,401.7bp)  .. (node_9);
\draw (187.16bp,404.59bp) node {$1$};
\draw [red,->] (node_4) ..controls (179.19bp,275.96bp) and (179.73bp,256.96bp)  .. (node_1);
\draw (189.16bp,259.4bp) node {$2$};
\draw [blue,->] (node_4) ..controls (145.64bp,274.43bp) and (107.99bp,251.32bp)  .. (node_10);
\draw (144.16bp,259.4bp) node {$3$};
\draw [red,->] (node_5) ..controls (121.66bp,130.55bp) and (121.66bp,111.63bp)  .. (node_7);
\draw (130.16bp,114.21bp) node {$2$};
\draw [green,->] (node_7) ..controls (121.66bp,60.654bp) and (121.66bp,41.582bp)  .. (node_6);
\draw (130.16bp,43.604bp) node {$1$};
\draw [blue,->] (node_8) ..controls (65.855bp,348.44bp) and (65.035bp,329.13bp)  .. (node_0);
\draw (74.157bp,332.0bp) node {$3$};
\draw [green,->] (node_8) ..controls (98.308bp,347.34bp) and (133.58bp,325.11bp)  .. (node_4);
\draw (145.16bp,332.0bp) node {$1$};
\draw [blue,->] (node_9) ..controls (178.66bp,349.66bp) and (178.66bp,330.34bp)  .. (node_4);
\draw (187.16bp,332.0bp) node {$3$};
\draw [red,->] (node_10) ..controls (78.657bp,203.84bp) and (95.944bp,182.8bp)  .. (node_5);
\draw (108.16bp,186.8bp) node {$2$};
\draw [green,->] (node_11) ..controls (121.66bp,637.44bp) and (121.66bp,618.37bp)  .. (node_12);
\draw (130.16bp,620.39bp) node {$1$};
\draw [red,->] (node_12) ..controls (121.66bp,566.7bp) and (121.66bp,547.66bp)  .. (node_13);
\draw (130.16bp,549.79bp) node {$2$};
\draw [blue,->] (node_13) ..controls (106.37bp,493.22bp) and (90.195bp,472.83bp)  .. (node_2);
\draw (109.16bp,477.19bp) node {$3$};
\draw [red,->] (node_13) ..controls (137.41bp,492.98bp) and (154.38bp,471.95bp)  .. (node_3);
\draw (166.16bp,477.19bp) node {$2$};
\end{tikzpicture}$}
.
\end{gather*}
The result follows because the consecutive $2$-strings have a sandwich dot.
Hence, if $\crystalgraph[{\fweight[i]}]$ is not entirely semisimple, then $\WAc(X)$ is not semisimple for some $n$.

\textit{(b).} If the components of $\bfweight$ are pairwise distinct, entirely semisimple crystal graphs,
then the argument of part (a) proves that $\WAc[n](X)$ is semisimple for all $n$.
Conversely, again by (a), having a component of $\bfweight$ that is not entirely semisimple implies that
$\WAc[n](X)$ is not semisimple. Having repeated components also implies
that $\WAc[n](X)$ is not semisimple since the block property fails.
\end{proof}

Except in the case of $\fweight[1]$ for type $\typec[3]$ and in type $\typeg[2]$, the proof of \autoref{P:ConsequencesEasyCrystals} shows that $\WAc(X)$ is semisimple for all $n$ if and only if $\crystalgraph[{\fweight[i]}]$ does not contain an octagon.

For the next statement we need the following lower and upper
bound formulas.

\begin{Definition}\label{D:ConsequencesSemisimple}
We define $l(i)$, for $i\in\vertices$ as follows:
\begin{center}
\rowcolors{2}{mygray!25}{}
\begin{tabular}{c|c}
\rowcolor{mygray!75}
Type & $l(i)$ for $i\in\vertices$
\\
$\typea[e]$ & $\infty$
\\
$\typeb[e>1]$ & $\infty$ for $i\in\{1,e\}$ and $2e-2i+3$ otherwise
\\
$\typec[e>2]$ & $8$ for $i=1$, $2i$ for $i\in\{2,\dots,e-1\}$, and $\infty$ if $i=e$
\\
$\typed[e>3]$ & $\infty$ for $i\in\{1,e,e-1\}$ and $2e-2i+1$ otherwise
\\
$\typee[6]$ & $\infty,7,5,11,7,\infty$
\\
$\typee[7]$ & $17,7,5,11,7,9,\infty$
\\
$\typee[8]$ & $17,7,5,11,7,9,11,29$
\\
$\typef[4]$ & $11,5,4,8$
\\
$\typeg[2]$ & $7,5$
\end{tabular}
,
\end{center}
where we list $l(i)$ starting from $i=1$ for the exceptional types.
We also set $u(i)=n(i)-l(i)$, where $n(i)$ is as in \autoref{D:ConsequencesStrings}.
\end{Definition}

By convention, $\set{l(i),l(i)+1,\dots,u(i)}$ is empty if $l(i)>u(i)$. This happens only if $i=3$ in type $\typec[3]$, or $i=2$ in type $\typeg[2]$. If $l(i)=\infty$, then the crystal graph $\crystalgraph[{\fweight[i]}]$ is entirely semisimple.

The following statement is in level one.

\begin{Proposition}\label{P:ConsequencesSemisimple}
Let $\ring$ be a field and $\ell=1$.
Let $\crystalgraph[{\fweight[i]}]$ be a crystal graph that is not entirely semisimple and suppose that $n\notin\set{l(i),l(i)+1,\dots,u(i)}$.
Then every block of $\WAc[n](X)$ contains a unique simple module.
\end{Proposition}

\begin{proof}
The simple
$\WAc[n](X)$-modules correspond to the vertices of $\crystalgraph[{\fweight[i]}]$
by \autoref{P:ConsequencesSimples}.  We need to show that
if $n\notin\set{l(i),l(i)+1,\dots,u(i)}$, then the vertices $\alpha\in\crystalgraph[{\fweight[i]}]$
with $\dist[\alpha]=n$ belong to distinct permutation classes of residue sequences.

By \autoref{L:FiniteTypesFlip} and \autoref{P:ConsequencesStrings}
we only need to prove that the lower bound holds.
Moreover, for $i=1$ the lower bound
follows from \autoref{E:FiniteTypesCrystalList}, while the
case $i=3$ in type $\typec[3]$ is easily verified, so we assume
that we are not in these cases.

To verify the lower bound in the classical types we will use the Young diagrams combinatorics as in the proof of \autoref{P:FiniteTypesSSR}.
We consider only type $\typeb[e>1]$, with types $\typec[e>2]$ and $\typed[e>3]$ being similar. The first situation where it is
possible to have two vertices $\lambda,\mu\in\crystalgraph$ with
$\lambda\neq\mu$ that have two paths $\bpath$ and $\cpath$ whose
residue sequences are permutations of one another is when the last two nodes in the column tableau of shape $(1^{i})$ are:
\begin{gather*}
i\neq 1,e\colon
\lambda\leftrightsquigarrow
\begin{ytableau}
i
\\
\scalebox{0.9}{$\overline{i}$}
\end{ytableau}
\text{ and }
\mu\leftrightsquigarrow
\begin{ytableau}
\scalebox{0.55}{$i{+}1$}
\\
\scalebox{0.55}{$\overline{i{+}1}$}
\end{ytableau}
.
\end{gather*}
Counting, using the tableau crystal graphs described in \autoref{Eq:ProofsCrystals},
shows that this happens at precisely $l(i)=2e-2i+3$ steps, for each $i\in\vertices$. The other classical cases are similar and are omitted.

The exceptional types were verified with the code that can be found in
\cite{MaTu-sagemath-finite-type-klrw}.
\end{proof}

\begin{Example}\label{E:ConsequencesSemisimple}
Consider the crystal graph $\crystalgraph[{\fweight[2]}]$ in
type $\typeb[4]$:
\begin{gather*}
\scalebox{0.6}{$\begin{tikzpicture}[anchorbase,>=latex,line join=bevel,xscale=0.75,yscale=0.4,every path/.style={very thick}]
\node (node_0) at (387.5bp,438.0bp) [draw,draw=none] {$\bullet$};
\node (node_11) at (295.5bp,366.0bp) [draw,draw=none] {$\bullet$};
\node (node_20) at (436.5bp,366.0bp) [draw,draw=none] {$\bullet$};
\node (node_1) at (207.5bp,798.0bp) [draw,draw=none] {$\bullet$};
\node (node_2) at (292.5bp,726.0bp) [draw,draw=none] {$\bullet$};
\node (node_4) at (119.5bp,726.0bp) [draw,draw=none] {$\bullet$};
\node (node_9) at (239.5bp,654.0bp) [draw,draw=none] {$\bullet$};
\node (node_33) at (436.5bp,654.0bp) [draw,draw=none] {$\bullet$};
\node (node_3) at (274.5bp,869.5bp) [draw,draw=none] {$\bullet$};
\node (node_10) at (358.5bp,798.0bp) [draw,draw=none] {$\bullet$};
\node (node_7) at (119.5bp,654.0bp) [draw,draw=none] {$\bullet$};
\node (node_5) at (392.5bp,150.5bp) [draw,draw=none] {$\bullet$};
\node (node_13) at (350.5bp,79.5bp) [draw,draw=none] {$\bullet$};
\node (node_6) at (387.5bp,582.0bp) [draw,draw=none] {$\bullet$};
\node (node_25) at (387.5bp,510.0bp) [draw,draw=none] {$\bullet$};
\node (node_29) at (105.5bp,582.0bp) [draw,draw=none] {$\bullet$};
\node (node_30) at (208.5bp,582.0bp) [draw,draw=none] {$\bullet$};
\node (node_8) at (208.5bp,510.0bp) [draw,draw=none] {$\bullet$};
\node (node_23) at (208.5bp,438.0bp) [draw,draw=none] {$\bullet$};
\node (node_19) at (457.5bp,726.0bp) [draw,draw=none] {$\bullet$};
\node (node_15) at (326.5bp,294.0bp) [draw,draw=none] {$\bullet$};
\node (node_31) at (157.5bp,294.0bp) [draw,draw=none] {$\bullet$};
\node (node_12) at (105.5bp,438.0bp) [draw,draw=none] {$\bullet$};
\node (node_14) at (157.5bp,366.0bp) [draw,draw=none] {$\bullet$};
\node (node_24) at (350.5bp,8.5bp) [draw,draw=none] {$\bullet$};
\node (node_28) at (237.5bp,222.0bp) [draw,draw=none] {$\bullet$};
\node (node_32) at (392.5bp,222.0bp) [draw,draw=none] {$\bullet$};
\node (node_16) at (480.5bp,294.0bp) [draw,draw=none] {$\bullet$};
\node (node_17) at (316.5bp,1011.5bp) [draw,draw=none] {$\bullet$};
\node (node_22) at (316.5bp,940.5bp) [draw,draw=none] {$\bullet$};
\node (node_18) at (525.5bp,582.0bp) [draw,draw=none] {$\bullet$};
\node (node_26) at (536.5bp,510.0bp) [draw,draw=none] {$\bullet$};
\node (node_21) at (358.5bp,869.5bp) [draw,draw=none] {$\bullet$};
\node (node_35) at (523.5bp,438.0bp) [draw,draw=none] {$\bullet$};
\node (node_27) at (53.5bp,510.0bp) [draw,draw=none] {$\bullet$};
\node (node_34) at (309.5bp,150.5bp) [draw,draw=none] {$\bullet$};
\draw [red,->] (node_0) ..controls (361.96bp,417.56bp) and (332.69bp,395.3bp)  .. (node_11);
\draw (361.0bp,402.0bp) node {$2$};
\draw [black,->] (node_0) ..controls (400.62bp,418.26bp) and (414.85bp,397.92bp)  .. (node_20);
\draw (426.0bp,402.0bp) node {$4$};
\draw [blue,->] (node_1) ..controls (231.27bp,777.43bp) and (257.34bp,755.96bp)  .. (node_2);
\draw (269.0bp,762.0bp) node {$1$};
\draw [black,->] (node_1) ..controls (182.66bp,777.24bp) and (155.03bp,755.27bp)  .. (node_4);
\draw (183.0bp,762.0bp) node {$4$};
\draw [black,->] (node_2) ..controls (277.94bp,705.78bp) and (262.4bp,685.24bp)  .. (node_9);
\draw (281.0bp,690.0bp) node {$4$};
\draw [red,->] (node_2) ..controls (333.84bp,704.9bp) and (380.99bp,681.98bp)  .. (node_33);
\draw (391.0bp,690.0bp) node {$2$};
\draw [black,->] (node_3) ..controls (256.38bp,849.7bp) and (236.41bp,828.99bp)  .. (node_1);
\draw (258.0bp,834.0bp) node {$4$};
\draw [blue,->] (node_3) ..controls (297.76bp,849.26bp) and (324.3bp,827.3bp)  .. (node_10);
\draw (335.0bp,834.0bp) node {$1$};
\draw [green,->] (node_4) ..controls (119.5bp,706.57bp) and (119.5bp,687.05bp)  .. (node_7);
\draw (128.0bp,690.0bp) node {$3$};
\draw [blue,->] (node_4) ..controls (153.36bp,705.25bp) and (193.07bp,682.08bp)  .. (node_9);
\draw (203.0bp,690.0bp) node {$1$};
\draw [blue,->] (node_5) ..controls (381.32bp,131.13bp) and (369.29bp,111.36bp)  .. (node_13);
\draw (385.0bp,115.0bp) node {$1$};
\draw [green,->] (node_6) ..controls (387.5bp,562.68bp) and (387.5bp,543.46bp)  .. (node_25);
\draw (396.0bp,546.0bp) node {$3$};
\draw [red,->] (node_7) ..controls (115.81bp,634.57bp) and (111.91bp,615.05bp)  .. (node_29);
\draw (123.0bp,618.0bp) node {$2$};
\draw [blue,->] (node_7) ..controls (144.21bp,633.56bp) and (172.52bp,611.3bp)  .. (node_30);
\draw (184.0bp,618.0bp) node {$1$};
\draw [red,->] (node_8) ..controls (208.5bp,490.2bp) and (208.5bp,470.79bp)  .. (node_23);
\draw (217.0bp,474.0bp) node {$2$};
\draw [red,->] (node_9) ..controls (281.31bp,633.23bp) and (331.74bp,609.37bp)  .. (node_6);
\draw (340.0bp,618.0bp) node {$2$};
\draw [green,->] (node_9) ..controls (231.24bp,634.36bp) and (222.35bp,614.27bp)  .. (node_30);
\draw (237.0bp,618.0bp) node {$3$};
\draw [black,->] (node_10) ..controls (340.64bp,778.05bp) and (320.92bp,757.15bp)  .. (node_2);
\draw (342.0bp,762.0bp) node {$4$};
\draw [red,->] (node_10) ..controls (386.14bp,777.46bp) and (418.05bp,754.9bp)  .. (node_19);
\draw (428.0bp,762.0bp) node {$2$};
\draw [black,->] (node_11) ..controls (301.22bp,348.54bp) and (306.75bp,333.5bp)  .. (312.5bp,321.0bp) .. controls (313.87bp,318.02bp) and (315.43bp,314.91bp)  .. (node_15);
\draw (321.0bp,330.0bp) node {$4$};
\draw [blue,->] (node_11) ..controls (256.72bp,345.33bp) and (210.31bp,321.79bp)  .. (node_31);
\draw (252.0bp,330.0bp) node {$1$};
\draw [red,->] (node_12) ..controls (119.58bp,418.04bp) and (135.14bp,397.1bp)  .. (node_14);
\draw (147.0bp,402.0bp) node {$2$};
\draw [red,->] (node_13) ..controls (350.5bp,60.442bp) and (350.5bp,41.496bp)  .. (node_24);
\draw (359.0bp,44.0bp) node {$2$};
\draw [green,->] (node_14) ..controls (157.5bp,346.57bp) and (157.5bp,327.05bp)  .. (node_31);
\draw (166.0bp,330.0bp) node {$3$};
\draw [blue,->] (node_15) ..controls (301.62bp,273.43bp) and (274.32bp,251.96bp)  .. (node_28);
\draw (302.0bp,258.0bp) node {$1$};
\draw [black,->] (node_15) ..controls (344.83bp,273.56bp) and (364.73bp,252.46bp)  .. (node_32);
\draw (376.0bp,258.0bp) node {$4$};
\draw [red,->] (node_16) ..controls (456.13bp,273.62bp) and (428.33bp,251.5bp)  .. (node_32);
\draw (455.0bp,258.0bp) node {$2$};
\draw [red,->] (node_17) ..controls (316.5bp,992.44bp) and (316.5bp,973.5bp)  .. (node_22);
\draw (325.0bp,976.0bp) node {$2$};
\draw [black,->] (node_18) ..controls (528.42bp,562.43bp) and (531.38bp,543.61bp)  .. (node_26);
\draw (540.0bp,546.0bp) node {$4$};
\draw [black,->] (node_19) ..controls (451.97bp,706.58bp) and (446.12bp,687.08bp)  .. (node_33);
\draw (457.0bp,690.0bp) node {$4$};
\draw [red,->] (node_20) ..controls (405.33bp,345.17bp) and (370.46bp,322.98bp)  .. (node_15);
\draw (404.0bp,330.0bp) node {$2$};
\draw [black,->] (node_20) ..controls (448.52bp,345.88bp) and (461.25bp,325.63bp)  .. (node_16);
\draw (472.0bp,330.0bp) node {$4$};
\draw [green,->] (node_21) ..controls (358.5bp,850.2bp) and (358.5bp,830.82bp)  .. (node_10);
\draw (367.0bp,834.0bp) node {$3$};
\draw [green,->] (node_22) ..controls (305.32bp,921.13bp) and (293.29bp,901.36bp)  .. (node_3);
\draw (309.0bp,905.0bp) node {$3$};
\draw [blue,->] (node_22) ..controls (327.68bp,921.13bp) and (339.71bp,901.36bp)  .. (node_21);
\draw (351.0bp,905.0bp) node {$1$};
\draw [green,->] (node_23) ..controls (232.59bp,417.62bp) and (260.08bp,395.5bp)  .. (node_11);
\draw (271.0bp,402.0bp) node {$3$};
\draw [blue,->] (node_23) ..controls (192.87bp,424.45bp) and (185.73bp,417.85bp)  .. (180.5bp,411.0bp) .. controls (174.18bp,402.72bp) and (168.7bp,392.45bp)  .. (node_14);
\draw (189.0bp,402.0bp) node {$1$};
\draw [green,->] (node_25) ..controls (387.5bp,490.2bp) and (387.5bp,470.79bp)  .. (node_0);
\draw (396.0bp,474.0bp) node {$3$};
\draw [black,->] (node_26) ..controls (533.05bp,490.43bp) and (529.56bp,471.61bp)  .. (node_35);
\draw (540.0bp,474.0bp) node {$4$};
\draw [blue,->] (node_27) ..controls (67.781bp,489.78bp) and (83.035bp,469.24bp)  .. (node_12);
\draw (95.0bp,474.0bp) node {$1$};
\draw [black,->] (node_28) ..controls (257.86bp,201.35bp) and (279.4bp,180.56bp)  .. (node_34);
\draw (291.0bp,186.0bp) node {$4$};
\draw [blue,->] (node_29) ..controls (91.503bp,562.16bp) and (76.183bp,541.53bp)  .. (node_27);
\draw (95.0bp,546.0bp) node {$1$};
\draw [red,->] (node_30) ..controls (208.5bp,562.68bp) and (208.5bp,543.46bp)  .. (node_8);
\draw (217.0bp,546.0bp) node {$2$};
\draw [black,->] (node_31) ..controls (179.39bp,273.84bp) and (203.95bp,252.36bp)  .. (node_28);
\draw (216.0bp,258.0bp) node {$4$};
\draw [green,->] (node_32) ..controls (392.5bp,202.7bp) and (392.5bp,183.32bp)  .. (node_5);
\draw (401.0bp,186.0bp) node {$3$};
\draw [blue,->] (node_32) ..controls (369.64bp,201.86bp) and (343.77bp,180.2bp)  .. (node_34);
\draw (370.0bp,186.0bp) node {$1$};
\draw [black,->] (node_33) ..controls (423.04bp,633.78bp) and (408.67bp,613.24bp)  .. (node_6);
\draw (426.0bp,618.0bp) node {$4$};
\draw [green,->] (node_33) ..controls (461.38bp,633.43bp) and (488.68bp,611.96bp)  .. (node_18);
\draw (500.0bp,618.0bp) node {$3$};
\draw [green,->] (node_34) ..controls (320.41bp,131.13bp) and (332.16bp,111.36bp)  .. (node_13);
\draw (344.0bp,115.0bp) node {$3$};
\draw [green,->] (node_35) ..controls (499.17bp,417.43bp) and (472.49bp,395.96bp)  .. (node_20);
\draw (498.0bp,402.0bp) node {$3$};
\end{tikzpicture}$}
.
\end{gather*}
In this case we have $l(2)=7=u(2)$ and indeed, only in the middle
of the crystal graph do we find different vertices with residue sequences that are related by face permutations. For example, there is a nontrivial face permutation between the paths with residue sequences $2344321$ and $2344312$. As a more complicated example, there is a nontrivial face permutation between the paths with residue sequences $2344321$ and $2132434$.

Note that the sandwiched algebras
are not necessarily trivial because having only one block does not
imply that $\WAc[n](X)$ is semisimple. Explicitly,
if $n=3$, then every vertex corresponds to a unique block, but
the idempotent $\idem[\ppath]$ for the path with residue sequence $\Res[\lambda]=234$ can get a sandwich dot on the solid $4$-string:
\begin{gather*}
\fsandwich{\lambda}{f}=\biggl\{
\begin{tikzpicture}[anchorbase,smallnodes,yscale=0.6,xscale=1.1]
\draw[redstring](0,0)node[below]{\,$2$}--+(0,1);
\draw[solid](-0.08,0)node[below]{$2$\,}--+(0,1);
\draw[ghost](0.92,0)--+(0,1)node[above]{$2$};
\draw[solid](0.84,0)node[below]{$3$}--+(0,1);
\draw[ghost](1.84,0)--+(0,1)node[above]{$3$};
\draw[solid](1.76,0)node[below]{$4$}--+(0,1);
\end{tikzpicture}
,\quad
\begin{tikzpicture}[anchorbase,smallnodes,yscale=0.6,xscale=1.1]
\draw[redstring](0,0)node[below]{\,$2$}--+(0,1);
\draw[solid](-0.08,0)node[below]{$2$\,}--+(0,1);
\draw[ghost](0.92,0)--+(0,1)node[above]{$2$};
\draw[solid](0.84,0)node[below]{$3$}--+(0,1);
\draw[ghost](1.84,0)--+(0,1)node[above]{$3$};
\draw[solid,dot=0.5](1.76,0)node[below]{$4$}--+(0,1);
\end{tikzpicture}
\biggr\}
.
\end{gather*}
In this case, the sandwiched algebra is isomorphic to $R[Y]/(Y^{2})$, so
$\WAc[3](X)$ is not semisimple.
\end{Example}

\begin{Example}\label{E:ConsequencesSemisimpleTwo}
Using the code from \cite{MaTu-sagemath-finite-type-klrw}, \autoref{P:ConsequencesSemisimple} for $\typee[6]$ was verified using the code:
\begin{sagemath}
sage: ' '.join(f"{CheckForRepeatedPositiveRoots(['E',6],[i], minlength=True)}"
....:          for i in range(1,7))
+Infinity 11 7 5 7 +Infinity
\end{sagemath}
Adjusting for the differences in different labeling conventions between \autoref{SS:FiniteTypesConventions} and SageMath, this agrees with \autoref{D:ConsequencesSemisimple}.
\end{Example}

%%%%%%%%%%%%%%%%%%%%%%%%%%%%%%%%%%%%%%%%%

\subsection{Projective modules in level one}\label{SS:ConsequencesProj}

%%%%%%%%%%%%%%%%%%%%%%%%%%%%%%%%%%%%%%%%%

Throughout this and the next section we work only with the cyclotomic quotients $\WAc(X)$ in level one such that $\CrystalVertices=\dCrystalVertices$.

For $\blam\in\CrystalVertices$ set $\RPlam=\WAc[n](X)\idem[\blam]$, which is a projective $\WAc[n](X)$-module, and let $\RDlam$ and $\RLlam$
be the corresponding cell module and graded simple with apex $\blam$, respectively.
%Similarly, let $\RPlam$, $\RDlam$ and $\RLlam$ be the corresponding $\WAc(X)$-modules.

% For $\blam\in\aCrystalVertices$ set $\WPlam=\WA[n](X)\idem[\blam]$, which is a projective $\WA[n](X)$-module, and let $\WDlam$ and $\WLlam$
% be the corresponding cell module and graded simple with apex $\blam$, respectively. Similarly, let $\RPlam$, $\RDlam$ and $\RLlam$ be the corresponding $\WAc(X)$-modules.

Since $\idem[\blam]$ is an idempotent, as graded algebras
$\End_{\WAc}(\RPlam)\cong\idem[\blam]\WAc\idem[\blam]$. Each choice of $\bS,\bT\in\Face(\blam)$ and $\bfi\in\Finch$ defines an endomorphism of these modules, given by
$\theta^{\bfi}_{\bS\bT}(\idem[\blam])=D^{\bfi}_{\bS\bT}$.

% Since $\idem[\blam]$ is an idempotent, as graded algebras
% $\End_{\aWA}(\WPlam)\cong\idem[\blam]\aWA\idem[\blam]$ and $\End_{\WAc(X)}(\RPlam)\cong\idem[\blam]\WAc(X)\idem[\blam]$. Each choice of $\bS,\bT\in\Face(\blam)$ and $\ba\in\Affch,\bfi\in\Finch$ defines an endomorphism of these modules, given by
% $\theta^{\ba,\bfi}_{\bS\bT}(\idem[\blam])=D^{\ba,\bfi}_{\bS\bT}$ and
% $\theta^{\bfi}_{\bS\bT}(\idem[\blam])=D^{\ba}_{\bS\bT}$.

\begin{Proposition}\label{P:ConsequencesInde}
Let $\ring$ be a field and set $\bfweight=(\Lambda_{i})$, for some $i$. Assume that $\CrystalVertices=\dCrystalVertices$. %Then:
%\begin{enumerate}

%\item If $\blam\in\aCrystalVertices$, then
%$\End_{\aWA}(\WPlam)\cong\idem[\blam]\aWA\idem[\blam]$, which has an $\ring$-basis
%given by $\set{D_{\bS\bT}^{\ba,\bfi}|\bS,\bT\in\Face(\blam),\ba\in\Affch,\bfi\in\Finch}$.

%\item
If $\blam\in\CrystalVertices$, then $\End_{\WAc(X)}(\RPlam)\cong\idem[\blam]\WAc(X)\idem[\blam]$, which has an $\ring$-basis
given by $\set{D_{\bS\bT}^{\bfi}|\bS,\bT\in\Face(\blam),\bfi\in\Finch}$.
%\end{enumerate}
Moreover, %$\WPlam$ is the projective cover of $\WLlam$ and
$\RPlam$ is the projective cover of $\RLlam$.
\end{Proposition}

\begin{proof}
%We only give the proof for $\aWA$ since the arguments for $\WAc(X)$ work verbatim.

First note that $\RPlam\supset\WAc\dotidem[\blam]$
by \autoref{Eq:ProofsSpanningPermutations}.
The endomorphism algebra of $\RPlam$ is, by construction,
given by all steady diagrams with bottom and top $\idem[\blam]$. Hence,
$\End_{\WAc}(\RPlam)\cong\idem[\blam]\WAc\idem[\blam]$.
Recall that we have classified all possible steady permutations as face permutations in \autoref{P:ProofsPermutationsSteady}. In particular, the detour permutations to $\blam$ and partial face permutations belong to the endomorphism algebra of $\WPlam$.
Moreover, the sandwiched algebras are contained in the endomorphism algebra by construction. By \autoref{T:ConstructionMain}, $P_{\blam}$ has no other endomorphisms, so $\End_{\WAc}(\WPlam)$ is generated as an $\ring$-module by the detour and face permutations, and the sandwich dots, as in the statement of the proposition.

As we are in level one, we can assume that all of the strings in the diagrams are to the left of the red strings by \autoref{L:ProofsKLRIsomorphismTwo}, so all detour permutations in $\End_{\WAc}(\RPlam)$ are invertible, modulo more dominant terms
by \autoref{P:ProofsDetour}. Hence, the endomorphism algebra of $\RPlam$ is a local ring by the explicit description of its basis given by \autoref{P:ProofsPermutationsSteady}. Therefore, the projective module $\RPlam=\WAc[n](X)\idem[\blam]$ is indecomposable. Hence, $\RPlam$ is the projective cover of $\RLlam$ by \autoref{T:ConstructionMain} and \autoref{P:ConsequencesSimples}.
\end{proof}

\begin{Corollary*}\label{C:ConsequencesIndeTwo}
Suppose that $\ring$ is a field, $\quiver$ is a quiver of type $\typea[e]$, $\bfweight=(\Lambda_{i})$ and $\CrystalVertices=\dCrystalVertices$.
Let $\blam\in\CrystalVertices$. Then
%$\End_{\aWA}(\WPlam)$ is positively graded and
$\End_{\WAc(X)}(\RPlam)$ is concentrated in degree $0$.
\end{Corollary*}

\begin{proof}
This follows from \autoref{D:ConstructionPermBasis} because all diagrammatic face permutations are compositions of (partial) face permutations around nonadjacent squares by \autoref{P:FiniteTypesSSR}, since $\quiver$ is of type $\typea[e]$. The associated permutation diagrams are thus of degree $0$. Hence, $\End_{\WAc(x)}(P_{\blam})$ is in degree $0$.
\end{proof}

Recall from \autoref{SS:MainExampleTwo} that $\vpar$ is our grading parameter. Let $[\RDlam:\RLlam[\bmu]]_{\vpar}$ be the \emph{graded decomposition multiplicity} of $\RLlam[\bmu]$ in $\RDlam$.

\begin{Corollary*}\label{C:ConsequencesDecompositionNumbers}
Suppose that $\ring$ is a field, $\CrystalVertices=\dCrystalVertices$, $\bfweight=(\Lambda_{i})$ and let $\blam,\bmu\in\CrystalVertices$.
Then
\begin{gather*}
[\RDlam:\RLlam[\bmu]]_{\vpar}
=\grdim(\idem[\bmu]\RDlam)
=\sum_{\bT\in\Face(\ppath,\ppath[\bmu])}\vpar^{\deg D_{\bT}}.
\end{gather*}
In particular,
$[\RDlam:\RLlam[\bmu]]_{\vpar}$ is independent of the characteristic of $\ring$.
\end{Corollary*}

\begin{proof}
By \autoref{P:ConsequencesInde}, $\RPlam[\bmu]=\WAc(X)\idem[\bmu]$ is the projective cover of $\RLlam[\bmu]$. Since $\idem[\bmu]$ is an idempotent,
\begin{gather*}
[\RDlam:\RLlam[\bmu]]_{\vpar}
=\dim_{\vpar}\Hom_{\WAc(X)}(\RPlam[\bmu],\RDlam)
=\dim_{\vpar}\idem[\bmu]\RDlam.
\end{gather*}
Hence, the result follows from \autoref{T:ConstructionGeneralMain}.
\end{proof}

Recall from \autoref{SS:ProofsGrothendieckRings} that
$L(\bfweight)\cong[\Proj\TAc[\oplus]]$ as $U_{\vpar}(\mathfrak{g})$-modules, where
$[\Proj\TAc[\oplus]]=\bigoplus_{n\geq 0}[\Proj\TAc]$. Hence, by \autoref{P:ConsequencesMorita}, $L(\bfweight)\cong[\Proj\WAc[\oplus]]$ as $U_{\vpar}(\mathfrak{g})$-modules. Henceforth, we identify these two modules. Note that $[\Rep\WAc[\oplus]]$ and $[\Proj\WAc[\oplus]]$ are finite dimensional $\Q(\vpar)$-modules since $\WAc(X)=0$ for $n\gg 0$ by \autoref{P:ConsequencesStrings}.

Let $[\Rep\WAc[\oplus]]=\bigoplus_{n\geq 0}[\Rep\WAc(X)]$. There is a natural inclusion $[\Proj\WAc[\oplus]]\hookrightarrow[\Rep\WAc[\oplus]]$, so we can consider $L(\bfweight)$ as a submodule of $[\Rep\WAc[\oplus]]$. If $M$ is a $\WAc(X)$-module, then let $[M]$ be its image in $[\Rep\WAc[\oplus]]$.

\begin{Lemma}\label{L:ConsequencesCrystal}
The Grothendieck group $[\Rep\WAc[\oplus]]$ has basis
$\set{[\RDlam]|\blam\in\dCrystalVertices}$.
\end{Lemma}

\begin{proof}
The simple $\WAc(X)$-modules, for $n\geq 0$, give a basis of $[\Rep\WAc[\oplus]]$ and the graded decomposition matrices of these algebras are unitriangular by \autoref{C:ConsequencesDecompositionNumbers}.
\end{proof}

Let $\mathscr{L}_{\bfweight}$ be the $\Z[\vpar]$-lattice of $[\Rep\WAc[\oplus]]$ spanned by $\set{[\RDlam]|\blam\in\CrystalVertices}$ and let $\mathscr{B}_{\bfweight}=\set{[\RPlam]+\vpar\mathscr{L}_{\bfweight}|\blam\in\CrystalVertices}$. Then $[\Rep\WAc[\oplus]]\cong\Q(\vpar)\otimes_{\Z[v]}\mathscr{L}_{\bfweight}$ and $\mathscr{B}_{\bfweight}$ is a $\Q$-basis of $\mathscr{L}_{\bfweight}/\vpar\mathscr{L}_{\bfweight}$ by \autoref{L:ConsequencesCrystal}. We will show that $(\mathscr{L}_{\bfweight},\mathscr{B}_{\bfweight})$ is a $0$-crystal base of $L(\bfweight)$ in the sense of \cite{Ka-crystal-bases}. Unlike \cite{Ka-crystal-bases}, we work with lattices over $\Z[\vpar]$. However, by base change, we obtain lattices over the ring of rational functions in $\Q(\vpar)$, which gives the setting used by Kashiwara.

For $i\in\vertices$ let $F^{\bfweight}_{i}$ be the $i$-induction functor, which is given by adding an $i$-string to the left of a diagram.

\begin{Theorem}\label{T:ConsequencesCanonicalBasis}
Suppose that $\ring$ is a field,  $\bfweight=(\Lambda_{i})$, $\CrystalVertices=\dCrystalVertices$, and that $\quiver$ is not of type $\typef[4]$. Let $\blam,\bmu\in\CrystalVertices$.
Then
\begin{gather*}
[\RDlam:\RLlam[\bmu]]_{\vpar}\in
\delta_{\blam\bmu}+\vpar\N[\vpar].
\end{gather*}
Consequently, under categorification $P_{\blam}$ corresponds to a canonical basis element in $L(\bfweight)$.
\end{Theorem}

\begin{proof}
If $\blam=\bmu$, then
$[\RDlam:\RLlam[\bmu]]_{\vpar}=1$ by \autoref{C:ConsequencesDecompositionNumbers}, so there is nothing to prove. If $\blam\neq\bmu$ and $D_{\bT}\in\Face(\ppath,\ppath[\bmu])$, then $D_{\bT}$ is a product of basic face permutations around adjacent squares and octagons together with at least one partial face permutation around an octagon.
This follows from \autoref{P:FiniteTypesSSR} and our assumption that $\quiver$
is not of type $\typef[4]$.
As all of the basic diagrammatic face permutations have nonnegative degree, and the partial diagrammatic face permutations have positive degree, it follows that
$[\RDlam:\RLlam[\bmu]]_{\vpar}\in\delta_{\blam\bmu}+\vpar\N[\vpar]$.

It remains to prove that $\RPlam$ corresponds to the canonical basis under Kang--Kashiwara's \cite{KaKa-categorification-via-klr} categorification of the highest weight module $L(\bfweight)$. By \autoref{P:ProofsKLRIsomorphismTwo} and \autoref{C:ConstructionMainTwo}, if $\bnu\in\CrystalVertices[n-1]$ and $i\in\vertices$, then
\begin{align*}
F^{\bfweight}_{i}[{\RPlam[\bnu]}]&=
\begin{cases*}
[\RPlam]& if $\bnu\xrightarrow{i}\blam$ is a path in $\crystalgraph[\bfweight]$,\\
0& if $\Res[{\ppath[\bnu]}]i$ is not a path in $\crystalgraph[\bfweight]$.
\end{cases*}
\intertext{%
It follows from the usual yoga that $\RPlam$ has a filtration with subquotients $\RDlam[\bmu]\otimes\sand[\mu]$ appearing with multiplicity $[\RDlam[\bmu]:\RLlam]_{\vpar}$ for $\bmu\in\CrystalVertices$. Therefore,
using the positivity property
$[\RDlam:\RLlam[\bmu]]_{\vpar}\in\delta_{\blam\bmu}+\vpar\N[\vpar]$ from the first paragraph, and working modulo $v\mathscr{L}_{\bfweight}$, }
F^{\bfweight}_{i}[\RPlam[\bnu]]&\equiv
\begin{cases*}
[\RDlam]& if $\bnu\xrightarrow{i}\blam$ is a path in $\crystalgraph[\bfweight]$,\\
0&otherwise.
\end{cases*}
\end{align*}
Therefore, $(\mathscr{L}_{\bfweight},\mathscr{B}_{\bfweight})$ is a $0$-crystal base for $L(\bfweight)$ in the sense of Kashiwara~\cite[Definition 2.3.1]{Ka-crystal-bases}. By \cite[Theorem~1]{Ka-global-crystal-bases}, this crystal base uniquely determines the lower global basis $\set{G_{\blam}|\blam\in\CrystalVertices, n\ge0}$, or canonical basis, of $L(\bfweight)$ as the unique bar invariant basis of $L(\bfweight)$ such that $G_{\blam}\equiv\RDlam\pmod{\vpar\mathscr{L}_{\bfweight}}$, for $\blam\in\CrystalVertices$. Hence, $\set{[\RPlam]|\blam\in\CrystalVertices,n\geq 0}$ is the canonical basis of $[\Rep\WAc[\oplus]]$, as claimed.
\end{proof}

\begin{Example}\label{E:ConsequencesPositivity}
The argument in the proof of \autoref{T:ConsequencesCanonicalBasis} does not work for type $\typef[4]$. To see this take the basic face permutation given by
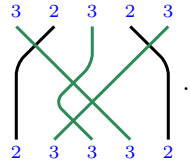
\begin{gather*}
\begin{tikzpicture}[anchorbase,smallnodes,rounded corners]
\draw[solid,black](0,0)node[below]{$2$}--++(0,1)--++(0.5,0.5) node[above,yshift=-1pt]{$2$};
\draw[solid,black](2,0)node[below]{$2$}--++(0,1)--++(-0.5,0.5)node[above,yshift=-1pt]{$2$};
\draw[solid,spinach](0.5,0)node[below]{$3$}--++(1,1)--++(0.5,0.5)node[above,yshift=-1pt]{$3$};
\draw[solid,spinach](1,0)node[below]{$3$}--++(-0.5,0.5)--++(0.5,0.5)--++(0,0.5)node[above,yshift=-1pt]{$3$};
\draw[solid,spinach](1.5,0)node[below]{$3$}--++(-1,1)--++(-0.5,0.5)node[above,yshift=-1pt]{$3$};
\end{tikzpicture}
.
\end{gather*}
This diagram is of degree $-2$. If the idempotent diagram is to the north of
this diagram, then \autoref{D:ConstructionRepeated} does not add extra dots, so the degree is negative.
\end{Example}

\begin{Remark}
The results in this section, and the next, follow from the fact that detour permutations in level one are invertible. This may still be true in level two, but this is difficult to check (compare \autoref{Ex:Octogon}). In general, even in type~$\typea[e]$, detour permutations are not invertible because this would imply that the decomposition matrices of $\WAc[n](X)$ are independent of the characters, contradicting \cite{Wi-analog-james} (see also \cite[Example 3.7.10]{Mathas:Singapore}).
\end{Remark}

\begin{Remark}\label{R:ConsequencesPositivity}
\autoref{E:ConsequencesPositivity} can be extended to show that the proof
\autoref{T:ConsequencesCanonicalBasis} only works in simply laced types for the more general wKLRW algebras that have red strings labeled by arbitrary dominant weights. More, precisely, according to
\cite{St-local-crystal-doubly-laced}, the analog \autoref{P:FiniteTypesSSR} for nonfundamental dominant weights has decagons appearing in types $\typeb[e>1]$, $\typec[e>2]$ and $\typeg[2]$, so the argument in the first paragraph of \autoref{T:ConsequencesCanonicalBasis} does not work for these types. This is consistent with \cite[Theorem 4.4]{VaVa-canonical-bases-klr} and \cite[Paragraph after Corollary B]{We-canonical-bases-higher-rep}.

Tsuchioka~\cite{Ts-non-positive-canonical-basis} shows that the structure constants of the canonical bases can have negative coefficients in nonsimply laced types, which implies that the projective indecomposable modules for the infinite KLR algebras cannot coincide with the Lusztig--Kashiwara canonical basis in general. This does not contradict \autoref{T:ConsequencesCanonicalBasis} because we are looking only at projective modules for wKLRW algebras attached to fundamental weights.
\end{Remark}

\begin{Example}\label{E:ConsequencesPositivityTwo}
Let $\quiver$ be a quiver of type $\typeg[2]$. Set
$\bfweight=(\Lambda_{1},\Lambda_{1},\Lambda_{2})$ and consider the tensor crystal $\Crystalgraph[\bfweight]$.
Given a path $\bpath$ in the crystal graph $\Crystalgraph[\bfweight]$, let $G(\bpath)$ be the corresponding canonical basis element. Tsuchioka~\cite{Ts-non-positive-canonical-basis} computes that
\begin{align*}
F_{2}G(121112211)=&\, G(1211122211)+[2]_{\vpar}\cdot G(1111222211)+G(2111112221)
+[2]_{\vpar}\cdot G(1211112221)\\
&+G(1111122221)
\mathop{\colorbox{tomato!50}{\mystrut -}}G(1112211122)+[2]_{\vpar}\cdot G(1122111122)
\end{align*}
Then $\Crystalgraph[\bfweight]$ has a path $\bpath$ with residue sequence $121112211$ but it does not have a path with residue sequence $1112211122$. Hence, the wKLRW algebras $\WAc(X)$ are not able to detect this negative multiplicity.
\end{Example}

%%%%%%%%%%%%%%%%%%%%%%%%%%%%%%%%%%%%%%%%%

\subsection{Simple modules in level one}\label{SS:ConsequencesDims}

%%%%%%%%%%%%%%%%%%%%%%%%%%%%%%%%%%%%%%%%%

We are still in level one.
For the next proposition recall
from \autoref{P:ConsequencesSimples} that
$\WAc[n](X)$ has simple modules $\RLlam$ indexed
by $\blam\in\CrystalVertices$.

\begin{Proposition}\label{P:ConsequencesDimensions}
Let $\ring$ be a field and supose that $\bfweight=(\Lambda_{i})$ and $\CrystalVertices=\dCrystalVertices$.
%Assume that $\WAc[n](X)$ is semisimple.
Then $\grdim\RLlam$ is equal to the number of detour permutations in $\Face(\blam)$.
\end{Proposition}

\begin{proof}
As noted already in the proof of \autoref{P:ConsequencesInde}, the action of every
a detour permutation $\bT\in\Face(\blam)$ on $\RDlam$ is invertible by \autoref{L:ProofsKLRIsomorphismTwo} and \autoref{P:ProofsDetour}, so its image in $\RLlam$ is nonzero. Conversely, any face permutation $\bT\in\Face(\ppath,\ppath[\bmu])$ for $\bmu\ne\blam$ is in the image of a map from $\RPlam[\bmu]\to\RDlam$, so the image of this permutation in $\RLlam$ is zero. Hence, the detour permutations index a basis of $\RLlam$. Finally, as in the proof of \autoref{T:ConsequencesCanonicalBasis}, the detour permutations have degree zero, so $\RLlam$ is concentrated in degree zero and $\grdim\RLlam$ is equal to the number of detour permutations.
\end{proof}

\begin{Example}\label{E:ConsequencesDimensions}
Let $\quiver$ be a quiver of type $\typea$. Then $\RLlam$ is one dimensional for all $\blam\in\CrystalVertices$. This follows because the only nontrivial faces of $\crystalgraph[\bfweight]$ are nonadjacent squares by \autoref{P:FiniteTypesSSR}, which are excluded by \autoref{D:ConstructionChoice}.(c). By \autoref{Eq:ConstructionTypeC3}, if $\quiver$ is a quiver of type $\typec[3]$ and $\bfweight=\fweight[1]$, then $\WAc[4]$ is semisimple and it has a two dimensional simple module with residue sequence $1232$.
\end{Example}

\begin{Example}
For a quiver of type $\typeb[2]$ and suppose that $\bfweight=\fweight[2]$ as in  \autoref{E:MainExampleTheBeastItself}. If $0\leq n\leq 10$, then the simple modules of $\WAc(X)$ are one dimensional unless $n=7$ in which case they are two dimensional.
\end{Example}

%%%%%%%%%%%%%%%%%%%%%%%%%%%%%%%%%%%%%%%%%

\subsection{Sandwich cellularity in affine types}\label{SS:Affine}

%%%%%%%%%%%%%%%%%%%%%%%%%%%%%%%%%%%%%%%%%

Fix $e\in\Z_{\geq 1}$. The \emph{affine types} are:
\begin{gather*}
\aonetypea:\quad
\begin{tikzpicture}[anchorbase]
\draw[directed=0.5](1,0)--(2,0);
\draw[directed=0.5](2,0)--(3,0);
\draw[directed=0.5](4,0)--(5,0);
\draw[directed=0.5](5,0)--(6,0);
\draw[directed=0.5](6,0) to (6,0.25) to (3.5,0.25);
\draw[directed=0.5](3.5,0.25) to (1,0.25) to (1,0);
\node at (3.5,0){$\cdots$};
\foreach \x in {1,...,3} {
\node[dynkin=\x] (\x) at (\x,0){};
}
\foreach \x [evaluate=\x as \c using {int(6-\x)}] in {4,5} {
\node[dynkin=e{-}\c] (\x) at (\x,0){};
}
\node[dynkin=e] (6) at (6,0){};
\node[tomato dynkin=above] (0) at (3.5,0.25){};
\end{tikzpicture}
,\quad
\atwotypea[{(e\ge1)}]:\quad
\begin{tikzpicture}[anchorbase]
\draw[directed=0.5,white](1,0)--(2,0);
\draw(1,0.025)--(2,0.025);
\draw(1,-0.025)--(2,-0.025);
\draw[directed=0.5](2,0)--(3,0);
\draw[directed=0.5](4,0)--(5,0);
\draw[directed=0.5,white](5,0)--(6,0);
\draw(5,0.025)--(6,0.025);
\draw(5,-0.025)--(6,-0.025);
\node at (3.5,0){$\cdots$};
\node[tomato dynkin] (0) at (1,0){};
\foreach \x [evaluate=\x as \c using {int(\x-1)}] in {2,...,3} {
\node[dynkin=\c] (\x) at (\x,0){};
}
\foreach \x [evaluate=\x as \c using {int(6-\x)}] in {4,5} {
\node[dynkin=e{-}\c] (\x) at (\x,0){};
}
\node[dynkin=e] (6) at (6,0){};
\end{tikzpicture}
,
\\
\athreetypea[{(e>3)}]:\quad
\begin{tikzpicture}[anchorbase,xscale=-1]
\draw[directed=0.5](2,0)--(1,0);
\draw[directed=0.5](2,0)--(2,1);
\draw[directed=0.5](3,0)--(2,0);
\draw[directed=0.5](5,0)--(4,0);
\draw[directed=0.5,white](6,0)--(5,0);
\draw(5,0.025)--(6,0.025);
\draw(5,-0.025)--(6,-0.025);
\node at (3.5,0){$\cdots$};
\foreach \x [evaluate=\x as \c using {int(\x-1)}] in {2,...,3} {
\node[dynkin=e{-}\c] (\x) at (\x,0){};
}
\node[dynkin=e] (1) at (1,0){};
\foreach \x [evaluate=\x as \c using {int(7-\x)}] in {4,5} {
\node[dynkin=\c] (\x) at (\x,0){};
}
\node[dynkin=1] (6) at (6,0){};
\node[tomato dynkin=right] (1) at (2,1){};
\end{tikzpicture}
,\quad
\atypeb[{e>1}]:\quad
\begin{tikzpicture}[anchorbase]
\draw[directed=0.5](1,0)--(2,0);
\draw[directed=0.5](2,1)--(2,0);
\draw[directed=0.5](2,0)--(3,0);
\draw[directed=0.5](4,0)--(5,0);
\draw[directed=0.5,white](5,0)--(6,0);
\draw(5,0.025)--(6,0.025);
\draw(5,-0.025)--(6,-0.025);
\node at (3.5,0){$\cdots$};
\foreach \x [evaluate=\x as \c using {int(\x+0)}] in {1,...,3} {
\node[dynkin=\c] (\x) at (\x,0){};
}
\foreach \x [evaluate=\x as \c using {int(6-\x)}] in {4,5} {
\node[dynkin=e{-}\c] (\x) at (\x,0){};
}
\node[dynkin=e] (6) at (6,0){};
\node at (6,-0.25){$e$};
\node[tomato dynkin=right] (1) at (2,1){};
\end{tikzpicture}
,\\
\atypec[{e>1}]:\quad
\begin{tikzpicture}[anchorbase]
\draw[directed=0.5,white](1,0)--(2,0);
\draw(1,0.025)--(2,0.025);
\draw(1,-0.025)--(2,-0.025);
\draw[directed=0.5](2,0)--(3,0);
\draw[directed=0.5](4,0)--(5,0);
\draw[directed=0.5,white](6,0)--(5,0);
\draw(5,0.025)--(6,0.025);
\draw(5,-0.025)--(6,-0.025);
\node at (3.5,0){$\cdots$};
\node[dynkin=1] (0) at (1,0){};
\node[dynkin=2] (0) at (2,0){};
\node[dynkin=3] (0) at (3,0){};
\node[dynkin=e-1] (0) at (4,0){};
\node[dynkin=e] (0) at (5,0){};
\node[tomato dynkin] (6) at (6,0){};
\end{tikzpicture}
,\quad
\aonetyped[{e>3}]:\quad
\begin{tikzpicture}[anchorbase]
\draw[directed=0.5](2,1)--(2,0);
\draw[directed=0.5](5,0)--(5,1);
\draw[directed=0.5](1,0)--(2,0);
\draw[directed=0.5](2,0)--(3,0);
\draw[directed=0.5](4,0)--(5,0);
\draw[directed=0.5](5,0)--(6,0);
\node at (3.5,0){$\cdots$};
\node[dynkin node={DarkBlue,right,e{-}1}] (1) at (5,1){};
\node[dynkin=e] (1) at (6,0){};
\foreach \x [evaluate=\x as \c using {int(\x+0)}] in {1,...,3} {
\node[dynkin=\c] (\x) at (\x,0){};
}
\foreach \x [evaluate=\x as \c using {int(6-\x)}] in {3,4} {
\node[dynkin=e{-}\c] (\x) at (\x+1,0){};
\node[tomato dynkin=right] (1) at (2,1){};
}
\end{tikzpicture}
,\\
\atwotyped[{(e>1)}]:\quad
\begin{tikzpicture}[anchorbase]
\draw[directed=0.5,white](2,0)--(1,0);
\draw(1,0.025)--(2,0.025);
\draw(1,-0.025)--(2,-0.025);
\draw[directed=0.5](2,0)--(3,0);
\draw[directed=0.5](4,0)--(5,0);
\draw[directed=0.5,white](5,0)--(6,0);
\draw(5,0.025)--(6,0.025);
\draw(5,-0.025)--(6,-0.025);
\node at (3.5,0){$\cdots$};
\node[tomato dynkin] (0) at (1,0){};
\foreach \x [evaluate=\x as \c using {int(\x-1)}] in {2,...,3} {
\node[dynkin=\c] (\x) at (\x,0){};
}
\foreach \x [evaluate=\x as \c using {int(6-\x)}] in {4,5} {
\node[dynkin=e{-}\c] (\x) at (\x,0){};
}
\node[dynkin=e] (6) at (6,0){};
\end{tikzpicture}
,\quad
\athreetyped:\quad
\begin{tikzpicture}[anchorbase]
\draw[directed=0.5](2,0)--(3,0);
\draw[directed=0.5](1,0)--(2,0);
\draw (1,0.05)--(2,0.05);
\draw (1,-0.05)--(2,-0.05);
\foreach \x in {1,...,2} {
\node[dynkin=\x] (\x) at (\x,0){};
}
\node[tomato dynkin] (44) at (3,0){};
\end{tikzpicture}
,
\\
\atypee:\quad
\begin{tikzpicture}[anchorbase]
\draw[directed=0.5](3,0)--(3,1);
\draw[directed=0.5](3,1)--(3,2);
\draw[directed=0.5](1,0)--(2,0);
\draw[directed=0.5](2,0)--(3,0);
\draw[directed=0.5](3,0)--(4,0);
\draw[directed=0.5](4,0)--(5,0);
\foreach \x in {1,...,3} {
\node[dynkin=\x] (\x) at (\x,0){};
}
\foreach \x in {5,...,6} {
\node[dynkin=\x] (\x) at (\x-1,0){};
}
\node[dynkin node={DarkBlue,right,4}] (6) at (3,1){};
\node[tomato dynkin=right] (6) at (3,2){};
\end{tikzpicture}
,\quad
\atwotypee:\quad
\begin{tikzpicture}[anchorbase]
\draw[directed=0.5](4,0)--(5,0);
\draw[directed=0.5,white](2,0)--(3,0);
\draw (2,0.025)--(3,0.025);
\draw (2,-0.025)--(3,-0.025);
\draw[directed=0.5](3,0)--(4,0);
\draw[directed=0.5](1,0)--(2,0);
\foreach \x in {1,...,4} {
\node[dynkin=\x] (\x) at (\x,0){};
}
\node[tomato dynkin] (22) at (5,0){};
\end{tikzpicture}
,
\\
\scalebox{0.95}{$\atypee[7]:\quad
\begin{tikzpicture}[anchorbase]
\draw[directed=0.5](0,0)--(1,0);
\draw[directed=0.5](3,0)--(3,1);
\draw[directed=0.5](1,0)--(2,0);
\draw[directed=0.5](2,0)--(3,0);
\draw[directed=0.5](3,0)--(4,0);
\draw[directed=0.5](4,0)--(5,0);
\draw[directed=0.5](5,0)--(6,0);
\foreach \x in {1,...,3} {
\node[dynkin=\x] (\x) at (\x,0){};
\node at (\x,-0.25){$\x$};
}
\foreach \x in {5,...,7} {
\node[dynkin=\x] (\x) at (\x-1,0){};
}
\node[dynkin node={DarkBlue,right,4}] (7) at (3,1){};
\node[tomato dynkin] (77) at (0,0){};
\end{tikzpicture}
,\quad
\atypee[8]:\quad
\begin{tikzpicture}[anchorbase]
\draw[directed=0.5](7,0)--(8,0);
\draw[directed=0.5](3,0)--(3,1);
\draw[directed=0.5](1,0)--(2,0);
\draw[directed=0.5](2,0)--(3,0);
\draw[directed=0.5](3,0)--(4,0);
\draw[directed=0.5](4,0)--(5,0);
\draw[directed=0.5](5,0)--(6,0);
\draw[directed=0.5](6,0)--(7,0);
\foreach \x in {1,...,3} {
\node[dynkin=\x] (\x) at (\x,0){};
}
\foreach \x in {5,...,8} {
\node[dynkin=\x] (\x) at (\x-1,0){};
}
\node[dynkin node={DarkBlue,right,4}] (8) at (3,1){};
\node[tomato dynkin] (77) at (8,0){};
\end{tikzpicture}$}
,
\\
\atypef:\quad
\begin{tikzpicture}[anchorbase]
\draw[directed=0.5](1,0)--(2,0);
\draw[directed=0.5,white](2,0)--(3,0);
\draw (2,0.025)--(3,0.025);
\draw (2,-0.025)--(3,-0.025);
\draw[directed=0.5](3,0)--(4,0);
\draw[directed=0.5](0,0)--(1,0);
\foreach \x in {1,...,4} {
\node[dynkin=\x] (\x) at (\x,0){};
}
\node[tomato dynkin] (22) at (0,0){};
\end{tikzpicture}
,\quad
\atypeg:\quad
\begin{tikzpicture}[anchorbase]
\draw[directed=0.5](2,0)--(3,0);
\draw[directed=0.5](1,0)--(2,0);
\draw (3,0.05)--(2,0.05);
\draw (3,-0.05)--(2,-0.05);
\node[dynkin=2] (11) at (3,0){};
\node[dynkin=1] (11) at (2,0){};
\node[tomato dynkin] (44) at (1,0){};
\end{tikzpicture}
.
\end{gather*}
The affine node is the red node with label $0$. The choices
of orientation for the simply laced edges are not needed below,
but we expect that these choices give the nicest ghost combinatorics.

\begin{Remark}\label{R:ConsequencesSageMath}
As in \autoref{R:FiniteTypesSageMath}, with a few exceptions our conventions above follow SageMath \cite{sage}.
\end{Remark}

The definitions in \autoref{SS:RecollectionRelations} gives a wKLRW algebra $\WAc(X)$ for each affine quiver.
We start with a slightly disappointing fact:

\begin{Proposition}\label{P:ConsequencesNotCellular}
Let $\ring$ be a field of characteristic $\neq 2$ and suppose that $\quiver$ of type $\aonetyped[e>3]$, $\atypee[6]$, $\atypee[7]$ or
$\atypee[8]$. Then $\WAc[n](X)$ is not cellular for some $n$.
\end{Proposition}

\begin{proof}
According to \cite[Section 5]{KlMu-affine-zigzag-klr},
the direct summand of the cyclotomic KLR algebras indexed by the null root
in these types is Morita equivalent to the so-called zigzag
algebra for the graph obtained by removing the affine node.
These zigzag algebras are not cellular by \cite[Theorem A]{EhTu-zigzag},
so the corresponding cyclotomic KLR algebra is not cellular by Morita invariance
of cellularity \cite[Section 5]{KoXi-cellular-inflation-morita}.
The cyclotomic KLR algebra is an idempotent truncation of $\WAc[n](X)$ so, because taking idempotent truncations
preserves cellularity by \autoref{T:SandwichMain}(a), the statement follows.
\end{proof}

\begin{Remark}\label{R:ConsequencesNotCellular}
We suspect that the restriction to odd characteristic
in \autoref{P:ConsequencesNotCellular} is not necessary. We add this condition only because it is required by \cite[Section 5]{KoXi-cellular-inflation-morita}. We expect that this condition is not needed because we have weakened the condition on the antiinvolution in \upshape(AC${}_{5}$\upshape).
\end{Remark}

\begin{Remark}\label{R:ConsequencesGraphsNotCellular}
% After we shared \autoref{Eq:ConsequencesList},
Robert Muth computed the zigzag-type algebras for $\atypeb[e>1]$
and $\atypec[e>1]$ for the null root, following the ideas of \cite[Section 5]{KlMu-affine-zigzag-klr}. It appears that these algebras are not cellular in type $\atypeb[e>2]$, but the zigzag algebras are probably cellular in types $\atypeb[2]$ and $\atypec[e>1]$.
Extending these ideas, Robert Muth's computations imply that finite types B and D are not cellular in general, see also \autoref{E:TypeDNotGood}.
\end{Remark}

In general, we hope that our methods give a sandwich cellular basis for all types. However, we expect that these bases are only honestly cellular in general in types A and C (affine and finite).

%%%%%%%%%%%%%%%%%%%%%%%%%%%%%%%%%%%%%%%%%%%%%
\bibliographystyle{alphaurl}
%\bibliography{references}
%%%%%%%%%%%%%%%%%%%%%%%%%%%%%%%%%%%%%%%%%%%%%

\end{document}